\pdfoutput=1
\documentclass[preprint,11pt]{elsarticle}
\usepackage[bookmarks=false,colorlinks=true,linkcolor=blue]{hyperref}

\biboptions{sort&compress}% compress bibtex entries

\usepackage{calc}
\usepackage[margin=.95in]{geometry}
\usepackage{mathtools, nccmath}
\usepackage{enumitem}

\newcommand{\citeCount}[1]{}

\usepackage[usenames,dvipsnames,svgnames,table]{xcolor}
\usepackage{amssymb,amsmath}
\usepackage{amsthm}
\usepackage{verbatim}

\usepackage{graphicx}    

\usepackage{tensor}

\usepackage{algorithm,algpseudocode}
\usepackage{algorithmicx}
\algrenewcommand\alglinenumber[1]{\footnotesize #1:} % Algorithm line number font size
\newcommand{\algFontSize}{\footnotesize}

\usepackage{caption}% http://ctan.org/pkg/caption

\usepackage{tikz}
\usetikzlibrary{shapes,arrows,decorations.markings}
\usetikzlibrary{arrows.meta}
\usetikzlibrary{calc}

\newlength{\tfwidth}
\newlength{\tfheight}
\newlength{\tfxa}
\newlength{\tfxb}
\newlength{\tfya}
\newlength{\tfyb}
\newcommand{\trimFigWithBox}[6]{%
\setlength\fboxsep{0pt}%
\setlength\fboxrule{1.0pt}% border thickness
\fbox{\includegraphics[width=#2, clip, trim=#3 #4 #5 #6]{#1}}%
}
\newcommand{\trimFigNoBox}[6]{%
\setlength\fboxsep{1pt}% note: make this 1pt and rule thickness zero so box size matches that below
\setlength\fboxrule{0.0pt}% border thickness
\fbox{\includegraphics[width=#2, clip, trim=#3 #4 #5 #6]{#1}}%
}
\newcommand{\trimFigHeightWithBox}[6]{%
\setlength\fboxsep{0pt}%
\setlength\fboxrule{1.0pt}% border thickness
\fbox{\includegraphics[height=#2, clip, trim=#3 #4 #5 #6]{#1}}%
}
\newcommand{\trimFigHeightNoBox}[6]{%
\setlength\fboxsep{1pt}% note: make this 1pt and rule thickness zero so box size matches that below
\setlength\fboxrule{0.0pt}% border thickness
\fbox{\includegraphics[height=#2, clip, trim=#3 #4 #5 #6]{#1}}%
}
\newsavebox\figBox

\newcommand{\trimw}[6]{%
\sbox\figBox{\includegraphics{#1}}
\setlength{\tfwidth}{\the\wd\figBox}
\setlength{\tfheight}{\the\ht\figBox}
\setlength{\tfxa}{\tfwidth*\real{#3}}%
\setlength{\tfxb}{\tfwidth*\real{#4}}%
\setlength{\tfya}{\tfheight*\real{#5}}%
\setlength{\tfyb}{\tfheight*\real{#6}}%
\trimFigNoBox{#1}{#2}{\tfxa}{\tfya}{\tfxb}{\tfyb}%
}

\newcommand{\trimwb}[6]{%
\sbox\figBox{\includegraphics{#1}}
\setlength{\tfwidth}{\the\wd\figBox}
\setlength{\tfheight}{\the\ht\figBox}
\setlength{\tfxa}{\tfwidth*\real{#3}}%
\setlength{\tfxb}{\tfwidth*\real{#4}}%
\setlength{\tfya}{\tfheight*\real{#5}}%
\setlength{\tfyb}{\tfheight*\real{#6}}%
\trimFigWithBox{#1}{#2}{\tfxa}{\tfya}{\tfxb}{\tfyb}%
}

\newcommand{\trimh}[6]{%
\sbox\figBox{\includegraphics{#1}}
\setlength{\tfwidth}{\the\wd\figBox}
\setlength{\tfheight}{\the\ht\figBox}
\setlength{\tfxa}{\tfwidth*\real{#3}}%
\setlength{\tfxb}{\tfwidth*\real{#4}}%
\setlength{\tfya}{\tfheight*\real{#5}}%
\setlength{\tfyb}{\tfheight*\real{#6}}%
\trimFigHeightNoBox{#1}{#2}{\tfxa}{\tfya}{\tfxb}{\tfyb}%
}

\newcommand{\trimhb}[6]{%
\sbox\figBox{\includegraphics{#1}}
\setlength{\tfwidth}{\the\wd\figBox}
\setlength{\tfheight}{\the\ht\figBox}
\setlength{\tfxa}{\tfwidth*\real{#3}}%
\setlength{\tfxb}{\tfwidth*\real{#4}}%
\setlength{\tfya}{\tfheight*\real{#5}}%
\setlength{\tfyb}{\tfheight*\real{#6}}%
\trimFigHeightWithBox{#1}{#2}{\tfxa}{\tfya}{\tfxb}{\tfyb}%
}
\newcommand{\red}{\color{red}}

\newcommand{\defeq}{\overset{{\rm def}}{=}}% definition = 
\newcommand{\eqdef}{\overset{{\rm def}}{=}}% definition = 

\def\ba#1\ea{\begin{align}#1\end{align}}

\def\bas#1\eas{\begin{align*}#1\end{align*}}

\def\bat#1\eat{\begin{alignat}{3}#1\end{alignat}}

\def\bats#1\eats{\begin{alignat*}{3}#1\end{alignat*}}

\newcommand{\bse}{\begin{subequations}}
\newcommand{\ese}{\end{subequations}}

\newcommand{\bogus}[1]{{}}
\newlength{\ycbTop}% For colour bar
\newlength{\ycbMid}%

\newcommand{\lmax}{{l_{\rm max}}}% denote total number of levels with \ml
\renewcommand{\ll}{{l}}% denote levels with \ll

\newcommand{\jv}{\mathbf{ j}}
\newcommand{\kv}{\mathbf{ k}}

\newcommand{\rv}{\mathbf{ r}}

\newcommand{\xv}{\mathbf{ x}}

\newcommand{\Mv}{\mathbf{ M}}

\newcommand{\Qv}{\mathbf{ Q}}

\newcommand{\Sv}{\mathbf{ S}}

\newcommand{\Uv}{\mathbf{ U}}

\newcommand{\Xv}{\mathbf{ X}}

\newcommand{\Gc}{{\mathcal G}}
\newcommand{\Integer}{\Bbb{Z}}

\newcommand{\eps}{\epsilon}
\def\a{\alpha }    \def\e{\epsilon }
    
\def\th{\theta }

\DeclareMathOperator{\sgn}{sgn}

\newcommand{\Oc}{{\mathcal O}}

\newcommand{\alv}{\boldsymbol{\a}}

\newcommand{\thetav}{\boldsymbol{\theta}}

\newcommand{\thv}{\thetav}
\newcommand{\thvBar}{\bar{\thetav}}
\newcommand{\thetaBar}{\bar{\th}}
\newcommand{\thbar}{\thetaBar}

\newcommand{\PhiVectorF}{%
 \begin{bmatrix}
  \phi_h (\cdot,\ \th_1,\ \th_2) & \phi_h (\cdot,\ \thbar_1,\ \thbar_2) & \phi_h (\cdot,\ \thbar_1,\ \th_2) &  \phi_h (\cdot,\ \th_1,\ \thbar_2)
  \end{bmatrix} 
}    
 
\newcommand{\PhiVectorEx}{%
	\begin{bmatrix}
		\phi_h (\xv,\ \thv) & \phi_h (\xv,\ \thvBar) 
	\end{bmatrix} 
}

\newcommand{\symE}{\mathbf }
\newcommand{\symF}{\hat }

\newcommand{\grid}{\mathcal }
\newcommand{\gridfunspace}{\mathbb }

\graphicspath{{./fig/}}

\usepackage[english]{babel}
\usepackage{amsthm}
\newtheorem{theorem}{Theorem}

\newtheorem{Corollary}{Corollary}

\usepackage{titlesec}
\usepackage{hyperref}

\titleclass{\subsubsubsection}{straight}[\subsection]

\newcounter{subsubsubsection}[subsubsection]
\renewcommand\thesubsubsubsection{\thesubsubsection.\arabic{subsubsubsection}}

\titleformat{\subsubsubsection}
{\itshape\normalsize}{\thesubsubsubsection.}{.5em}{}
\titlespacing{\subsubsubsection}
{0pt}{3.25ex plus 1ex minus .2ex}{0.5ex plus .2ex}

\makeatletter
\def\toclevel@subsubsubsection{4}
\def\l@subsubsubsection{\@dottedtocline{4}{7em}{4em}}
\makeatother

\mathchardef\mhyphen="2D % Define a "math hyphen"
\newcommand{\gcg}[3]{G_H^{#1}[L_H^{#2}\mhyphen I_h^{#3}]}% Galerkin coarse grid operator, order=#1 formed from L_h (order #2) and I_h order #3
\newcommand{\ngcg}[1]{L_H^{#1}}% non_Galerkin coarse grid operator of order #1

\newcommand{\LN}[2]{L_{#1}^{({#2})}[nG]} % non-Galerkin L
	\newcommand{\LNA}[1]{L_{#1}[nG]} % non-Galerkin L, same order as fine-level
	
\newcommand{\LG}[3]{L_{#1}[#2,\ #3]} % Galerkin L; in particular -- 
	\newcommand{\LGt}[2]{L_{#1}[#2]} % with the standard 2o transfers
	\newcommand{\LGG}[1]{L_{#1}[G]} % real Galerkin L
	\newcommand{\LGq}[2]{L_{#1}^{(#2)}[G]} % finest-level order q
	\newcommand{\LGA}[2]{L_{#1}^{({#2})}[G1]} % 1st-gen Galerkin L
		\newcommand{\LGAA}[1]{L_{#1}[G1]} % 1st-gen Galerkin L, same order as fine-level

\newcommand{\LGHh}[2]{\LGH{L_h^{(#1)}}{#2}}

\newcommand{\LNH}[1]{\LN{H}{#1}} % non-Galerkin L
\newcommand{\LNAH}{\LNA{H}} % non-Galerkin L, same order as fine-level

\newcommand{\LGH}[2]{\LG{H}{#1}{#2}} % Galerkin L; in particular -- 
\newcommand{\LGGH}{\LGG{H}} % real Galerkin L
\newcommand{\LGqH}[1]{\LGq{H}{#1}} % finest-level order q
\newcommand{\LGAH}[1]{\LGA{H}{#1}} % 1st-gen Galerkin L
\newcommand{\LGAAH}{\LGAA{H}} % 1st-gen Galerkin L, same order as fine-level

\newcommand{\commentB}[1]{}

\newcommand{\colourTwoD}{orange}% colour for title for 2D results

\newcommand{\rbDot}{node[fill,circle,scale=.6]{}}
\newcommand{\vdot}{node[fill,circle,scale=.45]{}}
\newcommand{\gdot}{node[fill,circle,scale=.30]{}}
\newcommand{\gdots}{node[fill,circle,scale=.45]{}}
\newcommand{\stencilColour}{blue}% colour of stencil in grid cartoons
\newcommand{\coarseColour}{orange}% colour of coarse grid in grid cartoons
\newcommand{\coarseColourRB}{red}% colour of coarse grid in grid cartoons for red-black coarsening

\newcommand{\standardGridOrderTwo}[1]{% 
  \begin{scope}[#1]%
    \draw[step=.25cm,black,thick] (-.5,-.5) grid (.5,.5);%
    \draw[\stencilColour,line width=2pt] (-.25,0) -- (.25,0); % stencil
    \draw[\stencilColour,line width=2pt] (0,-.25) -- (0,.25); % stencil
    \draw[\stencilColour] (-.25,0) \gdots (0,0) \gdots (.25,0) \gdots (0,-.25) \gdots (0,.25) \gdots ; 
    \draw[\coarseColour] (-.5,0) \gdot(0,0) \gdot (.5,0) \gdot (0,.5) \gdot(0,-.5) \gdot; 
    \draw[\coarseColour] (-.5,-.5) \gdot(0,0) \gdot (-.5,.5) \gdot (.5,-.5) \gdot(.5,.5) \gdot; 
  \end{scope}%
  }

\newcommand{\standardGridOrderFour}[1]{% 
  \begin{scope}[#1]%
    \draw[step=.25cm,black,thick] (-.5,-.5) grid (.5,.5);%
    \draw[\stencilColour,line width=2pt] (-.5,0) -- (.5,0); % stencil
    \draw[\stencilColour,line width=2pt] (0,-.5) -- (0,.5); % stencil
    \draw[\stencilColour] (-.5,0) \gdots (-.25,0) \gdots (0,0) \gdots (.25,0) \gdots (.5,.0) \gdots
                          (0,-.5) \gdots (0,-.25) \gdots (0,.25) \gdots (0,.5) \gdots; 
    \draw[\coarseColour] (-.5,0) \gdot(0,0) \gdot (.5,0) \gdot (0,.5) \gdot(0,-.5) \gdot; 
    \draw[\coarseColour] (-.5,-.5) \gdot(0,0) \gdot (-.5,.5) \gdot (.5,-.5) \gdot(.5,.5) \gdot; 
  \end{scope}%
  }

\newcommand{\standardGridOrderSix}[1]{% 
  \begin{scope}[#1]%
    \draw[step=.25cm,black,thick] (-.5,-.5) grid (.5,.5);%
    \draw[\stencilColour,line width=2pt] (-.75,0) -- (.75,0); % stencil
    \draw[\stencilColour,line width=2pt] (0,-.75) -- (0,.75); % stencil
    \draw[\stencilColour] (-.5,0) \gdots (-.25,0) \gdots (0,0) \gdots (.25,0) \gdots (.5,.0) \gdots
                          (0,-.5) \gdots (0,-.25) \gdots (0,.25) \gdots (0,.5) \gdots
                          (0,-.75) \gdots (0,.75) \gdots (-.75,0) \gdots (.75,.0) \gdots; 
    \draw[\coarseColour] (-.5,0) \gdot(0,0) \gdot (.5,0) \gdot (0,.5) \gdot(0,-.5) \gdot; 
    \draw[\coarseColour] (-.5,-.5) \gdot(0,0) \gdot (-.5,.5) \gdot (.5,-.5) \gdot(.5,.5) \gdot; 
  \end{scope}%
  }

\newcommand{\rotatedGrid}[1]{% 
  \begin{scope}[#1]%
    \draw[step=.25cm,black,thick] (-.5,-.5) grid (.5,.5);%
    \draw[\stencilColour,line width=2pt] (-.25,0) -- (.25,0); % stencil
    \draw[\stencilColour,line width=2pt] (0,-.25) -- (0,.25); % stencil
    \draw[\stencilColour] (-.25,0) \gdots (0,0) \gdots (.25,0) \gdots (0,-.25) \gdots (0,.25) \gdots;  % stencil dots 
    \draw[\coarseColourRB] (-.25,-.25) \gdot(0,0) \gdot (-.25,.25) \gdot (.25,-.25) \gdot(.25,.25) \gdot 
                           (-.5,-.5) \gdot(0,-.5) \gdot (.5,.-.5)  \gdot (-.5,.0) \gdot 
                           (-.5,+.5) \gdot(0,+.5) \gdot (.5,.+.5)  \gdot (+.5,.0) \gdot; 
    \draw[\coarseColourRB,thick] (-.5,-.5) -- (.5,.5);
    \draw[\coarseColourRB,thick] (-.5,.5) -- (.5,-.5);
    \draw[\coarseColourRB,thick] (-.5,0) -- (0,-.5) -- (.5,0) -- (0,.5) -- (-.5,0);

  \end{scope}%
  }

\newcommand{\VcycleTwoLevel}[2]{%  
   \begin{scope}[#1]
     \draw[-,thick] (-.5,.5)  \vdot --
                    (0,-.5)   \vdot --
                    (.5,.5)   \vdot;
     \draw (0,.5) node[anchor=north,yshift=6pt] {\scriptsize #2}; 
   \end{scope}
}

\newcommand{\VcycleTwoLevelBoxed}[2]{%  
   \begin{scope}[#1]
     \draw[-,thick,fill=white] (-.65,-.65) -- (.65,-.65) -- (.65,.65) -- (-.65,.65) -- (-.65,-.65); % border
     \draw[-,thick] (-.5,.5)  \vdot --
                    (0,-.5)   \vdot --
                    (.5,.5)   \vdot;
     \draw (0,.5) node[anchor=north,yshift=6pt] {\scriptsize #2}; 
   \end{scope}
}

\newcommand{\VcycleThreeLevel}[3]{%  
   \begin{scope}[#1]
     \draw[-,thick] (-.5,.5)  \vdot --
                    (-.25,0)   \vdot --
                    (0,-.5)   \vdot --
                    (.25,0)   \vdot --
                    (.5,.5)   \vdot;
     \draw (0,.5) node[anchor=north,yshift=6pt] {\scriptsize #2}; 
     \draw (0,-.5) node[anchor=west,xshift=0pt] {\scriptsize #3}; 
   \end{scope}
}

\newcommand{\WcycleThreeLevel}[3]{%  
   \begin{scope}[#1]
     \draw[-,thick] (   -.5,.50)  \vdot --
                    (-.3333, .0)  \vdot --
                    (-.1667,-.5)  \vdot --
                    (     0,  0)   \vdot --
                    ( .1667,-.5)   \vdot --
                    ( .3333, 0.)   \vdot --
                    (    .5,.5)   \vdot;
     \draw (0,.5) node[anchor=north,yshift=6pt] {\scriptsize #2}; 
     \draw ( .1667,-.5) node[anchor=west,xshift=0pt] {\scriptsize #3}; 
   \end{scope}
}

\newcommand{\orderTwoVcycleTwoLevel}[2]{%
 \begin{scope}[#1]
   \draw[-,thick,fill=white] (-.65,-.65) -- (1.9,-.65) -- (1.9,.65) -- (-.65,.65) -- (-.65,-.65); % border
   \standardGridOrderTwo{xshift=0cm,yshift=0cm}; 
   \VcycleTwoLevel{xshift=1.25cm,yshift=0cm}{#2}; 
 \end{scope}
}
\newcommand{\orderTwoVcycleThreeLevel}[3]{%
 \begin{scope}[#1]
   \draw[-,thick,fill=white] (-.65,-.65) -- (1.9,-.65) -- (1.9,.65) -- (-.65,.65) -- (-.65,-.65); % border
   \standardGridOrderTwo{xshift=0cm,yshift=0cm}; 
   \VcycleThreeLevel{xshift=1.25cm,yshift=0cm}{#2}{#3}; 
 \end{scope}
}

\newcommand{\orderFourVcycleTwoLevel}[2]{%
 \begin{scope}[#1]
   \draw[-,thick,fill=white] (-.65,-.65) -- (1.9,-.65) -- (1.9,.65) -- (-.65,.65) -- (-.65,-.65); % border
   \standardGridOrderFour{xshift=0cm,yshift=0cm}; 
   \VcycleTwoLevel{xshift=1.25cm,yshift=0cm}{#2}; 
 \end{scope}
}
\newcommand{\orderFourVcycleThreeLevel}[2]{%
 \begin{scope}[#1]
   \draw[-,thick,fill=white] (-.65,-.65) -- (1.9,-.65) -- (1.9,.65) -- (-.65,.65) -- (-.65,-.65); % border
   \standardGridOrderFour{xshift=0cm,yshift=0cm}; 
   \VcycleThreeLevel{xshift=1.25cm,yshift=0cm}{#2}; 
 \end{scope}
}

\newcommand{\orderSixVcycleTwoLevel}[2]{%
 \begin{scope}[#1]
   \draw[-,thick,fill=white] (-.65,-.65) -- (1.9,-.65) -- (1.9,.65) -- (-.65,.65) -- (-.65,-.65); % border
   \standardGridOrderSix{xshift=0cm,yshift=0cm}; 
   \VcycleTwoLevel{xshift=1.25cm,yshift=0cm}{#2}; 
 \end{scope}
}
\newcommand{\orderSixVcycleThreeLevel}[2]{%
 \begin{scope}[#1]
   \draw[-,thick,fill=white] (-.65,-.65) -- (1.9,-.65) -- (1.9,.65) -- (-.65,.65) -- (-.65,-.65); % border
   \standardGridOrderSix{xshift=0cm,yshift=0cm}; 
   \VcycleThreeLevel{xshift=1.25cm,yshift=0cm}{#2}; 
 \end{scope}
}

\newcommand{\orderTwoVcycleTwoLevelRotated}[2]{%
 \begin{scope}[#1]
   \draw[-,thick,fill=white] (-.65,-.65) -- (1.9,-.65) -- (1.9,.65) -- (-.65,.65) -- (-.65,-.65); % border
   \rotatedGrid{xshift=0cm,yshift=0cm}; 
   \VcycleTwoLevel{xshift=1.25cm,yshift=0cm}{#2}; 
 \end{scope}
}
\newcommand{\orderTwoVcycleThreeLevelRotated}[3]{%
 \begin{scope}[#1]
   \draw[-,thick,fill=white] (-.65,-.65) -- (1.9,-.65) -- (1.9,.65) -- (-.65,.65) -- (-.65,-.65); % border
   \rotatedGrid{xshift=0cm,yshift=0cm}; 
   \VcycleThreeLevel{xshift=1.25cm,yshift=0cm}{#2}{#3}; 
 \end{scope}
}

\newcommand{\orderTwoWcycleThreeLevelRotated}[3]{%
 \begin{scope}[#1]
   \draw[-,thick,fill=white] (-.65,-.65) -- (1.9,-.65) -- (1.9,.65) -- (-.65,.65) -- (-.65,-.65); % border
   \rotatedGrid{xshift=0cm,yshift=0cm}; 
   \WcycleThreeLevel{xshift=1.25cm,yshift=0cm}{#2}{#3}; 
 \end{scope}
}

\newcommand{\circleVortex}[3]{%
    \begin{scope}[xshift=#1cm,yshift=#2cm]%
      \draw[blue,xshift=0pt,yshift=0pt,thick] (0,0) circle (.3cm);%
      \draw[<-,thick,blue,xshift=0,yshift=.3cm] (0,.0) -- (0,.25) node[anchor=south,yshift=-2pt] {\scriptsize #3};
    \end{scope}%
}

\newcommand{\circleLeftLabel}[3]{%
    \begin{scope}[xshift=#1cm,yshift=#2cm]%
      \draw[blue,xshift=0pt,yshift=0pt,thick] (0,0) circle (.3cm);%
      \draw[<-,thick,blue,xshift=.3cm,yshift=0cm] (0,.0) -- (.25,0) node[anchor=west,xshift=-2pt] {\scriptsize #3};
    \end{scope}%
}

\newcommand{\circleRightLabel}[3]{%
    \begin{scope}[xshift=#1cm,yshift=#2cm]%
      \draw[blue,xshift=0pt,yshift=0pt,thick] (0,0) circle (.3cm);%
      \draw[<-,thick,blue,xshift=-.3cm,yshift=0cm] (0,.0) -- (-.25,0) node[anchor=east,xshift=+2pt] {\scriptsize #3};
    \end{scope}%
}

\newcommand{\circleLabel}[5]{%
    \begin{scope}[xshift=#1cm,yshift=#2cm]%
      \draw[#4,xshift=0pt,yshift=0pt,thick] (0,0) circle (#5cm);%
      \draw[<-,thick,#4,xshift=#5cm,yshift=0cm] (0,.0) -- (.5,.0)  node[anchor=west,xshift=-2pt,draw=#4,fill=#4!10] {\scriptsize #3};
    \end{scope}%
}

\newcommand{\circleLabelLeft}[5]{%
    \begin{scope}[xshift=#1cm,yshift=#2cm]%
      \draw[#4,xshift=0pt,yshift=0pt,thick] (0,0) circle (#5cm);%
      \draw[<-,thick,#4,xshift=-#5cm,yshift=0cm] (0,.0) -- (-.5,.0)  node[anchor=east,xshift=2pt,draw=#4,fill=#4!10] {\scriptsize #3};
    \end{scope}%
}

\newcommand{\ellipseLabel}[6]{%
    \begin{scope}[xshift=#1cm,yshift=#2cm]%
      \draw[#4,xshift=0pt,yshift=0pt,thick] (0,0) ellipse (#5cm and #6cm);%
      \draw[<-,thick,#4,xshift=#5cm,yshift=0cm] (0,.0) -- (.5,.0)  node[anchor=west,xshift=-2pt,draw=#4,fill=#4!10] {\scriptsize #3};
    \end{scope}%
}

\newcommand{\Etheta}{\mathbb{E}_{\th}}
\newcommand{\ThetaLowStd}{\Theta_{std}^{low}}
\newcommand{\ThetaLowRB}{\Theta_{rb}^{\text{low}}}
\newcommand{\Ebb}{\mathbb{E}}
\newcommand{\Fbb}{\mathbb{F}}
\newcommand{\Ethetav}{\mathbb{E}_{\thv}}
\newcommand{\Fthetav}{\mathbb{F}_{\thv}}

\newcommand{\rTarget}{r_{\text{target}}}% target grid spacing

\newcommand{\DIM}[1][2]{\textcolor{teal}{\textbf{#1}D}}
\newcommand{\ORD}[1][2]{\textcolor{Sepia}{order \textbf{#1}}}
\newcommand{\CC}[1][standard]{\textcolor{olive}{\textbf{#1} coarsening}}
\newcommand{\Cc}[1][S]{\textcolor{olive}{\textbf{#1}C}}
\newcommand{\LEV}[1][2]{\textcolor{Periwinkle}{\textbf{#1}-level}}
\newcommand{\CYC}[1]{\textcolor{Green}{\textbf{#1}}}

\newcommand{\old}[1]{{\color{OrangeRed}old: #1}}
\newcommand{\cut}[1]{{\color{black} #1}} % may be cut for short version
\newcommand{\edit}[1]{{\color{black} #1}} % need to be editted for short version

\newcommand{\wdh}[1]{{\color{blue}  #1}}
\newcommand{\kl}[1]{{\color{black} #1}}
\newcommand{\kP}[1]{{\color{black} #1}}

\begin{document}
	
	% space above and below equation:
	\setlength{\belowdisplayskip}{.5\baselineskip}
	\setlength{\belowdisplayshortskip}{.3\baselineskip} %short lines
	\setlength{\abovedisplayskip}{.5\baselineskip}
	\setlength{\abovedisplayshortskip}{.3\baselineskip}
	
	% space above and below figure:
	\setlength{\intextsep}{10pt plus 1.0pt minus 2.0pt}
	
	% space above and below itemize/enumerate:
	\setlist{noitemsep, topsep=0pt}

%---------- Title Page ----------------------------
\begin{frontmatter}
  \title{Multigrid with Nonstandard Coarsening
   % \today
  }
\author[rpi]{Kamala Liu\corref{cor1}\fnref{NSFgrantNew}}
\ead{liuc10@rpi.edu}

\author[rpi]{William D.~Henshaw\corref{cor1}\fnref{NSFgrantNew}}
\ead{henshw@rpi.edu}

\address[rpi]{Department of Mathematical Sciences, Rensselaer Polytechnic Institute, Troy, NY 12180, USA}

\cortext[cor1]{Department of Mathematical Sciences, Rensselaer Polytechnic Institute, 110 8th Street, Troy, NY 12180, USA.}

% \fntext[DOEThanks]{This work was performed under DOE contracts from the ASCR Applied Math Program.}

\fntext[NSFgrantNew]{Research supported by the National Science Foundation under grants DMS-1519934 and DMS-1818926.}
% \fntext[NSFgrantNew]{Research supported by the National Science Foundation under grant DMS-1519934.}

\begin{abstract}
We consider the numerical solution of Poisson's equation on structured grids using geometric multigrid
with nonstandard coarse grids and coarse level operators.
%% motivaion
We are motivated by the problem of developing high-order accurate numerical solvers for elliptic boundary value problems on complex geometry using overset grids. Overset grids are typically dominated by large Cartesian background
grids and thus fast solvers for Cartesian grids are highly desired. 
For flexibility in grid generation we would like to consider coarsening factors other than two, and lower-order accurate coarse-level approximations. 
% orderC
We show that second-order accurate coarse-level approximations are very effective for fourth- or sixth-order accurate fine-level finite difference discretizations. We study the use of different Galerkin and non-Galerkin coarse-level operators. 
We use red-black smoothers with a relaxation parameter $\omega$. Using local Fourier analysis we choose $\omega$ and the coarse-level operators to optimize the overall multigrid convergence rate.
% rbC
Motivated by the use of red-black smoothers in one dimension that can result in
a direct solver for the standard second-order accurate discretization to Poisson's equation, 
we show that this direct-solver property can be extended to two dimensions using a rotated grid that results from red-black coarsening.
We evaluate the use of red-black coarsening in more general settings.
% rC
We also study grid coarsening by a general factor and show that good convergence rates are retained for a range of coarsening factors near two. We ask the question of which coarsening factor leads to the most efficient algorithm.
\end{abstract}

\end{frontmatter}

\clearpage
\tableofcontents %65 pages

% ----------------------------------------------------------------------------------------------------
\clearpage
\section{Introduction}\label{Sec: intro}

%% motivation ----------------------------------------------------------------------------------------
We are motivated by the solution of elliptic partial differential equations (PDE) and boundary value problems in complex geometry using composite overlapping grids.
As shown in Figure~\ref{fig:introCompositeGrid}, a composite or ``overset'' grid consists of multiple overlapping structured grids used to cover a complex geometry. Solution values are interpolated at internal interpolation boundaries where two component grids overlap.
These grids have been shown to be very effective at efficiently solving a wide class of problems
including low speed flows~\cite{ICNS,ChandHenshawLundquistSinger2010,fins2020}, high-speed flows~\cite{reactamr2003b,mog2006,pog2008a,KapilaSchwendemanGambinoHenshaw2015},
fluid structure interactions~\cite{th2009,fis2014,fib2014,beamins2016,champ2016,flunsi2016,rbins3d2018,fibr2019}
and electromagnetics~\cite{max2006b,mxsosup2018,adegdm2020}. 
The approach is especially useful for problems with moving geometry.
Multigrid methods have been shown to be quite effective for solving problems on overset grids~\cite{CGMG,automg}, being fast and efficient and
having low startup costs as grids change in a moving geometry scenario.

{% -------------------
\newcommand{\labelsize}{\scriptsize}
\newcommand{\labelcol}{violet}
\newcommand{\figWidth}{5cm}%
\newcommand{\trimfig}[2]{\trimw{#1}{#2}{.05}{.05}{.0}{.0}}
\begin{figure}[htb]
\begin{center}
\begin{tikzpicture}[scale=1]
  \useasboundingbox (0,0.5) rectangle (15.5,5.5);  % set the bounding box (so we have less surrounding white space)
  \draw(0.0,0.0) node[anchor=south west,xshift=-15pt,yshift=-8pt] {\trimfig{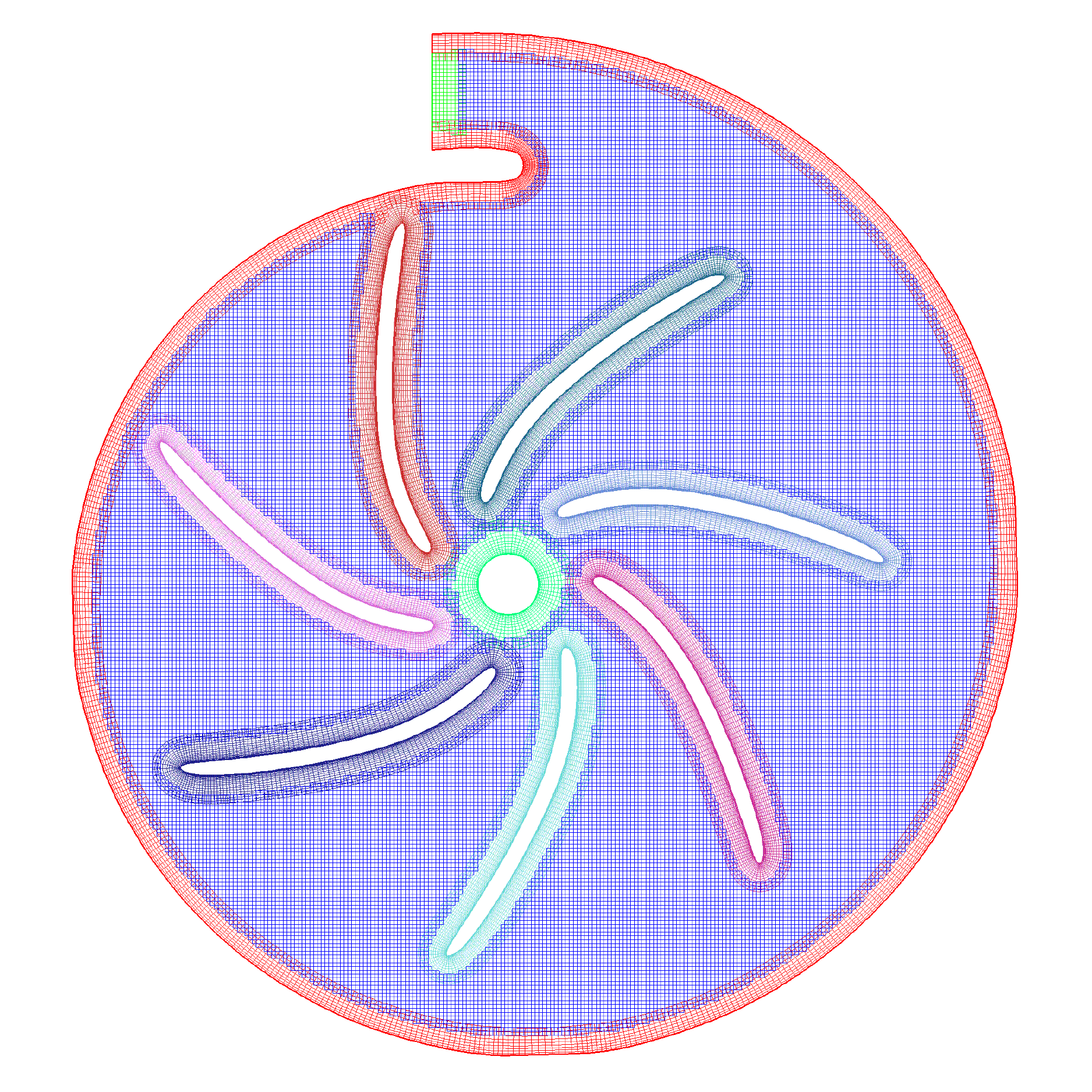}{\figWidth}};
  \draw(5.2,0.0) node[anchor=south west,xshift=-15pt,yshift=-8pt] {\trimfig{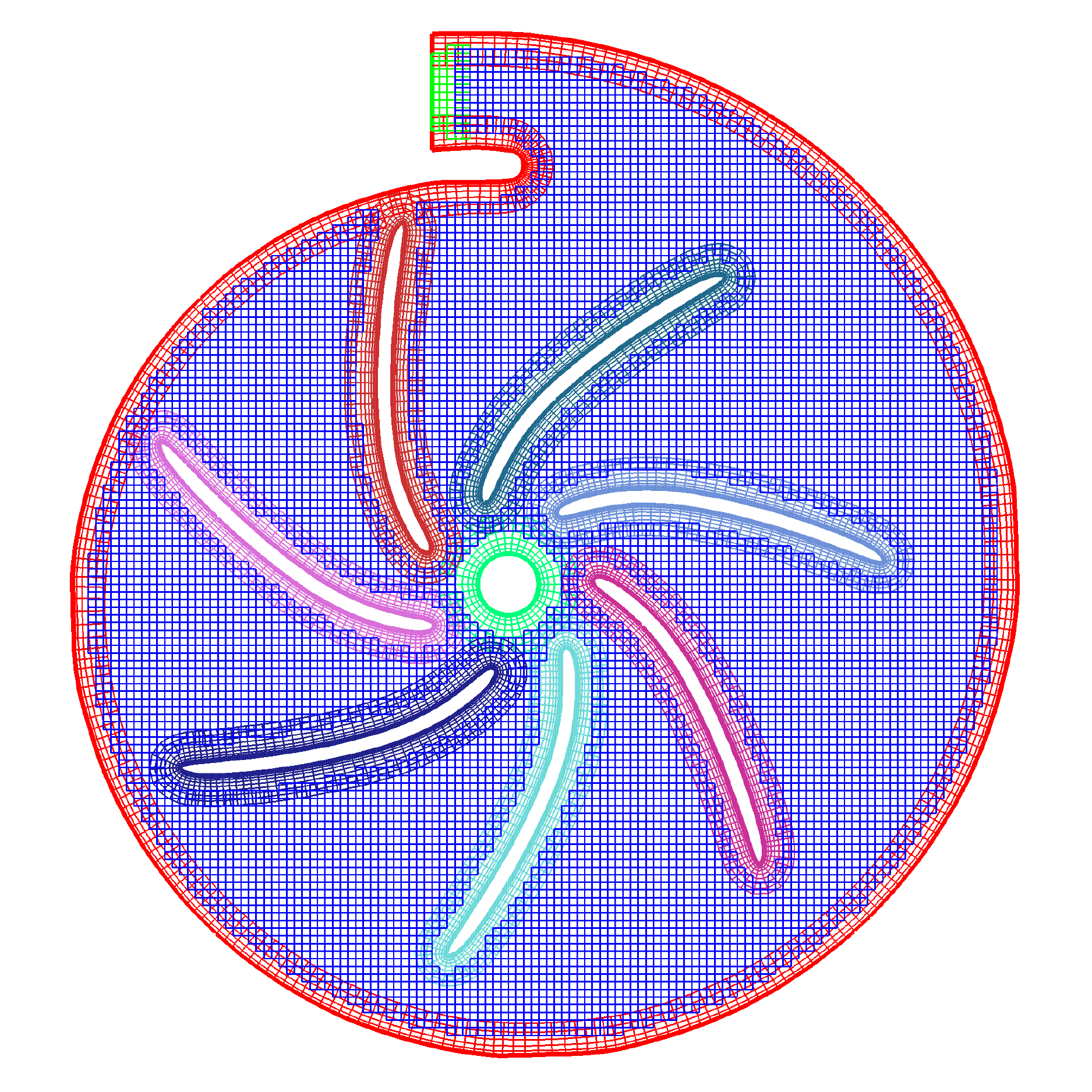}{\figWidth}};
  \draw(10.4,0.0) node[anchor=south west,xshift=-15pt,yshift=-8pt] {\trimfig{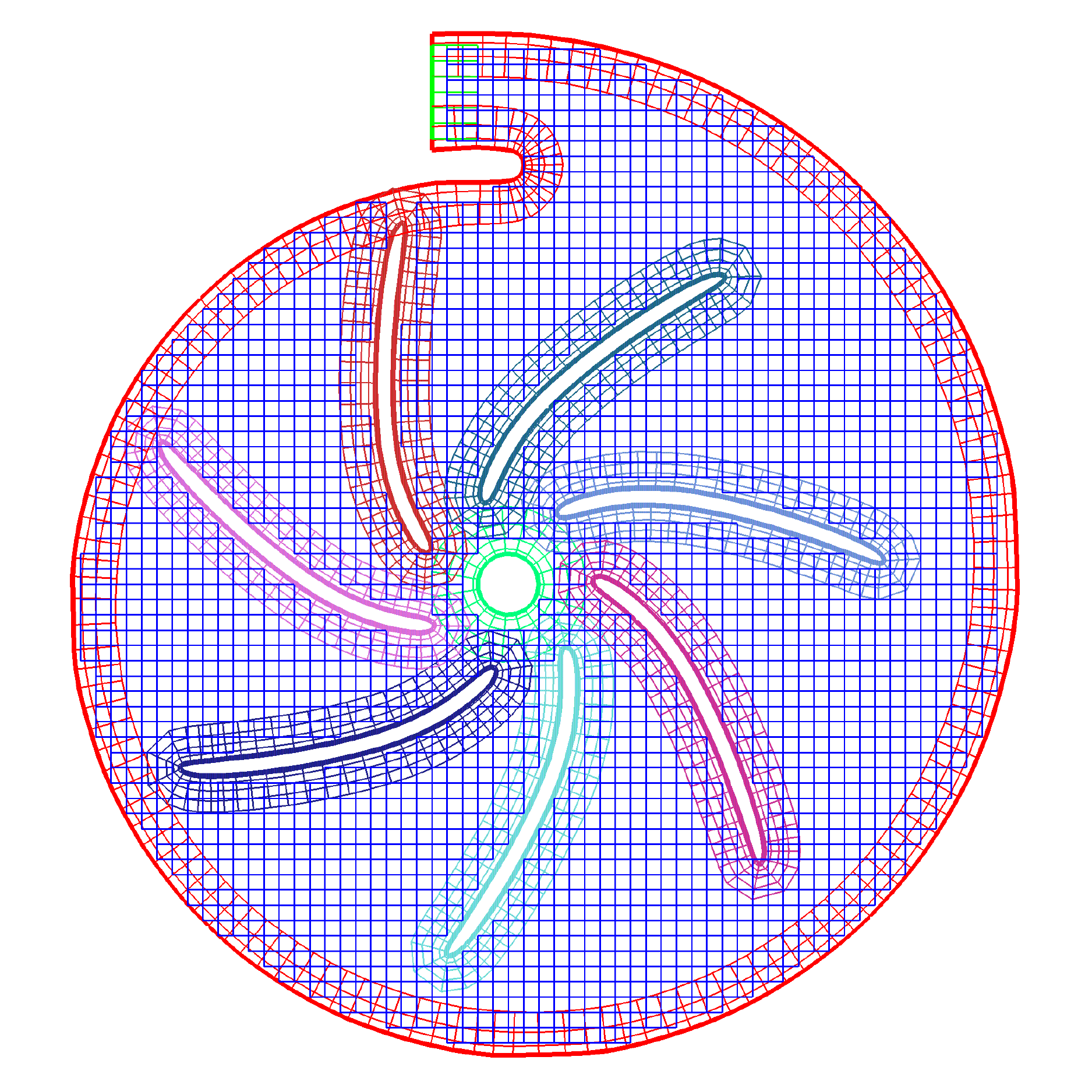}{\figWidth}};

  \draw (2.5,5) node[anchor=south,draw=\labelcol,fill=\labelcol!10] {\labelsize Fine grid $\ll=0$};

  \draw (7.7,5) node[anchor=south,draw=\labelcol,fill=\labelcol!10] {\labelsize Coarse level $\ll=1$};

  \draw (12.7,5) node[anchor=south,draw=\labelcol,fill=\labelcol!10] {\labelsize Coarse level $\ll=2$};

% grid:
%  \draw[step=1cm,gray] (0,0) grid (15.5,5.5);
\end{tikzpicture}
\end{center}
\caption{Overset grid for a centrifugal pump, fine grid and two multigrid coarsenings. As the mesh is refined further the
  majority of grid points will belong to the blue Cartesian background grid. Ideally,
  multigrid convergence rates for this overset grid should approach the convergence rates and efficiency of a single Cartesian grid.
 }
  \label{fig:introCompositeGrid}
\end{figure}
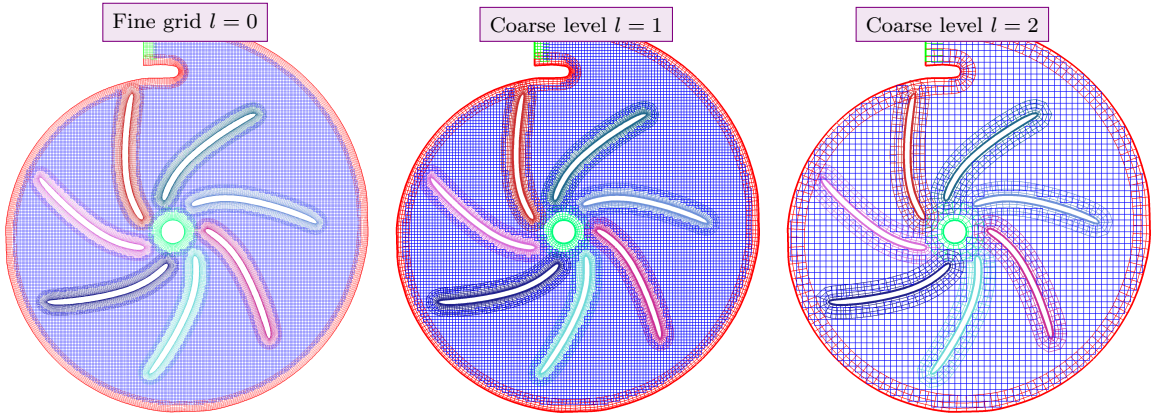
}

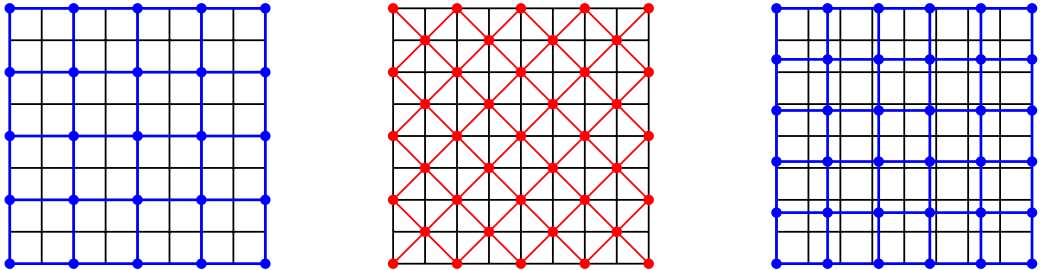
\begin{figure}[hbt]
  \begin{center}
  \resizebox{12cm}{!}{% START resize box
    \begin{tikzpicture}
    \useasboundingbox (0,.5) rectangle (14,4.25);
     % \draw[<-,ultra thick,red,dash pattern= on 8pt off 3pt] (0,3) -- (2,3);
     %

     % ---- standard coarsening --
     \begin{scope}[xshift=-1cm,yshift=0cm]  

        \draw[step=.5cm,black,thick] (0,0) grid (4,4);
        \draw[step=1cm,blue,very thick] (0,0) grid (4,4);
         \foreach \x in {0,1,...,4}
         \foreach \y in {0,1,...,4}
            \draw[blue] (\x,\y) \vdot;
     \end{scope}

     % ---- Red-Black coarsening --
     \begin{scope}[xshift=5cm,yshift=0cm]  

        \draw[step=.5cm,black,thick] (0,0) grid (4,4);
        % \draw[step=1cm,blue,very thick] (0,0) grid (4,4);
         \foreach \x in {0,1,...,4}
         \foreach \y in {0,1,...,4}
            \draw[red] (\x,\y) \vdot;
         \foreach \x in {.5,1.5,...,3.5}
         \foreach \y in {.5,1.5,...,3.5}
         {
            \draw[red] (\x,\y) \vdot;
            % -- lines on coarse grid ---
            \draw[red,thick] (\x-.5,\y-.5) -- (\x+.5,\y+.5); 
            \draw[red,thick] (\x-.5,\y+.5) -- (\x+.5,\y-.5); 
         }
         
     \end{scope}

     % ---- general coarsening --
     \begin{scope}[xshift=11cm,yshift=0cm]  

        \draw[step=.5cm,black,thick] (0,0) grid (4,4);
        \draw[step=.8cm,blue,very thick] (0,0) grid (4,4);
         \foreach \x in {0,.8,...,4.001}
         \foreach \y in {0,.8,...,4.001}
            \draw[blue] (\x,\y) \vdot;
     \end{scope}
% 
  % \draw[step=1cm,gray] (0,0) grid (14,4);
  \end{tikzpicture}
  }% --- END RESIZE BOX
  \end{center}
  \caption{Some coarsening strategies in two dimensions. 
  	   Left: standard coarsening, $H=2 h$.
  	   Middle: Red-Black coarsening results
  	   in a rotated grid with spacing $H=\sqrt{2} h$.
       Right: coarsening by a general factor $H = r h$. 
       }
  \label{fig:coarseningStrategies}
\end{figure}

In practical problems, the generation of high-quality grids is one of the most important steps in solving a PDE.
For a single Cartesian grid it is easy to construct a grid that can be coarsened many times by a factor of two, this is called standard coarsening (see Figure~\ref{fig:coarseningStrategies}, left). 
For flexibility in overset grid generation we do not want to place such constraints on the number of grid points in each individual component grids. % rC
Furthermore, a primary technical challenge in applying multigrid methods to overset grids is the automatic generation of coarser levels~\cite{CGMG}.
High-order accurate approximations require multiple layers of interpolation points to support wide stencils and sufficient overlap between component grids. 
We would like to relax these grid generation requirements on coarse levels and use lower-order accurate approximations instead so that more valid coarse levels can be generated. % orderC

{% -------------------
\newcommand{\labelsize}{\tiny\sffamily}
\newcommand{\labelcol}{violet}
\newcommand{\figWidth}{5cm}%
\newcommand{\trimfig}[2]{\trimw{#1}{#2}{.05}{.05}{.0}{.0}}
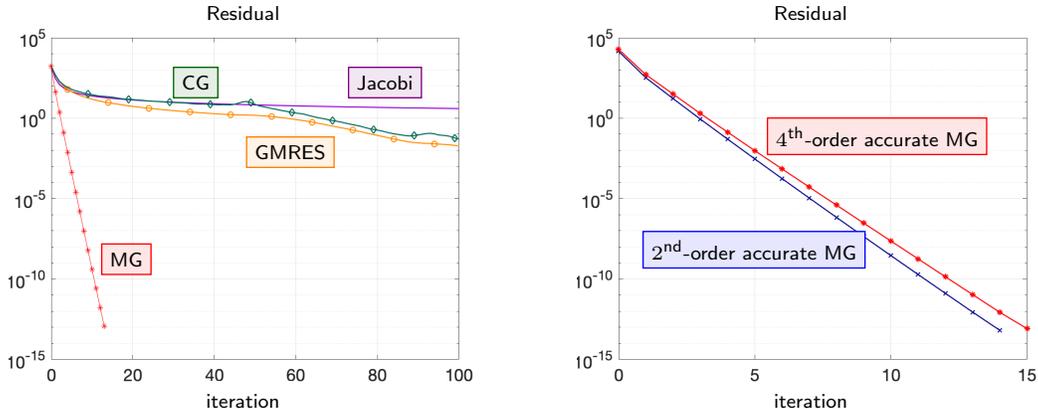
\begin{figure}[htb]
\begin{center}
\resizebox{14cm}{!}{% START resize box
\begin{tikzpicture}[scale=1]
  \useasboundingbox (0,0.25)  rectangle (11,4);  % set the bounding box (so we have less surrounding white space)

\begin{scope}[xshift=0cm]
   \draw(0,0) node[anchor=south west,xshift=-15pt,yshift=-8pt]{\trimfig{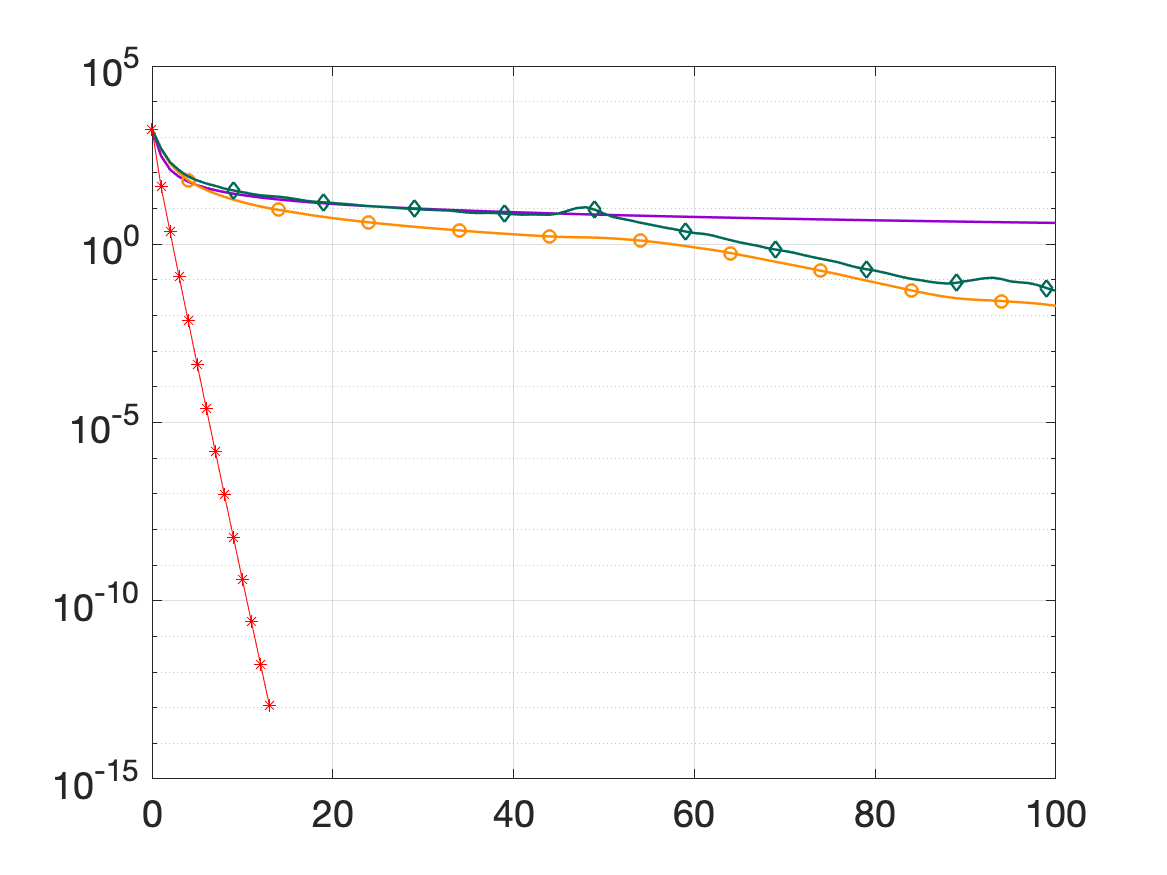}{\figWidth}};
   \draw (1,1.4) node[anchor=west,draw=red,fill=red!10,inner sep=2.5pt] {\labelsize MG};
   \draw (2,3.1) node[anchor=south,draw=DarkGreen,fill=DarkGreen!10,inner sep=2.5pt] {\labelsize CG};
   \draw (3,2.7) node[anchor=north,draw=orange,fill=orange!10,inner sep=2.5pt] {\labelsize GMRES};

   \draw (4,3.1) node[anchor=south,draw=violet,fill=violet!10,inner sep=2.5pt] {\labelsize Jacobi};
   \draw (2.5,0.2) node[anchor=north,fill=white,yshift=-2.5pt] {\labelsize iteration};
   \draw (2.5,3.9) node[anchor=south,fill=white,yshift=-3pt] {\labelsize Residual};
\end{scope}

   \begin{scope}[xshift=6cm]
     \draw(0,0) node[anchor=south west,xshift=-15pt,yshift=-8pt]{\trimfig{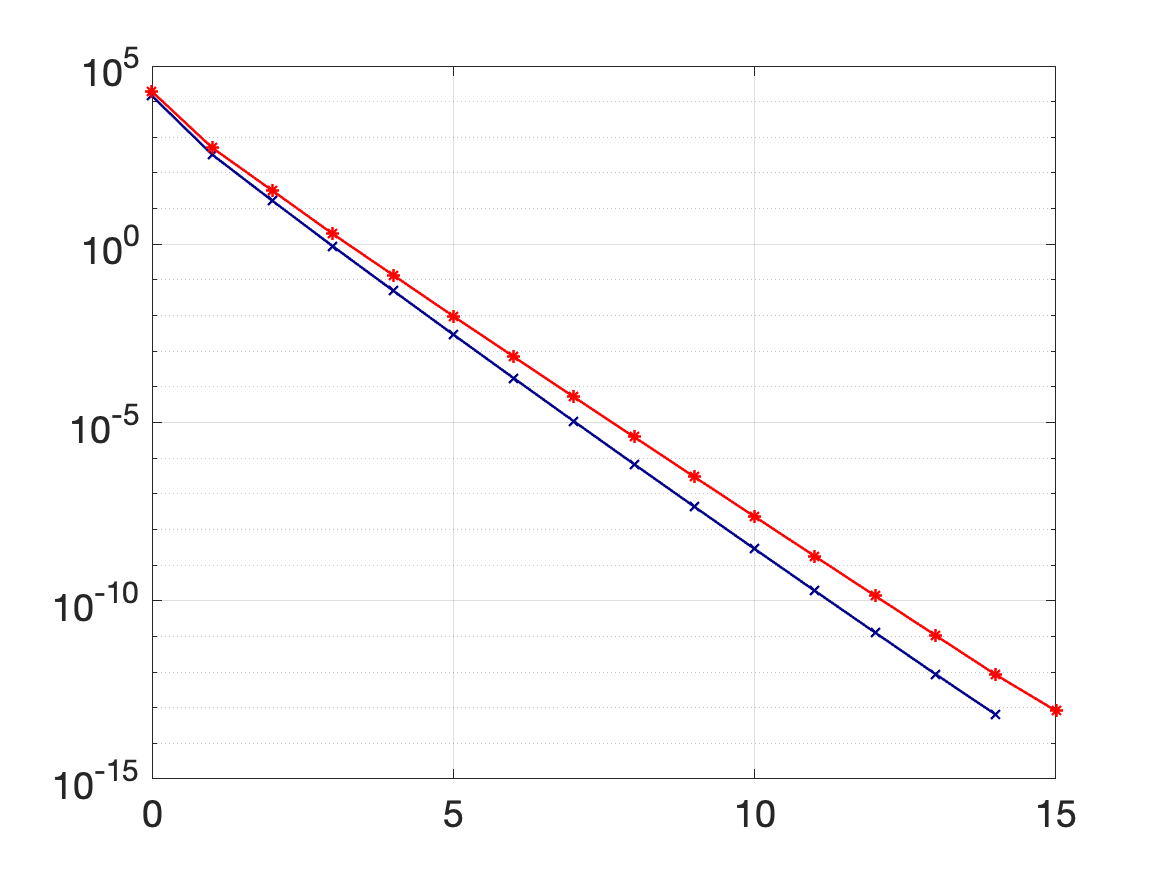}{\figWidth}};
     \draw (3.2,2.5) node[anchor=south,draw=red,fill=red!10,inner sep=2.5pt] {\labelsize $4^{\textsf{th}}$-order accurate MG};
     \draw (1.9,1.3) node[anchor=south,draw=blue,fill=blue!10,inner sep=2.5pt] {\labelsize $2^{\textsf{nd}}$-order accurate MG};
     \draw (2.5,0.2) node[anchor=north,fill=white,yshift=-2.5pt] {\labelsize iteration};
     \draw (2.5,3.9) node[anchor=south,fill=white,yshift=-3pt] {\labelsize Residual};
   \end{scope}

  % grid:
  % \draw[step=1cm,gray] (0,0) grid (11,4);

\end{tikzpicture}
}% END RESIZE BOX
\end{center}
\caption{Left: a comparison of multigrid to classical iterative methods, Jacobi, GMRES and Conjugate Gradient (CG), for solving the two-dimensional Poisson's equation.
  Right: A comparison of multigrid convergence for second- and fourth-order accurate discretizations to
  Poisson's equation.
  %{\red Fix colours of curves and labels to match. Increase font size. Remove legends?}
 }
  \label{fig:introMG}
\end{figure}
} % MG

Multigrid algorithms are well known to provide excellent iterative solvers for elliptic problems.
For example, the left graph of Figure~\ref{fig:introMG} compares the convergence of multigrid to some classical iterative solvers for the solution of Poisson's equation on a square using a second-order accurate discretization. While classical iterative solvers such as GMRES or conjugate gradient (CG)
have convergence rates that degrade as the mesh is refined, convergence rates for multigrid are essentially independent of the mesh spacing $h$, leading to so-called optimal algorithms.
% orderC
In this current work we wish to develop high-order accurate elliptic solvers that are nearly as efficient as second-order accurate schemes; sample results to this effect are shown in the right graph of Figure~\ref{fig:introMG} where the convergence of a fourth-order accurate multigrid scheme is compared to the convergence of a second-order accurate discretization.
The key point shown here is that it is possible to obtain convergence rates for the the fourth-order accurate scheme that are very similar to those from the second-order accurate scheme. Although a fourth-order accurate scheme is somewhat more expensive due to a larger stencil, it will generally be much more efficient than a second-order scheme in obtaining a solution to a given accuracy.
Another goal of the current work is to choose smoothing parameters and lower-order coarse level operators to optimize the overall multigrid convergence rate; this is illustrated in the left graph of Figure~\ref{fig:intro_topics} where an optimized algorithm is compared to a more standard approach.

%% main topics ------------------------------------------------------------------------------------
Overset grids are typically dominated by large Cartesian background grids. A well designed multigrid algorithm for overset grids might be expected to have convergence rates similar to those on a single Cartesian grid, and thus fast algorithms for Cartesian grids are highly desired.
% \commentB{On the other hand, good multigrid algorithms on a good overset grid should have convergence properties similar to those on a single Cartesian grid.}
Furthermore, many important applications
such as those involving incompressible fluids, incompressible elasticity, or electromagnetics often require the solution to Poisson's equation. 
In this paper we therefore focus on the solution of Poisson's equation on Cartesian grids where we can make use of local Fourier analysis to study properties of our multigrid algorithms.
% orderC (mainly)
We investigate the use of nonstandard coarsening and nonstandard coarse level operators.
We show that second-order accurate coarse level approximations are very effective for fourth- and sixth-order accurate discretizations on the fine grid. 
We study the use of different Galerkin and non-Galerkin coarse level approximations. 
A Galerkin coarse-level operator is constructed by averaging the fine-level operator using the coarse-to-fine interpolation operator and then transferring it to the coarser level using the restriction operator; % L_H = I_h^H L_h I_H^h
a non-Galerkin operator refers to one with the same stencil as the operator on the fine level. 
We assess the influence of these coarse-level operators on the multigrid convergence rate as well as the computation cost by studying the \textit{effective convergence rate} that
takes computational work in account when comparing different schemes.
Based on these results we propose more flexible and simpler Galerkin-like operators for coarse grids.
We study the use of over-relaxed red-black smoothers with a relaxation parameter $\omega$. Using local Fourier analysis we choose the value $\omega$ and the coarse-level operators to optimize the overall multigrid convergence rate, rather than the more common approach of choosing $\omega$ to optimize the smoothing rate in isolation from the full multigrid cycle.

% rbC
Motivated by the use of red-black smoothers in one dimension that can result in
a direct solver for a second-order accurate discrete Poisson's equation, 
we show that this direct-solver property of a two-level cycle can be carried over to two dimensions, 
using the rotated grid with coarsening factor $\sqrt{2}$ that is formed from red-black coarsening (see Figure~\ref{fig:coarseningStrategies}, middle).
We evaluate the use of red-black coarsening in general multigrid cycles with more than two levels  and discuss possible extensions to three dimensions.
% rC
Motivated by the needs of flexible grid generation for overset grids,
as well as the fast convergence obtained with red-black coarsening,
we also study grid coarsening by factors other than two (see Figure~\ref{fig:coarseningStrategies}, right).
We ask the question of which coarsening factor leads to the most efficient algorithm. 
For example, do the good convergence rates of red-black coarsening in two dimensions when coarsening by a factor $\sqrt{2}$ to a rotated grid carry over in any way to non-rotated grids with $\sqrt{2}$ coarsening?
Figure~\ref{fig:intro_topics} (right) shows that good effective convergence rates can be obtained when coarsening
by a general factor $r$ over a fairly large range of $r$. 
Local Fourier analysis is used to study the properties of general factor-$r$ coarsening.

  % coarsening strategies fig 

{% -------------------
\newcommand{\labelsize}{\tiny\sffamily}
\newcommand{\labelcol}{violet}
\newcommand{\figWidth}{5cm}%
\newcommand{\trimfig}[2]{\trimw{#1}{#2}{.05}{.05}{.0}{.0}}
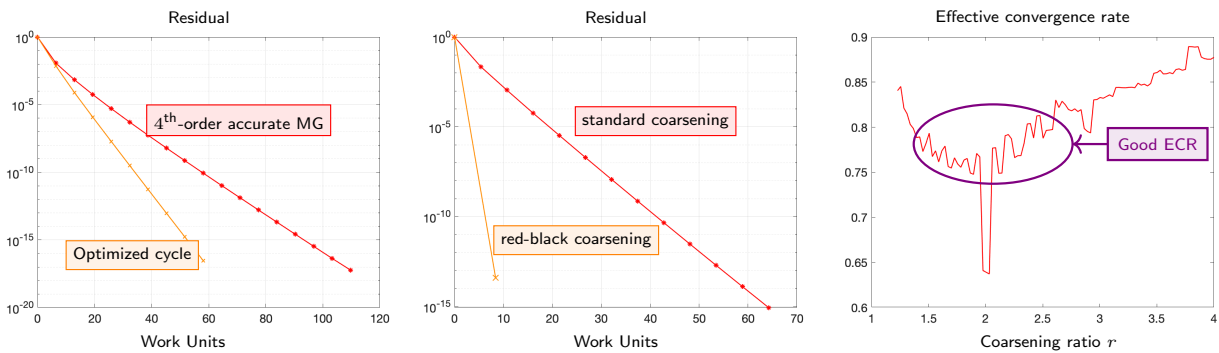
\begin{figure}[htb]
\begin{center}
\resizebox{16.5cm}{!}{% START resize box
\begin{tikzpicture}[scale=1]
  \useasboundingbox (0,0.4) rectangle (15.5,4);  % set the bounding box (so we have less surrounding white space)

  \begin{scope}[xshift=0cm] % orderC
    \draw(0,0) node[anchor=south west,xshift=-15pt,yshift=-8pt]{\trimfig{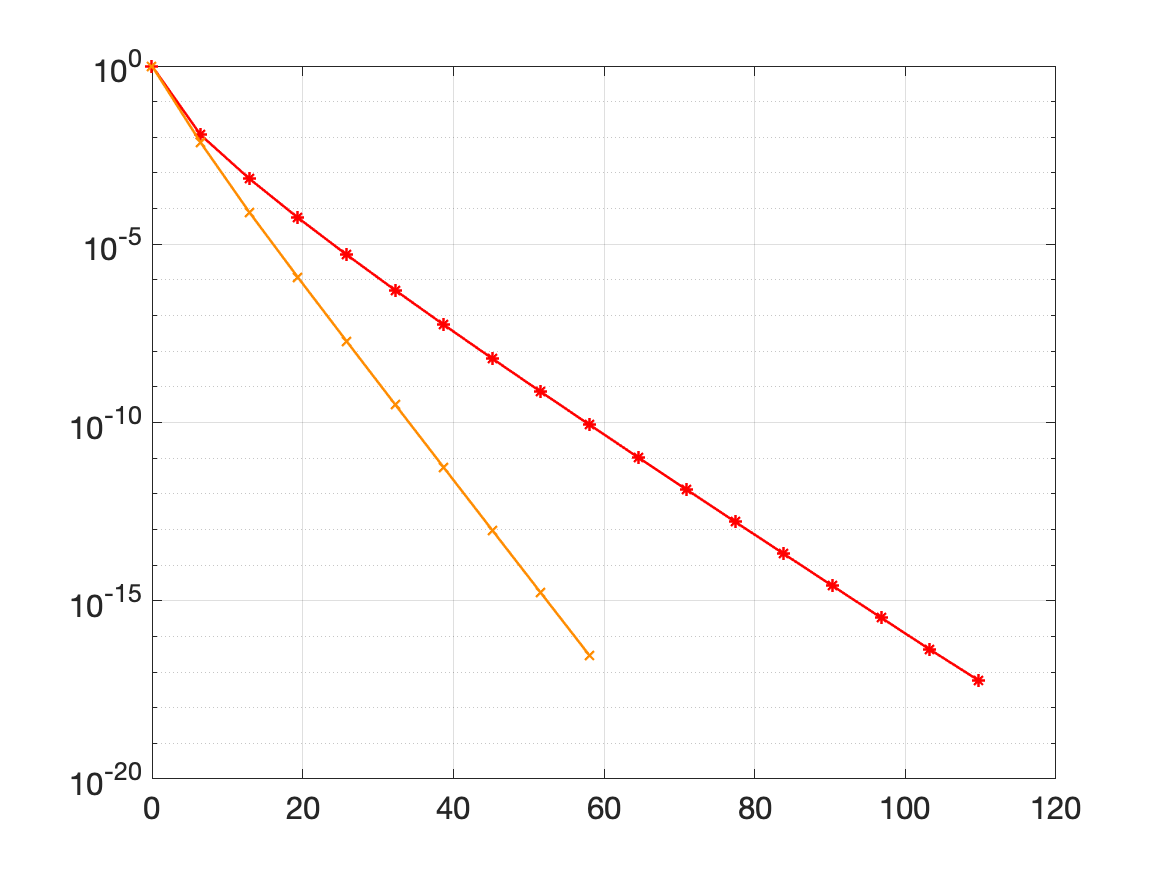}{\figWidth}};
    \draw (3,2.5) node[anchor=south,draw=red,fill=red!10,inner sep=2.5pt] {\labelsize $4^{\textsf{th}}$-order accurate MG};
    \draw (2.5,1.) node[anchor=east,draw=orange,fill=orange!10,inner sep=2.5pt] {\labelsize Optimized cycle};
    \draw (2.3,0.2) node[anchor=north,fill=white,yshift=-2pt] {\labelsize Work Units};
    \draw (2.5,3.9) node[anchor=south,fill=white,yshift=-3pt] {\labelsize Residual};
  \end{scope}

   \begin{scope}[xshift=5.25cm] % rbC
     \draw(0,0) node[anchor=south west,xshift=-15pt,yshift=-8pt]{\trimfig{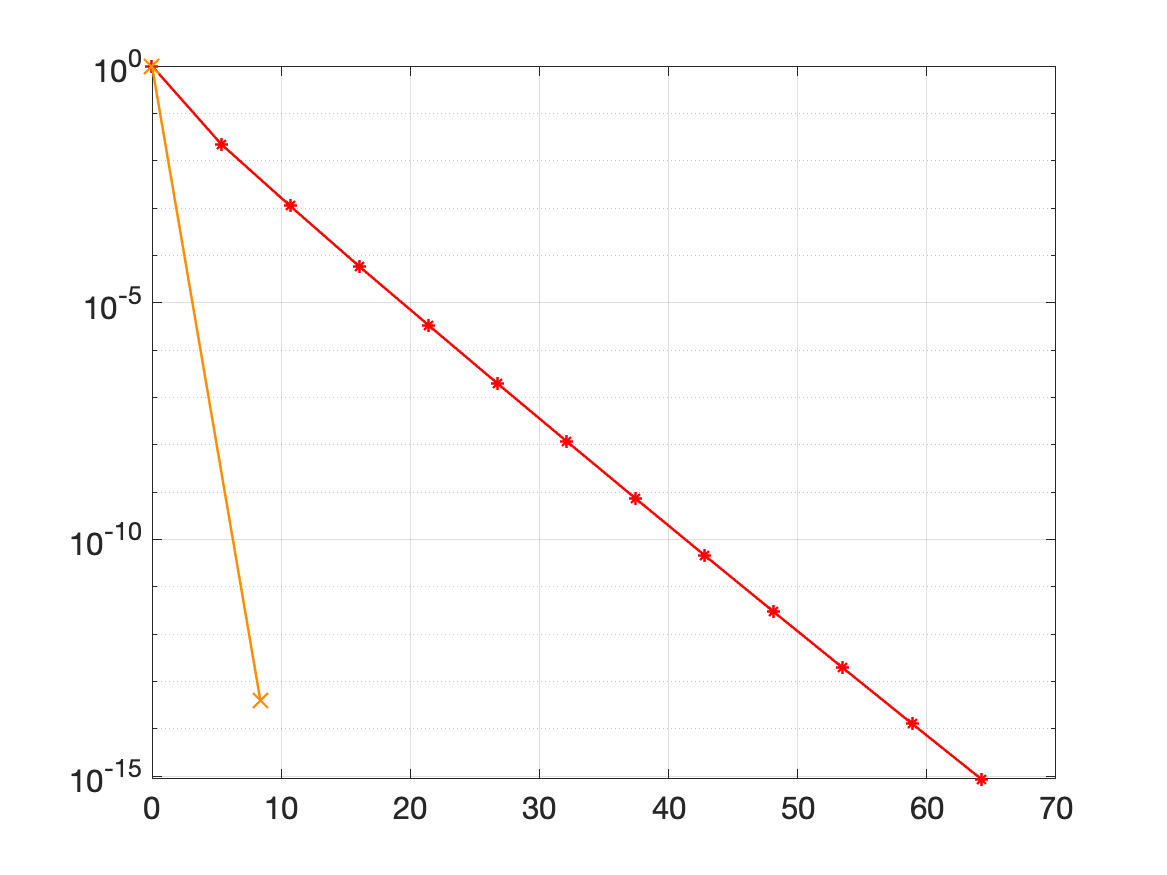}{\figWidth}};
     \draw (3,2.5) node[anchor=south,draw=red,fill=red!10,inner sep=2.5pt] {\labelsize standard coarsening};
     \draw (2,1) node[anchor=south,draw=orange,fill=orange!10,inner sep=2.5pt] {\labelsize red-black coarsening};
     \draw (2.5,0.2) node[anchor=north,fill=white,yshift=-2pt] {\labelsize Work Units};
     \draw (2.5,3.9) node[anchor=south,fill=white,yshift=-3pt] {\labelsize Residual};
   \end{scope}

   \begin{scope}[xshift=10.5cm] % rC
     \draw(0,0) node[anchor=south west,xshift=-15pt,yshift=-8pt]{\trimfig{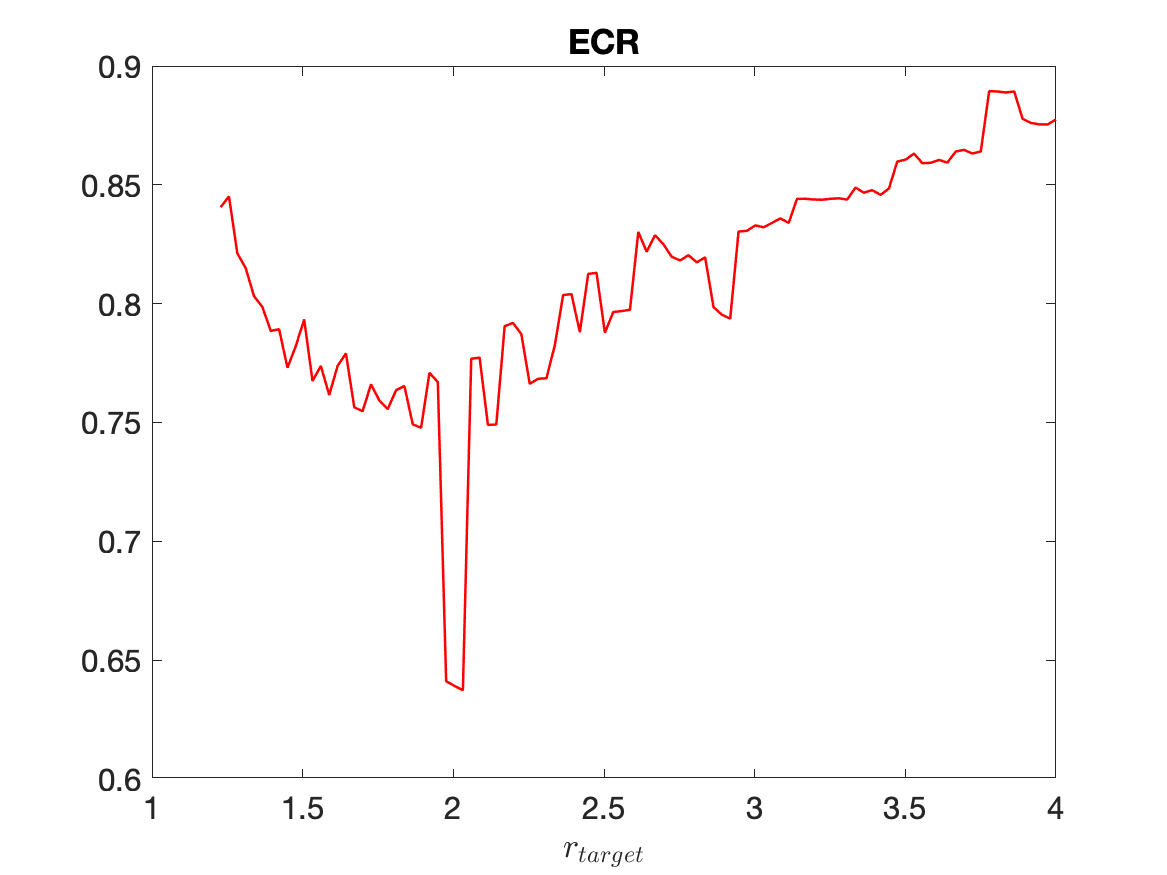}{\figWidth}};
     \draw (2.7,0.2) node[anchor=north,fill=white,yshift=-2pt] {\labelsize Coarsening ratio $r$};
     \draw (2.5,3.9) node[anchor=south,fill=white,yshift=-4pt] {\labelsize Effective convergence rate};
     \ellipseLabel{2}{2.4}{\labelsize Good ECR}{violet}{1}{.5};
   \end{scope}

  % grid:
 %  \draw[step=1cm,gray] (0,0) grid (15.5,4);

\end{tikzpicture}
}% END RESIZE BOX
\end{center}
\caption{Some main results (in two dimensions).
  Left: higher-order accurate multigrid convergence rates can be significantly improved 
  by optimizing the smoothing parameter and \textbf{lower-order} coarse-level operators.  % std-C
  % using Galerkin coarse grid operators and
  %        cycle-optimized relaxation parameters.
  Middle: \textbf{red-black coarsening} results in multigrid cycles with fast convergence for Poisson's equation to second-order accuracy. % 2-level 
  Right: good effective convergence rates can be obtained for a wide range of multigrid \textbf{coarsening factors $r$} near 2.
  %{\red maybve add grid to right fig ?}
 }
  \label{fig:intro_topics}
\end{figure}
}

%\bigskip
%% literature -------------------------------------------------------------------------------------
Multigrid is a well established field with a vast literature, see for example the textbooks
\cite{MGtext, Briggs2000, Trottenberg2001, Wesseling91, Brandt77} and the references therein. 
There has been much active research in algebraic multigrid methods~\cite{AMG, AMG_CW96, GGR97, AMG_KS98, PMJ01, Brezina06}. 
Red-black coarsening as proposed in this paper has also been mentioned 
for Poisson's equation (see~\cite{MGtext}, for example),
especially in the context of algebraic multigrid;
but this strategy differs from ours when generating multiple levels.
Multigrid for overset grids has been considered in~\cite{HinatsuFerziger91, Johnson95, Tu95, ZangStreet95}.
The use of over-relaxed red-black smoothers, in particular choosing the parameter $\omega$ to optimize the smoothing factor,
has been studied in~\cite{Yavneh96}.
Nonstandard coarsening for geometric multigrid including semi-coarsening (for anisotropy)
has been considered in for instance~\cite{semiC_NR93, Washio1998FlexibleMS}.
The development of multigrid algorithms for high-order accurate discretizations is often
accomplished by adding correction terms to a second-order accurate scheme~\cite{MGtext}.
For instance, a second-order accurate multigrid algorithm can be extended to higher order by adding a defect-correction to the right-hand-side on the finest level~\cite{Hackbusch82}; 
on the other hand, $\tau$-extrapolation~\cite{tau} avoids the direct use of a high-order
accurate operator altogether by adding a correction to the right-hand-side on the coarse level,
but this approach may be restricted to particular operators.
Our approach to high-order accuracy is to instead directly develop the multigrid algorithm, smoothers etc.,  for the high-order accurate discretization. Although not considered here, this approach works best when boundary conditions are discretized using compatibility so as to retain smooth behavior of the residual near boundaries~\cite{max2006b}.
The order of accuracy of the transfer operators has been studied in~\cite{Hemker90}, in terms of how it affects multigrid convergence.

\medskip
%% STRUCTURE: ---------------------------------------------------------------------------------------
The remainder of this paper is structured as follows.
In Section~\ref{Sec: defs}, we specify the problem and notations, give a brief overview of multigrid and introduce local Fourier analysis.
%% main -----
The main part of this paper consists of three topics, in Sections \ref{M1: orderC}, \ref{M2: rbC} and \ref{M3: rC}.
In Section~\ref{M1: orderC} we consider the use of lower-order accurate coarse-level operators for higher-order discretizations on the fine level. We argue through local Fourier analysis that multigrid cycles with lower-order accurate coarse-level operators can achieve comparable or better convergence rates than those with the same high order-of-accuracy on every level. The resulting multigrid algorithm will be less complex and more efficient.
In Section~\ref{M2: rbC} we develop and analyze multigrid with red-black coarsening in two dimensions, especially with the red-black smoother.
Analogous to the one-dimensional case, we construct a two-level cycle that is equivalent to a direct solver.
Section~\ref{Sec: rbC_multilevel} extends the algorithm to multiple levels.
In Section~\ref{Sec: rbC_3d}, we discuss the limitations of red-black coarsening and the possibilities of extending the algorithm to higher dimensions. 
In Section~\ref{M3: rC}, we give a simple multigrid algorithm that coarsens the grid by a general factor $r$ in each dimension, which is of practical importance. In the meantime we explore the optimal $r$ to aim at in terms of optimal convergence.
% conclusion
Finally, Section~\ref{Sec: conclusion} summarizes and concludes the paper.

%% theory
Most of the theoretical considerations are organized in~\ref{Sec: theory}. We record the local Fourier analysis for one and two dimensional multigrid, with standard and red-black coarsening (which are different in two dimensions).
In particular, \ref{Sec: rbC_LFA} corresponds to red-black coarsening in Section~\ref{M2: rbC},
with some generalization.
%{A three-level local Fourier analysis is given in \ref{Sec: rbC_3_LFA}, and the algorithm is  generalized to an anisotropic as well as a fourth-order accurate Laplacian in \ref{Sec: rbC_anisotropic} and \ref{Sec: rbC_4o}.}
In~\ref{Sec: cyclic_2d}, the algorithm in Section~\ref{Sec: rbC} is reformulated as red-black reduction; and in~\ref{Sec: connections} we discuss the connections of the proposed multigrid cycles with red-black coarsening with reduction methods and some of the ideas from algebraic multigrid.
\ref{Sec: orderC_LFA} corresponds to Section~\ref{M1: orderC}, and gives the local Fourier analysis for multigrid for the fourth-order accurate discretization, with discussion on the order of accuracy of the coarse-level correction.
In \ref{sec: WU} we give an estimate of the work-units of a general multi-level cycle.

\bigskip
%\clearpage
% ---------------- PROBLEM SPECIFICATION AND INTRODUCTION TO MULTIGRID ------------
\section{Problem specification and overview of multigrid}\label{Sec: defs}
In this section we specify the general and model problems to be considered in this paper along with some notations. We also give a brief introduction to the multigrid algorithm and analysis. 

\subsection{Problem specification and notation}
In general, we are interested in solving an elliptic boundary value problem in some domain $\Omega \subset \mathbb{R}^d$ (where $d = 1,\ 2,\ 3$ is the number of space dimensions):
\begin{align}
\begin{cases}
L u = f, \quad &\mathbf{x} \in \Omega\\
B u = g, &\mathbf{x} \in \partial \Omega,
\end{cases}
\end{align}
where $L$ is an elliptic operator and $B$ is some boundary operator. Consider the discretized problem 
\begin{align}
\begin{cases}
L_h u_h = f_h, \quad &\mathbf{x} \in \Omega_h\\
B_h u_h = g_h, &\mathbf{x} \in \partial \Omega_h,
\end{cases}
\end{align}
where $\Omega_h$ is the discretized domain, $L_h$ and $B_h$ are discrete approximations of the continuous operators, and $u_h$ is a grid function. 

It is known that the main properties of a well designed multigrid algorithm should not depend strongly on the particular domain or boundary conditions. Therefore, for the purposes of this paper, we consider $\Omega_h$ to be the infinite Cartesian grid in $d$ dimensions with a grid-spacing $h$:
\ba \label{G_h_d}
\grid G_h = \Big\{  \xv =\jv\, h:~ \jv \in \Integer^d \Big\}.
\ea
Let us denote the space of all the grid functions on $\grid G_h$ as 
$\gridfunspace G_h$.
%$\gridfunspace G_h = \{u_h:\ \grid G_h \to \mathbb{C}\}$. 
Multigrid algorithms involve multiple grids with different grid-spacings. We will need to define operators that map grid functions on some grid to those on another. Denote a linear operator from $\gridfunspace G_h$ to $\gridfunspace G_{h'}$ by $A_h^{h'}$.
For the case where one grid is embedded in another, we adopt stencil notation. If $\grid G_{h'} \subset \grid G_h$,
\ba
A_h^{h'} = [a_{\jv}]_h^{h'}:\ \gridfunspace G_h \to \gridfunspace G_{h'}
\ea
represents
\ba \label{bra}
A_h^{h'} u_h (\xv) = \sum_{\jv} a_{\jv}\ u_h(\xv + \jv h),
\quad \xv \in \grid G_{h'},
\ea
for any grid function $u_h \in \gridfunspace G_h$, where $a_{\jv}$ are the stencil coefficients.
If $\grid G_h$ and $\grid G_{h'}$ are the same, denote $A_h \equiv A_h^h$.
(For example, see equation (\ref{L_h_2}) of the 5-point stencil for the two-dimensional Laplacian.)
In particular, denote the identity operator as $I_h$.
On the other hand, if $\grid G_h \subset \grid G_{h'}$, the stencil notation with reversed brackets is used:
\ba
A_h^{h'} = ]a_{\jv}[_h^{h'}:\ \gridfunspace G_h \to \gridfunspace G_{h'},
\ea
representing
\ba 
A_h^{h'} u_h (\xv) = \sum_{\jv} \sum_{\xv - \jv h' \in \grid G_h} a_{\jv}\ u_h(\xv - \jv h'),
\quad \xv \in \grid G_{h'},
\ea
where the $a_{\jv}$ represent the weights with which the coarse grid values are distributed to the fine grid. (For example, see the linear interpolation operator $I_H^h$ in equation (\ref{transfer_1d}).)\\

Consider the inner product defined on $\gridfunspace G_h$ as 
%(infinite sum? l2)
\ba
<u_h,\ w_h>_h\ =\ h^d \sum_{\xv \in \grid G_h} u_h^*(\xv) w_h(\xv),\quad
u_h,\ w_h \in \gridfunspace G_h,
\ea
and the induced norm $\|\cdot\|_h$. The adjoint of an operator $A_h^{h'}$, denoted by $(A_h^{h'})^*$ or $(A^*)_{h'}^h$, is defined so that
\ba
<(A^*)_{h'}^h\ u_{h'},\ w_h>_h\ =\ <u_{h'},\ A_h^{h'} w_h>_{h'},\quad
\forall\ u_{h'} \in \gridfunspace G_{h'},\ w_h \in \gridfunspace G_h.
\ea

Finally, denote $\rho(A_h)$ or $\rho(\mathbf{A})$ as the spectral radius of an operator $A_h$ or a square matrix $\mathbf{A}$.

% -----------------------------------------------------------------------------------------------------------------------------------------

\subsection{Overview of multigrid and local Fourier analysis}\label{Sec: MG}
  Here we give an outline of a multi-level algorithm for a general discretized
  boundary value problem represented as 
  $L_h u_h = f_h$.
A geometric multigrid algorithm to solve this problem is often based on four key components:
\begin{enumerate}
  \item A sequence of grids (levels) of increasing coarseness. 
  \item an iterative procedure called a \textit{smoother} that is effective at reducing high-frequency components of the residual.
    For example, common smoothers are based on the Jacobi or Gauss-Seidel iteration.
  \item Fine-to-coarse (restriction) and coarse-to-fine (interpolation)
       operators that transfer a solution between a fine grid and a coarse grid.
  \item Coarse-level operators that approximate the fine-level operator $L_h$. 
\end{enumerate}    

In particular, consider an ($\lmax +1$)-level cycle with grids of grid-spacing $h_\ll\ (h_0 = h),\ \ll =
0,\ 1,\ \cdots,\ \lmax$ , on each level. The smoothing operator on the $\ll$-th level is denoted as
$S_{h_\ll}$, and restriction and interpolation operators between the $\ll$-th and $(\ll+1)$-st levels are
denoted as $I^{h_{\ll+1}}_{h_\ll}$ and $I_{h_{\ll+1}}^{h_\ll}$, respectively. The parameters $\nu_1$ and
$\nu_2$ indicate the number of pre- and post-smoothing sweeps per cycle. The multigrid cycle is
called a V cycle if $\gamma=1$, and a W cycle if $\gamma=2$.

% --------------------- ALGORITHM -----------------
\bigskip

\begin{comment}
\begin{algorithm}
\algFontSize % Text font size defined in r2p.tex. Line number font size also defined there.

\caption{$(\lmax+1)$-level algorithm for $L_h u_h = f_h$ }

\begin{algorithmic}[1]
\Function{$u_{h_\ll}  = \textcolor{violet}{MG}$}{$f_{h_\ll};\ \text{maxIters}$}

  \For {$n = 0,\ 1,\ \cdots,\ \text{maxIters}-1$}
      \While {$\| d_{h_\ll} \| / \| f_{h_\ll} \| < tol$}

         \State $\bar{u}_{h_\ll} ^{(n)}  \xleftarrow{S_{h_\ll} ^{\nu_1}}\ u_{h_\ll} ^{(n)} $ 
         \Comment Pre-smoothing
  
         \State $\bar{d}_{h_\ll} ^{(n)}  = f_{h_\ll}  -  L_{h_\ll}  \bar{u}_{h_\ll} ^{(n)}$
  
         \State $\bar{d}_{h_{\ll+1}}^{(n)} = I^{h_{\ll+1}}_{h_\ll} \bar{d}_{h_\ll}^{(n)} $ \Comment Restriction
  
         \If {$\ll+1 = \lmax$}
  	    \State $\bar{v}_{h_l}^{(n)} = L^{-1}_{h_l} \bar{d}_{h_l}^{(n)} $  \Comment Exact solve on the coarsest grid
   	\Else
  	    \State $\bar{v}_{h_{\ll+1}}^{(n)} = \textcolor{violet}{MG}(\bar{d}_{h_{\ll+1}}^{(n)};\ \gamma)$ \Comment Coarse-grid solve
  	\EndIf 
   
        \State $\tilde{v}_{h_\ll}^{(n)}  = I_{h_{\ll+1}}^{h_\ll} \bar{v}_{h_{\ll+1}}^{(n)}$  \Comment Interpolation
     
        \State $\bar{u}_{h_\ll}^{(n+1)} = \bar{u}_{h_\ll}^{(n)} +\tilde{v}_{h_\ll}^{(n)} $
  
       \State $u_{h_\ll}^{(n+1)} = S_{h_\ll}^{\nu_2}(\omega_2)\ \bar{u}_{h_\ll}^{(n+1)}$ \Comment Post-smoothing
  
       %\State $n \text{++}$ 
     \EndWhile  
   \EndFor
\EndFunction
 \end{algorithmic} 
 \label{alg:mg}
% Comments.
\end{algorithm}
\end{comment}

% ***TEST**** 
% \renewcommand{\algFontSize}{\normalsize}% *******

\noindent\begin{minipage}{.6\textwidth}
% ---- algorithm on left -----
% \captionof{algorithm}{\large $(l+1)$-level algorithm for $L_h u_h = f_h$}\label{myalg}
\begin{algorithmic}[1]
  \algFontSize % Text font size defined in adegdm.tex. Line number font size also defined there.
  \Function{$u_{h_\ll}  = \textcolor{violet}{MG}$}{$f_{h_\ll};\ \text{maxIters}$}

  \For {$n = 0,\ 1,\ \cdots,\ \text{maxIters}-1$}
      \While {$\| d_{h_\ll} \|_{h_\ll} / \| f_{h_\ll} \|_{h_\ll} < tol$}

         \State $\bar{u}_{h_\ll} ^{(n)} \xleftarrow{S_{h_\ll} ^{\nu_1}}\ u_{h_\ll} ^{(n)} $ \Comment Pre-smoothing
  
         \State $\bar{d}_{h_\ll} ^{(n)}  = f_{h_\ll}  -  L_{h_\ll}  \bar{u}_{h_\ll} ^{(n)}$
  
         \State $\bar{d}_{h_{\ll+1}}^{(n)} = I^{h_{\ll+1}}_{h_\ll} \bar{d}_{h_\ll}^{(n)} $ \Comment Restriction
  
         \If {$\ll+1 = \lmax$}
  	    \State $\bar{v}_{h_{\ll+1}}^{(n)} = L^{-1}_{h_{\ll+1}} \bar{d}_{h_{\ll+1}}^{(n)} $  \Comment Exact solve, coarsest level
   	\Else
  	    \State $\bar{v}_{h_{\ll+1}}^{(n)} = \textcolor{violet}{MG}(\bar{d}_{h_{\ll+1}}^{(n)};\ \gamma)$ \Comment Coarse-level solve
  	\EndIf 
   
        \State $\tilde{v}_{h_\ll}^{(n)}  = I_{h_{\ll+1}}^{h_\ll} \bar{v}_{h_{\ll+1}}^{(n)}$  \Comment Interpolation
     
        \State $\bar{u}_{h_\ll}^{(n+1)} = \bar{u}_{h_\ll}^{(n)} +\tilde{v}_{h_\ll}^{(n)} $
  
       \State $u_{h_\ll}^{(n+1)} \xleftarrow{ S_{h_\ll}^{\nu_2}}\ \bar{u}_{h_\ll}^{(n+1)}$ \Comment Post-smoothing
  
       %\State $n \text{++}$ 
     \EndWhile  
   \EndFor
\EndFunction
\end{algorithmic}
%
%  ----Figure on right ----
\end{minipage}%
\begin{minipage}{.4\textwidth}

\begin{center}
  \begin{tikzpicture}
  \useasboundingbox (0,1) rectangle (6,7.5);  % set the bounding box (so we have less surrounding white space)

   \begin{scope}[xshift=1cm,yshift=6cm]  

\begin{scope}[very thick,decoration={
    markings,
    mark=at position 0.5 with {\arrow{stealth}}}
    ] 
    \draw[postaction={decorate}] (0,0) -- (1,-2); 
    \draw[postaction={decorate}] (1,-2) -- (2,-4);
    \draw[postaction={decorate}, dotted] (2,-4) -- (2.5,-5); 
    \draw[postaction={decorate}, dotted] (2.5,-5) -- (3,-4); 
    \draw[postaction={decorate}] (3,-4) -- (4,-2);    
    \draw[postaction={decorate}] (4,-2) -- (5,0); 
\end{scope}

     % \draw[black,very thick,postaction={decorate}]] (0,0) -- (2,-4) -- (4,0);
%     \draw[->,black,very thick] (0,0) -- (1,-2); 
%     \draw[->,black,very thick] (1,-2) -- (0,-4); 
%     \draw[->,black,very thick] (0,-4) -- (1,-2); 
%     \draw[->,black,very thick] (1,-2) -- (2,0); 
    \draw[black] (0,0) \vdot node[anchor=east] {\small $S_{h_0}^{\nu_1}$};
    \draw[black] (5,0) \vdot node[anchor=west] {\small $S_{h_0}^{\nu_2}$};
    \draw[black] (.5,-1)  node[anchor=east] {\small $I_{h_0}^{h_1}$};
    \draw[black] (.5,-1)  node[anchor=west] {\small $\bar{d}_{h_1}^{(n)}$};

    \draw[black] (1,-2) \vdot  node[anchor=east] {\small $S_{h_1}^{\nu_1}$};
    \draw[black] (1.5,-3)  node[anchor=east] {\small $I_{h_1}^{h_2}$};
    \draw[black] (4.5,-1)  node[anchor=west] {\small $I_{h_1}^{h_0}$};
    \draw[black] (4.5,-1)  node[anchor=east] {\small $\bar{v}_{h_1}^{(n)}$};

    \draw[black] (4,-2) \vdot node[anchor=west] {\small $S_{h_0}^{\nu_2}$}; 
    \draw[black] (3.5,-3)  node[anchor=west] {\small $I_{h_2}^{h_1}$};
    %  coarse grid solve:
    \draw[black] (2.5,-5) \vdot node[anchor=west] {\small $\bar{v}_{h_\lmax}^{(n)} \!=\! L^{-1}_{h_\lmax} \bar{d}_{h_\lmax}^{(n)} $};
  \end{scope}  

%-  \draw[->,black,thick] (0,0) -- (7,0) node[right] {$x$};
%-  % \draw[->] (0,-3) -- (0,4.2) node[above] {$y$};
%-  % \draw[scale=0.5,domain=-3:3,smooth,variable=\x,blue] plot ({\x},{\x*\x});
%-  % \draw[scale=0.5,domain=-3:3,smooth,variable=\y,red]  plot ({\y*\y},{\y});
%-
%-  \begin{scope}[yshift=0cm]
%-    % plot sin(2*x)
%-    \draw[scale=1,domain=0:2*pi,smooth,variable=\x,blue,very thick]  plot ({\x},{sin(2*(180/pi)*\x)});
%-    \draw[scale=1,domain=0:2*pi,smooth,variable=\x,blue]  plot ({\x},{-sin(2*(180/pi)*\x)});
%-    % plot (-1)^j * sin(2*x_j)
%-    % (-1)^j = cos( (N/2)*x_j ) 
%-    % N=10, samples=N+1
%-    \draw[scale=1,domain=0:2*pi,samples=11,variable=\x,orange,very thick]  plot ({\x},{sin(2*(180/pi)*\x)*cos(5*(180/pi)*\x)});
%-  \end{scope}
%  draw grid temporarily while making figure 
%  \draw[step=1cm,gray] (0,0.5) grid (6,7.5);
% 
\end{tikzpicture}

\end{center}


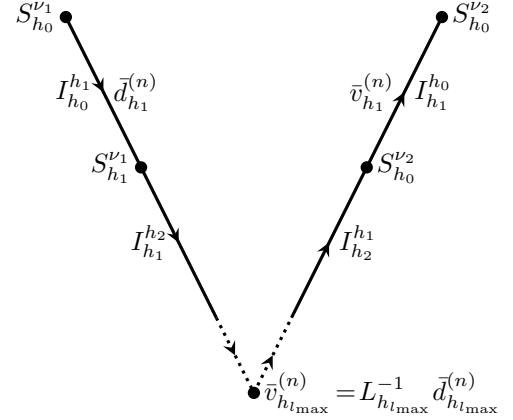
\captionof{figure}{Multigrid V Cycle.} \label{fig:mgCycle}
\end{minipage}

%- % ***TEST**** 
%- \noindent\begin{minipage}{.5\textwidth}
%- % ---- algorithm on left -----
%- \captionof{algorithm}{Euclid’s algorithm}\label{myalg}
%- \begin{algorithmic}[1]
%-   \Procedure{Euclid}{$a,b$}\Comment{The g.c.d. of a and b}
%-     \State $r\gets a\bmod b$
%-     \While{$r\not=0$}\Comment{We have the answer if r is 0}
%-       \State $a\gets b$
%-       \State $b\gets r$
%-       \State $r\gets a\bmod b$
%-     \EndWhile\label{euclidendwhile}
%-     \State \textbf{return} $b$\Comment{The gcd is b}
%-   \EndProcedure
%- \end{algorithmic}
%- %
%- %  ----Figure on right ----
%- \end{minipage}%
%- \begin{minipage}[t]{.5\textwidth}
%-   \centering
%-   \rule{0.3\textwidth}{50pt}
%-   \captionof{figure}{This is a figure caption.} \label{myfig1}
%-   \bigskip
%-   \rule{0.25\textwidth}{70pt}
%-   \captionof{figure}{This is another figure.} \label{myfig2}
%- \end{minipage}

\bigskip

Given any current solution $u_h^{(n)}$ (indexed by $n$), as an approximation of the solution $u_h$ to $L_h u_h = f_h$, we are going to consider the action of each component in a multigrid cycle on the error 
\begin{equation}
v_h^{(n)} \eqdef u_h - u_h^{(n)}
\end{equation}
and the residual
\begin{equation}
d_h^{(n)} \eqdef f_h - L_h u_h^{(n)}.
\end{equation}
Most importantly, the iteration operator of a multigrid cycle, denoted by $M_h$, relates the error between successive iterations as 
\ba
   v_h^{(n+1)} = M_h v_h^{(n)}.
\ea 
The convergence of multigrid can be analyzed by studying the operator $M_h$,
in particular its spectral radius and norms. The asymptotic convergence rate of the multigrid cycle is given by
% \commentB{
% \begin{equation}
% \lim_{n \to \infty}\frac{\|v_h^{(n+1)}\|_h}{\|v_h^{(n)}\|_h} =  \rho(M_h) \equiv \rho, 
% \end{equation}
% }
\begin{equation}
     \lim_{n \to \infty}\frac{\|v_h^{(n+1)}\|_h}{\|v_h^{(n)}\|_h} = \rho,  \label{eq:asymptoticConvergenceRate}
\end{equation}
where
\ba \label{rho}
   \rho \eqdef  \rho(M_h)
\ea
is the spectral radius of the multigrid iteration operator.
In practice, the convergence rate can be measured from the ratio of the norms of successive residuals since
\begin{equation}\label{eq: resReduction}
\lim_{n \to \infty}\frac{\|d_h^{(n+1)}\|_h}{\|d_h^{(n)}\|_h} = \rho(L_h M_h L_h^{-1}) = \rho( M_h ) = \rho.
\end{equation}

% -----------------------LFA -----------------------------------------
Local Fourier analysis (LFA) is a standard tool of analyzing multigrid algorithms, in particular their convergence properties \cite{MGtext}. When the boundary conditions are dealt with properly numerically, the most important properties of a multigrid cycle are ``local", in that its convergence is not affected much by the boundaries, as if on an infinite domain. This is why the infinite grid (\ref{G_h_d}) is introduced, and the discrete Fourier modes 
\begin{equation}\label{phi_d}
\phi_h (\mathbf{x},\ \thv) \eqdef e^{i \boldsymbol{\theta} \cdot \frac{\xv}{h}},\quad \xv \in \grid G_h,
\end{equation}
on the grid $\grid G_h$  are thus eigenfunctions for any linear constant-coefficient operator on $\gridfunspace G_h$, where $\boldsymbol{\theta} \in \Theta \eqdef [-\pi,\ \pi)^d$ is a parameter representing the frequency (see Figure \ref{fig:lfaModes} for an illustration).
We can analyze the behavior of the operators involved in a multigrid cycle (such as $S_h$) by analyzing their operation on the Fourier modes (or eigen-subspaces), based on the orthogonal decomposition of any grid function as
\begin{equation}
u_h(\xv) =\int_{\Theta} \hat{u}(\thv)\ \phi_h (\mathbf{x},\ \thv)\ d\thv,
\end{equation} 
where $\hat{u}(\thv)$ is the Fourier coefficient. For instance, by the Parseval's identity
$
\|u_h\|_h = \|\hat{u}\|_2,
$
the reduction of error in physical space can be analyzed by examining the reduction of each Fourier mode.
\cut{
% ----- Figure showing discrete Fourier modes ----
{% start scope of figure 
\newcommand*{\N}{10}% 
\newcommand*{\h}{2*pi/\N}
\newcommand{\mode}[2]{%
   \draw[scale=1,domain=-1.5*\h:2*pi+\h,smooth,variable=\x,#1,thick,densely dotted,samples=100]  plot ({\x},{cos(#2*(180/pi)*\x)});
        \foreach \i in {-1,...,\N} {
            \draw[#1, fill=#1] (\h*\i, {cos(#2*(180/pi)*\h*\i)}) circle (1.5pt);
        }
  }
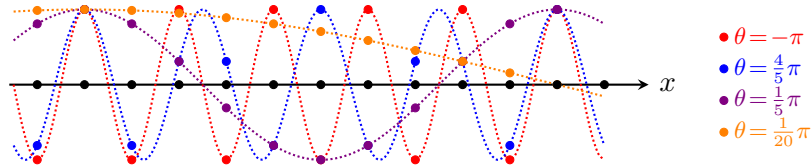
\begin{figure}[h!]
  \begin{center}
    \begin{tikzpicture}
      \useasboundingbox (-1,-.5) rectangle (10,1.2);  % set the bounding box (so we have less surrounding white space)

      \draw[-stealth,black,thick] (-1,0) -- (7.5,0) node[right] {$x$};

      % \draw[->] (0,-3) -- (0,4.2) node[above] {$y$};
      % \draw[scale=0.5,domain=-3:3,smooth,variable=\x,blue] plot ({\x},{\x*\x});
      % \draw[scale=0.5,domain=-3:3,smooth,variable=\y,red]  plot ({\y*\y},{\y});
    
      \begin{scope}[yshift=0cm]

        % plus minus mode plot sin(2*x)
        \mode{red}{5}; 
        \mode{blue}{4}; 
        \mode{violet}{1}; 
        \mode{orange}{.25}; 

        % grid points 
        \foreach \i in {-1,0,...,11} {
        	\draw[fill] (\h*\i, 0) circle (1.5pt);  % grid points 
        }

        \begin{scope}[xshift=8.5cm,yshift=-18pt]
           \draw[fill,red   ,yshift= 36pt] (0, 0) circle (1.5pt) node[anchor=west] {\footnotesize$\theta\!=\! -\pi$};
           \draw[fill,blue  ,yshift= 24pt] (0, 0) circle (1.5pt) node[anchor=west] {\footnotesize$\theta\!=\! \frac{4}{5}\pi$};
           \draw[fill,violet,yshift= 12pt] (0, 0) circle (1.5pt) node[anchor=west] {\footnotesize$\theta\!=\! \frac{1}{5}\pi$};
           \draw[fill,orange,yshift=  0pt] (0, 0) circle (1.5pt) node[anchor=west] {\footnotesize$\theta\!=\! \frac{1}{20}\pi$};
        \end{scope}
      \end{scope}
    %  draw grid temporarily while making figure 
   % \draw[step=1cm,gray] (-1,-1) grid (10,1);
    % 
    \end{tikzpicture}
  \end{center}
  \caption{There are a continuum of discrete Fourier modes $\phi_h(x,\ \theta)=e^{i\th x/h}$ represented on the infinite grid (1D).
    The highest frequency represented on the grid is the ``plus-minus'' mode (red) for $\th=-\pi$. Higher frequencies outside
    the primary interval $\th \in \Theta = [-\pi,\ \pi)$ are aliased onto modes for $\th \in \Theta$. 
     }
\label{fig:lfaModes}
\end{figure}
}% end scope of figure 
}

Since a multi-level cycle aims to have a convergence rate comparable to that of a two-level cycle, most of the analysis focus on two-level cycles. In particular, given the operator $L_h$ on the fine grid $\grid G_h$, denote the coarse grid as $\grid G_H$, and the coarse-level operator as $L_H$. 
For example, by \textsl{standard coarsening} (SC),
\ba \label{G_H_d}
\grid G_H = \Big\{  \xv =\jv\, H:~ \jv \in \Integer^d \Big\}
\ea
with a grid-spacing $H=2h$.
We say that the coarse-level operator $L_H$ has the \textsl{Galerkin} property if
\begin{equation}\label{L_H_G}
L_H =  I_h^H L_h I^h_H.
\end{equation}
Suppose the coarse-level solve is exact, with the restriction operator $I_h^H$ and interpolation operator $I_H^h$, we have the coarse-level correction operator
\begin{equation}\label{K}
  K_h^H = I_h - I^h_H L^{-1}_H I_h^H L_h.
\end{equation}
Then the iteration operator for a $V[\nu_1,\ \nu_2]$ cycle (with $\nu_1$ pre-smoothing and $\nu_2$ post-smoothing steps) is
\begin{equation}\label{M_2level}
M_h \equiv M_h^H = S_h^{\nu_2} K_h^H  S_h^{\nu_1}.
\end{equation}

According to whether the Fourier modes $\phi_h$ on $\grid G_h$ can be represented on the coarse grid  $\grid G_H$, they are categorized into low-frequency 
%modes if $\th \in \Theta^{low} \eqdef [-\frac{\pi}{2},\ \frac{\pi}{2})^d$, 
and high-frequency modes. 
%if $\th \in \Theta^{high} \eqdef \Theta \setminus \Theta^{low}$. 
In a multigrid cycle, the coarse-level correction aims to tackle the low-frequency components of the error, whereas the smoother aims to reduce the high-frequency ones.
To characterize the smoother's ability to reduce high-frequency components of the error, the smoothing factor is defined as
\begin{equation}\label{mu}
\mu \eqdef \rho^{\frac{1}{\nu}} \big( S_h^{\nu_2} Q_h^H S_h^{\nu_1}\big),
\end{equation}
where $\nu = \nu_1 + \nu_2,$
and $Q_h^H$ is the (theoretical) idealized coarse-level correction which is 0 on the low-frequency modes and identity on the high-frequency modes. If the coarse-level correction $K_h^H$ were to be replaced by the ideal one $Q_h^H$, the convergence rate $\rho$ would become $\mu^{\nu}$. Thus $\mu^{\nu}$ is used as a reference convergence rate for a multigrid cycle. (Note that $\rho$ can often be better than $\mu^{\nu}$.)
The smoothing parameter $\omega$ is often chosen in textbooks to optimize the smoothing factor $\mu$.
Instead, the value of $\omega$ can be chosen to optimize the actual multigrid cycle ($\rho$) which in general gives better results as can be seen from the various convergence graphs presented in our results.

In a general multi-level cycle with grid-spacings $(h_0,\ h_1,\ \cdots,\ h_{\lmax})$ 
($h_0=h$, $\lmax>1$), the coarse-level correction  operator on the finest level is given by
\begin{equation}\label{K_l}
    K_h = I_h - I^h_{h_1} \tilde{L}^{-1}_{h_1} I_h^{h_1} L_h,
\end{equation}
where the coarse level solve $\tilde{L}^{-1}_{h_1}$ is in turn obtained recursively by an $\lmax$-level cycle.
The multi-level iteration operator is then given by
\begin{equation}\label{M}
M_h = S_h^{\nu_2}\ K_h\  S_h^{\nu_1}.
\end{equation}

% ------------ computation -------------------------------
For computations we restrict ourselves to finite domains.
For the computed results presented in this paper, the domain is taken to be $[0,\ 1]^d$ with periodic boundaries, discretized to a fine grid  
\begin{equation}\label{G_h_N}
\grid G_h = \Big\{  \xv_{\jv} =\jv\, h:~ \jv \in \{0,\ 1,\ \cdots,\ N \}^d \Big\},
\end{equation}
with $h=\frac{1}{N}$, $N \in \Integer$. 
(In most cases $N$ is assumed to be an even number for the red-black smoother.)
%so that there will be an equal number of red and black points for the smoother.) 
The computed convergence rate (CR) is given by
\ba \label{CR}
\text{CR} = \frac{\|d_h^{(n+1)}\|_h}{\|d_h^{(n)}\|_h},
\ea
where the iteration number $n$ is assumed to be sufficiently large so that the reduction of error and residual is close to the asymptotic convergence rate $\rho$.
In order to compare the relative cost of different algorithms we define a normalized convergence rate called the
\textsl{effective convergence rate} (ECR) which takes into account the computational cost of an algorithm.
Define one work unit to be the number of floating-point operations (FLOPS)\footnote{We acknowledge that FLOPS do not tell the whole story on modern computer architectures and so the results presented here are more of a rough guideline.} needed for a single Jacobi smoothing step 
for the difference equation on the finest level.
The total computational cost in each cycle is measured in work-units (WU).
Define the effective convergence rate of a multigrid cycle to be 
\begin{equation}\label{ECR}
\text{ECR} \eqdef \text{CR}^{\frac{1}{\text{WU}}}.
\end{equation}
In particular, note that for the Jacobi iteration, its ECR is equal to its CR.
The ECR can be interpreted as the convergence rate
that a Jacobi (or Gauss-Seidel) method would need to achieve, per iteration,
in order to match the convergence of the multigrid cycle. 
A good multigrid algorithm may have an ECR in the range 0.5 to 0.8. In comparison,
the ECR for Jacobi is of order $1-\Oc(h^2)$ and thus approaches one as $h\rightarrow 0$.
In this paper, we often show both the CR, obtained from numerical computations on a finite grid, and the %spectral radius of the multigrid iteration operator, 
asymptotic convergence rate,
$\rho = \rho(M_h)$, obtained from local Fourier analysis; the agreement between these two, in accordance with \eqref{eq: resReduction}, provides some evidence for the correctness of the results. We also show how the CR and ECR depend on the smoothing parameter $\omega$ in the results. It is of particular interest to see how far away the optimal $\omega$ (to minimize the convergence rate) is from the default value $\omega=1$.  

Note that even though for a finite domain we do not have a continuous spectrum of frequency $\Theta$ as assumed in local Fourier analysis, the CR should be close to $\rho$ regardless of the size of the problem $N$ (as in \eqref{G_h_N}).
On the other hand, the WU per cycle, which consists of the work for smoothing, restriction and interpolation, as well as the coarse-level solve, may depend on $N$.
For a two-level cycle, the exact solve on the coarse grid can be the dominant cost of the the algorithm, for the size of the coarse problem can still be quite large. Hence the WU and so the ECR could potentially be dependent on $N$. 
(We assume on the coarsest grid an exact or iterative solver for a banded system is used.)
%In one dimension, however, the bandwidth is fixed from the stencil, which results in the WU of the coarse-level solve independent of $N$.)
In practice, the number of levels chosen for a multigrid cycle is determined from the size of problem $N$, in which case the WU of the coarsest-level solve would no longer be dominant and the dependence of WU (and hence ECR) per cycle from $N$ would be negligible. 
An estimate of the WU for a general multi-level cycle is given in \ref{sec: WU}.

% -----------------------------------------------------------------------------------------------------------------------------------------
\subsection{The model problem of Poisson's equation}
In this paper we focus on the important model problem of Poisson's equation in two dimensions (2D)
\begin{equation}\label{MP}
L u = f,\quad L = -\Delta = -(\partial_x^2 + \partial_y^2)
\end{equation}
and its discretization 
\begin{equation}\label{MP_h}
L_h u_h = f_h
\end{equation}
on the infinite grid (\ref{G_h_d}) with $d=2$.
%\ba 
%\grid G_h = \Big\{  \xv =\jv\, h:~ \jv=[j_1,\ j_2]^T \in \Integer^2 \Big\}.
%\ea
In particular, we consider the standard finite difference discretization of the negative Laplacian to second-order accuracy 
\begin{equation}\label{L_h_2}
L_h = \frac{1}{h^2}
\begin{bmatrix}
& -1 \\
-1 & 4 & -1\\
& -1
\end{bmatrix}_h,
\end{equation}
and to fourth-order accuracy
\begin{equation}\label{L_h_2d_4o}
L_h = \frac{1}{12h^2}
\begin{bmatrix}
& & 1 \\  & & -16 \\
1 & -16 & 60 & -16 & 1\\
& & -16 \\ & & 1 
\end{bmatrix}_h.
\end{equation}

%\medskip
\subsection{Notation for schemes}
Consider a multigrid scheme for equation (\ref{MP_h}) with ($\lmax+1$) levels, on each level the grid-spacing being
$h_\ll,\ \ll = 0,\ 1,\ \cdots,\ \lmax$, $h_0 = h$. The operators on the finest level ($\ll=0$) is chosen to be \kl{$L_h = \LN{h}{p}$}, the standard finite-difference discretizations with order of accuracy $p\ (\equiv p_0)$.
%such as~\eqref{L_h_2} and~\eqref{L_h_2d_4o}. 
The scheme consists of a number of components, between two levels $\ll$ (fine) and $\ll+1$ (coarse), $\ll = 0,\ 1,\ \cdots,\ \lmax-1$:
\begin{enumerate}
	\item[-] a fine-level operator $L_{h_\ll}$,
	\item[-] a fine-level smoother $S_{h_\ll}$ based on $L_{h_\ll}$,
	\item[-] transfer operators: restriction $I_{h_\ll}^{h_{\ll+1}}(t)$ and interpolation $I_{h_{\ll+1}}^{h_\ll}(t) = (I_{h_\ll}^{h_{\ll+1}}(t))^*$ (order of accuracy $t$), 
	\item[-] together with a coarse-level operator $L_{h_{\ll+1}}$.
\end{enumerate}
In most cases, we will choose the red-black Gauss-Seidel (rb-GS) smoother with a relaxation parameter $\omega$. Different choices, however, will be made for the coarse-level operators and we now define some notation to clarify the presentation of the results. Either
\begin{itemize}
	\item $L_{h_{\ll+1}} = \kP{ \LN{h_{\ll+1}}{p_{\ll+1}} }$: non-Galerkin coarse-level operator (order of accuracy $p_{\ll+1} \leq p_\ll$)
	(in particular, denote $\kP{ \LNA{h_{\ll+1}} } \eqdef \LN{h_{\ll+1}}{p}$ with the same stencil as the finest-level operator); or
	\item $L_{h_{\ll+1}} = \kP{ \LG{h_{\ll+1}}{\tilde{L}_{h_\ll}}{t} } \eqdef \ I_{h_\ll}^{h_{\ll+1}}(t)\ \tilde{L}_{h_\ll}\ I_{h_{\ll+1}}^{h_\ll}(t)  $: Galerkin coarse-level operator, constructed from a fine-level operator $\tilde{L}_{h_\ll}$ and the transfer operators. We use the term `Galerkin' quite generally here; some specific examples we are going to consider are
	\begin{itemize}
		\item $\kP{ \LGG{h_{\ll+1}}} \eqdef \LG{h_{\ll+1}}{L_{h_\ll}}{t}$: the real ($(\ll+1)$-st-generation) Galerkin operator, constructed from the actual operator $L_h$ that is used on the fine level; 
		\item $\kP{ \LGq{h_{\ll+1}}{q} } \eqdef \LG{h_{\ll+1}}{\LGq{h_{\ll}}{q}}{t}$: the ($(\ll+1)$-st-generation) Galerkin operator constructed from the operator $\LN{h}{q}$ with order of accuracy $q\ (\leq p_\ll)$ on the finest level
		(in particular, $\LGq{h_{\ll+1}}{p} \equiv \LGG{h_{\ll+1}}$);
		\item $\kP{ \LGA{h_{\ll+1}}{q} } \eqdef \LG{h_{\ll+1}}{\LN{h_\ll}{q}}{t}$: the first-generation Galerkin operator from the operator $\LN{h_\ll}{q}$ with order of accuracy $q \leq p_\ll$.
		(In particular, denote $\kP{ \LGAA{h_{\ll+1}} } \eqdef \LGA{h_{\ll+1}}{p}$, constructed from the operator $\LN{h_\ll}{p}$ with the same stencil as the finest-level operator.)
	\end{itemize}
\end{itemize}

We assume the transfers to be standard linear interpolation and its adjoint restriction with $t=2$, unless otherwise specified. Hence we also adopt the shorthand notation for a Galerkin operator as
$\kP{ \LGt{h_{\ll+1}}{\tilde{L}_{h_\ll}} } \eqdef \LG{h_{\ll+1}}{\tilde{L}_{h_\ll}}{2}$.
(Only in Section \ref{Sec: orderC_t} we also consider $t=4$.)

In Section \ref{M1: orderC}, we consider $p = 4$ (or $p = 6$ in Section \ref{Sec: orderC_6o}), we choose the order of accuracy on each level $p_\ll \in \{p,\ q\}$, with $q = 2 < p$. 
%In most cases, $t = 2$, but we also consider $t=4$ in Section \ref{Sec: orderC_t}.
In Sections \ref{M2: rbC} and \ref{M3: rC}, on the other hand, we consider $p=q=2$. (For \ref{Sec: rbC_4o}, $p=4$, $q=2$.)

\bigskip
% MAIN---------------------------------------------------------------------------------------------------
%-----------------------------------------------------------------------------------------------
% topic-1: HIGH-ORDER ACCURACY WITH LOW ORDER COARSE-LEVEL-------------------------------------------
\section{High-order accurate discretizations with second-order accurate coarse-level operators}\label{M1: orderC}

In this section we consider using lower-order accurate coarse-level operators when solving the model
problem for Poisson's equation to fourth- and sixth-order accuracy.  This will be a desirable
approach to use in practice, since it is likely to be less expensive and importantly, it eases the
process of generating coarser level grids for the case of overset grids.  We will show that use of
lower-order accurate coarse level operators is very effective and generally results in convergence
rates comparable to, or better than, those obtained using high-order accurate coarse-level solves.

There is a simple heuristic argument that suggests why the coarse-level solves can be obtained using lower-order accurate approximations.
The argument is as follows.
In a typical (good)  multigrid algorithm, the residual reduction per cycle is approximately the asymptotic
convergence rate $\rho$ given by~\eqref{rho}, where $\rho$ is  roughly in the range $10^{-2}$ to $10^{-1}$
(corresponding to an ECR in the range 0.5 to 0.7, for example).
This range for $\rho$ is roughly independent of the order of accuracy on the fine grid.
Thus a typical coarse-level solve at level $\ll=1$
(obtained recursively with a multigrid cycle if $\lmax>1$) provides a correction $\tilde{v}_{h}^{(n)}$ to the fine grid at level $\ll=0$ ($h_0 = h$) that has a relative error $\e$ of roughly the size $10^{-1}$ to $10^{-2}$:
\ba
  \tilde{v}_{h}^{(n)} = \bar{v}_{h}^{(n)} - K_h v_{h}^{(n)}, \qquad \| K_h v_{h}^{(n)} \|_h = \Oc(\eps)\ \| \bar{v}_{h}^{(n)}\|_h , 
 \ea
where $\bar{v}_{h}^{(n)} = \bar{u}_{h} - u_{h}^{(n)}$ is the actual current error on the fine grid. 
It therefore seems reasonable to suppose that a correction with such a relative error can be provided instead 
by a lower-order (e.g. second-order) accurate coarse-level approximation, instead of using a high-order accurate one. 
  
Looking more closely, the target rate of a multigrid cycle (with $\nu$ smoothing sweeps) is $\rho
\approx \mu^{\nu}$, based on the smoothing rate with an ideal coarse-level correction $Q_h^H$.  As
illustrated in Figure \ref{fig: SK_HL}, a current error $v_h^{(n)}$ can be decomposed into its
high-frequency and low-frequency components.  In a multigrid cycle, the smoother mainly reduces the
high-frequency component of the error, while the coarse-level correction mainly reduces the
low-frequency component. Suppose in a multigrid cycle, via smoothing the high-frequency component of the error is reduced, from $v_h^{(n)}$ to $\bar{v}_h^{(n)} = S_h^{\nu} v_h^{(n)}$, by a factor of $\mu^{\nu}$. To achieve a multigrid convergence rate $\rho \approx \mu^{\nu}$, the coarse-level correction needs to
reduce the low-frequency component of the error by the same factor; further reduction of
low-frequency error does not improve the overall convergence rate.  For reference, the ideal
coarse-level correction $Q_h^H$, corresponding to a multigrid convergence rate of exactly
$\mu^{\nu}$, actually removes all low-frequency error.  In reality, the actual coarse level
correction $K_h$ is a perturbation of $Q_h^H$, and reduces the error along a subspace that is
closely aligned with the low-frequency component: $v_h^{(n+1)} = \bar{v}_h^{(n)} -
\tilde{v}_{h}^{(n)} = K_h v_{h}^{(n)}$. But this perturbation does not affect the main idea of our
heuristic argument.  (In \ref{Sec: orderC_heuristic} we give a more detailed argument through
local Fourier analysis.)

This heuristic argument suggests that second-order accurate coarse level operators can be productively used for fine level discretizations to arbitrarily high order of accuracy.
Here we show this to be true for fourth- and sixth-order accurate fine level discretizations.   

\newcommand{\drawaline}[4]{
	\draw [extended line=1cm,stealth-stealth] (#1,#2)--(#3,#4);
}
\begin{figure}[hbt]
	\centering
\begin{tikzpicture}[scale=1,extended line/.style={shorten >=-#1,shorten <=-#1},]
  \useasboundingbox (0, 0) rectangle (7,5.5);  % set the bounding box (so we have less surrounding white space)
  \begin{scope}[xshift=1cm,yshift=.3cm] 
  %\draw [help lines] (0,0) grid (6,5);
  % Euclidean
  \draw [extended line=.2cm,stealth-stealth](0,0)--(0,5) node[right]{$high$};
  \draw [extended line=.2cm,stealth-stealth](0,0)--(5,0) node[below]{$low$};
  % polar coordinate
  \draw [thin,dashed] (0,0)--(4,4) node[anchor=south west]{$v_h^{(n)}$};
  \draw [thin,dashed] (4,4)--(4,0);
  \draw [thin,dashed] (4,4)--(0,4);
  %\draw (0.8,1)node[anchor=west,rotate=45]{$r=\sqrt{x^2+y^2}$};
  %\draw [fill=green](0,0) -- (0.75,0) arc (0:45:0.75cm);
  %\draw (1,0.3) node[rotate=10]{$\theta=\tan^{-1} \frac{y}{x}$};
  \fill [red](4,4) circle(2pt);
  \draw [red](0,4) circle(2pt);
  %\draw [red](4,0) circle(2pt);
  \draw [orange](0,1) circle(2pt);
  \fill [orange](4,1) circle(2pt); 
  \draw [thin,dashed] (0,1)--(4,1) node[anchor=south west]{$\bar{v}_h^{(n)}$};
  \draw [-stealth, orange, very thick] (4,4) -- (4,1);
  %\draw [-stealth, orange, dashed, very thick] (0,4) -- (0,1);
  \draw [orange] (4,2.5) node[anchor=west] {$S_h^{\nu}$};
  \draw [orange](4,0) circle(2pt);
  
  \fill [orange](1,1) circle(3pt) node[anchor=south]{\textbf{target}}; 
  \draw [orange](1,0) circle(2pt);
  \draw [thin,dashed] (1,1)--(1,0); 
  
  %\draw [-stealth, DarkGreen, dashed, very thick] (4,0) -- (0,0);
  \draw [-stealth, DarkGreen, very thick] (4,1) -- (0,1) node[anchor=south east] {$Q_h^H \bar{v}_h^{(n)}$};
  \draw [DarkGreen] (2,1) node[anchor=south] {$Q_h^H$};
  \fill [DarkGreen](0,1) circle(2pt);
  \draw [DarkGreen](0,0) circle(2pt);
  
  \draw [orange](410/101, 41/101) circle(2pt);
  \draw [brown, thin,dashed] (4,1)--(410/101, 41/101) node[anchor= west]{$\tilde{v}_h^{(n)}$};
  %\draw [thin,dashed] (4,1)--(-6/101, 60/101);
  %\draw [orange](-6/101, 60/101) circle(2pt);
  \draw [brown] (2.2,0.5) node[anchor=south] {$K_h$};
  \fill [brown](-6/101, 60/101) circle(2pt);
  \draw [-stealth, brown, very thick] (4,1)--(-6/101, 60/101) node[anchor=east] {$v_h^{(n+1)}= K_h^H \bar{v}_h^{(n)}$};
  %\draw [-stealth, dashed, brown, very thick] (410/101, 41/101)--(0, 0);

  \draw [brown,extended line=.2cm,stealth-stealth] (0,0)--(5,.5);
  %node[above]{$\mathcal{C}(I_H^h)$}; 
  \draw [brown,extended line=.2cm,stealth-stealth] (0,0)--(-.5,5);
  %node[above]{$\mathcal{N}(I_h^H L_h)$};  
                                                
  \end{scope}
\end{tikzpicture}
\caption{Illustration: reduction of error in a two-level $[\nu,\ 0]$ cycle. The smoother mainly reduces the high-frequency component of the error, while the coarse-level correction mainly reduces the low-frequency component. 
If the smoother reduces the high-frequency component of the error by a factor of $\mu^{\nu}$, the coarse-level correction only needs to reduce the low-frequency component by the same factor (to \textcolor{orange}{\textbf{target}}).}
\label{fig: SK_HL}
\end{figure} 
%\vspace{-10pt}

% ----------------------------------------------------------------------
\subsection{Fourth-order accurate discretizations and second-order accurate coarse-level operators}
Consider the model problem (\ref{MP}) discretized to fourth-order accuracy with the operator given by (\ref{L_h_2d_4o}).                                     
We consider the red-black Gauss-Seidel smoother with an over-relaxation parameter $\omega$.
Separating the grid points in $\grid G_h$ into red and black points as
\begin{equation}
\begin{cases}
\grid G_h^R = \{\mathbf{x}_{\mathbf{j}} = 
h\ [j_1,\ j_2]^T :\ j_1+j_2 \in \mathbb{Z}_1\},\\
\grid G_h^B = \{\mathbf{x}_{\mathbf{j}} = 
h\ [j_1,\ j_2]^T :\ j_1+j_2 \in  \mathbb{Z}_0\},
\end{cases}
\end{equation}
(where $\Integer_1$ and $\Integer_0$ denote odd and even integers respectively), the red-black smoother consists of two partial Gauss-Seidel steps
\begin{align}\label{S_rb}
S_h = S^B_h\ S^R_h,\qquad
S^R_h = \begin{cases}
S_h^{GS}  & \text{on}\quad 
\grid G_h^R\\
I_h & \text{on}\quad
\grid G_h^B
\end{cases},\quad
S^B_h = \begin{cases}
S_h^{GS}  & \text{on}\quad \grid G_h^B\\
I_h & \text{on}\quad \grid G_h^R
\end{cases},
\end{align}
in which 
%the Gauss-Seidel smoothing operator is given by
\begin{equation}\label{S_GS}
S_h^{GS}(\omega) =  -(L^+_h)^{-1}L^-_h,
\end{equation}
with 
\begin{equation}
L^+_h = \frac{1}{12h^2}
\begin{bmatrix}
& & 1 \\  & & 0 \\
1 & 0 & 60\frac{1}{\omega} & 0 & 0\\
& & 0 \\ & & 0 
\end{bmatrix}_h,\quad L^-_h = L_h - L^+_h.
\end{equation}

Consider an $(h,\ H)$ two-level cycle with the coarse grid \eqref{G_H_d}, $d=2$, by standard coarsening with $H = 2h$, and the full-weighting restriction and linear interpolation operators 
\begin{equation}\label{transfer_2d}
I_h^H = \frac{1}{16}
\begin{bmatrix}
1 & 2 & 1 \\
2 & 4 & 2\\
1 & 2 & 1
\end{bmatrix}_h^{2h},\quad
I^h_H = \frac{1}{4} \left]
\begin{matrix}
1 & 2 & 1 \\
2 & 4 & 2\\
1 & 2 & 1
\end{matrix}\right[^h_{2h}.
\end{equation}
Note that these transfer operators are second-order accurate.
On the coarse grid $\grid G_H$, it is natural to use the operator with the same stencil as the fine-level operator
\begin{equation} \label{L_H_4n_2d}
L_H = \frac{1}{12H^2} \begin{bmatrix}
& & 1 \\  & & -16 \\
1 & -16 & 60 & -16 & 1\\
& & -16 \\ & & 1 
\end{bmatrix}_H,
\end{equation}
%with eigenvalue
%\begin{equation}
%\symF L_H({2\th}) = \frac{4}{H^2} \big(\sin^2 \th + \frac{1}{3}\sin^4 \th \big);
%end{equation}
or the Galerkin coarse-level operator 
\begin{equation}\label{L_H_4G}
L_H =  I_h^H L_h I^h_H.
\end{equation}
\begin{comment}
with eigenvalue 
\begin{equation}
\symF L_H(2\th) = \symE I_h^H(\th)\ \symE L(\th)\ \symE I^h_H (\th),
\end{equation}
where
\begin{equation}
\symE L(\th) = \begin{bmatrix}
\symF L_h(\th)\\ & \symF L_h(\bar{\th})
\end{bmatrix}
\end{equation}
\end{comment}
With the coarse-level operator and the transfer operators we have the coarse-level correction operator $K_h^H$ as in (\ref{K}), and from the smoother and the coarse-level correction we have the two-level iteration operator $M_h$, as given in (\ref{M_2level}).

The main idea of this section is to consider coarse-level operators induced from the second-order accurate fine-level operator
\begin{equation}\label{L_h_2_bra}
L_h^{(2)} = \frac{1}{h^2} \begin{bmatrix}
& -1 \\
-1 & 4 & -1\\
& -1
\end{bmatrix}_h
\end{equation}
instead of (\ref{L_h_2d_4o}). That is, the non-Galerkin 
\begin{equation}\label{L_H_2d}
L_H = \frac{1}{H^2} \begin{bmatrix}
& -1 \\
-1 & 4 & -1\\
& -1
\end{bmatrix}_H,
\end{equation}
and the ``Galerkin" 
\begin{equation}\label{L_H_2G_2d}
L_H =  I_h^H L_h^{(2)} I^h_H.
\end{equation}
Here ``Galerkin" is in quotation since the coarse-level operator (\ref{L_H_2G_2d}) does not come from the actual operator on the fine level (\ref{L_h_2d_4o}). (Thus for instance the resulting coarse-level correction operator $K_h^H$ will not be a projector in this case.)
Note that (\ref{L_H_2d}) or (\ref{L_H_2G_2d}) makes the coarse-level problem $L_H\ v_H = d_H$ to be of second-order accuracy, which requires less complexity in terms of both the grid and the solver on the coarse level. At the same time, we claim that reducing the coarse-level operator to lower-order accuracy does not hinder the multigrid convergence. 

{
\newcommand{\figWidth}{6cm}% height 
\newcommand{\trimfig}[2]{\trimh{#1}{#2}{.0}{.0}{.0}{.0}}
\begin{figure}[h!]
\begin{center}
\begin{tikzpicture}[scale=1]
  \useasboundingbox (0,0.6) rectangle (15,12.1);  % set the bounding box (so we have less surrounding white space)

   \draw(0.0,0.0) node[anchor=south west,xshift=-15pt,yshift=-8pt]{\trimfig{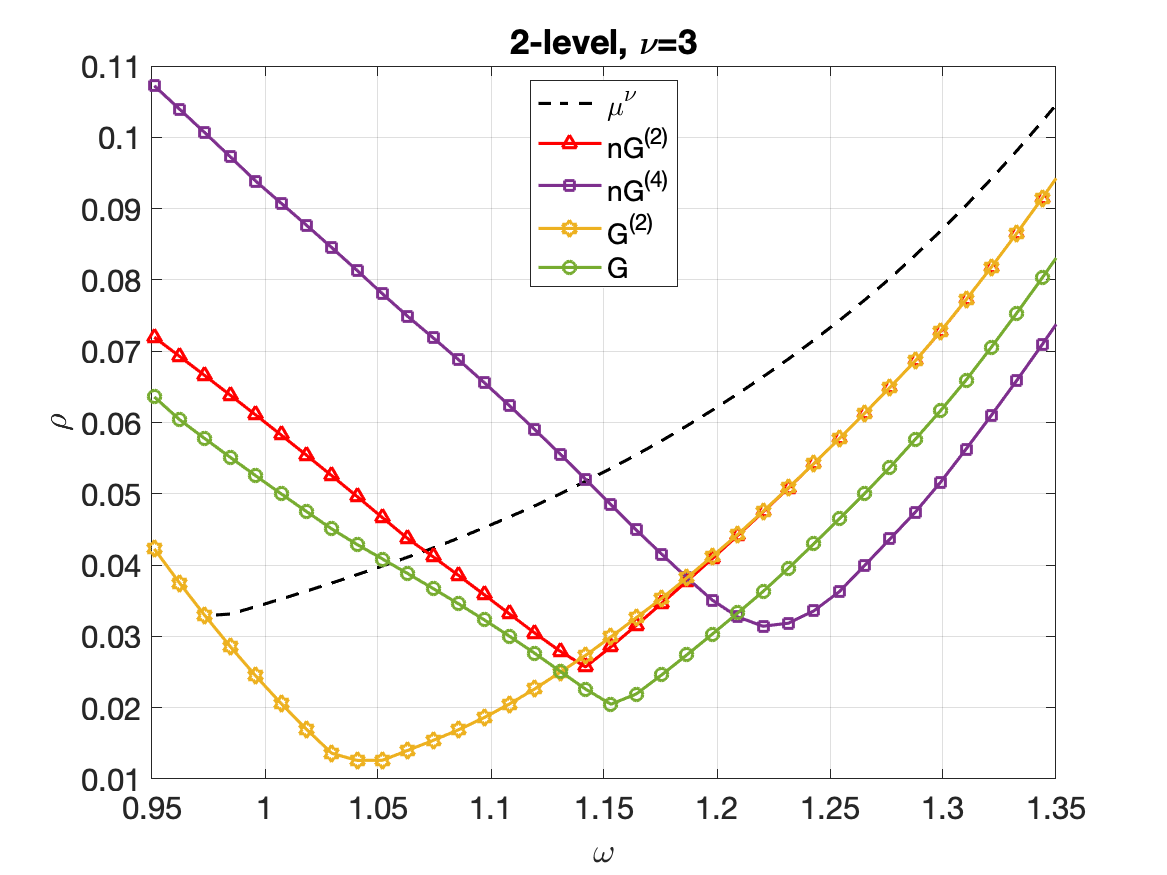}{\figWidth}};
   \draw(7.5,0.0) node[anchor=south west,xshift=-15pt,yshift=-8pt]{\trimfig{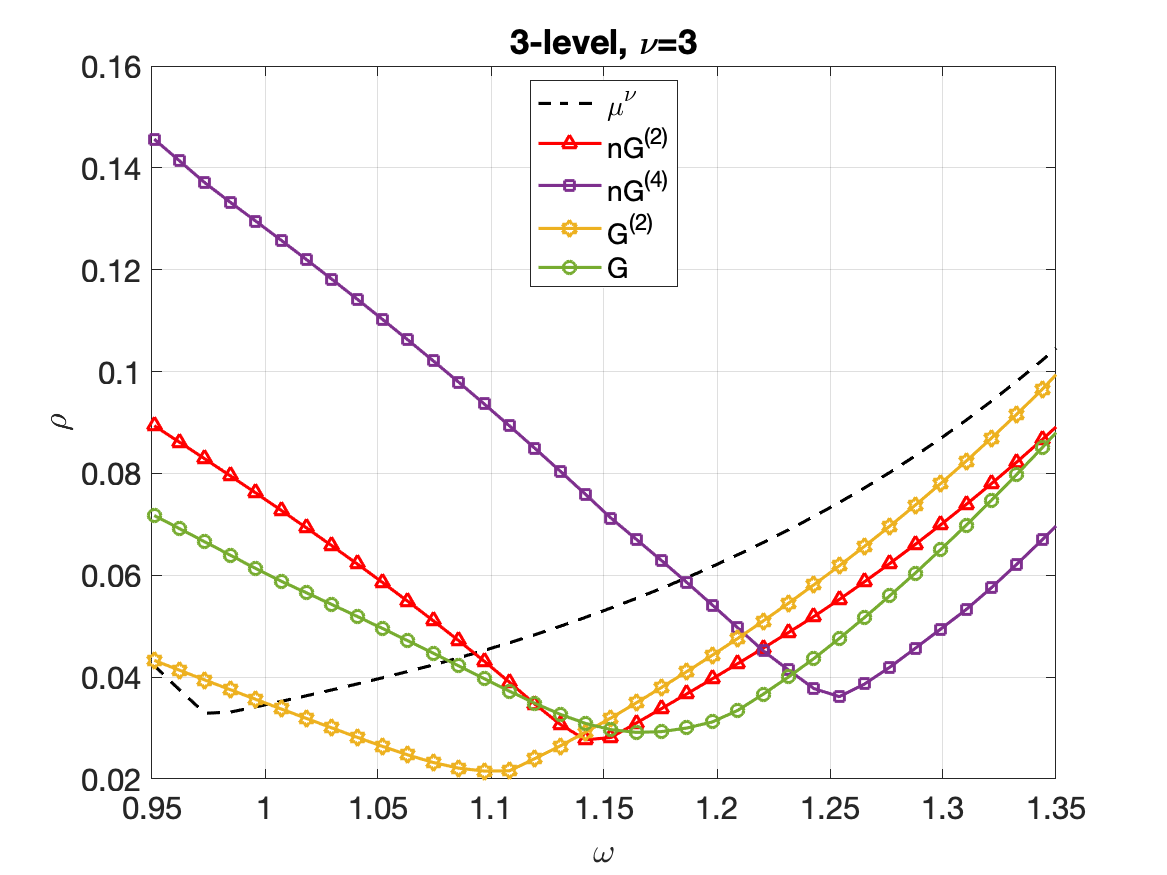}{\figWidth}};

   %  circles to highlight best results
   \circleLeftLabel{2.55}{.8}{opt CR};
   \circleLeftLabel{10.9}{.7}{opt CR};

   \circleLabel{1.5}{1.7}{opt $\mu^\nu$}{violet}{.25};
   \circleLabelLeft{9}{1.}{opt $\mu^\nu$}{violet}{.25};

   \begin{scope}[yshift=6cm]
     \draw(0.0,0)     node[anchor=south west,xshift=-15pt,yshift=-8pt]{\trimfig{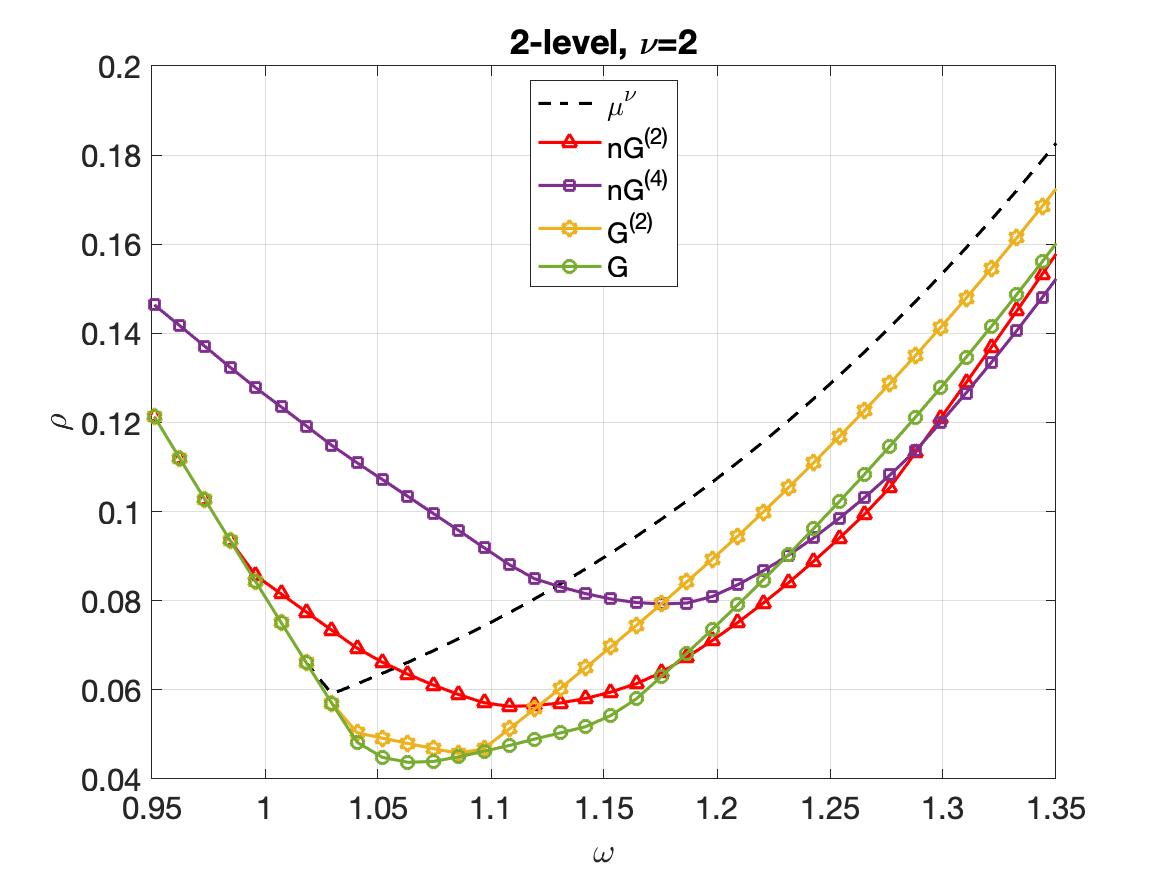}{\figWidth}};
     \draw(7.5,0)   node[anchor=south west,xshift=-15pt,yshift=-8pt]{\trimfig{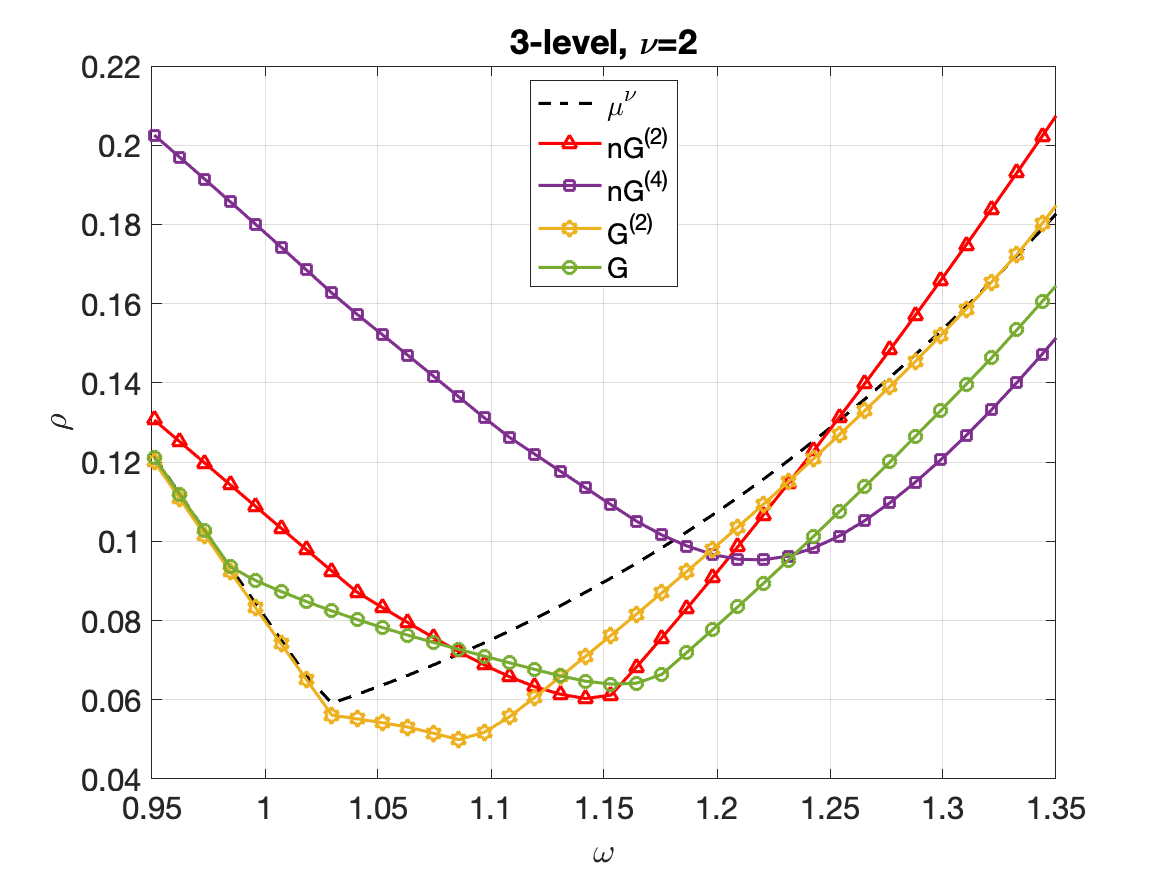}{\figWidth}};

     %  circles to highlight best results
     \circleRightLabel{2.9}{.7}{opt CR};
     \circleLeftLabel{10.7}{.9}{opt CR};

     \circleLabelLeft{2.3}{1.15}{opt $\mu^\nu$}{violet}{.25};
     \circleLabelLeft{9.8}{1.1}{opt $\mu^\nu$}{violet}{.25};

   \end{scope}   
   % title
   \draw(7.5,11.8)  node[draw=\colourTwoD,very thick,fill=\colourTwoD!20,anchor=south,inner sep=2.5pt] {\small \DIM, \ORD[4], 2nd-order transfers, \Cc};

   % small grid cartoons
   % \standardGridOrderFour{xshift= 5.75cm,yshift=4.5cm}; 
   %\standardGridOrderFour{xshift=13.00cm,yshift=4.5cm}; 

   %\standardGridOrderFour{xshift= 6cm,yshift=10.5cm}; 
   %\standardGridOrderFour{xshift=13.5cm,yshift=10.5cm}; 

   % V cycle
   %\VcycleTwoLevel{xshift= 2.5cm,yshift=10.5cm}; 
   %\VcycleThreeLevel{xshift=10.5cm,yshift=10.5cm}; 
   % \VcycleTwoLevel{xshift= 2.5cm,yshift= 4.5cm}; 
   %\VcycleThreeLevel{xshift=10.5cm,yshift= 4.5cm}; 

   % Combined cartoons
   \orderFourVcycleTwoLevel{xshift=1.7cm,yshift=4.75cm,scale=1}{$\nu\!=\!3$}
   \begin{scope}[xshift=0cm,yshift=6cm]
     \orderFourVcycleTwoLevel{xshift=1.7cm,yshift=4.75cm,scale=1}{$\nu\!=\!2$}
   \end{scope}
     
   \begin{scope}[xshift=7.5cm]
     \orderFourVcycleThreeLevel{xshift=1.7cm,yshift=4.75cm,scale=1}{$\nu\!=\!3$}
   \end{scope}
   \begin{scope}[xshift=7.5cm,yshift=6cm]
     \orderFourVcycleThreeLevel{xshift=1.7cm,yshift=4.75cm,scale=1}{$\nu\!=\!2$}
   \end{scope}

   % \begin{scope}[xshift=1.7cm,yshift=4.75cm,scale=1]
   %   \draw[-,thick,fill=white] (-.65,-.65) -- (1.9,-.65) -- (1.9,.65) -- (-.65,.65) -- (-.65,-.65); % border
   %   \standardGridOrderFour{xshift=0cm,yshift=0cm}; 
   %   \VcycleTwoLevel{xshift=1.25cm,yshift=0cm}; 
   % \end{scope}

% grid:
% \draw[step=1cm,gray] (0,0) grid (15,13);
\end{tikzpicture}
\end{center}
\caption{$\rho$ versus $\omega$. \DIM, \ORD[4], 2-level and 3-level V cycles with red-black Gauss-Seidel smoothing ($\nu$ sweeps per cycle) and \CC,
\kl{
	2nd-order accurate transfer operators;
	4th-order accurate non-Galerkin (nG$^{(4)}=\LNH{4} \equiv \LNAH$),
	4th-order accurate Galerkin (G = $\LGGH \equiv \LGHh{4}{2}$),
	2nd-order accurate non-Galerkin (nG$^{(2)}=\LNH{2}$),
	and
	2nd-order accurate Galerkin (G$^{(2)}= \LGqH{2} \equiv \LGHh{2}{2}$) coarse-level operators.
}
  }
  \label{fig:rho_V_2d_4o_std}
\end{figure}
}

In Figure \ref{fig:rho_V_2d_4o_std} we plot the asymptotic convergence rate $\rho$ of the two- and three-level V cycles verses the smoothing parameter $\omega$. 
(For three-level cycles the same $\omega$ is used on both the fine level and the first coarse level.) 
We compare the choice of the various coarse-level operators discussed: the fourth-order accurate non-Galerkin operator (nG$^{(4)}$) given by (\ref{L_H_4n_2d}), the (``fourth-order accurate"\footnote{Note that by calling a Galerkin coarse-level operator ``fourth-order accurate", we mean that it is constructed from a \textsl{fine}-level operator of fourth-order accuracy. This does not necessarily mean that the resulting coarse-level operator $L_H$ is a fourth-order accurate approximation, in terms of $\mathcal{O}(H^4)$, to the continuous operator $L$, particularly when the transfer operators are only of second-order accuracy.}) Galerkin operator (G) given by (\ref{L_H_4G}); and the second-order accurate non-Galerkin and ``Galerkin" operators given by (\ref{L_H_2d}) and (\ref{L_H_2G_2d}) respectively.  
We also show the smoothing rate $\mu^\nu$ for comparison, since the smoothing rate is an estimate for the expected convergence rate.
We can see that the second-order accurate coarse-level operators do no worse than the fourth-order accurate ones. In fact, by adjusting the value of the over-relaxation parameter $\omega$ of the smoother, we can usually achieve better
convergence rates with second-order accurate coarse-level operators.

A numerical example has been given in Figure \ref{fig:intro_topics} (left), which shows computed reduction of residual against WU for a multi-level $V[2,\ 1]$ cycle with the second-order accurate Galerkin coarse-level operator (\ref{L_H_2G_2d}), with the smoothing parameter tuned at $\omega = 1.1$. For reference, we compare this against the basic choice of the fourth-order accurate non-Galerkin operator (\ref{L_H_4n_2d}) without over-relaxation ($\omega = 1$), and substantial improvement of convergence can be observed.

% ---------------------------------------------------------------------------------------
% ------- 4TH ORDER, 4TH-ORDER COARSE GRID  V-CYCLES -----
{
\newcommand{\figWidth}{6cm}% height 
\newcommand{\trimfig}[2]{\trimh{#1}{#2}{.0}{.0}{.0}{.0}}
\begin{figure}[h!]
\begin{center}
\begin{tikzpicture}[scale=1]
  \useasboundingbox (0,0.5) rectangle (15,12);  % set the bounding box (so we have less surrounding white space)

  % \draw(0,6) node[anchor=south west,xshift=-15pt,yshift=-8pt]{\trimfig{W_rho_2_4o_2g_some}{\figWidth}};
  % \draw(7.5,6) node[anchor=south west,xshift=-15pt,yshift=-8pt]{\trimfig{W_rho_2_4o_3g_some}{\figWidth}};
  \draw(0.0,0.0) node[anchor=south west,xshift=-15pt,yshift=-8pt]{\trimfig{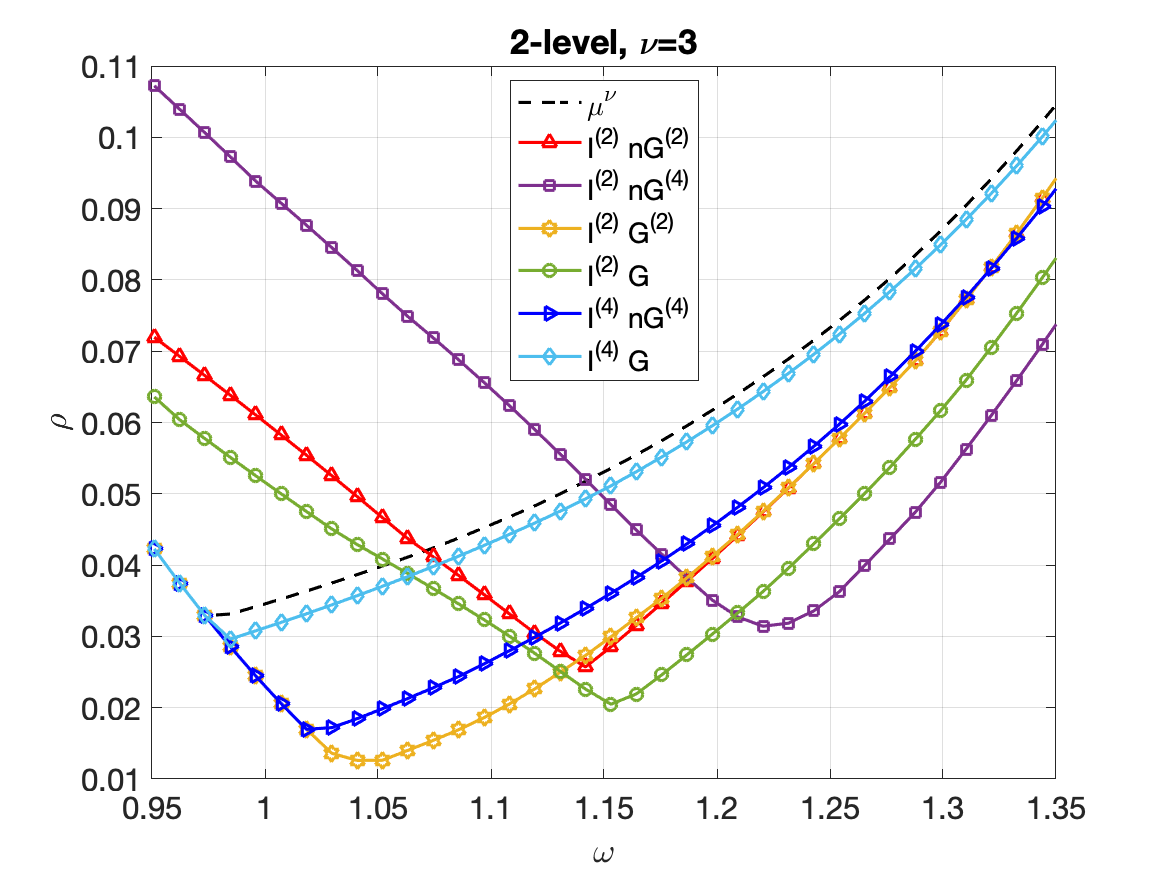}{\figWidth}};
  \draw(7.5,0.0) node[anchor=south west,xshift=-15pt,yshift=-8pt]{\trimfig{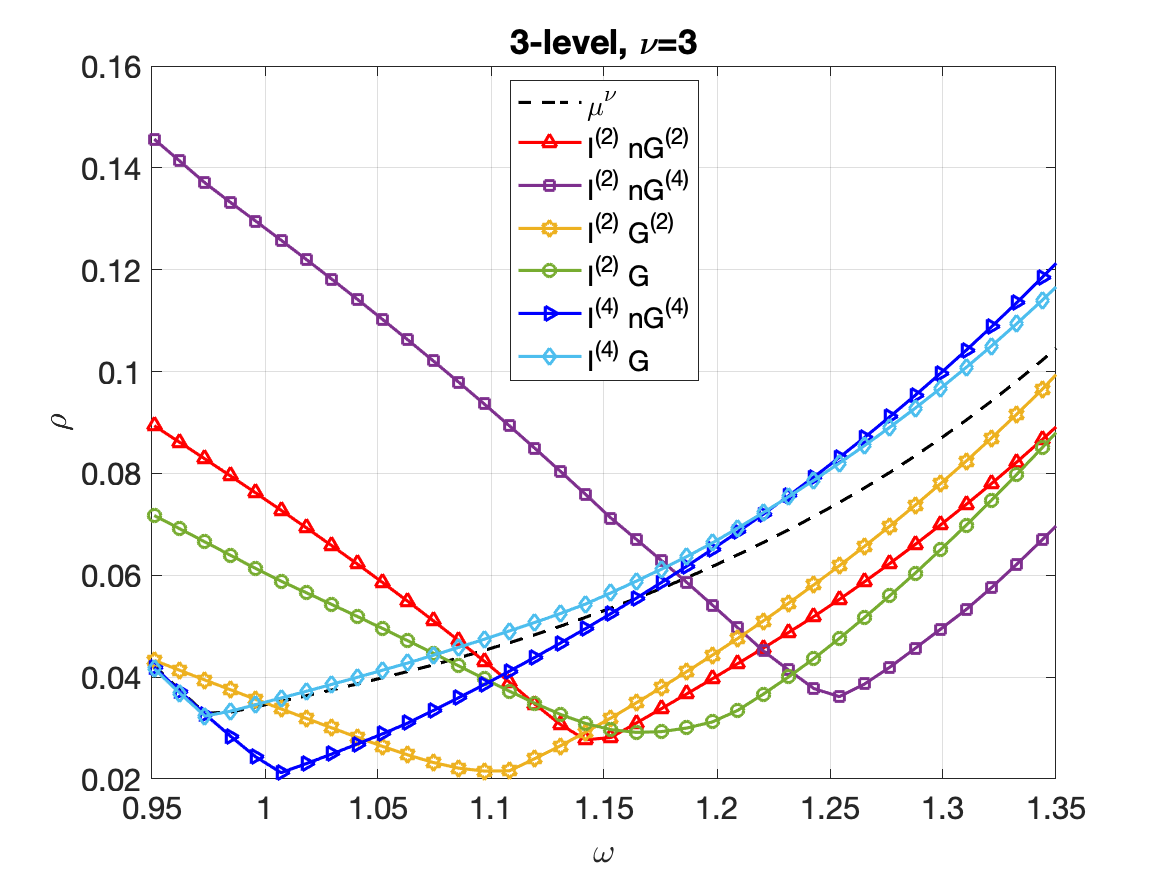}{\figWidth}};

  %  circles to highlight best results
   \circleLeftLabel{2.55}{.8}{opt CR};
   \circleLeftLabel{10.9}{.7}{opt CR};

   \circleLabel{1.5}{1.7}{opt $\mu^\nu$}{violet}{.25};
   \circleLabelLeft{9}{1.}{opt $\mu^\nu$}{violet}{.25};
   
   \begin{scope}[yshift=6cm]
     \draw(  0,0) node[anchor=south west,xshift=-15pt,yshift=-8pt]{\trimfig{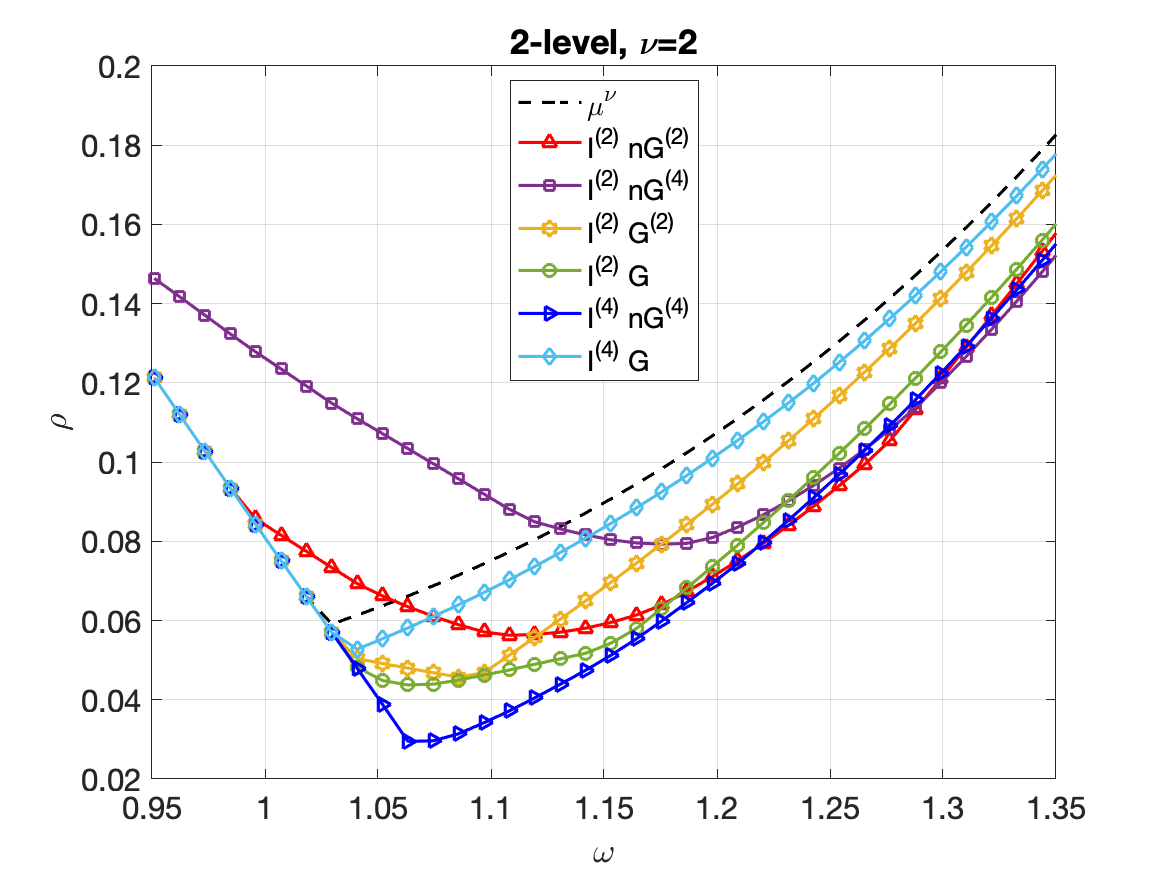}{\figWidth}};
     \draw(7.5,0) node[anchor=south west,xshift=-15pt,yshift=-8pt]{\trimfig{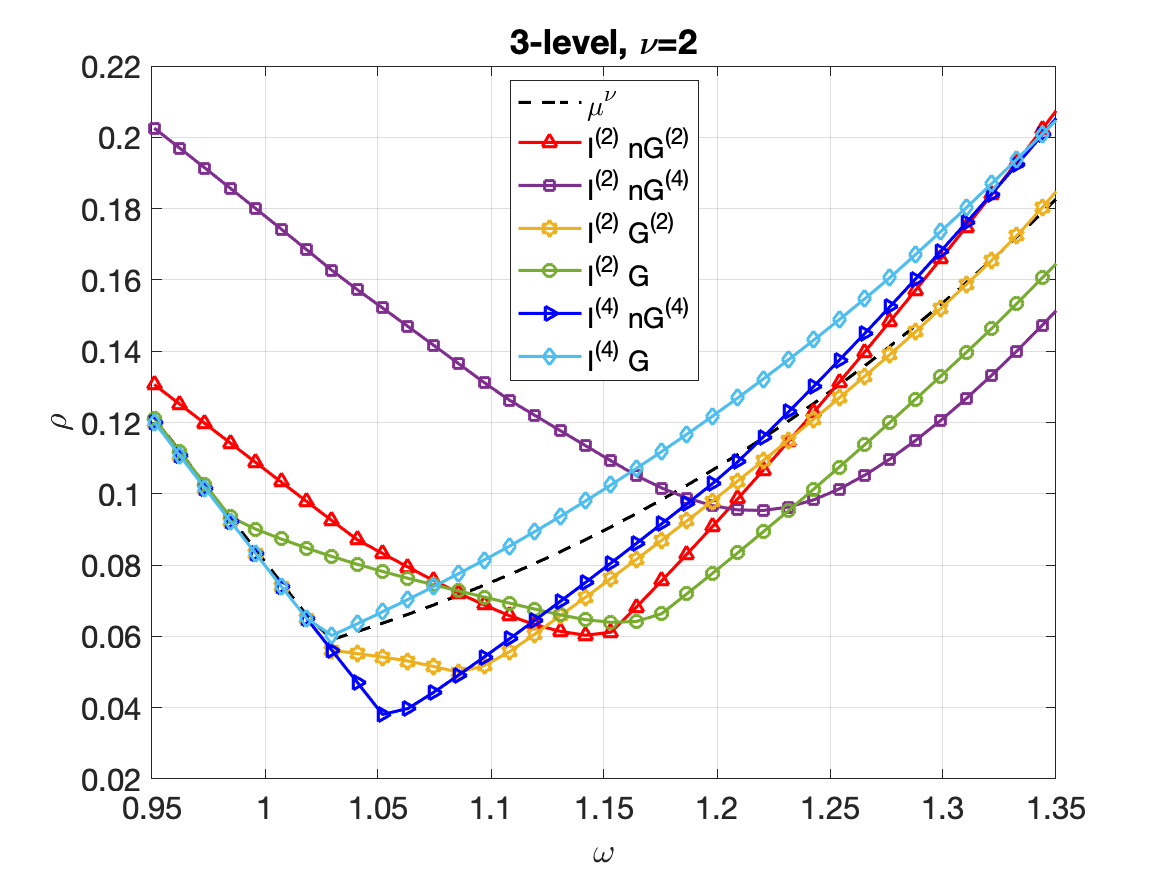}{\figWidth}};

     %  circles to highlight best results
     \circleRightLabel{2.9}{.8}{opt CR};
     \circleLeftLabel{10.2}{1}{opt CR};

     \circleLabelLeft{2.3}{1.65}{opt $\mu^\nu$}{violet}{.25};
     \circleLabelLeft{9.8}{1.5}{opt $\mu^\nu$}{violet}{.25};

   \end{scope}
   
  % title
  \draw(7.5,11.7)  node[draw=\colourTwoD,very thick,fill=\colourTwoD!20,anchor=south,inner sep=2.5pt] {\small \DIM, \ORD[4], \Cc}; 

%   % small grid cartoons
%   \standardGridOrderFour{xshift= 5.75cm,yshift=4.5cm}; 
%   \standardGridOrderFour{xshift=13.50cm,yshift=4.5cm}; 
%   \standardGridOrderFour{xshift=5.50cm,yshift=10.5cm}; 
%   \standardGridOrderFour{xshift=13.0cm,yshift=10.5cm}; 
%
%   % V cycle
%   \VcycleThreeLevel{xshift= 2.5cm,yshift=10.5cm}; 
%   \VcycleThreeLevel{xshift=10.5cm,yshift=10.5cm}; 
%   \VcycleThreeLevel{xshift= 2.5cm,yshift= 4.5cm}; 
%   \VcycleThreeLevel{xshift=10.5cm,yshift= 4.5cm}; 

   % Combined cartoons
   \orderFourVcycleTwoLevel{xshift=1.7cm,yshift=4.75cm,scale=1}{$\nu\!=\!3$}
   \begin{scope}[xshift=0cm,yshift=6cm]
     \orderFourVcycleTwoLevel{xshift=1.7cm,yshift=4.75cm,scale=1}{$\nu\!=\!2$}
   \end{scope}
     
   \begin{scope}[xshift=7.5cm]
     \orderFourVcycleThreeLevel{xshift=1.7cm,yshift=4.75cm,scale=1}{$\nu\!=\!3$}
   \end{scope}
   \begin{scope}[xshift=7.5cm,yshift=6cm]
     \orderFourVcycleThreeLevel{xshift=1.7cm,yshift=4.75cm,scale=1}{$\nu\!=\!2$}
   \end{scope}

   % grid:
% \draw[step=1cm,gray] (0,0) grid (15,12.5);
\end{tikzpicture}
\end{center}
\caption{$\rho$ versus $\omega$. \DIM, \ORD[4],
  V cycles with red-black Gauss-Seidel smoothing ($\nu$ sweeps per cycle) and \CC;  
	2nd-order (I$^{(2)}$) and 4th-order accurate transfer operators (I$^{(4)}$);
	4th-order accurate non-Galerkin (nG$^{(4)}=\LNH{4} \equiv \LNAH$),
	4th-order accurate Galerkin (G=$\LGGH \equiv \LGHh{4}{t}$, \kP{$t=2$ for I$^{(2)}$ and $t=4$ for I$^{(4)}$}),
	2nd-order accurate non-Galerkin (nG$^{(2)}=\LNH{2}$),
	and
	2nd-order accurate Galerkin (G$^{(2)}=\LGqH{2} \equiv \LGHh{2}{2}$) coarse-level operators.}
  \label{fig:rho_V_2d_4o_more_std}
\end{figure}
}
\subsection{High-order accurate transfer operators}\label{Sec: orderC_t}
Besides the order of the coarse-level operators, it may also be a natural consideration to use higher-order accurate transfer operators for higher-order accurate discretizations. For completeness we also consider cubic interpolation \eqref{interp_2d4o} and its adjoint restriction \eqref{restr_2d4o} operators, which are of fourth-order accuracy, in place of the second-order accurate transfer operators (\ref{transfer_2d}). In this case, we still consider the fourth-order accurate non-Galerkin coarse-level operator (\ref{L_H_4n_2d}), as well as the Galerkin operator constructed from the fourth-order accurate transfer operators. 

In Figure \ref{fig:rho_V_2d_4o_more_std} we show the same results with the second-order accurate transfer operators (I$^{(2)}$) as given in Figure \ref{fig:rho_V_2d_4o_std}, and compare them with the multigrid cycles with the fourth-order accurate transfer operators (I$^{(4)}$) and the fourth-order accurate non-Galerkin operator (nG$^{(4)}$) or the (fourth-order accurate) Galerkin operator (G).
We can observe that fourth-order accurate transfer operators along with fourth-order accurate coarse-level solves give convergence rates that are closer to the smoothing rate, which reflects the fact that the resulting coarse-level correction operators are higher-order accurate approximations to the ideal $Q_h^H$. On the other hand, this is the most computationally expensive combination.

% ---- summary fig------
{
\newcommand{\figWidth}{6cm}% height 
\newcommand{\trimfig}[2]{\trimh{#1}{#2}{.0}{.0}{.0}{.0}}
\begin{figure}[h!]
\begin{center}
\begin{tikzpicture}[scale=1]
  \useasboundingbox (0,0.5) rectangle (16,5.9);  % set the bounding box (so we have less surrounding white space)

     \draw(0,0) node[anchor=south west,xshift=-15pt,yshift=-8pt]{\trimfig{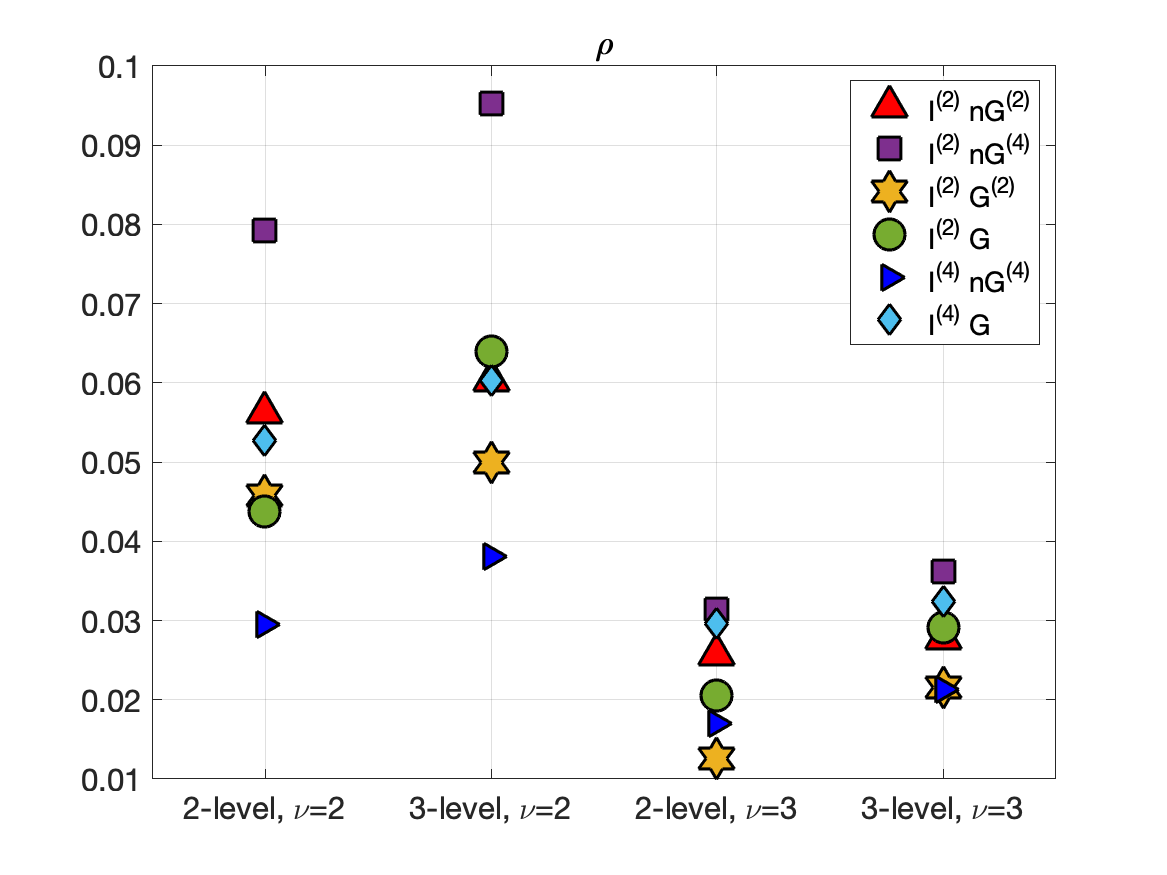}{\figWidth}};
     \draw(8,0) node[anchor=south west,xshift=-15pt,yshift=-8pt]{\trimfig{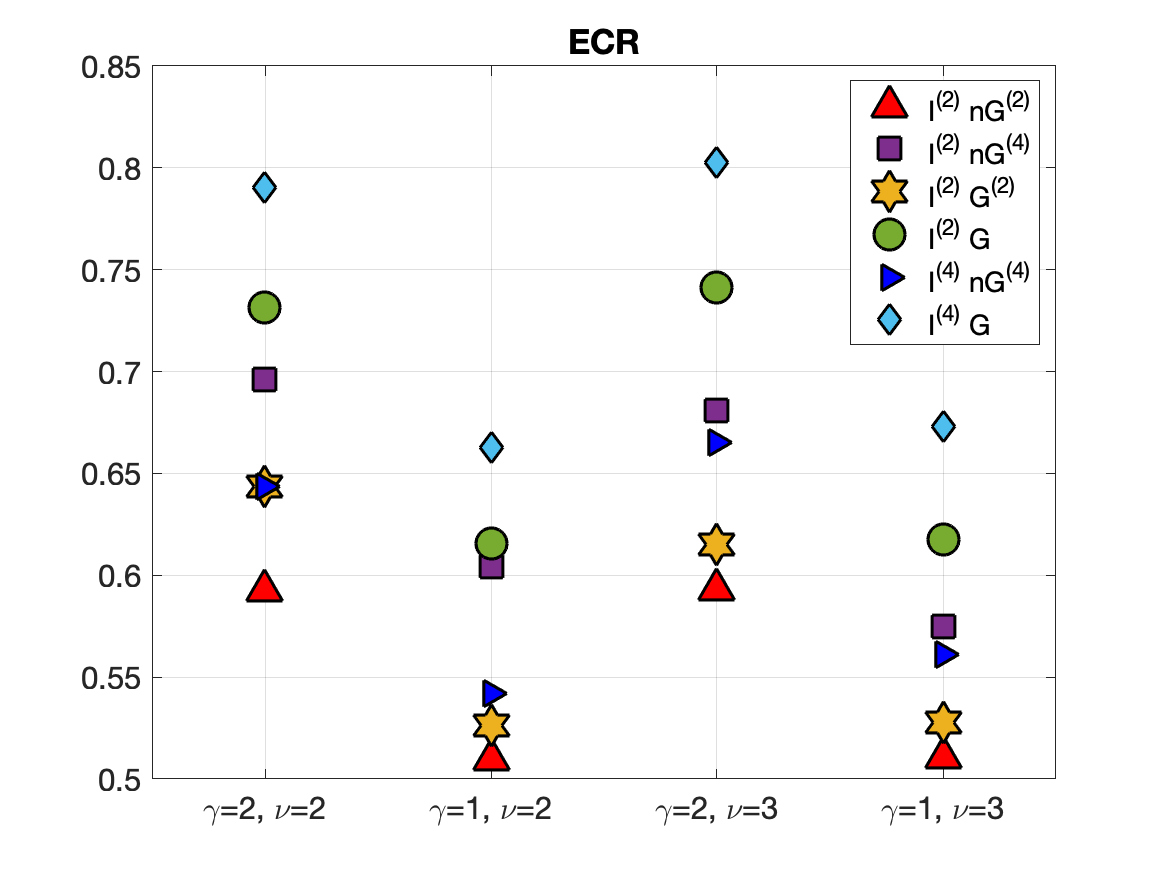}{\figWidth}};
   % title
     \draw(8.2,5.7)  node[draw=\colourTwoD,very thick,fill=\colourTwoD!20,anchor=south] {\small Summary: \DIM, \ORD[4], \Cc}; 
   
% grid:
% \draw[step=1cm,gray] (0,0) grid (8,6.5);
\end{tikzpicture}
\end{center}
\caption{Optimal $\rho$ and ECR. \DIM, \ORD[4], V ($\gamma = 1$) and W ($\gamma = 2$) cycles with red-black Gauss-Seidel smoothing ($\nu$ sweeps per cycle) and \CC;
  	2nd-order (I$^{(2)}$) and 4th-order accurate transfer operators (I$^{(4)}$);
  	4th-order accurate non-Galerkin (nG$^{(4)}=\LNH{4} \equiv \LNAH$),
  	4th-order accurate Galerkin (G=$\LGGH \equiv \LGHh{4}{t}$, \kP{$t=2$ for I$^{(2)}$ and $t=4$ for I$^{(4)}$}),
  	2nd-order accurate non-Galerkin (nG$^{(2)}=\LNH{2}$),
  	and
  	2nd-order accurate Galerkin (G$^{(2)}=\LGqH{2} \equiv \LGHh{2}{2}$) coarse-level operators;
  	at their respective optimal $\omega$.
}
  \label{fig:rho_opt}
\end{figure}
}

From Figure \ref{fig:rho_V_2d_4o_more_std} we can make out the optimal $\omega$ for each choice of
coarse-level correction for a multigrid $\gamma$-cycle. Note that the three-level CRs are representative
of multi-level V cycles ($\gamma=1$) in practice, while the two-level CRs are representative of W
cycles ($\gamma=2$).  In Figure \ref{fig:rho_opt} (left) we record the optimal asymptotic
convergence rates $\rho$ for the various choices of transfer and coarse-level operators at their
respective optimal $\omega$.
It can be seen that the use of second-order accurate coarse level operators gives convergence rates comparable to those with the fourth-order accurate ones. 
A more fair comparison between the schemes will also take into account the computational cost and therefore in Figure \ref{fig:rho_opt} (right) we provide estimates of the ECR. Work-units for a W cycle are associated with CRs of two-level cycles while work-units of a V cycle are associated with CRs of the three-level cycles.
Note that one WU is defined as the number of FLOPS in a single Jacobi smoothing step for the
operator $L_h$ on the finest level, which in this case is (\ref{L_h_2d_4o}), a sparse stencil of
\textsl{fourth}-order accuracy.
We estimate the work-units for a multigrid cycle based on using as many levels as possible
(see the discussion in \ref{sec: WU} for the details).  
Non-Galerkin coarse-level operators have some advantage in terms of WU over Galerkin operators because of their sparser stencils, especially for the simple example of the standard discrete Laplacians.
In terms of effective convergence rates, the schemes using second-order accurate coarse-level corrections are found to be always as good as or better than those using higher-order accurate ones. These results provide some good justification for the use of second-order accurate coarse-level correction operators for fourth-order accurate fine-level discretizations, especially given the added benefits for grid generation that are associated with the use of lower-order accurate coarse-level operators.

% ------- SIXTH ORDER ----------------------------------------------------------
{
\newcommand{\figWidth}{6cm}% height 
\newcommand{\trimfig}[2]{\trimh{#1}{#2}{.0}{.0}{.0}{.0}}
\begin{figure}[h!]
\begin{center}
\begin{tikzpicture}[scale=1]
  \useasboundingbox (0,0.6) rectangle (15,12.2);  % set the bounding box (so we have less surrounding white space)

   \begin{scope}[yshift=6cm]
     \draw(0.0,0)     node[anchor=south west,xshift=-15pt,yshift=-8pt]{\trimfig{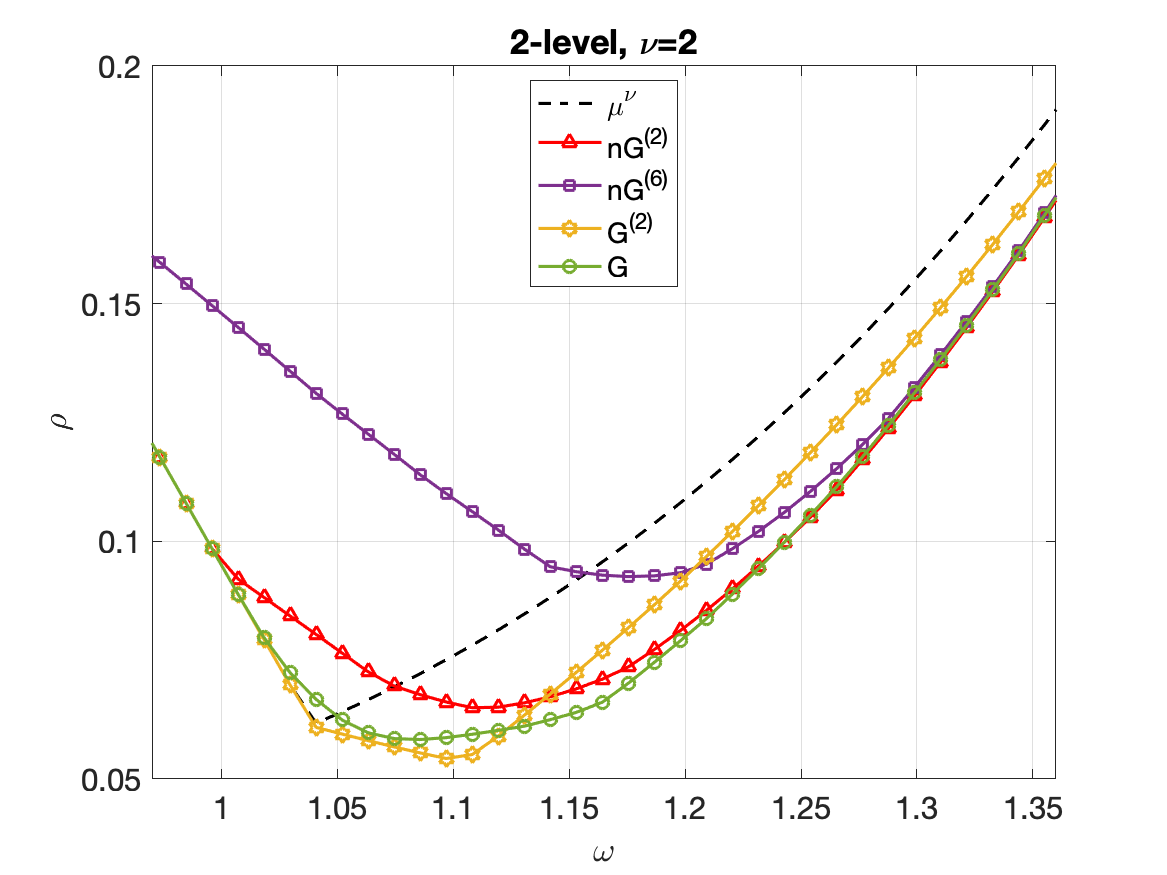}{\figWidth}};
     \draw(7.5,0)   node[anchor=south west,xshift=-15pt,yshift=-8pt]{\trimfig{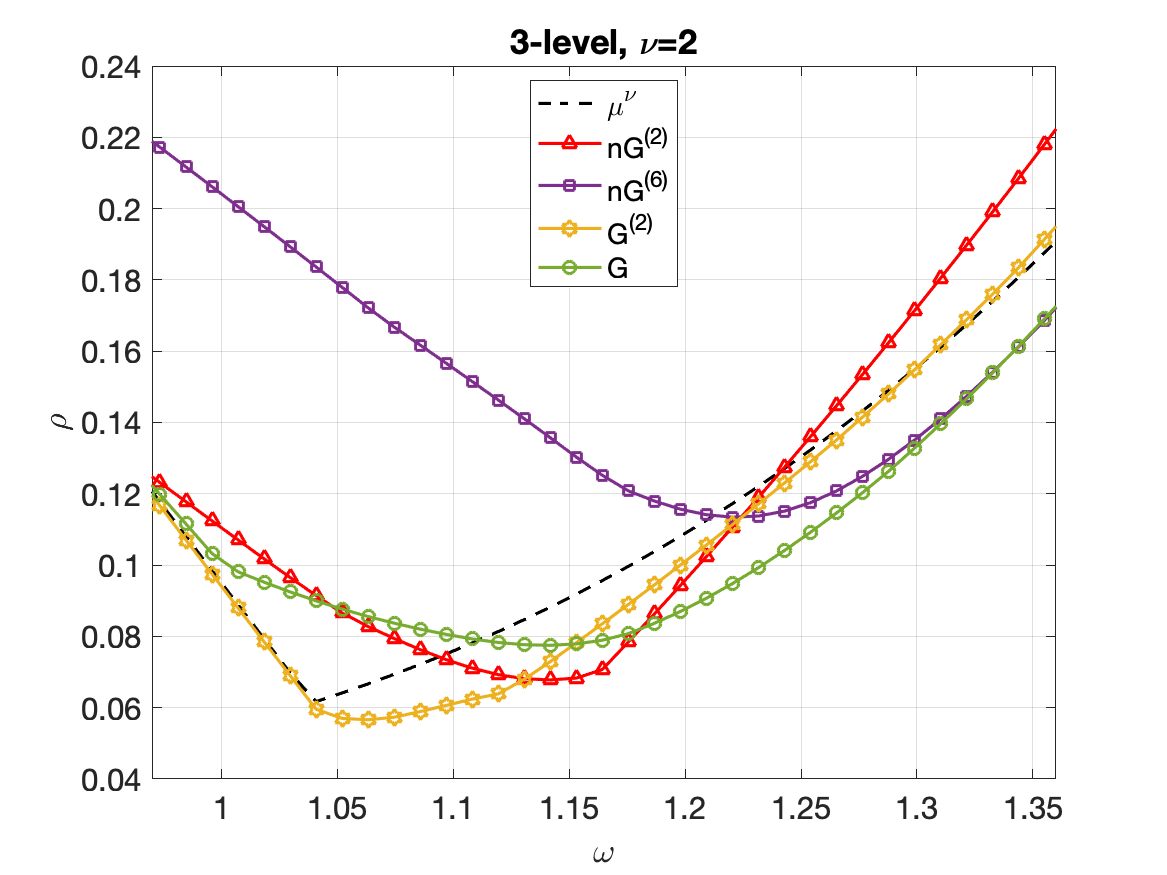}{\figWidth}};

     %  circles to highlight best results
     \circleLeftLabel{3.2}{.8}{opt CR};
     \circleLeftLabel{10.0}{1}{opt CR};

     \circleLabelLeft{2.3}{.9}{opt $\mu^\nu$}{violet}{.25};
     \circleLabelLeft{9.7}{1.2}{opt $\mu^\nu$}{violet}{.25};

   \end{scope}

   % \draw(0,6)     node[anchor=south west,xshift=-15pt,yshift=-8pt]{\trimfig{W_rho_2_6o_2g}{\figWidth}};
   % \draw(7.5,6)   node[anchor=south west,xshift=-15pt,yshift=-8pt]{\trimfig{W_rho_2_6o_3g}{\figWidth}};
   \draw(0.0,0.0) node[anchor=south west,xshift=-15pt,yshift=-8pt]{\trimfig{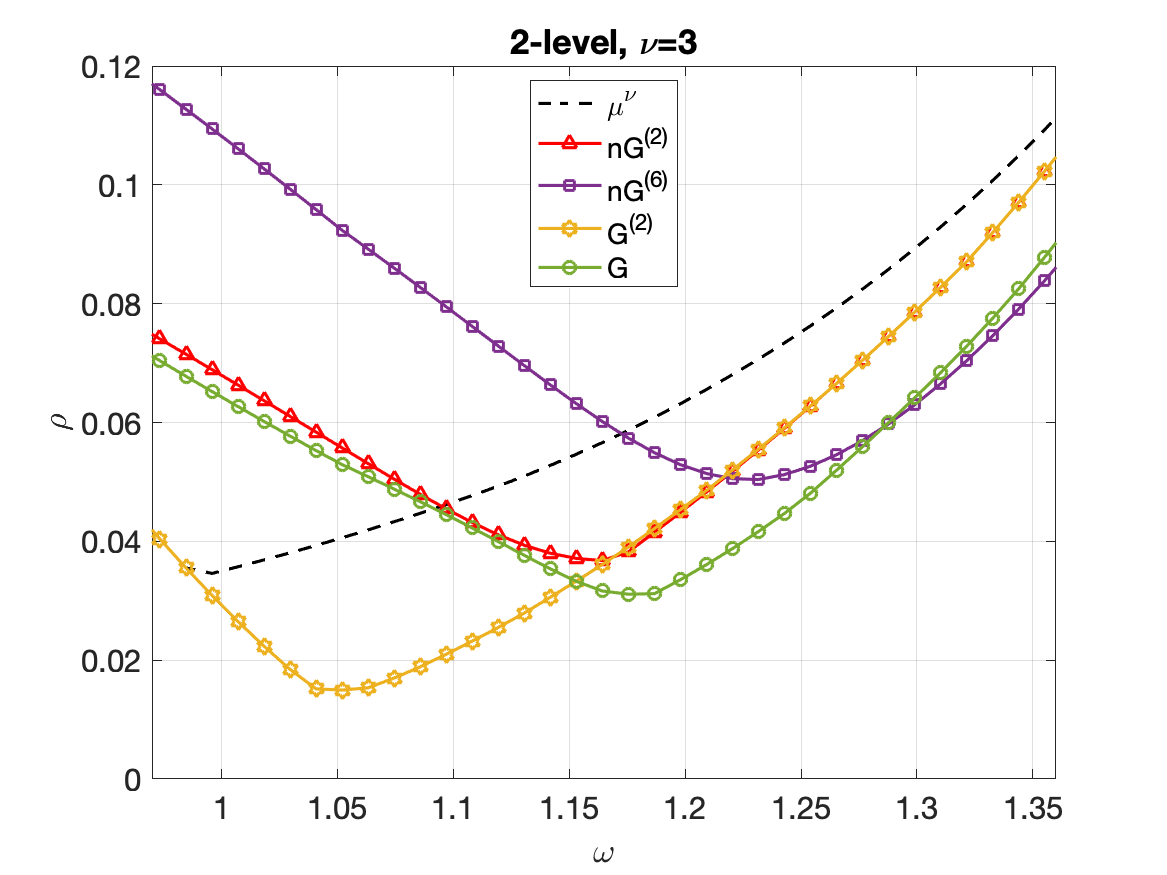}{\figWidth}};
   \draw(7.5,0.0) node[anchor=south west,xshift=-15pt,yshift=-8pt]{\trimfig{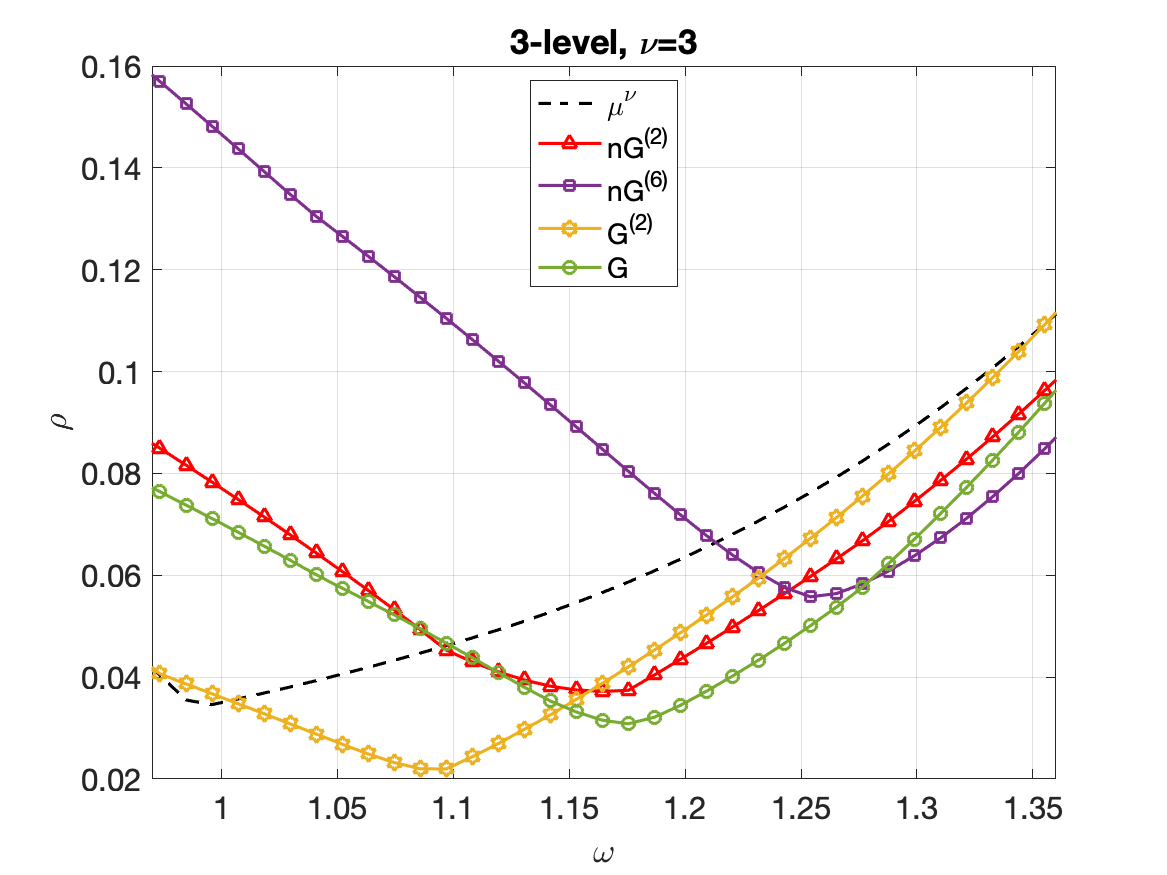}{\figWidth}};
   
  %  circles to highlight best results
   \circleLeftLabel{2.4}{1.1}{opt CR};
   \circleLeftLabel{10.6}{.7}{opt CR};

   \circleLabel{1.43}{1.9}{opt $\mu^\nu$}{violet}{.25};
   \circleLabelLeft{9}{1.}{opt $\mu^\nu$}{violet}{.25};

   % Combined cartoons
   \orderSixVcycleTwoLevel{xshift=1.7cm,yshift=4.75cm,scale=1}{$\nu\!=\!3$}
   \begin{scope}[xshift=0cm,yshift=6cm]
     \orderSixVcycleTwoLevel{xshift=1.7cm,yshift=4.75cm,scale=1}{$\nu\!=\!2$}
   \end{scope}
     
   \begin{scope}[xshift=7.5cm]
     \orderSixVcycleThreeLevel{xshift=1.7cm,yshift=4.75cm,scale=1}{$\nu\!=\!3$}
   \end{scope}
   \begin{scope}[xshift=7.5cm,yshift=6cm]
     \orderSixVcycleThreeLevel{xshift=1.7cm,yshift=4.75cm,scale=1}{$\nu\!=\!2$}
   \end{scope}

   % title
   \draw(7.5,11.7)  node[draw=\colourTwoD,very thick,fill=\colourTwoD!20,anchor=south,inner sep=2.5pt] {\small \DIM, \ORD[6], \Cc}; 
   
 % grid:
%  \draw[step=1cm,gray] (0,0) grid (15,12.5);
\end{tikzpicture}
\end{center}
\caption{$\rho$ versus $\omega$. \DIM, \ORD[6], 2-level and 3-level V cycles with red-black Gauss-Seidel smoothing ($\nu$ sweeps per cycle) and \CC,
\kl{
	2nd-order accurate transfer operators;
	6th-order accurate non-Galerkin (nG$^{(6)}=\LNH{6} \equiv \LNAH$),
	6th-order accurate Galerkin (G=$\LGGH \equiv \LGHh{6}{2}$),
	2nd-order accurate non-Galerkin (nG$^{(2)}=\LNH{2}$),
	and
	2nd-order accurate Galerkin (G$^{(2)}= \LGqH{2} \equiv \LGHh{2}{2}$) coarse-level operators.
}
}
  \label{fig:rho_V_2d_6o_std}
\end{figure}
}
\subsection{Sixth-order accurate discretizations with second-order accurate coarse-level operators} \label{Sec: orderC_6o}
In this section we present results for sixth-order accurate fine-level discretizations.
%and using second-order accurate coarse-level operators. 
We wish to see if second-order accurate coarse-level operators can also be used even for sixth-order accurate discretizations without any significant effect on the convergence rates. 

We consider the standard sixth-order accurate finite difference stencil 
\begin{equation}
L_h = \frac{1}{180h^2}\begin{bmatrix}
   &	&	   & -2 \\
   &	&	   &  27  \\
   &    &      & -270 \\
-2 & 27 & -270 &  980 & -270 & 27 -2\\
   &    &      & -270 \\
   &	&	   &  27  \\
   &	&	   & -2    
\end{bmatrix}_h  
\label{eq:LhOrder6}
\end{equation} 
for the two-dimensional (negative) Laplacian.
We use the second-order accurate transfer operators~\eqref{transfer_2d}.
On the coarse levels, we consider Galerkin and non-Galerkin operators constructed from
the second-order accurate fine-level operator~\eqref{L_h_2_bra},
in comparison with coarse-level operators constructed
from the actual (sixth-order accurate) fine-level operator~\eqref{eq:LhOrder6}.
Convergence rates 
%(as a function of the red-black smoothing parameter $\omega$) 
are given in Figure~\ref{fig:rho_V_2d_6o_std} for two-level and three-level $V[1,1]$ and $V[2,1]$ cycles.
These results, computed using local Fourier analysis,
are similar to those obtained for the fourth-order accurate fine-level operator (Figure~\ref{fig:rho_V_2d_4o_std}). 
Overall the best convergence rates are obtained using the second-order accurate Galerkin coarse-level operators. 

The results from this and the previous sections 
show that second-order accurate coarse-level operators are effective for
both fourth- and sixth-order accurate fine-level discretizations. 
The same conclusion will likely hold for eighth- and higher-order accurate discretizations.

% topic-2: RED BLACK COARSEING ----------------------------------------------------------------------
\bigskip
\section{Red-black coarsening and red-black smoothing} \label{M2: rbC}

Since a majority of grid points on an overset grid often belong to Cartesian background grids,
we are interested in finding the very best multigrid algorithms for Cartesian grids. Moreover, there
are many important application areas that require the solution to Poisson's equation (e.g. incompressible fluid flow and electromagnetics).
Therefore there is a strong incentive to find the best multigrid algorithms for solving Poisson's equation on Cartesian grids.
It is well known~\cite{MGtext} that, in one dimension, a particular multigrid algorithm
based on the red-black smoother becomes a direct solver for a second-order order accurate
discretization to Poisson's equation. This special multigrid algorithm converges in one iteration after
a post-smoothing step.
The question now arises as to whether this one-dimensional result can be extended to two or three space dimensions, perhaps using a nonstandard coarsening strategy. 
In this section we show that there is a special two-level multigrid algorithm in two dimensions based on the red-black smoother and \textit{red-black coarsening} that becomes a direct solver for the standard second-order accurate discretization to Poisson's equation.

{
\newcommand{\figWidth}{6cm}% height 
\newcommand{\trimfig}[2]{\trimh{#1}{#2}{.0}{.0}{.0}{.0}}
\begin{figure}[h!]
\begin{center}
\begin{tikzpicture}[scale=1]
  \useasboundingbox (0,0.4) rectangle (15,12);  % set the bounding box (so we have less surrounding white space)

   \begin{scope}[yshift=6cm]
     \draw(0.0,0.0) node[anchor=south west,xshift=-15pt,yshift=-8pt] {\trimfig{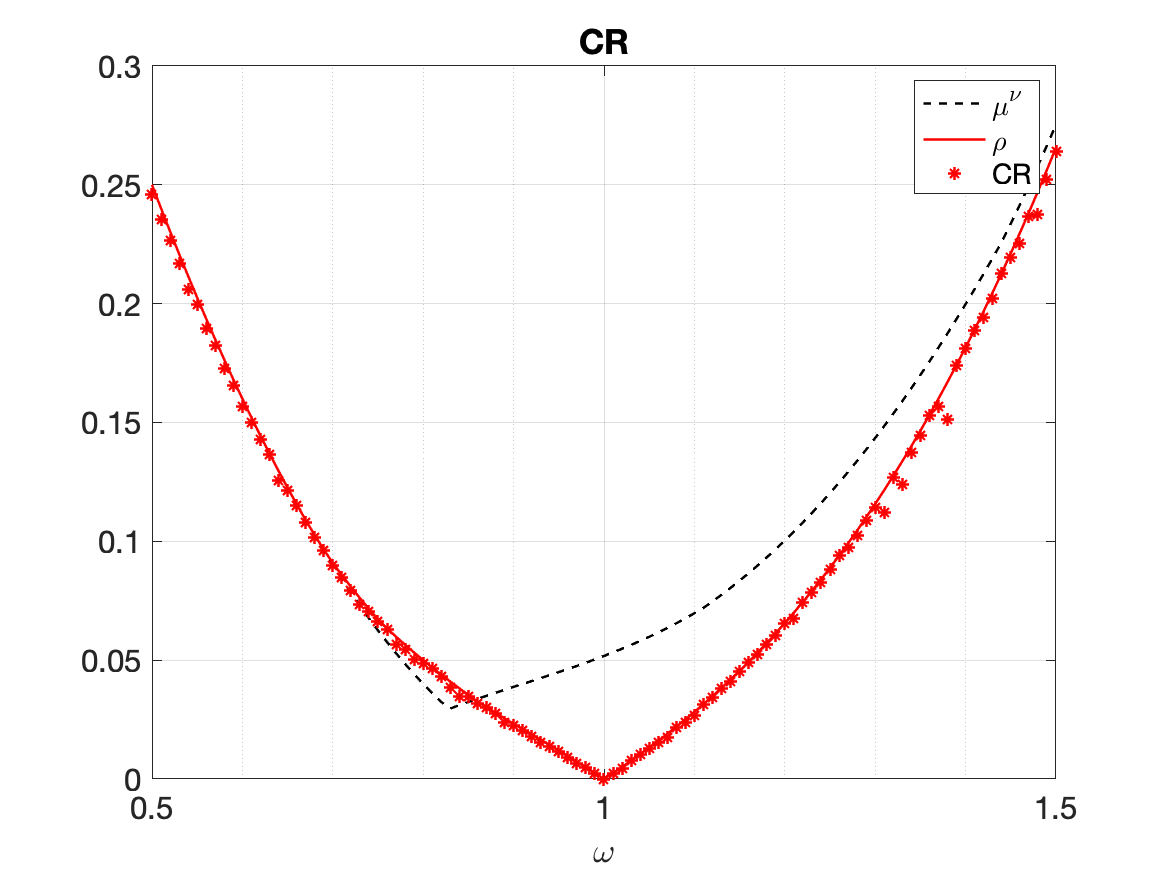}{\figWidth}};
     \draw(7.5,0.0) node[anchor=south west,xshift=-15pt,yshift=-8pt] {\trimfig{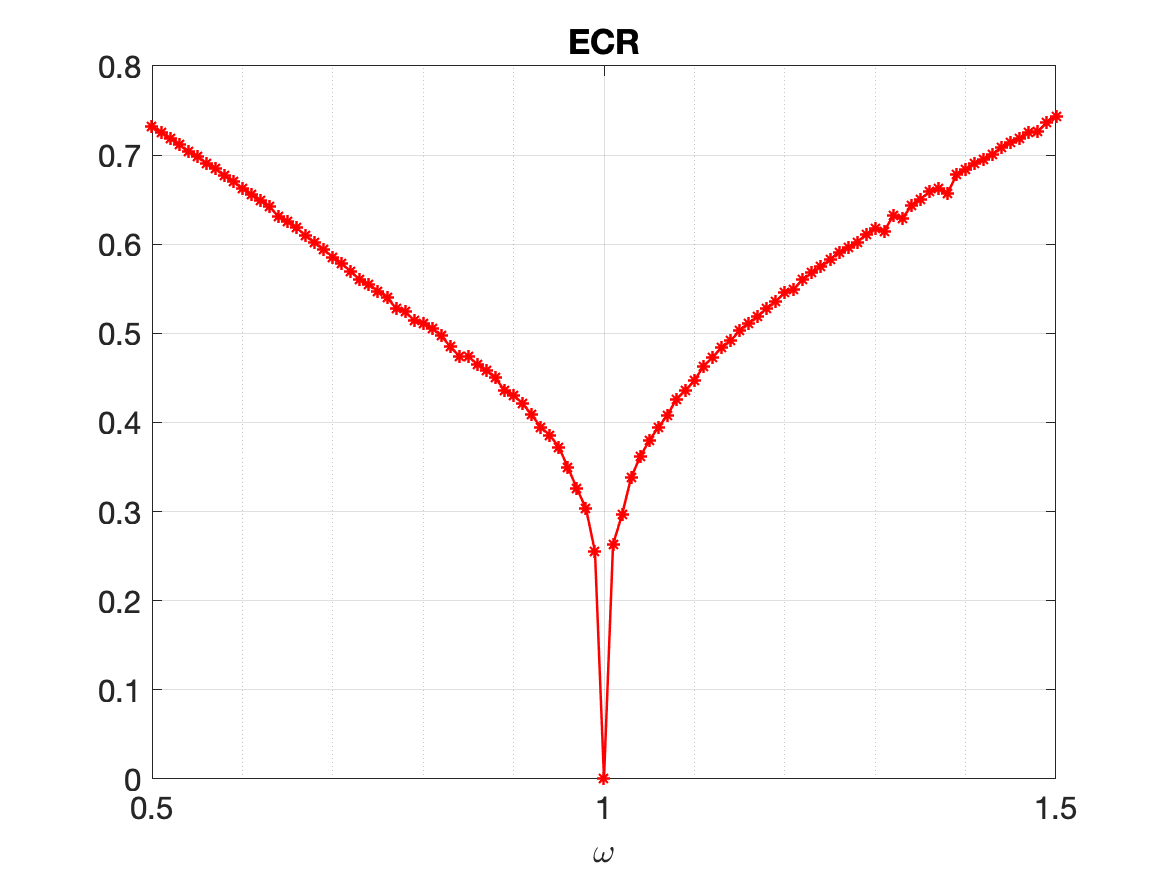}{\figWidth}};

     \circleLabel{11.675}{.7}{direct solver, $\omega\!=\!1$}{DarkGreen}{.25};

     \circleLabelLeft{3.1}{1.05}{opt $\mu^\nu$}{violet}{.25};

     \draw(7.5,3.5)  node[draw=orange,very thick,fill=orange!20,anchor=south] {\small Red-Black smoothing, \DIM[1], \LEV, \Cc}; 

     \VcycleTwoLevelBoxed{xshift=7.5cm,yshift=2.5cm,scale=1}{$\nu\!=\!2$}
   \end{scope}
   
   \draw(0.0,0.0) node[anchor=south west,xshift=-15pt,yshift=-8pt] {\trimfig{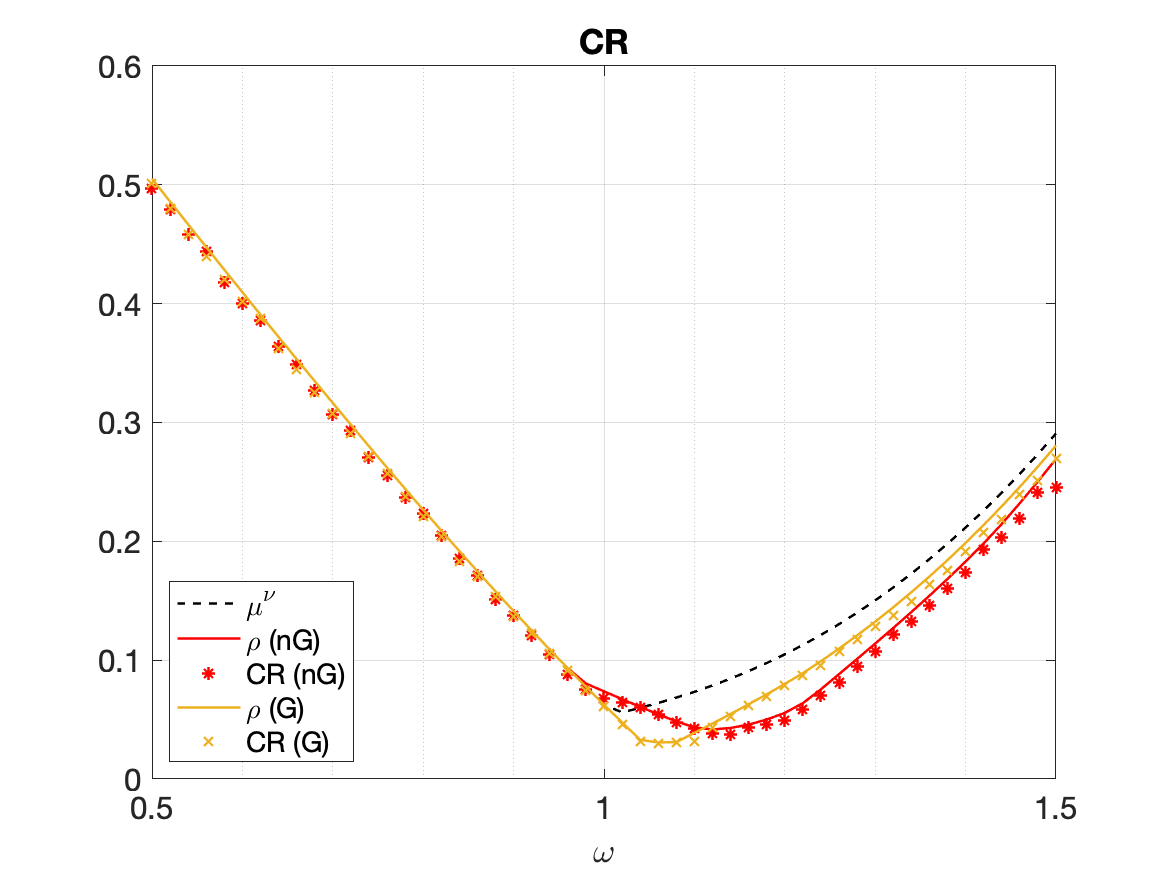}{\figWidth}};
   \draw(7.5,0.0) node[anchor=south west,xshift=-15pt,yshift=-8pt] {\trimfig{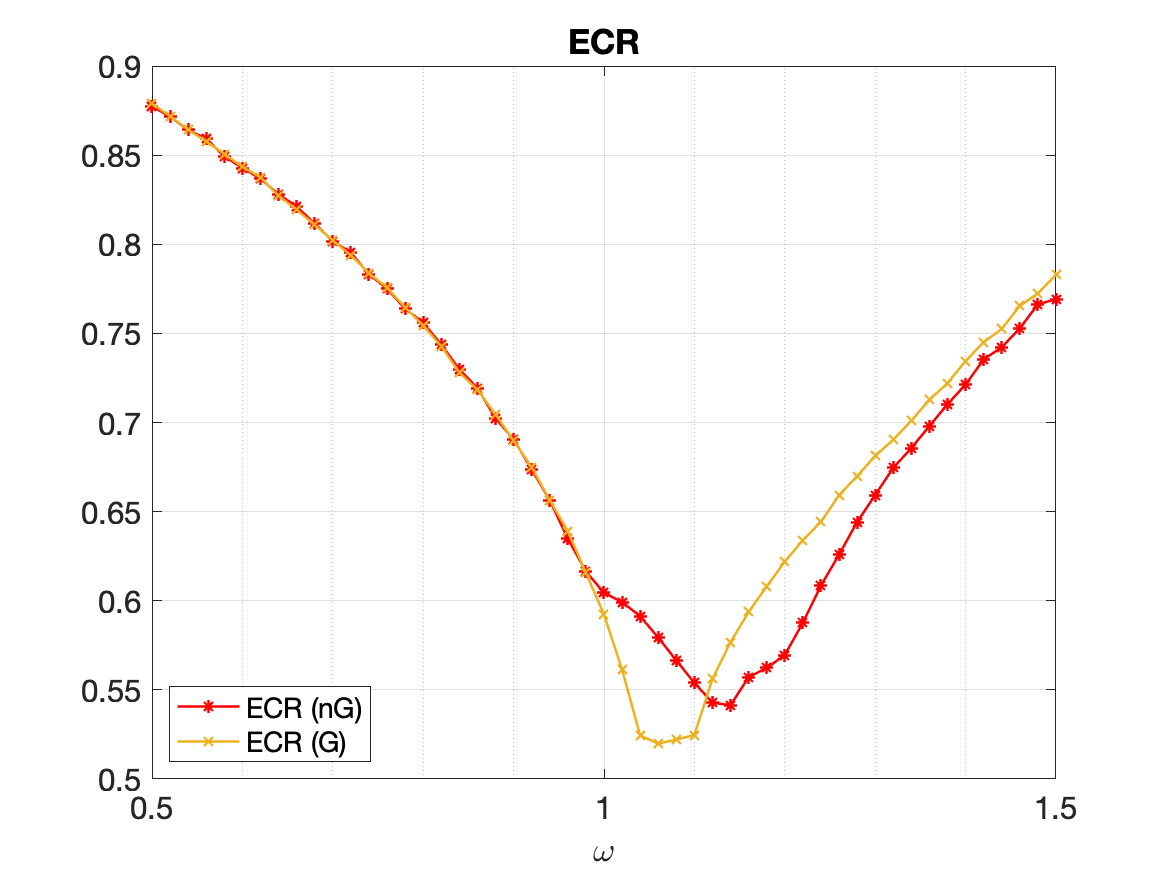}{\figWidth}};

   \draw(7.75,4.75)  node[draw=\colourTwoD,very thick,fill=\colourTwoD!20,anchor=south,yshift=-4pt] {\small Red-Black smoothing, \DIM, \LEV, \Cc}; 

   \circleLabel{12.1}{1.}{optimal $\omega>1$}{DarkGreen}{.35};
   \circleLabelLeft{4.3}{1.15}{opt $\mu^\nu$}{violet}{.25};

   %  label:
   % \circleVortex{11.675}{.8}{direct solver for $\omega=1$};
   % Combined Cartoons
   \orderTwoVcycleTwoLevel{xshift=7.2cm,yshift=3.2cm,scale=1}{$\nu\!=\!2$}

   % \draw[<-,thick,xshift=1pt,yshift=2pt](11.7,.5) -- (12.5,1.6) node[anchor=west,xshift=-2pt,yshift=2pt] {\scriptsize exact solver}; 
 
% grid:
%  \draw[step=1cm,gray] (0,0) grid (15,12);
\end{tikzpicture}
\end{center}
\caption{CR and ECR versus $\omega$, comparing one and two dimensions;
	\ORD, \LEV, $V[1, 1]$
	cycle with red-black smoothing and \CC ($N = 32$).
  Top: \DIM[1];
  bottom: \DIM, 
  %non-Galerkin (nG=$\ngcg{2}$) and Galerkin (G=$\gcg{2}{2}{2}$) 
  %non-Galerkin (nG=$\LNAH \equiv \LNH{2}$) and Galerkin (G=$\LGGH \equiv \LGHh{2}{2}$)
  with non-Galerkin (nG=$\LNAH$) and Galerkin (G=$\LGGH$) 
  coarse-level operators. 
  %{\red Explain solid-lines versus dots here and in text -- do the same for other figures}
}
  \label{fig:redBlackSmootherStandardCoarsening1d2d}
\end{figure}
}

We consider the solution of Poisson's equation~\eqref{MP_h} in $d=2$ space dimensions,  discretized with the
second-order accurate five-point star~\eqref{L_h_2}.
We choose our smoother to be red-black with a relaxation parameter $\omega$.
Figure~\ref{fig:redBlackSmootherStandardCoarsening1d2d} shows two-level results in one and two dimensions using standard coarsening and a V[1,1] cycle. 
In this and other figures, the computed CRs are plotted along with
the theoretical asymptotic convergence rates $\rho$ (solid lines)
to confirm the consistency of the computations with local Fourier analysis.
The transfer operators are taken as linear interpolation and full weighting restriction.
The top figures show the convergence rate (CR) and effective convergence rate (ECR) in
one dimension. The CR (and ECR) is zero at $\omega=1$ where the multigrid algorithm converges in one cycle.
This rate is much better than that predicted by the smoothing rate $\mu^\nu$ with an ideal coarse-level correction. The coarse-level correction operator thus interacts in a special way with the red-black smoother so as to achieve this ideal convergence rate. 
It is noted that the ECR increases rapidly as $\omega$ moves away from the optimal value $\omega=1$. 

Figure~\ref{fig:redBlackSmootherStandardCoarsening1d2d} (bottom) shows results in two dimensions using
Galerkin and non-Galerkin coarse-level operators (note that in one dimension the Galerkin and non-Galerkin coarse-level operators are identical). 
The CR (and ECR) in two dimensions are still very good but there is no $\omega$ for which the rates are zero. The Galerkin coarse-level operator gives nearly the same convergence rates as the non-Galerkin one.
It is noted that the optimal CR occurs with an over-relaxed smoother ($\omega >1$) and that the
optimal CR is close to the smoothing rate $\mu^\nu$. Our intuition is that an $\omega >1$ is required for optimal CR since the red-black smoother and the coarse-level correction operator are not perfectly balanced (unlike the one-dimensional case where $\omega=1$ gives optimal convergence),
and thus the properties of the smoother in two dimensions need to be adjusted to obtain an optimal CR. 

It is therefore seen that, unlike the case in one dimension, in two dimensions the coarse-level correction does not conspire with the red-black smoother to give a direct solver.
It is natural to ask if there is some other coarse-level correction operator in two or more space dimensions that does lead to a direct solver.
To answer this question we first review the theory for the one-dimensional case in Section~\ref{sec:redBlack1d},
and then move on to two dimensions 
in Section~\ref{Sec: rb_2d_std}, where we analyze the discrepancy between the red-black smoother and standard coarsening in two dimensions.
In Section~\ref{Sec: rbC} we introduce red-black coarsening in two dimensions, and recover the direct-solver property for a two-level cycle.
Extensions to more than two levels are provided in Section~\ref{Sec: rbC_multilevel}, where we show fast convergence with LFA and numerical results.
Remarks on the three-dimensional case are made in Section~\ref{Sec: rbC_3d}.
The details of the local Fourier analysis will be left to \ref{Sec: theory}, where we also consider a few generalizations to the operator (\ref{L_h_2}).

%----------------------------------------------------------------------------------------------------------------------
%\subsection{\old{Multigrid with red-black smoothing and standard coarsening}}\label{Sec: rbC_problem}
% motivation: rbC
\subsection{Standard coarsening and red-black smoothing in one dimension}\label{sec:redBlack1d}
We begin by reviewing the one-dimensional case and showing a well-known result~\cite{MGtext} 
summarized in the following theorem. 
% -------------------------------- THEOREM -------------------------------------
\begin{theorem}[Red-black smoothing and standard coarsening in one dimension]\label{Th_1d}
  Consider the one-dimensional discrete Poisson's equation $L_h u_h = f_h$ with 	
  \begin{equation}\label{L_h_1}
    L_h = \frac{1}{h^2} \begin{bmatrix}
      -1 & 2 & -1
    \end{bmatrix}_h
  \end{equation}
  on the infinite grid (\ref{G_h_d}) with $d=1$.
  An $(h,\ H)$ two-level multigrid cycle consisting of a red-black smoother ($\omega=1$)
  and standard coarsening ((\ref{G_H_d}), $H=2h$), with the transfer operators
  \begin{equation}\label{transfer_1d}
    I_h^H = \frac{1}{4}\begin{bmatrix}
      1 & 2 & 1
    \end{bmatrix}_h^{2h},\quad
    I_H^h =  \frac{1}{2}\left]\begin{matrix}
	1 & 2 & 1
      \end{matrix}\right[_{2h}^h,
  \end{equation}
  as well as a coarse-level operator with the same stencil as~\eqref{L_h_1}, 
  gives rise to a direct solver. That is, the cycle converges in a single step (with post-smoothing).
\end{theorem}  
  Note that since the coarse-level operator is the same as the fine-level, it follows that the direct-solver
  property holds for any number of levels.
\begin{Corollary}
   When extended to more than two levels, the algorithm in Theorem~\ref{Th_1d} remains a direct solver.
\end{Corollary}
This result can be shown through local Fourier analysis.
\edit{
We review the key points here and leave some of the details to \ref{Sec: theory}, where we also show the equivalence of this multigrid cycle to cyclic reduction in \ref{Sec: cyclic_1d}.}

\begin{figure}[h!]
	\begin{center}
		\begin{tikzpicture}[scale=.9]
		\useasboundingbox (1, -.1) rectangle (7,.8); 
		\draw (0,0) -- (8,0);
		\draw[thick, red]  (2,0) -- (6,0);
		\draw[fill] (0,0) circle (1pt) node[align=left,   below] {$-\pi$};
		\draw[red,fill=red] (2,0) circle (1pt) node[align=left, red, below] {$-\frac{\pi}{2}$};
		\draw[fill] (4,0) circle (1pt) node[align=center, below] {$0$} 
		node[align=center, red, above=5pt] {$\Theta^{\text{low}}$};
		\draw[red,fill=red] (6,0) circle (1pt) node[align=left, red, below] {$\frac{\pi}{2}$};
		\draw[fill] (8,0) circle (1pt) node[align=right, below] {$\pi$};
		\end{tikzpicture} 		
	\end{center}
	\caption{$\Theta$, 1D, 2-level}
	%	\vspace{-10pt}
	\label{fig: Theta_1d}
\end{figure}
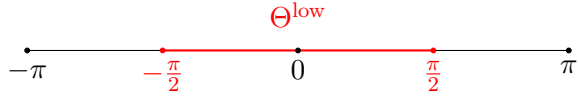

As discussed in Section \ref{Sec: defs}, the operator $L_h$ has the discrete Fourier modes
\begin{equation}
  \phi_h(x,\ \th) = e^{i \th \frac{x}{h}},\ x \in \grid G_h,\quad 
  \th \in \Theta \equiv [-\pi,\ \pi)
\end{equation}
as its eigenfunctions, with corresponding eigenvalues 
\ba
\symF L_h({\th}) = \frac{4}{h^2} \xi_{\th},\quad \xi_{\th} \eqdef \sin^2 \frac{\th}{2}.
\ea
$\symF L_h({\th})$ is also called the Fourier \textit{symbol} of $L_h$.
By standard coarsening with $H = 2h$, on the coarse grid (\ref{G_H_d}) we have Fourier modes
\begin{equation}
\phi_H(x,\ 2\th) = e^{i 2\th \frac{x}{H}},\quad x \in \grid G_H,\quad 
\th \in [-\frac{\pi}{2},\ \frac{\pi}{2}).
\end{equation}
According to whether the Fourier modes $\phi_h(\cdot,\ \th)$ on the fine grid $\grid G_h$ can be represented on the coarse grid $\grid G_H$,
they can be separated into low- and high-frequency modes, for which
$\Theta^{low} \equiv [-\frac{\pi}{2},\ \frac{\pi}{2})$, and $\Theta^{high} \equiv \Theta \setminus \Theta^{low}$.
An illustration is shown in Figure \ref{fig: Theta_1d}.
Note that this separation gives an equal amount of low-frequency and high-frequency modes.
Define the two-dimensional subspace $\Etheta$, parameterized by $\th \in \Theta^{low}$, to be that spanned by Fourier modes $\phi_h (\cdot,\ \th)$ and $\phi_h (\cdot,\ \bar{\th})$:
\begin{equation}\label{E_1d}
  \Etheta = \text{span} \{\phi_h (\cdot,\ \th),\quad \phi_h (\cdot,\ \bar{\th})\},
\end{equation}
where 
\begin{align}\label{th_bar}
  \bar{\th} \eqdef \th - \sgn(\th) \pi = \th + \pi\!\!\! \mod{\Theta}.
\end{align}
\cut{
The modes in $\Etheta$ are indistinguishable on the coarse grid, in particular,
\begin{equation}
  \phi_h(x,\ \th) = \phi_h(x,\ \bar{\th}) = \phi_H(x,\ 2\th),\quad \forall\ x \in \grid G_H,\quad \th \in \Theta^{low},
\end{equation}
which comes from the simple fact that
\begin{equation}\label{eqn_minus}
  \phi_h(jh,\ \bar{\th}) = (-1)^j \phi_h(jh,\ \th),\quad j \in \Integer,
\end{equation}
as illustrated in Figure \ref{fig:sine}.
\begin{figure}[h!]
	\begin{center}
		\begin{tikzpicture}
  \useasboundingbox (0,-.5) rectangle (7,1);  % set the bounding box (so we have less surrounding white space)
  \draw[->,black] (0,0) -- (7,0) node[right] {$x$};
  % \draw[->] (0,-3) -- (0,4.2) node[above] {$y$};
  % \draw[scale=0.5,domain=-3:3,smooth,variable=\x,blue] plot ({\x},{\x*\x});
  % \draw[scale=0.5,domain=-3:3,smooth,variable=\y,red]  plot ({\y*\y},{\y});

  \begin{scope}[yshift=0cm]
    % plot sin(2*x)
    \draw[scale=1,domain=0:2*pi,smooth,variable=\x,blue,very thick]  plot ({\x},{sin(2*(180/pi)*\x)});
    \draw[scale=1,domain=0:2*pi,smooth,variable=\x,blue,densely dotted, thick]  plot ({\x},{-sin(2*(180/pi)*\x)});
    % plot (-1)^j * sin(2*x_j)
    % (-1)^j = cos( (N/2)*x_j ) 
    %N=10, samples=N+1
    \newcommand*{\N}{10}
    %\draw[scale=1,domain=0:2*pi,samples=11,variable=\x,orange,very thick]  plot ({\x},{sin(2*(180/pi)*\x)*cos(5*(180/pi)*\x)});
    \draw[scale=1,domain=0:2*pi,smooth, variable=\x,orange,very thick]  plot ({\x},{sin((2-\N/2)*(180/pi)*\x)});
    
    % grid points:
    \newcommand*{\h}{2*pi/\N}
    \foreach \i in {0,...,\N} {
    	\draw[fill] (\h*\i, 0) circle (1pt);
    	\draw[blue, fill=blue] (\h*\i, {sin(2*(180/pi)*\h*\i)}) circle (1.5pt);
    	\draw[orange, fill=orange] (\h*\i, {sin((2-\N/2)*(180/pi)*\h*\i)}) circle (1.5pt);
    }
  \end{scope}
%  draw grid temporarily while making figure 
 %\draw[step=1cm,gray] (0,-1) grid (7,1);
% 
\end{tikzpicture}
	\end{center}
	\caption{Representative modes of $\Etheta$. $\Im \phi_h(x,\ \theta)=\sin {\frac{\th}{h} x}$ (\textcolor{blue}{solid blue}) and $\Im \phi_h(x,\ \thetaBar)= \sin{\frac{\th-\pi}{h} x}$ (\textcolor{orange}{orange}) satisfy (\ref{eqn_minus}) at grid points $x = jh$. $\theta=2 h$.}
	\label{fig:sine}
\end{figure}
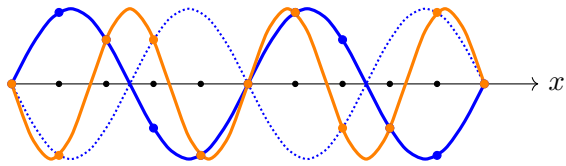

}

The key point we shall establish in the following is that both the red-black smoothing operator $S_h$ and the coarse-level correction operator $K_h^H$ by standard coarsening ($H=2h$) are invariant on the same two-dimensional subspaces $\mathbb{E}_{\th}$.
The Red-black smoother separates the grid points in $\grid G_h$ into red and black points as
\begin{equation}
\begin{cases}
\grid G_h^R = \Big\{  x =j h:~ j \in \Integer_1 \Big\}, \\
\grid G_h^B = \Big\{  x =j h:~ j \in \Integer_0 \Big\}.
\end{cases}
\end{equation}
Note that the coarse grid with standard coarsening coincides exactly with the black points in the fine grid, that is,
\begin{equation}
\grid G_H \equiv \grid G_h^B.
\end{equation}
The red-black smoothing operator $S_h$ is given by (\ref{S_rb}), and
for the second-order accurate operator~\eqref{L_h_1} (with the three-point stencil) a Gauss-Seidel step is equivalent to a Jacobi step with 
\begin{gather}
S_h^{GS} \equiv S_h^J = I_h - \frac{\omega h^2}{2} L_h.
\end{gather}
It can be shown (see \ref{Sec: rb_1d}) that $S_h$ is invariant on each subspace $\Etheta$, and on this subspace can be represented by the $2\times 2$ \textit{matrix symbol}
%\old{with the matrix representation}
\begin{equation}\label{S_hat_rb}
\symE S(\th) = \symE S^B (\th)\ \symE S^R(\th),
\end{equation}
where
\begin{equation}\label{S_RB}
  \begin{aligned}
  \symE S^R(\th) 
  \commentB{&= 
  \frac{1}{2}\begin{bmatrix}
  1+\symF S^J_h({\th}) & 1-\symF S^J_h(\bar{\th}) \\
  1-\symF S^J_h({\th}) & 1+\symF S^J_h(\bar{\th})
  \end{bmatrix}}
  = \mathbf{I} - \omega \begin{bmatrix}
  1 & -1\\ -1 & 1
  \end{bmatrix} \begin{bmatrix}
  \xi_{\th} \\ & 1-  \xi_{\th}
  \end{bmatrix},\quad
  \symE S^B (\th)  
  \commentB{&=
  \frac{1}{2}\begin{bmatrix}
  1+\symF S^J_h({\th}) & -(1-\symF S^J_h(\bar{\th})) \\
  -(1-\symF S^J_h({\th})) & 1+\symF S^J_h(\bar{\th})
  \end{bmatrix}}
  = \mathbf{I} - \omega \begin{bmatrix}
  1 & 1\\ 1 & 1
  \end{bmatrix} \begin{bmatrix}
  \xi_{\th} \\ & 1-  \xi_{\th}
  \end{bmatrix}.
  \end{aligned}
\end{equation}

With the full-weighting restriction and linear interpolation operators between $\grid G_h$ and $\grid G_H$ given by (\ref{transfer_1d}), for each $\th \in \Theta^{low}$, we have, for $x \in \grid G_H$,
\begin{equation}
  \begin{cases}
  I_h^{H} \begin{bmatrix}
  \phi_h (x,\ \th) & \phi_h (x,\ \bar \th) 
  \end{bmatrix}  
  = \phi_{H}(x,\ 2\th)\ \symE I_h^{H}(\th),\\
  I^h_{H}\ \phi_{H}(x,\ 2\th) =  \begin{bmatrix}
  \phi_h (x,\ \th) & \phi_h (x,\ \bar \th) 
  \end{bmatrix}  \symE I^h_{H} (\th),
  \end{cases}
\end{equation}
where 
\begin{equation}\label{transfer_hat_1d}
   \symE I_h^{H}(\th) = \begin{bmatrix}
  1-\xi_{\th} & \xi_{\th}
  \end{bmatrix},\quad \symE I^h_{H}(\th) = \symE I_h^{H}(\th)^T
\end{equation}
\cut{are matrix representations of the restriction and interpolation operators.}
%{(\blue this seems ok to be -- they are matrix representations aren't they?)} 
% \kP{(I want to avoid this, since the matrices are represented for transfers between $\Etheta$ and $\phi_H$, rather than on some invariant subspace ($\Etheta$) like $\mathbf{S}$ and $\mathbf{K}$; if you want to call them specifically matrix representations, we need to state more clear the domain and codomain spaces for the transfer maps.)}}
% 
On the coarse grid $\grid G_H$, we consider the coarse-level operator 
\ba
L_H = \frac{1}{H^2} \begin{bmatrix}
-1 & 2 & -1
\end{bmatrix}_H 
\ea
with the same stencil as the fine-grid operator $L_h$. We have 
%$\phi_H (\cdot,\ 2\th)$ being an eigenfunction of $L_H$,
\ba
  L_H\ \phi_H (\cdot,\ 2\th) = \symF L_H(2\th)\ \phi_H (\cdot,\ 2\th),
\ea
with symbol
\ba
  \symF L_H(2\th) = \frac{4}{H^2} \sin^2 \th.
\ea
Note that in one dimension, the (non-Galerkin) operator $L_H$ has the special property that it
the same as the Galerkin operator~\eqref{L_H_G},
or equivalently,
\begin{equation}
   \symF L_H(2\th) = \symE I_h^H(\th)\ \symE L(\th)\ \symE I^h_H (\th) =\frac{4}{h^2} \xi_{\th} (1-\xi_{\th}),
\end{equation}
where $\symE L(\th)$ is the matrix symbol of $L_h$ on $\mathbb{E}_{\th}$ given by
\begin{equation}
  \symE L(\th) =
   \begin{bmatrix}
    \symF L_h(\th)\\ & \symF L_h(\bar{\th})
  \end{bmatrix} = \frac{4}{h^2}
  \begin{bmatrix}
    \xi_{\th}\\ & 1-\xi_{\th}
  \end{bmatrix}.
\end{equation}
In general, however, the Galerkin and non-Galerkin coarse-level operators are not the same.
The coarse-level operator along with the transfer operators gives the coarse-level correction operator
$
    K_h^H = I_h - I^h_H L^{-1}_H I_h^H L_h.
$
Thus we can derive its matrix representation on each $\Etheta$ ($\th \neq 0$) to be
\begin{equation}\label{K_hat_1d}
  \symE K (\th) = \mathbf{I} - \symF L_H^{-1}(2\th)\
  \symE I^h_H (\th)\ \symE I_h^H(\th)\ \symE L(\th) 
  = \begin{bmatrix}
  1 \\ -1
  \end{bmatrix}\begin{bmatrix}
  \xi_{\th} & -(1-\xi_{\th})
  \end{bmatrix}.
\end{equation}
\commentB{
With the Galerkin property of $L_H$, $K_h^H$ is a projector ($(K_h^H)^2=K_h^H$). The coarse-level correction reduces the current error $v_h$ as $K_h^H v_h = v_h - \tilde{v}_h$, where the estimated error 
\begin{equation}
\tilde{v}_h = I^h_H L^{-1}_H I_h^H L_h v_h
\end{equation}
is the projection of the error along the range of the interpolation operator $I_H^h$.
}

We can now present the following theorem, which shows that in one space dimension the two-level multigrid algorithm of Theorem~\ref{Th_1d} (with post-smoothing $\nu_2>0$)
%, consisting of a coarse-level correction and one post-smoothing sweep, 
is indeed a direct solver.
\begin{theorem}[Red-black smoothing and standard coarsening in one dimension - LFA]\label{Th_1d_LFA}
  The $(h,\ 2h)$ two-level iteration operator $M_h^H$ of a $V[\nu_1,\ \nu_2]$ cycle as described in Theorem \ref{Th_1d}
  has its matrix representation on each eigen-subspace $\mathbb{E}_{\th}$ (\ref{E_1d}) ($\th \neq 0$)  given by
  \begin{equation}
    \symE M(\th) = \symE S^{\nu_2}(\th) \ \symE K(\th)\ \symE S^{\nu_1}(\th),
  \end{equation}
  where $\symE S(\th)$ is given by (\ref{S_hat_rb}), and $\symE K(\th)$ is given by (\ref{K_hat_1d}).
  Thus we have the asymptotic convergence rate of the two-level cycle as
  \begin{equation}\label{rho_1d}
    \rho = \rho \big(M_h^H\big)
    = \sup_{\th \in \Theta^{\text{low}}} { \rho \big(\symE M(\th)\big)}.
  \end{equation}
\end{theorem}  
%\noindent 
Therefore in the special case of $\omega=1$ we have 
%{\red note: I shifted the row vector up -- maybe fix other instances if we like this.}
\begin{equation}\label{eqn: zero_1d}
  \symE S^R(\th)\ \symE K (\th) =
        \begin{bmatrix}
  1-\xi_{\th} \\ \xi_{\th}
  \end{bmatrix}
  \begin{matrix} 
    \begin{bmatrix}
      1 & 1
    \end{bmatrix}
      \\ \mbox{}
  \end{matrix} 
  \begin{bmatrix}
      1 \\ -1
  \end{bmatrix}
  \begin{matrix} 
    \begin{bmatrix}
      \xi_{\th} & -(1-\xi_{\th})
    \end{bmatrix}
  \\ \mbox{}
  \end{matrix} 
  =
  \begin{bmatrix}
    0 & 0 \\
    0 & 0 
  \end{bmatrix}
  =   \mathbf{O},
\end{equation}
that is, the method always converges after one post-smoothing with $\omega=1$, establishing Theorem~\ref{Th_1d}.

%In fact, $\symE S^R = \symE I - \symE K$, so that $\symE S^R \symE K = \symE K - \symE K^2 = \mathbf{O}$.

% ===============================================================================================================================
\subsection{Standard coarsening and red-black smoothing in two dimensions}\label{Sec: rb_2d_std}
In two dimensions, we consider the solution to Poisson's equation~\eqref{MP_h} on $\grid G_h$ given by~\eqref{G_h_d} with $d=2$,
with the operator $L_h$ being the standard five-point stencil~\eqref{L_h_2} of second-order accuracy.
\commentB{
\ba
L_h = 
 \frac{1}{h^2}
\begin{bmatrix}
  & -1 \\
  -1 & 4 & -1\\
  & -1
\end{bmatrix}_h . 
\ea
}
We present some results of local Fourier analysis for multigrid with red-black smoothing and standard coarsening in two dimensions \edit{(for further details see \ref{Sec: rb_2d})}.
We discuss how standard coarsening leads to a misalignment between the coarse-level correction and the red-black smoother, in terms of their eigen-subspaces,
in contrast to the one-dimensional case studied in Section~\ref{sec:redBlack1d}. 
In the subsequent Section~\ref{Sec: rbC} we consider an alternative coarsening, red-black coarsening, which
more naturally aligns with the red-black smoother. With red-black coarsening,
the smoothing and coarse-level correction operators are invariant on the same two-dimensional subspaces, 
and in this case a two-level direct solver can be constructed with an appropriate choice of transfer and coarse-level operators.

As discussed in Section~\ref{Sec: defs}, the operator $L_h$ has Fourier modes 
\begin{equation}
\phi_h (\mathbf{x},\ \thv) = e^{i \boldsymbol{\theta} \cdot \frac{\xv}{h}},\quad \xv \in \grid G_h,\quad \boldsymbol{\theta} = [\theta_1,\  \theta_2]^T \in \Theta \eqdef [-\pi,\ \pi)^2,
\end{equation}
as eigenfunctions, with corresponding symbol
\begin{equation}
\symF L_h (\thv) = \frac{8}{h^2} \xi_{\thv}, 
\end{equation}
where we define
\begin{equation}\label{xi_r2}
\xi_{\thv} \eqdef \frac{1}{2} \left(\sin^2 \frac{\th_1}{2} + \sin^2 \frac{\th_2}{2}\right).
\end{equation}
This definition for $\xi_{\thv}$ is chosen so that the subsequent LFA results in two dimensions closely resemble the corresponding one-dimensional results in Section~\ref{sec:redBlack1d}. 
% This notation of $\xi_{\thv}$ is for the purpose of making the two-dimensional analysis formally analogous to the one-dimensional case. 
The red-black smoothing operator is given as before in~\eqref{S_rb}, with
% the \kP{Gauss-Seidel operator (would let to avoid calling it that since it's different from the operator for an actually GS iteration)} being given by 
\begin{equation}
S_h^{GS} \equiv S_h^J = I_h - \frac{\omega h^2}{4} L_h.
\end{equation}
When operating on a Fourier mode $\phi_h (\cdot,\ \thv)$,
%\old{Considering the smoother acting on each Fourier mode $\phi_h (\cdot,\ \thv)$, }
the red-black smoother, as in any dimension, only introduces one aliasing mode with
$\bar{\thv} \eqdef [\bar{\th}_1,\ \bar{\th}_2]^T$
for each $\thv \in \Theta$. 
That is, $S_h$ is invariant on the two-dimensional subspaces 
%({\blue should we use Theta Low as in 1D?} \kP{No. This works for any $\thv \in \Theta$, only with repetitions. There is no Theta Low at this point, since there is no coarse grid yet, only concerning the smoother})
\begin{equation}\label{E_2d}
    \mathbb{E}_{\thv} = \text{span} \Big\{\, \phi_h (\cdot,\ \thv),~ \phi_h (\cdot,\ \bar{\thv}) \, \Big\},\quad \thv \in \Theta,
\end{equation}
on which it has the $2\times 2$ matrix symbol (for details see \ref{Sec: rb_2d})
\begin{equation}\label{S_hat_2d}
\symE S(\thv) = 
% \symE S_\Ebb (\thv)\eqdef 
\symE S^B (\thv)\ \symE S^R(\thv),
\end{equation}
with
\begin{equation}\label{S_hat_RB_2d}
\begin{aligned}
\symE S^R(\thv) = 
\mathbf{I} - \omega \begin{bmatrix}
1 & -1\\ -1 & 1
\end{bmatrix} \begin{bmatrix}
\xi_{\thv} \\ & 1-  \xi_{\thv}
\end{bmatrix},\quad
\symE S^B (\thv)  = \mathbf{I} - \omega \begin{bmatrix}
1 & 1\\ 1 & 1
\end{bmatrix} \begin{bmatrix}
\xi_{\thv} \\ & 1-  \xi_{\thv}
\end{bmatrix}.
\end{aligned}
\end{equation}
This is of the exact same form as the one-dimensional case (\ref{S_RB}) (with only the definition of $\xi_{\thv}$ being different).

%\old{On the other hand,}
Now consider the coarse grid $\grid G_H$ \eqref{G_H_d} (with $d=2$)
by standard coarsening with $H = 2h$. 
Note that $\grid G_H$ is embedded in $\grid G_h$ and consists of one quarter of the grid points. 
The eigenfunctions of any linear operator $L_H$ on $\grid G_H$ are
\begin{equation}
  \phi_H (\mathbf{x},\ 2\thv) = e^{i 2\boldsymbol{\theta} \cdot \frac{\xv}{H}},\quad \xv \in \grid G_H,\quad \theta\in\ThetaLowStd, 
\end{equation}
where have introduced the set of low modes associated with standard coarsening corresponding with
\ba
    \ThetaLowStd \eqdef [-\frac{\pi}{2},\ \frac{\pi}{2})^2.
\ea
This gives the separation of $\Theta_{std}^{low}$ and $\Theta_{std}^{high} \eqdef \Theta \setminus \Theta_{std}^{low}$ of the Fourier modes $\phi_h (\cdot,\ \thv)$, as shown in Figure \ref{fig: theta_2d} (left). Note that only one quarter of the modes are low-frequency modes, that is, only one-quarter of the modes can be represented on the coarse grid $\grid G_H$.
Given the coarse grid by standard coarsening, the resulting coarse-level correction operator $K_h^H$ has the four-dimensional eigen-subspaces 
\begin{equation}\label{F}
\Fthetav \eqdef \text{span} \Big\{ \, \phi_h (\cdot,\ \th_1,\ \th_2),~ \phi_h (\cdot,\ \thbar_1,\ \thbar_2) ,~ \phi_h (\cdot,\ \thbar_1,\ \th_2),~  \phi_h (\cdot,\ \th_1,\ \thbar_2)\, \Big \},
\quad \thv = [\theta_1,\  \theta_2]^T   \in \ThetaLowStd.
% \quad \thv = [\theta_1,\  \theta_2]^T \in \Theta^{low}.
\end{equation}
These eigenspaces come from the fact that, for each 
$\phi_h (\cdot,\ \thv)$, 
%$\thv \in \Theta$,
the coarse grid $\grid G_H$ can not distinguish 
%the modes in $\Fthetav$;
it from the three other aliasing modes in the basis of $\Fthetav$;
thus for each $\thv = [\th_1,\ \th_2]^T$, the coarse-level correction introduces three other modes corresponding to $[\thbar_1,\ \thbar_2]^T$, $[\thbar_1,\ \th_2]^T$ and $[\th_1,\ \thbar_2]^T$.
% corresponding to $\thv = [\th_1,\ \th_2]^T$.
%  
Comparing $\mathbb{E}_{\thv}$ with $\mathbb{F}_{\thv}$, we note that 
\begin{equation}\label{FE}
    \mathbb{F}_{\th_1,\ \th_2} =  \mathbb{E}_{\th_1,\ \th_2} + \mathbb{E}_{\bar{\th}_1,\ \th_2},\quad [\theta_1,\  \theta_2]^T \in \ThetaLowStd.
\end{equation}
Thus the red-black smoothing operator is also invariant on the four-dimensional subspaces $\Fthetav$. 
With respect to $\Fthetav$, the red-black smoother has the $4\times4$ matrix symbol
\begin{equation}\label{S_hat_F}
\symE S_{\mathbb{F}}(\th_1,\ \th_2) 
%= \symE S_{\mathbb{F}}^B (\thv)\ \symE S_{\mathbb{F}}^R(\thv) 
= \begin{bmatrix}
\symE S(\th_1,\ \th_2)\\
& \symE S(\bar \th_1,\ \th_2)
\end{bmatrix}.
\end{equation}
The smoothing factor of the red-black smoother for a multigrid cycle with standard coarsening is therefore 
\begin{equation}\label{mu_2_2}
\mu = \sup_{\thv \in \Theta_{std}^{\text{low}}} { \rho^{\frac{1}{\nu}} \big(\symE Q_{\mathbb{F}}\ \symE S_{\mathbb{F}}^{\nu}(\thv) \big)},
\end{equation}
where the ideal coarse-level correction operator $Q_h^H$ has the matrix symbol on  $\Fthetav$ as
\ba
 \symE Q_{\mathbb{F}} = \begin{bmatrix}  0\\ & 1\\ && 1\\ &&&1 \end{bmatrix}.  
\ea
Meanwhile, the two-level multigrid iteration operator $M_h^H$ has its matrix symbol on $\mathbb{F}_{\thv}$ as 
\begin{equation}\label{M_2}
\symE M_{\mathbb{F}}(\thv) = \symE S_{\mathbb{F}}^{\nu_2}(\thv) \ \symE K(\thv)\ \symE S_{\mathbb{F}}^{\nu_1}(\thv),
\end{equation}
where the matrix symbol $\symE K(\thv)$ of $K_h^H$ is given in (\ref{K_hat_2}) (for details see \ref{Sec: rb_2d}).
Based on $\symE M_{\mathbb{F}}(\thv)$, the asymptotic convergence rate of the two-level cycle with standard coarsening is
\begin{equation}\label{rho_2d_std}
\rho = \rho \big(M_h) = \sup_{\thv \in \Theta_{std}^{\text{low}}} { \rho \big(\symE M_{\mathbb{F}}(\thv)\big)}.
\end{equation}	

Note that in the one-dimensional case, both $S_h$ and $K_h^H$ introduce the same aliasing frequency $\thbar$ for each $\th$, and thus the two-level iteration operator $M_h$ has eigenspaces $\Ebb_\theta$~\eqref{E_1d} of dimension two.
In two dimensions, however, while the coarse-level correction $K_h^H$ has a 
representation $\symE K(\thv)$ that is a full $4\times 4$ matrix, the smoothing operator has a matrix representation $S_{\mathbb{F}}$ that is block diagonal, consistent with the subspace $\Fthetav$ being a sum of two $\Ethetav$ subspaces~\eqref{FE}.
%{is separable, that is $S_{\mathbb{F}}$ is block-diagonal, consistent with the relationship of thesubspaces (\ref{FE}).}
Thus, while the coarse-level operator $K_h^H$ mixes all four modes in the basis of $\Fthetav$, the smoother only mixes modes in the separate two-dimensional subspaces $\Ebb_{\theta_1,\theta_2}$ and $\Ebb_{\bar\theta_1,\theta_2}$.
This inconsistency between the invariant subspaces for the smoother and coarse-level operators 
is likely to be the reason why the smoother and the coarse-level correction
do not cooperate to give as good convergence rates as in one dimension. 
\commentB{Heuristically, while $K_h^H$ considers each subspace
$\mathbb{F}_{\boldsymbol{\theta}}$ as a single item (and acts on it by $\symE K(\thv)$), $S_h$
actually discriminates(\textcolor{red}{?}) the 4-dimensional subspaces further into two
2-dimensional subspaces (and acts on them separately by $\symE S(\th_1,\ \th_2)$ and $\symE S(\bar
\th_1,\ \th_2)$).
}

% -----------------  Low and High Modes Figure, Standard coarsening and red-black coarsening 
\input ./fig/theta_2d

We seek to construct a coarse grid such that the coarse-level correction operator $K_h^H$ is indeed
invariant on the sames subspaces $\mathbb{E}_{\thv}$ as the smoothing operator.
In Figure~\ref{fig: theta_2d} (middle) we illustrate the behavior of the red-black smoother on the Fourier modes $\phi_h (\cdot,\ \thv)$, $\thv \in \Theta$,
by plotting contours of the spectral radius of the $2 \times 2$ matrix $\symE S(\thv)$ as a function of $\thetav$\footnote{Since we can not plot the matrix symbol $\symE S(\thv)$, we represent it by its spectral radius. But note that $\rho(\symE S(\thv))$ does not tell the whole story of $S_h$.
%the operator is fully represented by its matrix symbol, but can not be fully represented by a mere spectral radius of the matrix; 
In particular, when two operators composite, their matrix symbols multiply, but the spectral radii do not.}.
(Note that some of the symmetries in the figure arise from $\mathbb{E}_{\bar{\thv}} \equiv \mathbb{E}_{\thv}$.)
\kP{These contours roughly represent how effective the smoother is in reducing each frequency $\thv \in \Theta$; for example, for $\thv$ close to $\mathbf{0}$, the smoother is very ineffective.}
We observe that besides the periodic property resulting from $\bar{\thv}$, there are also symmetry properties of $\symE S(\thv)$ as a
result of the symmetries in (\ref{xi_r2}):
$
\xi_{\pm\th_1,\ \pm \th_2} =  \xi_{\pm\th_2,\ \pm \th_1}
$
(the symmetry in $\th_1$ and $\th_2$ comes from the symmetry of the Laplacian).
In particular, for $\th_1,\ \th_2 > 0$, we have
\begin{equation}
\symE S(\th_1,\ \th_2) \equiv \symE S(\th_1-\pi,\ \th_2-\pi) = \symE S(\pi-\th_1,\ \pi-\th_2) = \symE S(\pi-\th_2,\ \pi-\th_1),
\end{equation}
thus the symmetry of $\symE S$ about $\th_1+\th_2 =\pi$ (shown as the red line in the first quadrant); and similarly for the other three quadrants. 
From these symmetries of $\symE S(\thv)$ we can derive a partition of $\Theta$ in two regions of equal area
as shown in Figure \ref{fig: theta_2d} (middle, red lines),
so that the smoother is fully represented by $\thv$ inside the red box. 
%(i.e. $\thv \in \ThetaLowRB)$, to be discussed in Section~\ref{Sec: rbC}} (\wdh{true?}). 
%\old{so that $\thv$ inside the red box can exhaust all $\mathbb{E}_{\thv}$ for $\thv \in \Theta$} 
With standard coarsening, unfortunately,
the smoothing operator mixes all the modes in both $\mathbb{E}_{\thv}$ and $\mathbb{E}_{\bar{\th}_1,\ \th_2}$ so that its spectral radius over $\Fthetav$ is the maximum of two spectral radii
\commentB{ 
\ba
    \symE S_{\mathbb{F}}(\th_1,\ \th_2) 
%\equiv \symE S_{\mathbb{F}}(\bar\th_1,\ \bar\th_2) 
\equiv \symE S_{\mathbb{F}}(\bar\th_1,\ \th_2)
%\equiv \symE S_{\mathbb{F}}(\th_1,\ \bar\th_2)
\ea
so that
}
\ba
  \rho \big(\symE S_{\mathbb{F}}(\thv) \big) = \max \Big\{\, \rho \big(\symE S(\th_1,\ \th_2) \big),~ \rho \big(\symE S(\bar \th_1,\ \th_2)\big)\, \Big\},\quad \thv = [\theta_1,\  \theta_2]^T   \in \Theta.
\ea
The behavior of $S_h$ on $\mathbb{F}_{\thv}$, represented as contours of the spectral radius of the $4\times4$ matrix $S_\Fbb(\thetav)$,
is illustrated in Figure \ref{fig: theta_2d} (right) for comparison. 
On the other hand, if we could find some coarse grid whose low-frequency $\Theta^{low}$
corresponds to the red box in Figure~\ref{fig: theta_2d} (middle), then
we would be able to recover a similar Fourier structure between the smoother and coarse-level operator as the one-dimensional case.
In the next section, we will show that red-black coarsening (in two dimensions) results in the desired set of low-frequency modes.

 %1D2D std-coarsening
\subsection{Red-black coarsening and red-black smoothing in two dimensions}\label{Sec: rbC}

Red-black coarsening (RBC) is a well known coarsening strategy~\cite{Trottenberg2001}
that consists of choosing the coarse-level grid points to be the black (or red) points in a red-black (two-color, odd-even) partition of the fine grid points, $\grid G_H \equiv \grid G_h^B$, as shown in Figure~\ref{fig:redBlackCoarseLevelGrids}.
These grid points can be viewed as a rotated grid with spacing $H=\sqrt{2} h$. 
The coarse grid thus has one half of the grid points in the fine grid and 
consequently one half of the number of Fourier modes on the fine grid in local Fourier analysis.
As we shall see, this choice of coarse grid
means that the coarse-level correction operator $K_h^H$
will be invariant on the same eigen-subspaces $\mathbb{E}_{\thv}$,
as the red-black smoother. This property leads us to the construction of a special
two-level multigrid algorithm in two dimensions
that becomes a direct solver.  

% --------------------- RED BLACK COARSE GRIDS -------------------
{
\newcommand{\labelFont}{\footnotesize}
\begin{figure}[hbt]
  \begin{center}
  \resizebox{14cm}{!}{% START resize box
    \begin{tikzpicture}
    \useasboundingbox (0,.5) rectangle (18,4.75);
     % \draw[<-,ultra thick,red,dash pattern= on 8pt off 3pt] (0,3) -- (2,3);
     %

     % ---- level 0 fine grid -------
     \begin{scope}[xshift=0cm,yshift=0cm]  

        \draw (2,4.25) node[anchor=south,fill=white,draw=black,inner sep=2.5pt] {\labelFont level $l=0$};

        \draw[step=.5cm,black,thick] (0,0) grid (4,4);

        % \draw[step=1cm,blue,very thick] (0,0) grid (4,4);
         \foreach \x in {0,1,...,4}
         \foreach \y in {0,1,...,4}
            \draw[black] (\x,\y) \rbDot;
         \foreach \x in {.5,1.5,...,4}
         \foreach \y in {.5,1.5,...,4}
            \draw[black] (\x,\y) \rbDot;

         \foreach \x in {.5,1.5,...,4}
         \foreach \y in {0,1,...,4}
            \draw[red] (\x,\y) \rbDot;
         \foreach \x in {0,1,...,4}
         \foreach \y in {.5,1.5,...,4}
            \draw[red] (\x,\y) \rbDot;
     \end{scope}

     % ---- level 1 --
     \begin{scope}[xshift=4.75cm,yshift=0cm]  

        \draw (2,4.25) node[anchor=south,fill=white,draw=black,inner sep=2.5pt] {\labelFont level $\ll\!=\!1$};

         \foreach \x in {0,1,...,4}
         \foreach \y in {0,1,...,4}
            \draw[black] (\x,\y) \rbDot;
         \foreach \x in {.5,1.5,...,3.5}
         \foreach \y in {.5,1.5,...,3.5}
         {
            \draw[black] (\x,\y) \rbDot;
            % -- lines on coarse grid ---
            \draw[black,thick] (\x-.5,\y-.5) -- (\x+.5,\y+.5); 
            \draw[black,thick] (\x-.5,\y+.5) -- (\x+.5,\y-.5); 
         }
         
     \end{scope}

     % ---- level 2 -------
     \begin{scope}[xshift=9.5cm,yshift=0cm]  

        \draw (2,4.25) node[anchor=south,fill=white,draw=black,inner sep=2.5pt] {\labelFont level $\ll\!=\!2$};

        \draw[step=1cm,black,thick] (0,0) grid (4,4);

        % \draw[step=1cm,blue,very thick] (0,0) grid (4,4);
         \foreach \x in {0,1,...,4}
         \foreach \y in {0,1,...,4}
            \draw[black] (\x,\y) \rbDot;

         \foreach \x in {1,3}
         \foreach \y in {0,2,4}
            \draw[red] (\x,\y) \rbDot;
         \foreach \x in {0,2,4}
         \foreach \y in {1,3}
            \draw[red] (\x,\y) \rbDot;

     \end{scope}

     % ---- level 3 --
     \begin{scope}[xshift=14.25cm,yshift=0cm]  

        \draw (2,4.25) node[anchor=south,fill=white,draw=black,inner sep=2.5pt] {\labelFont level $\ll\!=\!3$};

         \foreach \x in {0,2,...,4}
         \foreach \y in {0,2,...,4}
            \draw[black] (\x,\y) \rbDot;

         \foreach \x in {1,3}
         \foreach \y in {1,3}
         {
            \draw[black] (\x,\y) \rbDot;
            % -- lines on coarse grid ---
            \draw[black,thick] (\x-1,\y-1) -- (\x+1,\y+1); 
            \draw[black,thick] (\x-1,\y+1) -- (\x+1,\y-1); 
         }
         
     \end{scope}
% 
%   \draw[step=1cm,gray] (0,0) grid (18,4);
  \end{tikzpicture}
  }% --- END RESIZE BOX
  \end{center}
  \caption{Red-black coarsening (2D). Coarse level $\ll\!=\!1$ is a rotated grid with spacing $h_1=\sqrt{2}h_0$ formed from the black points of level $\ll\!=\!0$.
   Coarse level $\ll\!=\!1$ has spacing $h_2=2 h_0$, while level $\ll\!=\!3$ is a rotated grid with spacing $h_3=2\sqrt{2} h_0$. 
  }
  \label{fig:redBlackCoarseLevelGrids}
\end{figure}
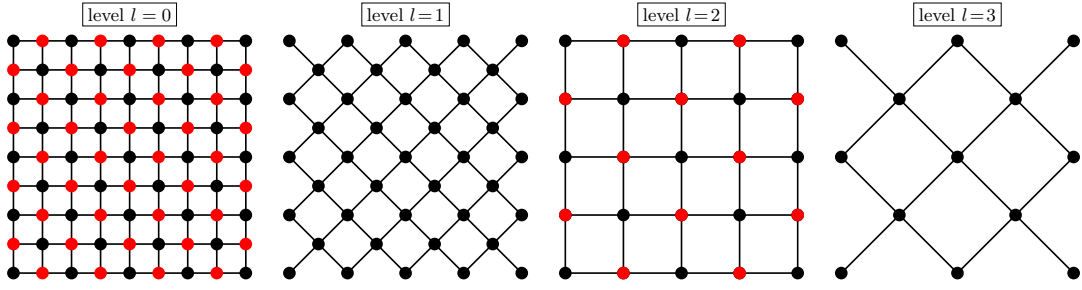 
}

Introduce the ``parameter space'' coarse grid, using coordinates $\rv$ and with grid spacing $H=\sqrt{2} h$, 
\begin{gather}\label{G_r2}
\grid G_H = \Big\{  \rv =\Uv \xv = \kv\, H ~:~ \kv \in \mathbb{Z}^2 \Big\},
\end{gather}
where
the corresponding physical space coordinates $\xv$ are  
defined through the unitary mapping (rotation) 
\begin{equation}
\xv=\Xv(\rv)=\Uv^{T}\rv, \quad  \Uv \eqdef \frac{1}{ \sqrt{2}}\begin{bmatrix} 1 & -1 \\ 1 & 1 \end{bmatrix}.
\end{equation}
The grid points $\xv=\Xv(\rv)$, for $\rv\in\grid G_H$, correspond exactly to the black grid points in $\grid G_h^B$,
  but here we consider these to be from a rotated Cartesian
  grid with grid-spacing $H = \sqrt{2}h$, as can be seen in Figure~\ref{fig:redBlackCoarseLevelGrids} (second from left).
Note that the Laplacian $\Delta$ is invariant under rotation:
\begin{equation}
\Delta_{\mathbf{x}} = \partial_{\mathbf{x}}^T \partial_{\mathbf{x}} = 
\partial_{\mathbf{r}}^T \mathbf{U U}^T \partial_{\mathbf{r}} = 
\Delta_{\mathbf{r}},
\end{equation}
where $\partial_{\mathbf{x}} \equiv [\partial_{x_1},\ \partial_{x_2}]^T$, $\partial_{\mathbf{r}} \equiv [\partial_{r_1},\ \partial_{r_2}]^T$.
Thus any discretization stencil for the Laplacian on an equally spaced Cartesian grid takes the same form
on the rotated grid. 

We shall show that by using red-black smoothing and red-black coarsening, analogous to Theorem~\ref{Th_1d}, we have the following result:
% -------------------------------- THEOREM --------------------
\begin{theorem}[red-black smoothing and red-black coarsening in two dimensions]
\label{Th_2d}
Consider the two-dimensional discrete Poisson's equation (\ref{MP_h}) with (\ref{L_h_2})
on the infinite grid (\ref{G_h_d}) with $d=2$.
An $(h,\ H)$ two-level multigrid cycle consisting of the red-black smoother ($\omega=1$) and red-black coarsening ((\ref{G_r2}), $H=\sqrt{2}h$), with the transfer operators
\begin{equation}\label{transfer_r2}
I_h^H = \frac{1}{8}
\begin{bmatrix}
& 1 &  \\
1 & 4 & 1\\
& 1 & 
\end{bmatrix}_h^{H},\quad
I^h_{H} = \frac{1}{4} \left]
\begin{matrix}
& 1 &  \\
1 & 4 & 1\\
& 1 & 
\end{matrix}\right[^h_{H},
\end{equation}
as well as the Galerkin coarse-level operator $L_{H} =  I_h^{H} L_h I^h_{H}$,
gives rise to a direct solver.
%consisting of a coarse-level correction and a post-smoothing step.
That is, the cycle converges in a single step (as long as it includes a post-smoothing step after a coarse-level correction).
\end{theorem}  
% --------------------------- END THEOREM ---------------
The recovery of the direct-solver property in two dimensions will be shown through local Fourier analysis, and
 will rely on the spectral properties of the red-black smoothing and red-black coarse-level correction operators being invariant on the same eigenspaces.

The Fourier modes on the parameter-space coarse grid~\eqref{G_r2} are given by
\begin{equation}
\phi_H (\rv,\ \alv) = e^{i \alv \cdot \frac{\rv}{H}},\quad \rv \in \grid G_{H},\quad  \alv\in \Theta.
%= [\a_1,\  \a_2]^T 
\end{equation}
Then the coarse grid separates low- and high-frequency modes $\phi_h (\cdot,\ \thetav)$ on $\grid G_h$, based on whether they can be represented on $\grid G_H$, as  
\begin{align}\label{Th_L_2d_rb}
\Theta_{rb}^{\text{low}} \eqdef \{\thv = [\theta_1,\  \theta_2]^T \in \Theta:\quad \sqrt{2} \mathbf{U} \thv = [\th_1-\th_2,\  \th_1+\th_2]^T \in \Theta\},
%&=\{[\th_1,\  \th_2]^T \in \Theta:\quad \th_1-\th_2,\  \th_1+\th_2 \in [-\pi,\ \pi)\},
\end{align}
and $\Theta_{rb}^{\text{high}} = \Theta \setminus \Theta_{rb}^{\text{low}}$, as shown in Figure \ref{fig: theta_2d} (left).
We see that red-black coarsening partitions the frequencies into two sets of equal size, in exactly the same way as proposed in Figure \ref{fig: theta_2d} (middle).
\commentB{
\wdh{Comment: I am not sure if bisection is the correct word here -- to me it means split along the middle.}
\old{We can see that this gives a bisection of $\Theta$ into low and high-frequency modes
exactly the same way as proposed in Figure \ref{fig: theta_2d} (middle).}}
The $(h,\ H)$ two-level cycle then has the subspaces $\mathbb{E}_{\thv},\ \thv \in \ThetaLowRB$ as eigenspaces.
Modes in each of the subspaces are indistinguishable on the coarse grid $\grid G_H$, in particular,  
\begin{align}
    \phi_h (\xv,\ \thv) &= \phi_h (\xv,\ \bar \thv)
    = \phi_H (\rv(\xv),\ \sqrt{2} \mathbf{U} \thv),\quad \forall\ \rv(\xv) \in \grid G_H,
    \quad \thv \in \Theta_{rb}^{\text{low}} . 
\end{align}
This is from the facts that for $\mathbf{x} = [j_1,\ j_2]^T h$ such that $\rv =\Uv \xv \in \grid G_H$, $e^{i\pi (j_1 \pm j_2)} = 1$, and
$\thv \cdot \xv\frac{1}{h} = \thv \cdot \mathbf{U}^T\rv\frac{1}{h}  = \sqrt{2} \mathbf{U} \thv \cdot \rv\frac{1}{H}.$

The smoothing factor of the red-black smoother then becomes
\begin{equation}\label{mu_2d_rbC}
\mu = \rho^{\frac{1}{\nu}} \big( Q_h^H S_h^{\nu}\big)
= \sup_{\thv \in \Theta_{rb}^{\text{low}}} { \rho^{\frac{1}{\nu}} \big(\symE Q\ \symE S^{\nu}(\thv) \big)},
\qquad   \symE Q = \begin{bmatrix} 0  \\ & 1 \end{bmatrix}, 
\end{equation}
where $\symE S(\thv)$ is given by \eqref{S_hat_2d}, and 
$\symE Q$ is the matrix representation of the ideal coarse-level correction $Q_h^H$ on $\mathbb{E}_{\th}$.
\edit{The smoothing factor \eqref{mu_2d_rbC} is completely analogous to the one-dimensional case \eqref{mu_1d}.}
Consider the smoothing factor with standard coarsening \eqref{mu_2_2}, and
note that according to the block-diagonal structure,
\begin{equation}
  \rho^{\frac{1}{\nu}} \big(\symE Q_{\mathbb{F}}\ \symE S_{\mathbb{F}}^{\nu}(\thv) \big) =
  \max \Big\{ \, \rho^{\frac{1}{\nu}} \big(\symE Q\ \symE S^{\nu}(\th_1,\ \th_2) \big),~ \rho \big(\symE S(\bar \th_1,\ \th_2)\big) \, \Big\}.
\end{equation}
That is, for $\thv \in \Theta_{std}^{low}$,
not only $\phi_h(\cdot,\ \bar{\thv})$, but also the modes in $\mathbb{E}_{\bar{\th}_1,\ \th_2}$ are high-frequency according to standard coarsening and have to be dealt with by the smoother. 
As a result,  the smoothing factor $\mu$  with standard coarsening (\ref{mu_2_2}) will inevitably be worse than that with red-black coarsening \eqref{mu_2d_rbC}.
% \old{\wdh{Are these next comments still true?} 
We suspect the imbalance of $\Theta_{std}^{\text{low}}$ and $\Theta_{std}^{\text{high}}$ (in contrast to the balance of $\Theta_{rb}^{\text{low}}$ and $\Theta_{rb }^{\text{high}}$) to be one of the reasons why multigrid by standard coarsening in two dimensions does not yield convergence rates as good as those in one dimension, and
we are going to examine this suggestion in Section \ref{Sec: rbC_3d} and Section \ref{Sec: rC_J}.
Contours of $\rho(\Qv\ \Sv(\thv))$, corresponding to red-black coarsening, 
are compared with contours of $\rho(\Qv_\Fbb\Sv_\Fbb(\thv))$, corresponding to standard coarsening, 
in Figure~\ref{fig:rho_T_2o}. 
As noted previously, the worst value of $\rho(\Qv\ \Sv(\thv))$ 
is smaller than the worst value of  $\rho(\Qv_\Fbb\Sv_\Fbb(\thv))$, for $\thetav\in\Theta$.

\medskip   
We construct the bilinear restriction and interpolation operators between $\grid G_h$ and $\grid G_H$ as given by (\ref{transfer_r2}),
with  $I_h^H = (I^h_H)^* = \frac{1}{2} (I^h_H)^T$.
For each $\boldsymbol{\theta} \in \Theta_{rb}^{\text{low}}$, we have
\begin{equation}
\begin{cases}
I_h^{H} \PhiVectorEx = \phi_{H}\big(\rv(\xv),\ \th_1-\th_2,\ \th_1+\th_2\big)\ \symE I_h^{H}(\boldsymbol{\theta}),\\
I^h_{H}\ \phi_{H}\big(\rv(\xv),\ \th_1-\th_2,\ \th_1+\th_2\big)=  \PhiVectorEx \symE I^h_{H} (\boldsymbol{\theta}),
\end{cases}
\end{equation}
where
\begin{comment}
\begin{equation}
\symE I_h^{H}(\thv) = \begin{bmatrix}
\symF I_h^H(\thv) &
\symF I_h^H(\bar{\thv}) 
\end{bmatrix},\quad \symE I^h_{H}(\thv) = \symE I_h^{H}(\thv)^T,
\end{equation}
with 
\begin{equation}
\symF I_h^H(\th_1,\ \th_2) =1- \frac{1}{2}(\sin^2 \frac{\th_1}{2} + \sin^2 \frac{\th_2}{2}) 
,\quad [\th_1,\ \th_2]^T \in \Theta.
\end{equation}
Or more specifically, 
\end{comment}
\begin{equation}
\symE I_h^{H}(\boldsymbol{\theta}) = \begin{bmatrix}
1-\xi_{\thv} & \xi_{\thv}
\end{bmatrix}, \quad \symE I^h_{H}(\thv) = \symE I_h^{H}(\thv)^T,
\end{equation}
which is completely analogous to the one-dimensional case (\ref{transfer_hat_1d}).
On the coarse grid $\grid G_H$, as usual, we consider both the non-Galerkin operator with the same stencil as $L_h$, that is,
% Illustration of the smoothing factor and iteration operator ...
{
  \newcommand{\labelFont}{\scriptsize}
  \newcommand{\figWidth}{5.2cm}% height 
  \newcommand{\trimfig}[2]{\trimh{#1}{#2}{.0}{.0}{.0}{.0}}
  \begin{figure}[t!]
    \begin{center}
      \begin{tikzpicture}[scale=1]
	\useasboundingbox (0,0.65) rectangle (15,5.5);  % set the bounding box (so we have less surrounding white space)
	
	\draw(-0.5,0.0) node[anchor=south west,xshift=-15pt,yshift=-8pt] {\trimfig{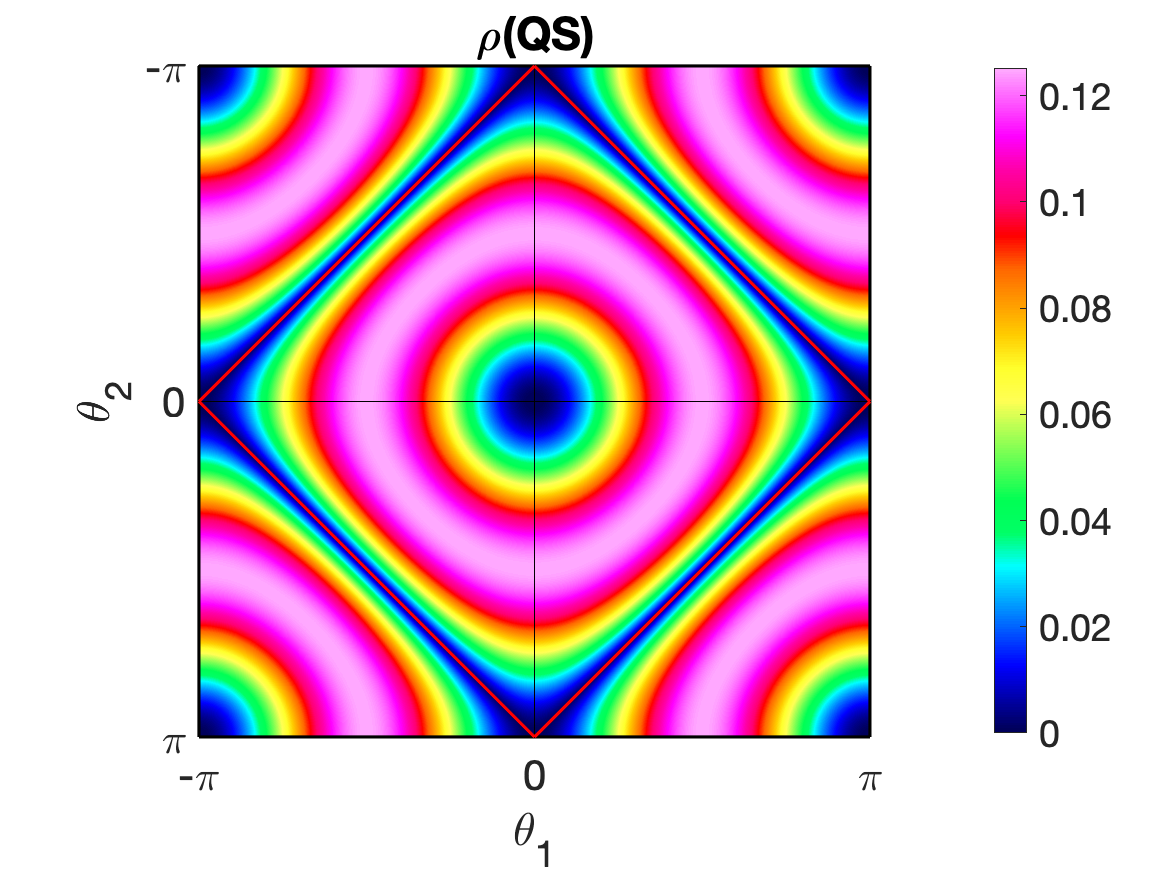}{\figWidth}};
        \draw  (2.8,4.75) node[anchor=south,fill=white,yshift=-1pt] {\labelFont %$\rho(\Qv_\Ebb\Sv_\Ebb(\thetav))$}; 
        $\rho(\Qv\ \Sv(\thetav))$}; 
        \draw (.5,4.25) node[anchor=north west,fill=white,draw=black,inner sep=1.5pt] {\labelFont red-black coarsening};

        \begin{scope}[xshift=7.5cm]
	  \draw(-0.5,0.0) node[anchor=south west,xshift=-15pt,yshift=-8pt] {\trimfig{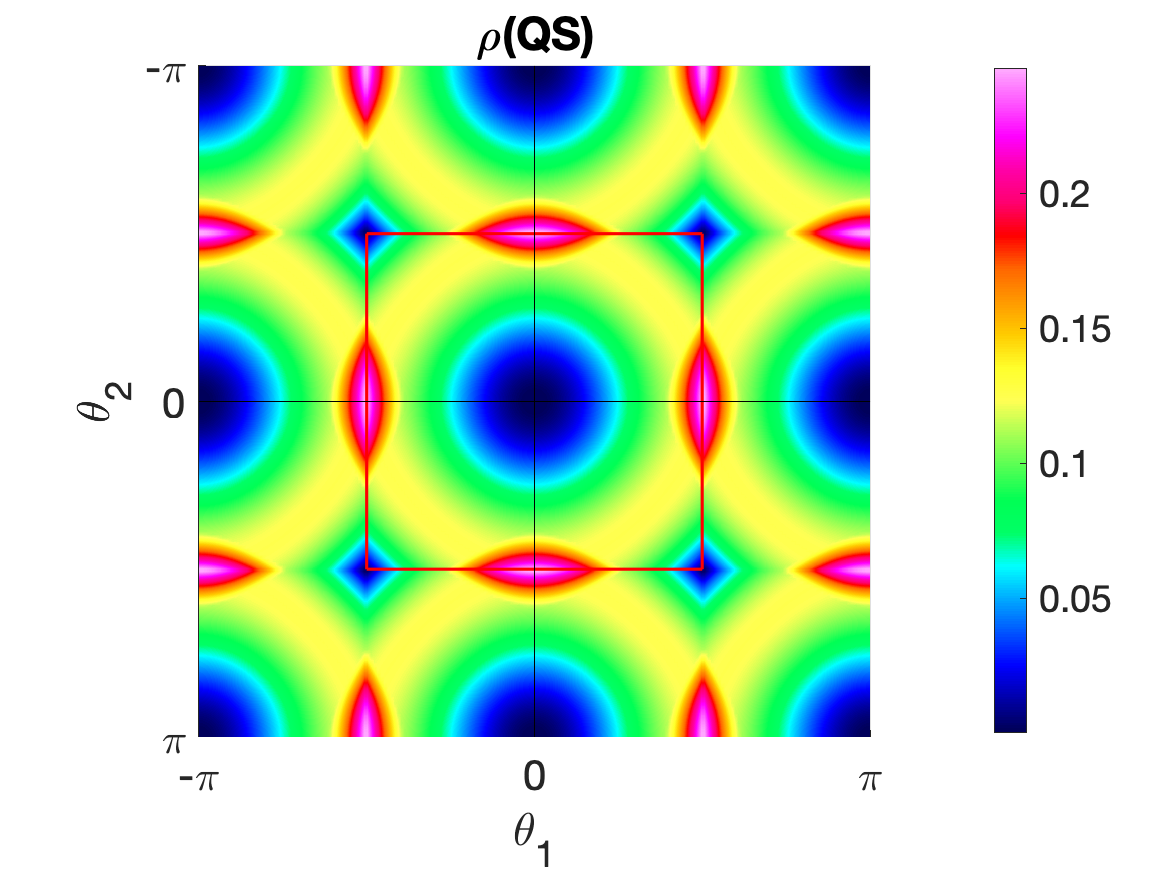}{\figWidth}};
          \draw (2.8,4.75) node[anchor=south,fill=white,yshift=-1pt] {\labelFont $\rho(\Qv_\Fbb\Sv_\Fbb(\thetav))$}; 
          \draw (.5,4.25) node[anchor=north west,fill=white,draw=black,inner sep=1pt] {\labelFont standard coarsening};
        \end{scope}

	% grid:
	% \draw[step=1cm,gray] (0,0) grid (16,6);
      \end{tikzpicture}
    \end{center}
    \captionof{figure}{Spectral radii of the matrix symbols of $Q_h^H S_h$, corresponding to a cycle with a single red-black smoothing step ($\nu=1$, $\omega=1$) and the ideal coarse-level correction, versus $\thv \in \Theta.$ 
    The maxima of these contours represent the \textbf{smoothing factor}s $\mu$ of the multigrid cycles.
  Left: 
  %$\rho(\Qv_\Ebb\Sv_\Ebb(\thv))$ 
  $\rho(\Qv\ \Sv(\thetav))$ for the $2\times2$ matrix symbol of $Q_h^H S_h$ on the 
  %two-dimensional 
  subspace $\mathbb{E}_{\thv}$
  associated with red-black coarsening.
  Right: $\rho(\Qv_\Fbb\Sv_\Fbb(\thv))$ for the $4\times 4$ matrix symbol of $Q_h^H S_h$ 
     on the 
     %four-dimensional 
     subspace $\mathbb{F}_{\thv}$ associated with  standard coarsening.
\commentB{ 
      Illustration of the operator $Q_h^H S_h$ (with $\omega=1$) on $\Theta$.
      Left: $\rho(\symE Q\ \symE S(\thv))$, right: $\rho(\symE Q_{\mathbb{F}}\ \symE S_{\mathbb{F}}(\thv))$; $\thv \in \Theta.$ }
  }
   \label{fig:rho_T_2o}
  \end{figure}
}

{
  \newcommand{\labelFont}{\scriptsize}
  \newcommand{\figWidth}{5.2cm}% height 
  \newcommand{\trimfig}[2]{\trimh{#1}{#2}{.0}{.0}{.0}{.0}}
  \begin{figure}[ht!]
    \begin{center}
      \begin{tikzpicture}[scale=1]
	\useasboundingbox (0,0.5) rectangle (15,5.5);  % set the bounding box (so we have less surrounding white space)
	
	\draw(-0.5,0.0) node[anchor=south west,xshift=-15pt,yshift=-8pt] {\trimfig{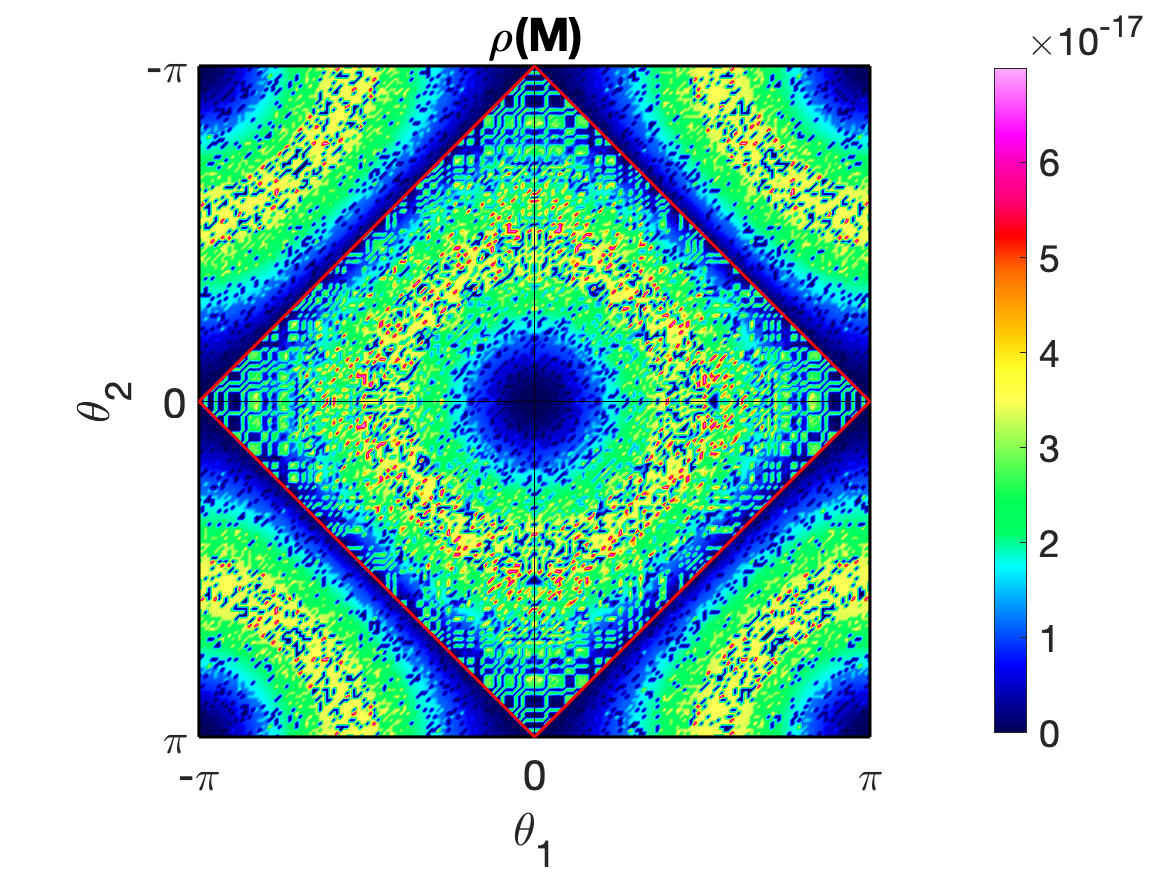}{\figWidth}};
        \draw  (2.8,4.75) node[anchor=south,fill=white,yshift=-1.5pt] {\labelFont %$\rho(\Mv_\Ebb(\thetav))$}; 
        $\rho(\Mv(\thetav))$}; 
        \draw (.5,4.25) node[anchor=north west,fill=white,draw=black,inner sep=1.5pt] {\labelFont red-black coarsening};
        \begin{scope}[xshift=7.5cm]
    	  \draw(-.5,0.0) node[anchor=south west,xshift=-15pt,yshift=-8pt] {\trimfig{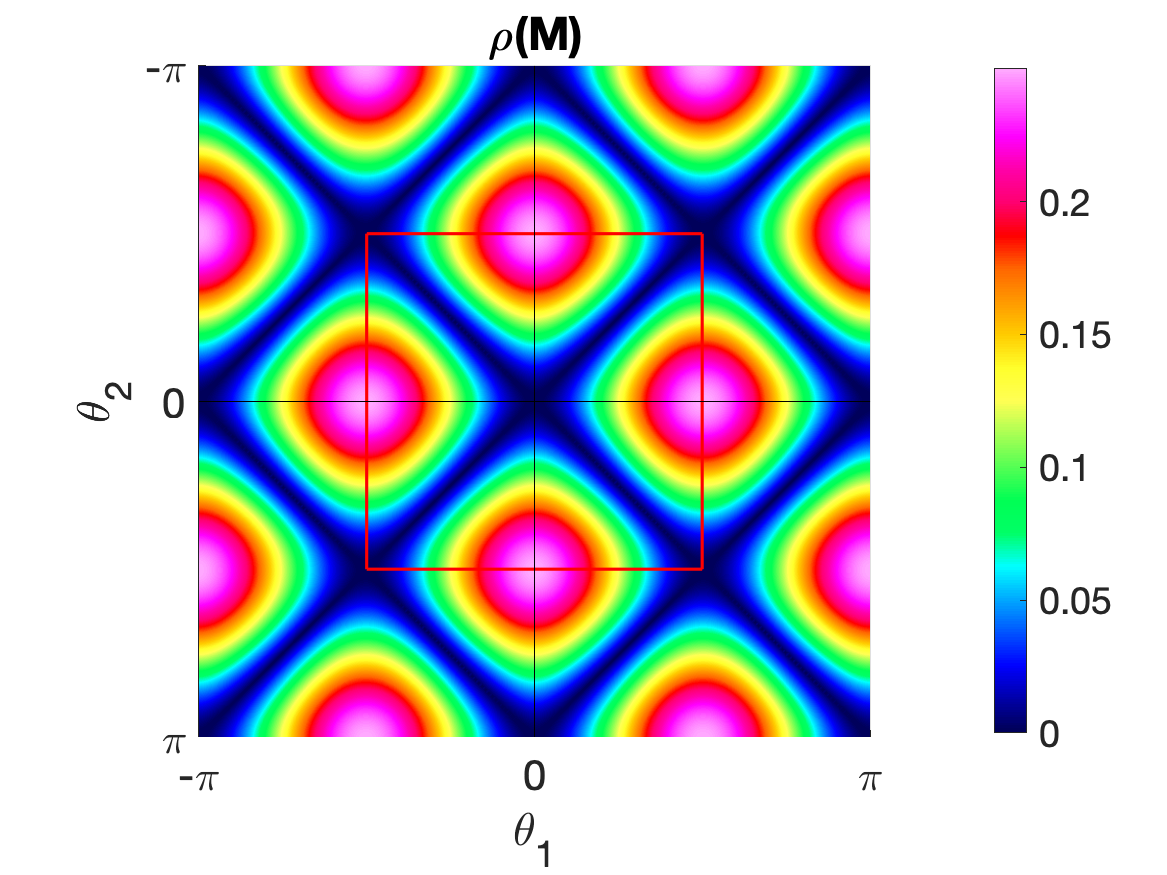}{\figWidth}};
          \draw (2.8,4.75) node[anchor=south,fill=white,yshift=-1pt] {\labelFont $\rho(\Mv_\Fbb(\thetav))$}; 
          \draw (.5,4.25) node[anchor=north west,fill=white,draw=black,inner sep=1pt] {\labelFont standard coarsening};
        \end{scope}
	% g     rid:
	% \draw[step=1cm,gray] (0,0) grid (16,6);
      \end{tikzpicture}
    \end{center}
    \captionof{figure}{
      (Computed) spectral radii of the matrix symbols of the multigrid iteration operator
      	$M_h = K_h^H S_h $, with $\omega=1$, $\nu=1$, and Galerkin coarse-level operators, versus $\thv \in \Theta.$ 
      	\kP{The maxima of these contours represent the \textbf{convergence rate}s $\rho$ of the multigrid cycles.}
          Left: $\rho(\Mv(\thetav))$
%          $\rho(\Mv_\Ebb(\thetav))=\rho(\Sv_\Ebb \Kv_\Ebb)$           
          for the $2\times2$ matrix symbol of
          $M_h$ on the 
          %two-dimensional 
          subspace $\mathbb{E}_{\thv}$ associated with red-black coarsening (zero to machine epsilon indicating a direct solver).
          Right: $\rho(\Mv_\Fbb(\thetav))$
          %$\rho(\Mv_\Fbb(\thetav))=\rho(\Sv_\Fbb \Kv_\Fbb)$ 
          for the $4\times 4$ matrix symbol of $M_h$ on the
%four dimensional 
	subspace $\mathbb{F}_{\thv}$ associated with standard coarsening.
      \commentB{Illustration of the 2-level iteration operator $M_h$ (with $\omega=1$, $\nu=1$, and Galerkin coarse-level operators) on $\Theta$.
      Left: red-black coarsening, right: standard coarsening.} 
  }
   \label{fig:rho_M_2o}
  \end{figure}
}

\begin{equation}
L_H = \frac{1}{H^2} \begin{bmatrix}
& -1 \\
-1 & 4 & -1\\
& -1
\end{bmatrix}_H,
\end{equation}
with Fourier symbol 
\begin{equation}
\symF L_H (\thv) = \frac{4}{H^2} (\sin^2 \frac{\th_1-\th_2}{2} + \sin^2 \frac{\th_1+\th_2}{2}) = \frac{8}{h^2} (\xi_{\thv} - \sin^2 \frac{\th_1}{2} \sin^2 \frac{\th_2}{2})
\end{equation} 
corresponding with the eigen-mode $\phi_{H}\big(\rv(\xv),\ \th_1-\th_2,\ \th_1+\th_2\big)$;
and the Galerkin operator
\begin{equation}\label{L_G1}
  L_{H} =  I_h^{H} L_h I^h_{H} = \frac{1}{4H^2}
  \begin{bmatrix}
     -1 & -2 & -1\\
     -2 & 12 & -2\\
     -1 & -2 & -1
  \end{bmatrix}_H,
\end{equation}
with symbol
\begin{equation}
    \symF L_H(\thv) = \symE I_h^{H}(\boldsymbol{\theta})\ \symE L(\boldsymbol{\theta})\ \symE I^h_{H} (\boldsymbol{\theta}) = \frac{8}{h^2}\ \xi_{\thv} (1-\xi_{\thv} ),
\qquad 
\symE L(\thv)  = \frac{8}{h^2} \begin{bmatrix} \xi_{\thv}\\ & 1-\xi_{\thv} \end{bmatrix} , 
\end{equation}
where $\symE L(\thv)$  is the matrix symbol of $L_h$ on $\mathbb{E}_{\thv}$. 
\begin{comment}
We have
\begin{equation}
L_{H}\ \phi_{H}\big(\rv(\xv),\ \th_1-\th_2,\ \th_1+\th_2\big) = \symF L_H(\thv)\ \phi_{H}\big(\rv(\xv),\ \th_1-\th_2,\ \th_1+\th_2\big).
\end{equation}
\end{comment}
Note that the stencils are written on the rotated grid.
The coarse-level correction is then constructed from $L_H$ together with the transfer operators (\ref{transfer_r2}) as
\begin{equation}\label{K_rb}
K_h^H = I_h - I^h_H L^{-1}_H I_h^H L_h, 
\end{equation}
whose matrix symbol on each $\mathbb{E}_{\thv}$ ($\thv \neq \mathbf{0})$ is 
\begin{equation}\label{K_hat_2d_rb}
\symE K (\thv) = \mathbf{I} - \symF L_H^{-1}(\thv)\
\symE I^h_H (\thv)\ \symE I_h^H(\thv)\ \symE L(\thv).
\end{equation}
In particular, if $L_H$ is the Galerkin operator~\eqref{L_G1}, then the coarse-level correction operator is represented by
\begin{equation}
  \symE K (\thv)  =
    \begin{bmatrix}
        1 \\ -1
    \end{bmatrix}
  \begin{matrix}
    \begin{bmatrix}
      \xi_{\thv} & -(1-\xi_{\thv})
     \end{bmatrix}
   \\ \mbox{}
  \end{matrix},
\end{equation}
which is a projector, and of the exact same form as the one-dimensional case (\ref{K_hat_1d}).

Thus we have the following main result for red-black coarsening analogous to Theorem \ref{Th_1d_LFA} in one dimension:
%----------------------- THEOREM 2D DIRECT SOLVER LFA ----------------------------
\begin{theorem}[red-black smoothing and red-black coarsening in two dimensions - LFA] \label{Th_rbC_LFA}
  The $(h,\ H)$ two-level iteration operator
  $
  M_h^H = S_h^{\nu_2} K_h^H  S_h^{\nu_1}
  $
  of a $V[\nu_1,\ \nu_2]$ cycle as described in Theorem \ref{Th_2d}
  %with red-black smoothing operator $S_h$ and coarse-level correction operator $K_h^H$ with red-black coarsening as in (\ref{K_rb}),
  has the matrix representation on each two-dimensional eigenspace $\mathbb{E}_{\thv}$ (\ref{E_2d}) ($\thv \neq \mathbf{0}$) as
  \begin{equation}
    \symE M(\thv) = \symE S^{\nu_2}(\thv) \ \symE K(\thv)\ \symE S^{\nu_1}(\thv),
  \end{equation}
  where $\symE S (\thv)$ is given by (\ref{S_hat_2d}), and $\symE K(\thv)$ is given by (\ref{K_hat_2d_rb}).
  Thus the asymptotic convergence rate of the two-level cycle is given by
  \begin{equation}\label{rho_2d_rb}
    \rho = \rho \big(M_h^H\big)
    = \sup_{\thv \in \Theta_{rb}^{\text{low}}} { \rho \big(\symE M(\thv)\big)},
  \end{equation}
  where $\Theta_{rb}^{\text{low}}$ is given by (\ref{Th_L_2d_rb}).
\end{theorem}
%----------------------- END THEOREM ------------------------------------
In particular, if $\omega=1$, and if the coarse-level operator $L_H$ is Galerkin (\ref{L_G1}), we have the same situation as the one-dimensional case, and the two-level cycle converges after one post red-smoothing, for
\begin{equation}
\symE S^R(\thv)\ \symE K (\thv) = \begin{bmatrix}
1-\xi_{\thv} \\ \xi_{\thv}
\end{bmatrix}\begin{bmatrix}
1 & 1
\end{bmatrix}
\begin{bmatrix}
1 \\ -1
\end{bmatrix}\begin{bmatrix}
\xi_{\thv} & -(1-\xi_{\thv})
\end{bmatrix} = \mathbf{O},
\end{equation}
analogous to (\ref{eqn: zero_1d}).
The iteration operator $M_h^H$ (with $\nu=1$) is illustrated by $\rho(\symE M(\thv)),\ \thv \in \Theta$
in Figure~\ref{fig:rho_M_2o}, compared with the standard coarsening case \eqref{M_2}. 
(Figure~\ref{fig:rho_M_2o} should also be compared with Figure~\ref{fig:rho_T_2o} with ideal coarse-level correction.)

In summary, given our definition (\ref{xi_r2}) of the averaged $\xi_{\thv}$, the local Fourier analysis in two dimensions with red-black coarsening on the subspaces $\mathbb{E}_{\thv}$ is reduced to very much the same form as the one-dimensional case, especially when the Galerkin coarse-level operator is considered. As with the smoothing factor $\mu$, we would expect the asymptotic convergence rate (\ref{rho_2d_rb}) to be smaller than (\ref{rho_2d_std}) with standard coarsening, and comparable with the convergence rates (\ref{rho_1d}) of the multigrid cycles in one dimension.

%\subsection{two-level Results}\label{Sec: results} --------------------------------------------
In Figure~\ref{fig:CR_2level_2d_rb}, we present the results from both local Fourier analysis and
numerical computations for two-level cycles with red-black coarsening in two dimensions.  To compare
red-black coarsening with standard coarsening, we assume the same size $N$, which is small enough so
that two-level cycles are reasonable to use.  The results confirm that the multigrid algorithm with the Galerkin
coarse-level operator converges in one cycle
%\old{has exact convergence} 
with $\omega = 1$. In fact, the overall curves of
convergence rates with the ideal coarse-level correction ($\mu^{\nu}$) and the Galerkin coarse-level
correction in Figure \ref{fig:CR_2level_2d_rb} are very similar to those of the one-dimensional case
in Figure \ref{fig:redBlackSmootherStandardCoarsening1d2d} (top). These observations are consistent
with our analysis.

{
\newcommand{\figWidth}{6cm}% height 
\newcommand{\trimfig}[2]{\trimh{#1}{#2}{.0}{.0}{.0}{.0}}
\begin{figure}[h!]
\begin{center}
\begin{tikzpicture}[scale=1]
  \useasboundingbox (0,0.75) rectangle (15,6.5);  % set the bounding box (so we have less surrounding white space)

   \draw(0.0,0.0) node[anchor=south west,xshift=-15pt,yshift=-8pt] {\trimfig{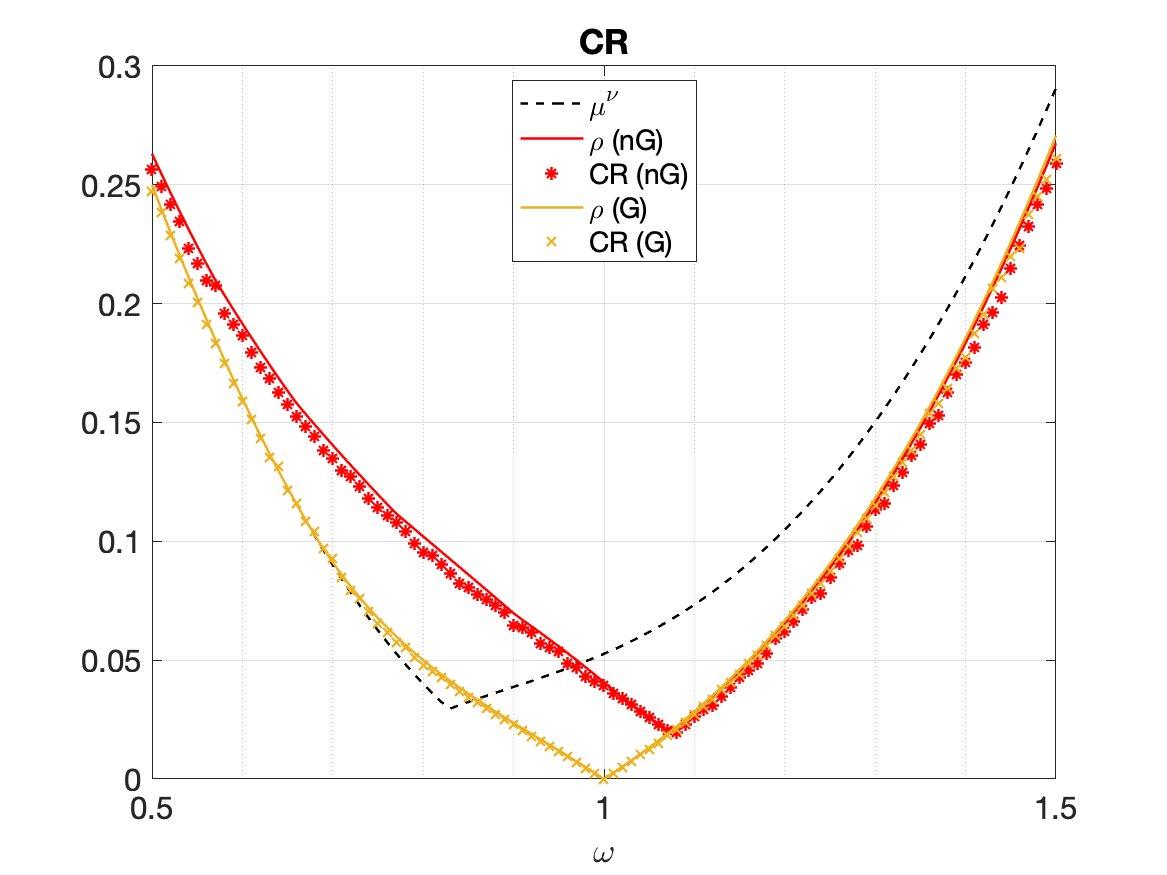}{\figWidth}};
   \draw(7.5,0.0) node[anchor=south west,xshift=-15pt,yshift=-8pt] {\trimfig{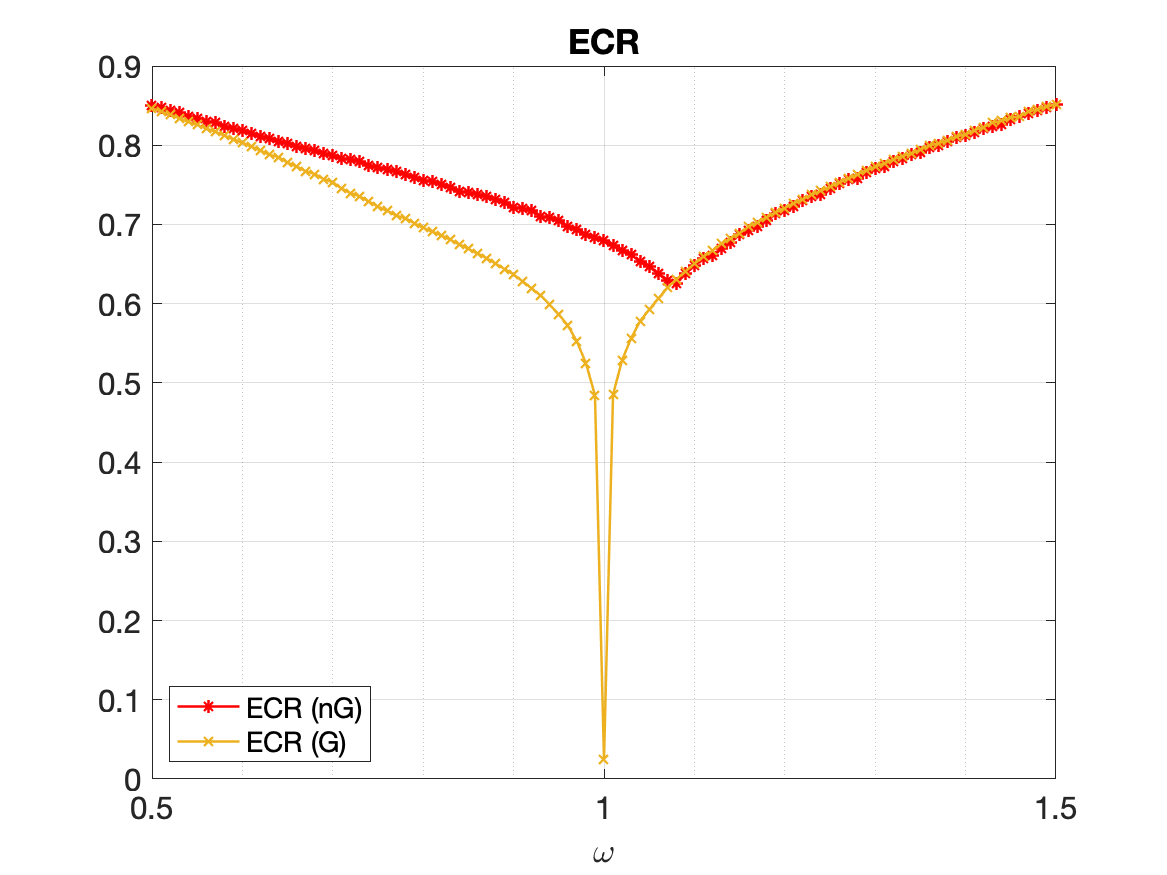}{\figWidth}};
   % title: 
   \draw(7.5,6)  node[draw=\colourTwoD,very thick,fill=\colourTwoD!20,anchor=south,yshift=-4pt,inner sep=2.5pt] {\small \DIM, \ORD, \LEV, \Cc[RB]}; 

   % \circleVortex{11.675}{.8}{direct solver with $\omega=1$};
   \circleLabel{11.675}{.8}{direct solver with $\omega=1$}{DarkGreen}{.25};

%   % small grids 
%   \rotatedGrid{xshift=6cm,yshift=4.5cm}; 
%   \rotatedGrid{xshift=14cm,yshift=2.5cm}; 
%
%   % V cycle
%   \VcycleTwoLevel{xshift= 2.5cm,yshift=4.5cm}; 
%   \VcycleTwoLevel{xshift= 9.5cm,yshift=2.5cm}; 

   % Combined cartoons
   \orderTwoVcycleTwoLevelRotated{xshift=7.5cm,yshift=2.5cm,scale=1}{$\nu\!=\!2$}
   %\orderTwoVcycleTwoLevelRotated{xshift=12.8cm,yshift=2.5cm,scale=1}{$\nu\!=\!2$}

% grid:
% \draw[step=1cm,gray] (0,0) grid (15,6.5);
\end{tikzpicture}
\end{center}
\caption{CR and ECR versus $\omega$. \DIM, \ORD, \LEV\ $V[1, 1]$ cycle with red-black smoothing and \CC[red-black],  
	%(nG=$\ngcg{2}$) and Galerkin (G=$\gcg{2}{2}{2}$)
	non-Galerkin (nG=$\LNAH$) and Galerkin (G=$\LGGH$) 
  coarse-level operators ($N = 32$).}
  \label{fig:CR_2level_2d_rb}
\end{figure}
}

 %2-level
\subsection{Multi-level red-black coarsening}\label{Sec: rbC_multilevel}

We now extend our red-black coarsening algorithm in two dimensions to more than two levels.
In most large scale problems 
%on overset grids 
we will want more than two multigrid levels so that at the coarsest level the equations can be approximately solved by, for example, a Krylov space method without a significant impact
on the overall performance\footnote{Note that for standard coarsening in $d$ dimensions there are $2^{d\, \lmax}$ times fewer points on the coarsest level while for red-black coarsening there are $\sqrt{2}^{d\, \lmax}$ times fewer points so that the cost of solving the coarse level equations decreases rapidly with $\lmax$.}.  
We shall show in this section that a general $(\lmax+1)$-level W cycle with red-black
coarsening can yield convergence rates that are quite comparable to that of a two-level cycle, although we lose the direct-solver property.
Unfortunately W cycles are expensive and we therefore consider variations with much less WU. 

Figure~\ref{fig:redBlackCoarseLevelGrids} shows the coarse-level grids $\ll=1,2,3$ obtained with red-black coarsening.
\cut{If the finest level has grid spacing $h_0$, level one is a rotated grid with spacing $h_1=\sqrt{2} h_0$, level two is a non-rotated grid with $h_2=2 h_0$ while
level three is again a rotated grid with spacing $h_3 = 2^{3/2} h_0$. }
Thus the grid spacing on level $\ll$ is $h_{\ll}= (\sqrt{2})^\ll \, h_0$ and the
number of grid points decreases by a factor of two on successive levels. The number of Fourier modes similarly decreases by a factor of two on successive levels: on level $\ll$ we have $\phi_h(\cdot,\ \thv),\ \thv \in \Theta_{rb}^{(\ll)}$, where
%{\blue We don't seem to use the sets $\Theta^{(\ll)}$. No just show some meaning by itself.} 
%- Letting $\ThetaRB{\ll}$ denote the set of Fourier frequencies $\thv$ on level $\ll=0,1,2,\ldots$, then
\begin{equation}
 	\Theta_{rb}^{(\ll)} = \{\thv \in \Theta:\quad 
 	 	2^{\frac{\ll}{2}} \mathbf{U}^\ll\ \thv \in \Theta\}.
 	\end{equation}
The set $\Theta_{rb}^{(\ll)}$ consists of $\frac{1}{2^{\ll}}$ of $\Theta$ and rotates by 45 degrees as $\ll$ increases (see Figure~\ref{fig: Theta_2d_rb2} for a three-level cycle).

On coarse levels in a general multigrid cycle, we consider both the non-Galerkin and Galerkin coarse-level operators,
the latter are defined recursively using the transfer operators~\eqref{transfer_r2} as
\ba \label{L_G_l}
L_{h_{\ll+1}} =  I_{h_{\ll}}^{h_{\ll+1}}\ L_{h_{\ll}}\ I^{h_{\ll}}_{h_{\ll+1}},\quad
h_{\ll+1} = \sqrt{2} \, h_{\ll},\quad
\ll = 0,\ 1,\ \cdots,\ \lmax-1.
\ea
In practice, however, the Galerkin recursion process can be rather tiresome and unnecessary. As we
go down more levels, the exact coarse-level operator becomes less important, since we are just
solving the error equation on the coarse levels approximately anyway. So besides the real Galerkin
operators (\ref{L_G_l}), we can also use more flexible `Galerkin' operators. In particular, we
consider: 
\begin{enumerate}[label=\roman*.]
	\item(`G1') using the same stencil as the first-generation Galerkin operator (\ref{L_G1}) on every
	coarse level ;
	\item(`Gn') using the non-Galerkin operator on every coarse levels except for the first
	one (the second finest level).
\end{enumerate}

\cut{
We consider a general multi-level $\gamma$-cycle, where recall that $\gamma=1$ corresponds to a V cycle and $\gamma=2$ to a W cycle.
Note that that as $\gamma \rightarrow \infty$, the $\ll=0$ coarse-level correction operator $K_h$ approaches the two-level correction operator $K_h^H$ and
the convergence rate of the multi-level cycle would approach that of the two-level cycle.
Thus if we were to obtain a sufficiently accurate
coarse-level correction with large enough $\gamma$, we would approximately recover the direct-solver property when using the Galerkin operator on level $\ll=1$.  
That is, only one cycle would be needed if $\gamma$ is large enough. 
However, the computational cost of a $\gamma$-cycle increases rapidly with increasing $\gamma$ and so taking $\gamma$ large is not generally useful in practice.
}

\newcommand{\Nmin}{N_{\rm min}}
\newcommand{\SC}{\text{SC}}% std coarsening
\newcommand{\RBC}{\text{RBC}}% red-black coarsening
\newcommand{\ECR}{\text{ECR}}

In practice the number of multigrid coarse levels $\lmax$ is often chosen based on the number of grid points on the finest level so that the coarsest level has roughly some given number of grids points, say $\Nmin$, in each direction.
In order to compare results for red-black coarsening (\RBC) with those for standard coarsening (\SC),
we therefore choose $\lmax$ according to 
\footnote{Assuming $\lmax$ is determined by the size of the fine-level problem so that the work of the direct (order iterative) solve on the coarsest grid is negligible compared to the work on the other levels.}
\ba
&    \lmax[\SC] \approx \log_2\Big(\frac{N}{\Nmin}\Big) + 1, \\
&    \lmax[\RBC] \approx 2\log_2\Big(\frac{N}{\Nmin}\Big) + 1, 
\ea
so that RBC uses approximately twice as many coarse levels as SC in two dimensions. 
Figure~\ref{fig:CR_7level_2d_rb} shows multi-level convergence results with red-black coarsening using V and W cycles, given a fine grid with a fixed size.
The W cycle gives convergence rates very close to zero at $\omega=1$ with Galerkin coarse-level operators, consistent with the direct-solver property of the two-level RBC scheme. 
The V cycle RBC scheme has good convergence rates, but not as good as those from a two-level scheme,
and the V Cycle rates are optimized with over-relaxation $\omega >1$. This indicates that the coarse-level solves in the V cycle are not accurate enough to obtain rates that are close to a direct solver.
However, with respect to the effective convergence rates, the V cycle ($\ECR \approx .6$)
is better than the W cycle ($\ECR \approx .75$) with the Galerkin coarse-level operators, since the 
W cycle requires too many work units.

Work-units estimates for a general $(\lmax+1)$-level
$\gamma$-cycle with red-black coarsening or standard coarsening
in two dimensions
(see \ref{sec: WU}) are 
\ba
    \text{WU}[\SC;\ \gamma,\ \nu] &\approx \frac{4}{4-\gamma}(\nu+2) \qquad \text{(for $\gamma < 4$)},  \\
    \text{WU}[\RBC;\ \gamma,\ \nu] &\approx
    \begin{cases}
        %  \frac{2}{5} N_{min} \frac{1}{2^{\lmax}} +  2(\nu+2)(1-\frac{1}{2^{\lmax}} ),  & \text{$\gamma=1$: V-cycle},  \\
        2(\nu+2) ,  & \text{$\gamma=1$: V cycle},  \\
        {\lmax}(\nu+2) ,    & \text{$\gamma=2$: W cycle} .
    \end{cases}
\ea
It is apparent that a W cycle with red-black coarsening involves too much computational work to be generally useful for large values of levels since the work-units scale with ${\lmax}$.

% 2D, 2nd-order, 7-level, RB coarsening:
{
\newcommand{\figWidth}{6cm}% height 
\newcommand{\trimfig}[2]{\trimh{#1}{#2}{.0}{.0}{.0}{.0}}
\begin{figure}[h!]
\begin{center}
\begin{tikzpicture}[scale=1]
  \useasboundingbox (0,0.5) rectangle (15,6.);  % set the bounding box (so we have less surrounding white space)
% V ----------------------
   \draw(-1,0.0) node[anchor=south west,xshift=-15pt,yshift=-8pt] {\trimfig{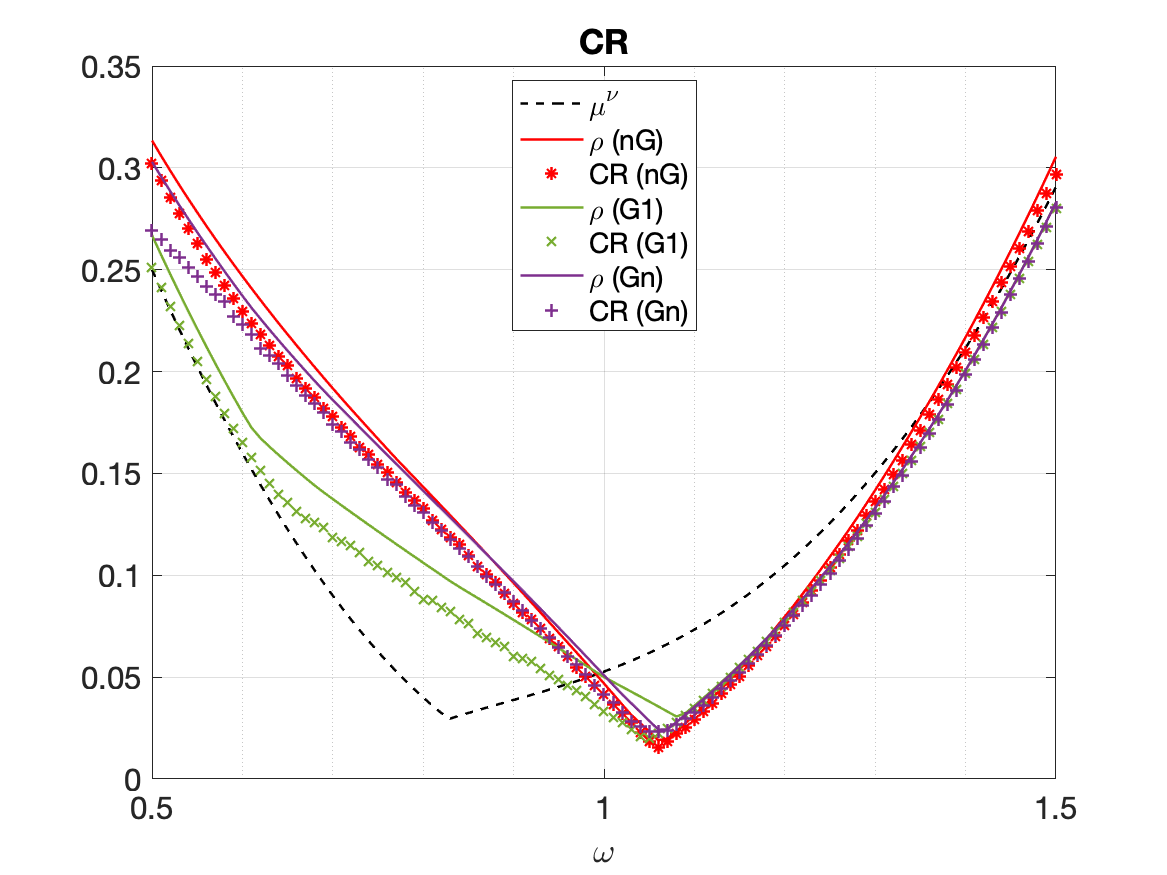}{\figWidth}};
   %\draw(7.5,7.0) node[anchor=south west,xshift=-15pt,yshift=-8pt] {\trimfig{ECR_2drt_l6gamma10_nu2}{\figWidth}};   

   % \circleVortex{4.5}{0.7}{$ECR=0.59$};
   \circleLabelLeft{3.5}{0.7}{$ECR=0.59$}{purple}{.3};

   % Combined cartoons
   \orderTwoVcycleThreeLevelRotated{xshift=4.6cm,yshift=4.6cm,scale=1}{$\nu\!=\!2$}{$l6$};
   %\orderTwoVcycleThreeLevelRotated{xshift=5.5cm,yshift=4.8cm,scale=1}{$\nu\!=\!2$}{$l6$};
   %\orderTwoVcycleThreeLevelRotated{xshift=11.25cm,yshift=4.5cm,scale=1}{$\nu\!=\!2$}{$l6$};
% W ----------------------
   \draw(7.5,0.0) node[anchor=south west,xshift=-15pt,yshift=-8pt] {\trimfig{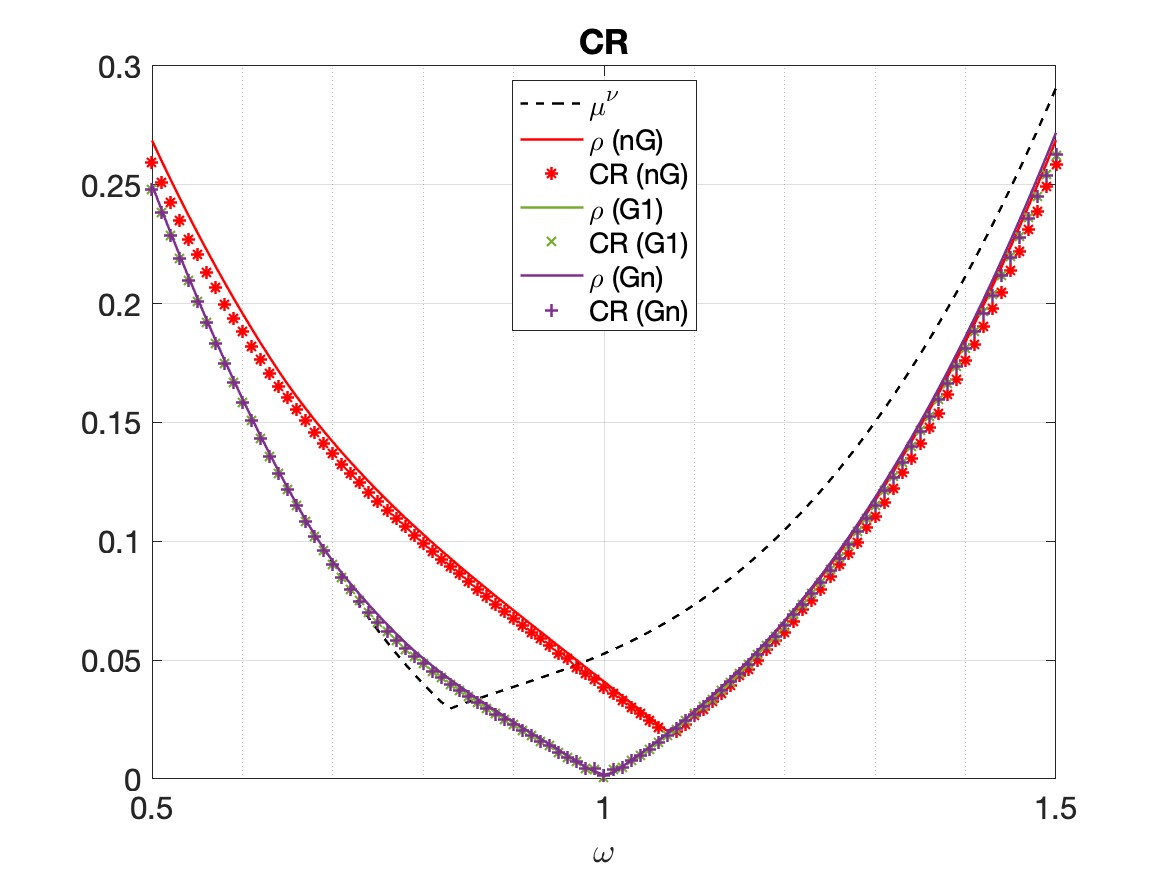}{\figWidth}};
   %\draw(7.5,0.0) node[anchor=south west,xshift=-15pt,yshift=-8pt] {\trimfig{ECR_2drt_l6gamma20_nu2}{\figWidth}};

   % \circleVortex{11.7}{0.5}{$ECR=0.75\ (\omega^*=1)$};
   % \circleLabelLeft[x][y][label][colour][radius]
   % \circleLabelLeft{11.7}{0.6}{$ECR=0.75\ (\omega^*=1)$}{purple}{.3};
   \circleLabel{11.7}{0.6}{$ECR=0.75\ (\omega^*=1)$}{purple}{.3};
   % Combined cartoons
   \orderTwoWcycleThreeLevelRotated{xshift=13.1cm,yshift=4.6cm,scale=1}{$\nu\!=\!2$}{$l6$};
   %\orderTwoWcycleThreeLevelRotated{xshift=5.5cm,yshift=4.8cm,scale=1}{$\nu\!=\!2$}{$l6$};
   %\orderTwoWcycleThreeLevelRotated{xshift=11.25cm,yshift=4.5cm,scale=1}{$\nu\!=\!2$}{$l6$};
   
   % title:
   \draw(7.5,5.8)  node[draw=\colourTwoD,very thick,fill=\colourTwoD!20,anchor=south,yshift=-4pt] {\small \DIM, \ORD, \LEV[7], \Cc[RB]}; 

% grid:
  %\draw[step=1cm,gray] (0,0) grid (15,6.5);
\end{tikzpicture}
\end{center}
  \caption{CR versus $\omega$. \DIM, \LEV[7]\ cycles with red-black smoothing and \textcolor{olive}{\textbf{red-black} coarsening}, 
  	non-Galerkin (nG=$\LNAH$) and Galerkin (G1=$\LGAAH$, Gn) 
  	coarse-level operators ($N = 128$).
  Left: $\CYC{V}[1, 1]$ cycle, right: $\CYC{W}[1, 1]$ cycle.}
  \label{fig:CR_7level_2d_rb}
\end{figure}
}
% 2D, W-cycles, RB coarsening:
{
\newcommand{\figWidth}{6cm}% height 
\newcommand{\trimfig}[2]{\trimh{#1}{#2}{.0}{.0}{.0}{.0}}
\begin{figure}[h!]
\begin{center}
\begin{tikzpicture}[scale=1]
  \useasboundingbox (0,0.5) rectangle (15,6);  % set the bounding box (so we have less surrounding white space)

   \draw(-1,0.0) node[anchor=south west,xshift=-15pt,yshift=-8pt] {\trimfig{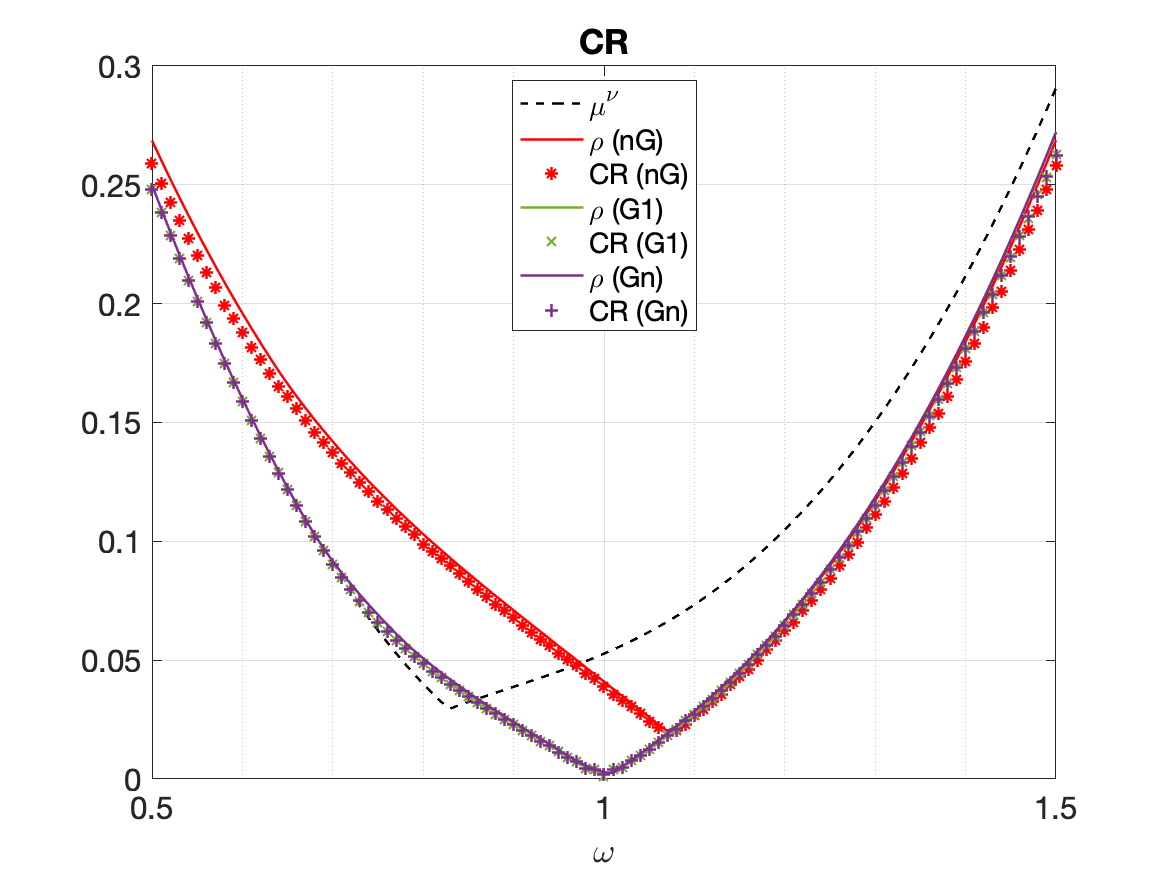}{\figWidth}};
   %\draw(7.5,0.0) node[anchor=south west,xshift=-15pt,yshift=-8pt] {\trimfig{ECR_2drt_l6gamma21_nu2}{\figWidth}};
   %\circleVortex{4.2}{0.5}{$ECR=0.56\ (\omega^*=1)$};
   \circleLabel{3.2}{0.65}{$ECR=0.56\ (\omega^*=1)$}{purple}{.3};
   
    \draw(7.5,0.0) node[anchor=south west,xshift=-15pt,yshift=-8pt] {\trimfig{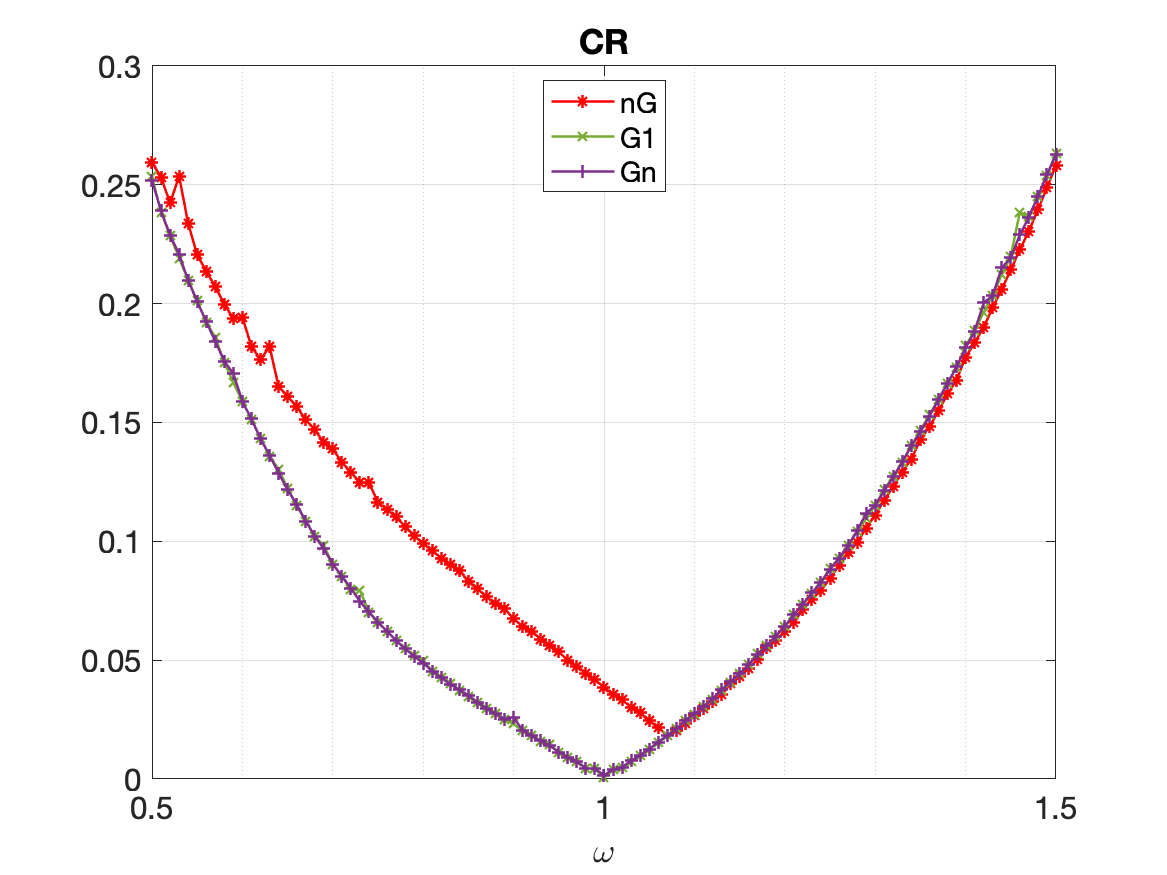}{\figWidth}};

   %\circleVortex{11.6}{0.35}{$ECR=0.59\ (\omega^*=1)$};
   \circleLabel{11.6}{0.65}{$ECR=0.59\ (\omega^*=1)$}{purple}{.3};
    
    \rotatedGrid{xshift=7.5cm,yshift=2.5cm}; %2o   
% title:
\draw(7.5,5.8)  node[draw=\colourTwoD,very thick,fill=\colourTwoD!20,anchor=south,yshift=-4pt] {\small \DIM, variant \CYC{`W'} cycles, \Cc[RB]}; 
% grid:
% \draw[step=1cm,gray] (0,0) grid (16,6.5);
\end{tikzpicture}
\end{center}
  	\caption{CR versus $\omega$. \DIM, 
  		$\CYC{`W'}[1, 1]$ cycles with red-black smoothing and \CC[red-black], non-Galerkin (nG=$\LNAH$) and Galerkin (G1=$\LGAAH$, Gn) coarse-level operators ($N = 128$). 
  Left: \LEV[7] \CYC{Wn} cycle with \CC[red-black];
  right: \LEV[5] \CYC{W} cycle with \CC[variable]. 
  %(\CC[red-black] only for the first two coarse levels). 
  %and standard coarsening after the third level)
  }
  \label{fig:CR_WW_2d_rb}
\end{figure}
}
{
\newcommand{\figWidth}{6cm}% height 
\newcommand{\trimfig}[2]{\trimh{#1}{#2}{.0}{.0}{.0}{.0}}
\begin{figure}[h!]
\begin{center}
\begin{tikzpicture}[scale=1]
  \useasboundingbox (0,0.5) rectangle (15,5.65);  % set the bounding box (so we have less surrounding white space)

   \draw(4.0,0.0) node[anchor=south west,xshift=-15pt,yshift=-8pt] {\trimfig{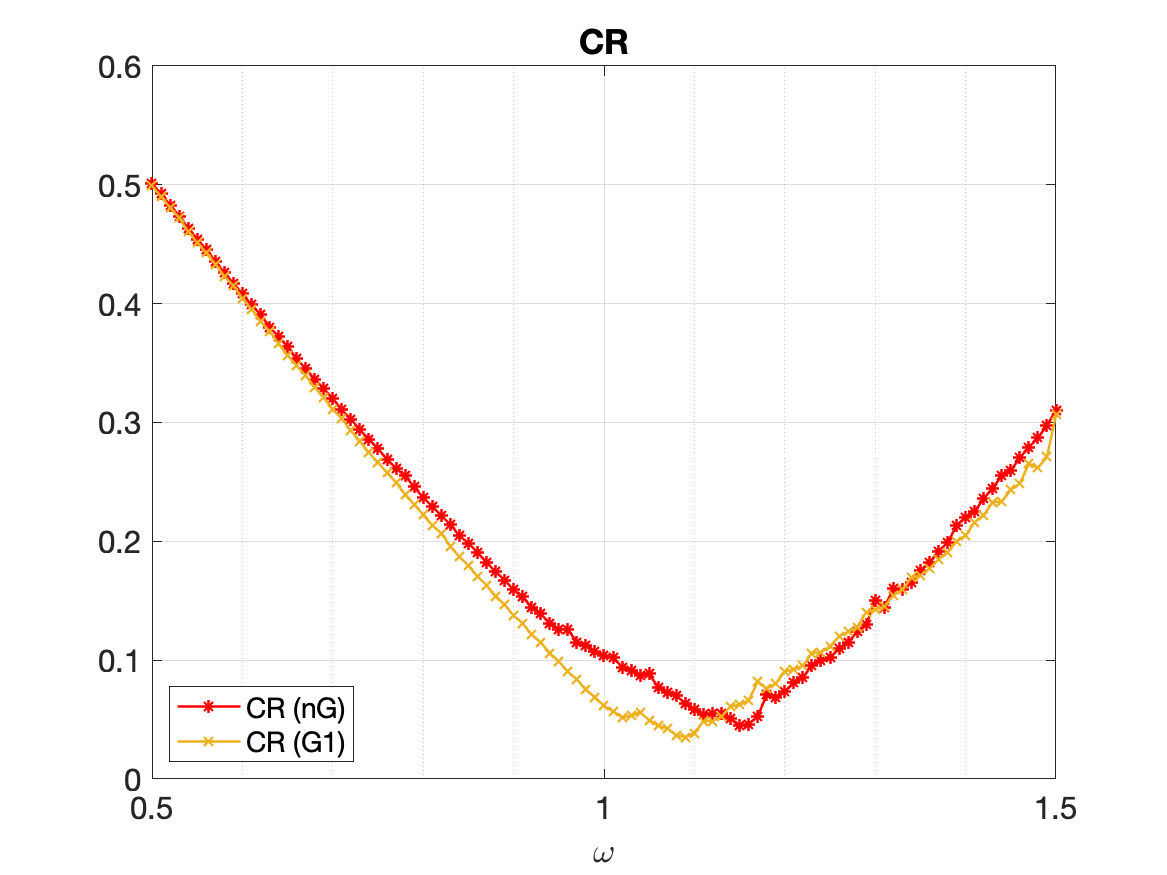}{\figWidth}};
   %\draw(7.5,0.0) node[anchor=south west,xshift=-15pt,yshift=-8pt] {\trimfig{ECR_2d_l3gamma1_nu2}{\figWidth}};
   %\circleVortex{8.7}{0.7}{$ECR=0.53$};   
   \circleLabel{8.7}{0.7}{$ECR=0.53$}{purple}{.3};
    \orderTwoVcycleThreeLevel{xshift=9cm,yshift=4.5cm,scale=1}{$\nu\!=\!2$}{$l4$};
    
    \draw(7.75,2.75)  node[draw=\colourTwoD,very thick,fill=\colourTwoD!20,anchor=south,yshift=-4pt] {\small \DIM, \LEV[4] \CYC{V} cycle, \Cc}; 
% grid:
% \draw[step=1cm,gray] (0,0) grid (16,6);
\end{tikzpicture}
\end{center}
  \caption{CR versus $\omega$. \DIM, \LEV[4] $\CYC{V}[1, 1]$ cycle with red-black smoothing and \CC, non-Galerkin (nG=$\LNAH$) and Galerkin (G1=$\LGAAH$) coarse-level operators ($N = 128$).}
  \label{fig:CR_4levelV_2d}
\end{figure}
}

As a compromise, we consider more flexible ``W''-like cycles that do not require as much
work as a standard W cycle, but
nevertheless provide approximate coarse-level solves accurate enough to ensure CRs comparable to a full W cycle.
Define the ``Wn'' cycle to be one where the cycle parameter 
$\gamma$ depends on the level 
as $\gamma = 2$ for $\ll = 0,\ 1$ and $\gamma = 1$ for $\ll \geq 2$.
As shown in Figure~\ref{fig:CR_WW_2d_rb} (left),
the convergence rates for the Wn cycle are almost the same as those of full W cycles but the ECRs are significantly better.
As an alternative approach to reduce the cost of red-black coarsening and its many levels,
we evaluate a ``variable coarsening'' (VC) algorithm that uses
red-black coarsening to form the first two coarse levels, $\ll=1,2$, and then standard coarsening for the remaining levels, 
since the idea is just to approximate the two-level cycle by solving the coarse-level equation recursively with coarser levels. 
Figure~\ref{fig:CR_WW_2d_rb} (right), shows results for a $(\log_2{\frac{N}{N_{min}}} + 2)$-level W cycle with variable coarsening.
It is seen that this VC approach gives comparable effective convergence rates to a Wn cycle.
The RBC results in Figures~\ref{fig:CR_7level_2d_rb} and~\ref{fig:CR_WW_2d_rb}
can be compared to the multi-level results with standard coarsening for the same fine-level problem in Figure~\ref{fig:CR_4levelV_2d}, which shows a best ECR that is similar to the best ECRs from standard coarsening.
In particular, with red-black coarsening, we do not need to adjust the over-relaxation parameter since the optimal $\omega$ is approximately 1; with standard coarsening, in contrast, over-relaxation (and tuning the parameter $\omega$) is needed for faster convergence and comparable ECRs.

%---------------------------------------------------------------------------------------------------------------------------------------
\subsection{Red-black coarsening in three dimensions}\label{Sec: rbC_3d}
Now that we've extended red-black coarsening to two dimensions, it is natural to ask whether it can be further extended to higher dimensions. As seen before, in one and two dimensions, by red-black coarsening the coarse grid consists half of the points in the fine grid, and in the mean time the Fourier modes on the coarse grid corresponding to $\Theta^{\text{low}}$ consist half of the modes on the fine grid corresponding to $\Theta$. 
In one dimension we have $H = 2h$, and $H = \sqrt{2}h$ in two dimensions. 
For red-black smoothing in particular, the structure of red-black coarsening is special in that the smoother and coarse-level correction operator have the same eigenspaces.
We consider two factors contributing to fast convergence with red-black coarsening: 
\begin{enumerate}[label=\roman*.]
	\item the section of $\Theta$ into two halves with equal size and thus the balance of low- and high- frequency modes, regardless of the smoother;
	\item the alignment of the coarse-level correction operator with the red-black smoother.
\end{enumerate}
This gives us two perspectives to try to extend to three dimensions (3D). 

The first idea is to look directly at the Fourier space and aim to construct a cube with half the volume of $\Theta$ to be $\Theta^{\text{low}}$, which will require a coarse grid with $H = \sqrt[3]{2}h$. 
However, there will be no construction like $\ThetaLowRB$ in Figure \ref{fig: theta_2d} (left) that gives such special properties. Then a natural way to try would be to construct the coarse grid 
$
\grid G_H = \Big\{  \xv =\jv\, H:~ \jv \in \Integer^3 \Big\}
$
with spacing $H \approx \sqrt[3]{2}h$. 
\kl{This idea will be pursued in Section~\ref{M3: rC}. But unfortunately, as we shall see, this method does not seem to bring any advantage over standard coarsening. 
We will give an explanation on why this balancing property is not the main factor to improve convergence in Section~\ref{Sec: rC_J} via an analysis with the simple example of a Jacobi smoother.}
%This may be explained partially by the fact that $H = \sqrt[3]{2}h$ can not be reached exactly in this manner. 

The second idea is to continue with red-black coarsening and set the coarse grid as the black points of the fine grid $\grid G_H = \grid G_h^B$. 
This corresponds to red-black reduction for the matrices of the 3D problem. 
However, the resulting coarse grid still has $H = \sqrt{2}h$ (instead of $\sqrt[3]{2}h$), thus 
red-black coarsening in three dimensions would no longer have the balancing property; meanwhile, the coarse grid would no longer be orthogonal.
(In any dimension $d$, the red-black section of the grid points, both for red-black smoothing and red-black coarsening, will mix up $\thv \in \Theta$
%have eigenspaces$\mathbb{E}_{\thv} = \text{span} \{\phi_h (\cdot,\ \thv),\quad \phi_h (\cdot,\ \bar{\thv})\},\ \thv \in \Theta,$ 
with $\bar{\thv} \defeq [\bar{\th}_1,\ \bar{\th}_2,\ \cdots,\ \bar{\th}_d]^T$.)
However, since the set of black points as a grid is non-orthogonal, 
it is not straightforward to correspond $\thv \in \Theta^{\text{low}}$ to
%this section of $\Theta$ will not directly correspond to 
the eigen-modes on the coarse grid (in such a manner as in \eqref{Th_L_2d_rb}). 
%[\textcolor{red}{What are the eigenmodes on a scattered grid?}])
A non-orthogonal grid would also bring many more complications: in particular, the difference stencil would consists too many points; meanwhile the transfer operators between the fine and coarse grids would be hard to keep track of. But nevertheless red-black coarsening can be done in any dimension, at least in theory, by means of red-black reduction as in \ref{Sec: cyclic_2d}, where the transfer operators as well as the Galerkin coarse-level operators can be constructed automatically.

As we've discussed, these two ideas lead to the same red-black coarsening in one and two dimensions, but they do not coincide in three and more dimensions. It would be interesting to ask what are the special geometric properties of the one and two dimensional spaces, that leads to such a distinction.

 % multi-level; 3D

\bigskip
% topic-3: COARSENING BY A GENERAL FACTOR -----------------------------------------------------------
\section{Coarsening by a general factor $r$}\label{M3: rC}
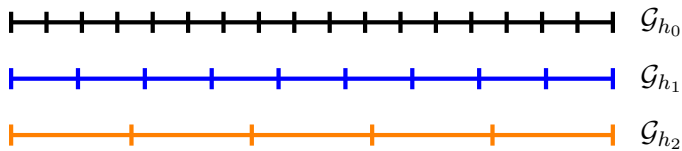
\begin{figure}[hbt]
  \centering
  \begin{tikzpicture}
    \useasboundingbox (0,.8) rectangle (8,3);
     % \draw[<-,ultra thick,red,dash pattern= on 8pt off 3pt] (0,3) -- (2,3);
     %
    \begin{scope}[yshift=-0.5cm]
      \draw[ultra thick ,black] (0,3) -- (8,3)  node [black,anchor=west,xshift=6pt] {$\Gc_{h_0}$};
      \foreach \x in {0,0.470588235294118,...,8}% 18 pts, 17 cells 
        \draw[ultra thick,black] (\x,3cm-4pt) -- (\x,3cm +4pt); 
    \end{scope}

    \begin{scope}[yshift=.75cm]
      \draw[ultra thick,blue] (0,1) -- (8,1)  node [black,anchor=west,xshift=6pt] {$\Gc_{h_1}$};
      % \foreach \x in {0,1,2,3,4,5,6,7,8}
      \foreach \x in {0,.8888888,...,8}% 9 cells
        \draw[ultra thick,blue] (\x,1cm-4pt) -- (\x,1cm +4pt); 
    \end{scope}

    \begin{scope}[yshift=0cm]
      \draw[ultra thick,orange] (0,1) -- (8,1)  node [black,anchor=west,xshift=6pt] {$\Gc_{h_2}$};
      % \foreach \x in {0,1,2,3,4,5,6,7,8}
      \foreach \x in {0,1.5999,...,8}% 5 cells 
        \draw[ultra thick,orange] (\x,1cm-4pt) -- (\x,1cm +4pt); 
    \end{scope}

  \end{tikzpicture}
  \caption{Coarsening by a general factor r. The target grid-spacing ratio is $\rTarget = 2$.
  	  %The number of grid points on each level are $N_0+1=18,\ N_1+1=10,\ N_2+1=6$.
      The grid-spacings are $h_0=1/17$, $h_1=1/9$, $h_2=1/5$. 
  }
  \label{fig:generalCoarsening}
\end{figure}

As noted in the introduction, for flexibility in the generation of high-quality overset grids,
we do not wish to place constraints on the number of grid points of each component grid to force it to support coarsening by sufficient factors of two to allow some number of levels for multigrid.
We thus are interested in allowing coarsening by a factor $r$, for instance in a neighborhood of two other than exactly two (see Figure \ref{fig:generalCoarsening} for a 1D example), and in this section we investigate multigrid algorithms based on factor-$r$ coarsening. 
There is a second reason to investigate general factor-$r$ coarsening, which arises
from our results of good convergence rates from red-black coarsening with a coarsening factor $r=\sqrt{2}$ in two dimensions.
It is thus of interest, in general, to see if there is a value for $r$ 
%different from two 
that leads to better multigrid convergence rates than those for $r=2$. 
We consider the model problem (\ref{MP_h}) with second-order accuracy, on the fine grid $\grid G_h$ as given in (\ref{G_h_N}).
We start with the one-dimensional case in Section~\ref{sec:generalCoarseningOneDimension}, and then extend the construction to multiple dimensions in a straightforward fashion.
A local Fourier analysis of a multigrid cycle with general factor-$r$ coarsening is then studied in Section~\ref{Sec: rC_J}
in order to explain the results. 

% -----------------------------------------------------------------------------------------------
\subsection{General factor-$r$ coarsening in one dimension} \label{sec:generalCoarseningOneDimension}

Let $\rTarget\ (>1)$ denote the target grid coarsening factor, and consider the fine grid $\grid
G_h$ as given in (\ref{G_h_N}) with $d=1$ and size $N$.  The goal is to construct a coarse grid with
grid-spacing $H = r h$, where $r \approx \rTarget$ is the actual coarsening factor.
Aiming at a ratio $\frac{H}{h} \approx \rTarget$, set 
\ba
   M \eqdef \left\lfloor \frac{N}{\rTarget} \right\rfloor,
   \quad H = \frac{1}{M},
   \quad r \eqdef \frac{H}{h} = \frac{N}{M}, 
\ea 
where the floor
function $\lfloor a \rfloor$ denotes the largest integer less than or equal to $a \in \mathbb{R}$.
Given the fine grid $\grid G_h$ and coarse grid $\grid G_H = \{x_k^H = kH:\ k = 0,\ 1,\ \cdots,\ M \}$, we need to construct the restriction operator $I_h^H$ and interpolation operator $I_H^h$, as well as the coarse level operator $L_H$.
Given each $x^h_j \in \grid G_h$, we find out where it lies between the grid points on $\grid G_H$,
\begin{equation}
 x^h_j \in \Big[ x^H_k,\ x^H_{k+1} \Big],
     \quad k = \Big\lfloor \frac{x^h_j}{H} \Big\rfloor =  \left\lfloor  \frac{j}{r} \right\rfloor,
\end{equation}
and construct the interpolation operator $I_H^h$ with coefficients $r_{jk}$ from the linear interpolant,
% ; next apply linear interpolation, from which we construct the coefficients $r_{jk}$ in the interpolation operator $I_H^h$: 
\begin{equation}
  I_H^h u_H (x^h_j) = (1-r_{jk})\ u_h(x^H_k) + r_{jk}\ u_h(x^H_{k+1}),\quad
      r_{jk} \eqdef \frac{x_j^h - x^H_k}{H} = \frac{j}{r} - k,
\end{equation}
for any grid function $u_H$ on the coarse grid.
% Then we construct its adjoint as the restriction operator: 
The restriction operator $I_h^H$ is defined as adjoint of $I_H^h$, 
\begin{equation} \label{interp_r}
I_h^H = (I_H^h)^* = \frac{1}{r} (I_H^h)^T.
\end{equation}
(This is the full-weighting operator when $r=2$.)
% On $\grid G_H$,
The non-Galerkin coarse-level operator $L_H$ uses the same stencil as $L_h$, while
the Galerkin coarse-level operator is based on the transfer operators just constructed,
$L_H = I_h^H L_h I_H^h$.
Note that when $r = 2$, these two would be the same in one dimension.

Figure~\ref{fig:CR_r_1d} presents convergence rates and effective convergence rates
of the $(h,\ rh)$ two-level multigrid cycles with a red-black smoother in one dimension for the model problem with $L_h$ given in~\eqref{L_h_1}.
Results are given with Galerkin and non-Galerkin coarse-level operators using a $V[1,1]$ cycle. 
The first thing to note is the convergence rate tends to generally increase with increasing r.
The convergence rate goes to zero at $r=1$ since the coarse grid approaches the fine grid and
there is an exact solve on the coarse level. There is a singular behavior in the convergence rate
at $r=2$ where the one-dimensional multigrid algorithm becomes a direct solver. 
  The fact that values of $r$ near two do not exhibit \kP{nearly} the optimal convergence behavior seen exactly at $r=2$ indicates that 
  it is not only the coarsening factor two, but
  the special structure of the red-black coarse-level correction that results in the optimal convergence at $r=2$.
The ECR also shows the trend of generally increasing with increasing $r$. The ECR with the Galerkin
coarse-level operator is quite flat for $r\in[1.2,\ 3]$ apart from the singular behavior at $r=2$ where the ECR is zero.
One important conclusion to draw from these first results is that there is a range $r$ near $r=2$
where the multigrid scheme behaves well.
{
\newcommand{\figWidth}{6cm}% height 
\newcommand{\trimfig}[2]{\trimh{#1}{#2}{.0}{.0}{.0}{.0}}
\begin{figure}[h!]
\begin{center}
\begin{tikzpicture}[scale=1]
  \useasboundingbox (0,0.5) rectangle (15,6.4);  % set the bounding box (so we have less surrounding white space)

   \draw(0.0,0.0) node[anchor=south west,xshift=-15pt,yshift=-8pt] {\trimfig{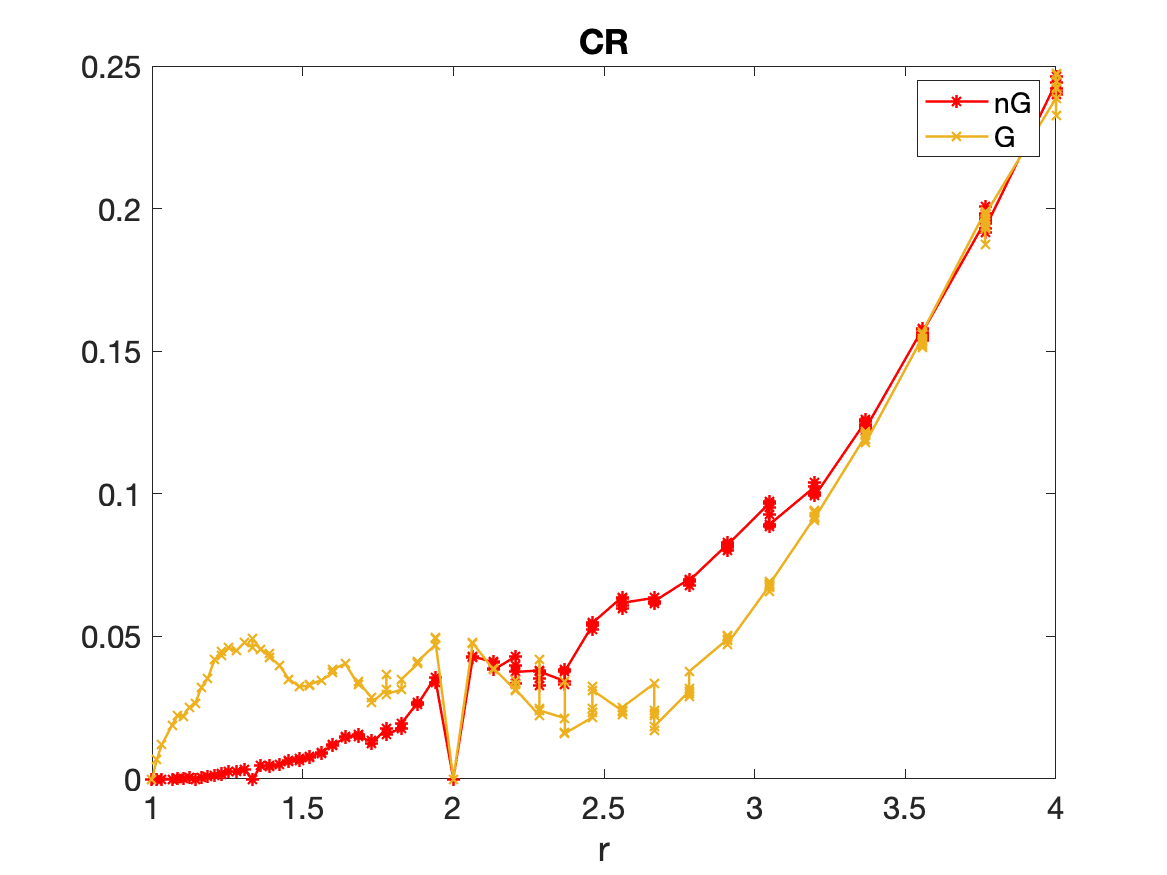}{\figWidth}};
   \draw(7.5,0.0) node[anchor=south west,xshift=-15pt,yshift=-8pt] {\trimfig{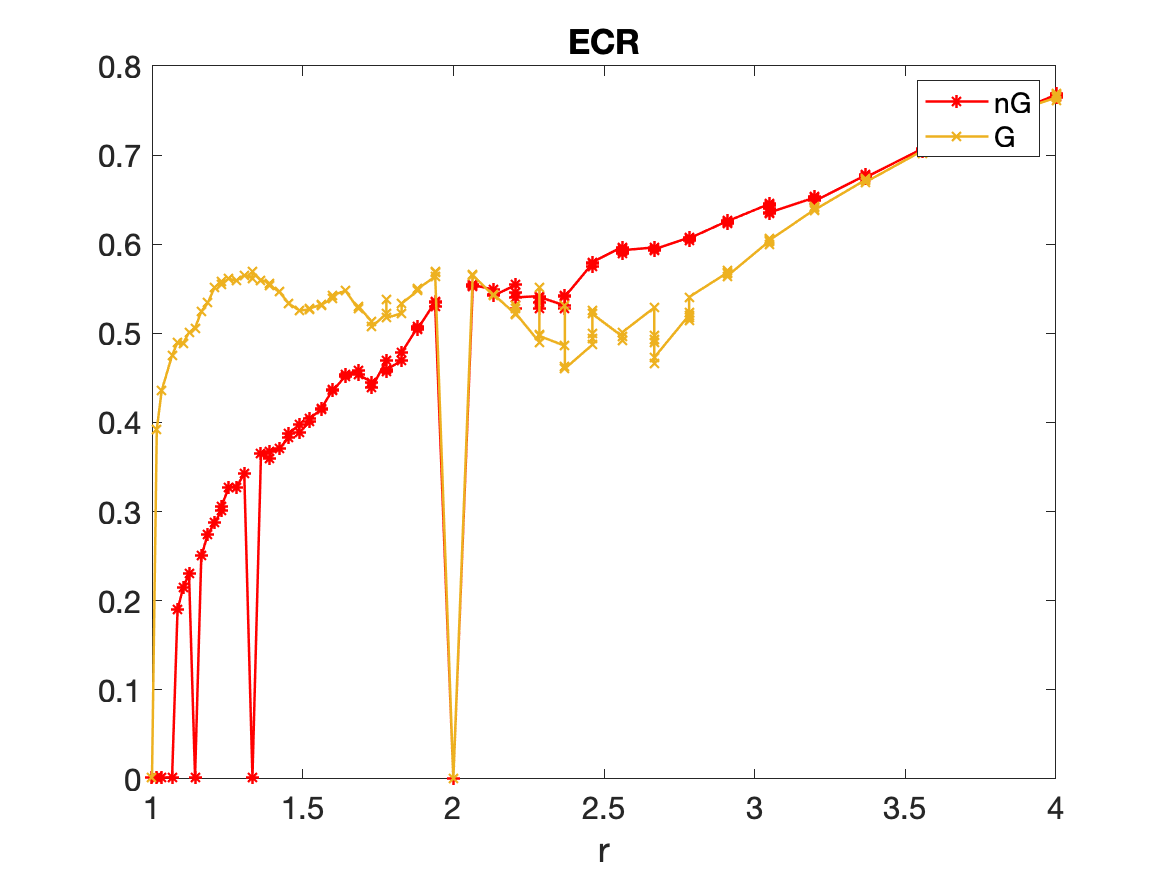}{\figWidth}};

   %\orderTwoVcycleTwoLevel{xshift=1.7cm,yshift=4.75cm,scale=1}{$\nu\!=\!2$}; % wrong: std-C
   \VcycleTwoLevel{xshift=9.5cm,yshift=4.75cm,scale=1}{$\nu\!=\!2$};

  % \circleVortex{10.6}{3.8}{*fix me* Good ECR over a range of $\tilde{r}$}; 

   %\circleVortex{10.65}{.7}{direct solver at $r=2$}; 
   \circleLabel{10.65}{.7}{direct solver at $r=2$}{DarkGreen}{.3};

    % title:
   \draw(7.5,5.9)  node[draw=\colourTwoD,very thick,fill=\colourTwoD!20,anchor=south,yshift=-4pt,inner sep=2.5pt]
             {\small \CC[factor-$r$], \DIM[1], \LEV}; 

% grid:
% \draw[step=1cm,gray] (0,0) grid (16,6);
\end{tikzpicture}
\end{center}
\caption{CR and ECR versus $r$.
  \DIM[1],
  $(h,\ H)$ \LEV\ $V[1, 1]$ cycle with red-black smoothing
  ($\omega = 1$), $\frac{H}{h} = r$ ($N = 64$).
  (nG=$\LNAH$, G=$\LGGH$.)
  %{\red make x-axis labels consisent between different factor r results.}
}

  \label{fig:CR_r_1d}
\end{figure}
}

Based on the two-level construction, % including the coarse level grid and operator, as well as transfers
it is straightforward to construct a general multi-level algorithm with factor-$r$ coarsening. Given the target grid coarsening factor $\rTarget$, let $r_\ll$ denote the actual coarsening factor that is used coarsening from level $\ll$ to level $\ll+1$. If there are $N_0 \equiv N = \frac{1}{h}$ grid points on the finest level, then the number of grid points $N_\ll$ at level $\ll$, and the corresponding coarsening factor $r_\ll$, are given by the formulae 
\ba
       N_{\ll+1} = \left\lfloor \frac{N_{\ll}}{\rTarget}  \right\rfloor, \quad
      r_{\ll} = \frac{N_{\ll}}{ N_{\ll+1} } = \frac{h_{\ll+1}}{h_{\ll}},
      \quad \ll=0,\ 1,\ \cdots,\ \lmax-1.
\ea
Multi-level results in one dimension are shown in Figure~\ref{fig:CR_r_l_1d}, for  $V[1,1]$ cycles with  non-Galerkin coarse-level operators.
We allow the number of multigrid levels to vary with $r$ since 
fixing the number of levels would result in more expensive solves
on the coarsest level for $r$ closer to one.
Therefore, for a more fair comparison in terms of the optimal $\rTarget$, the number of levels is chosen so that the coarsest level grid
has an approximately fixed number of grid points, $\Nmin$, 
\ba
   \lmax(\rTarget) \approx \log_{\rTarget} \Big( \frac{N_0}{\Nmin} \Big) \propto \frac{1}{\log \rTarget}. \label{eq:numberOfLevelsVersusR}
\ea
This means that the number of levels increases as $r$ decreases as shown in the insert in Figure~\ref{fig:CR_r_l_1d}. 
%The multi-level results with the non-Galerkin operator are shown in Figure~\ref{fig:CR_r_l_1d} for  $V[1,1]$ cycles.
As with the two-level results, the convergence rate tends to increase as r increases. The convergence rate is again zero at $r=2$ where this particular scheme converges in one cycle.
The effective convergence rate\footnote{Work-units estimates for general multi-level cycles with factor-$r$ coarsening are given in Section~\ref{Sec: rC_J}.}
is quite flat for a large range of $r\in[1.5,3.5]$, which is the behavior we desire.
% (A WU estimate for general multilevel cycles is given in Section~\ref{Sec: rC_J}.)
It is not surprising that the optimal choice for $\rTarget$ would be 2 in one dimension.
Neglecting the singular behavior at $r=2$, however, there is no obvious optimal value of the coarsening factor based on the ECR. 
{
\newcommand{\figWidth}{6cm}% height 
\newcommand{\trimfig}[2]{\trimh{#1}{#2}{.0}{.0}{.0}{.0}}
\newcommand{\figWidtha}{2.5cm}% height 
\newcommand{\trimfiga}[2]{\trimhb{#1}{#2}{.0}{.0}{.0}{.0}}
\begin{figure}[h!]
\begin{center}
\begin{tikzpicture}[scale=1]
  \useasboundingbox (0,0.5) rectangle (16,6.2);  % set the bounding box (so we have less surrounding white space)

   \draw(0.0,0.0) node[anchor=south west,xshift=-15pt,yshift=-8pt] {\trimfig{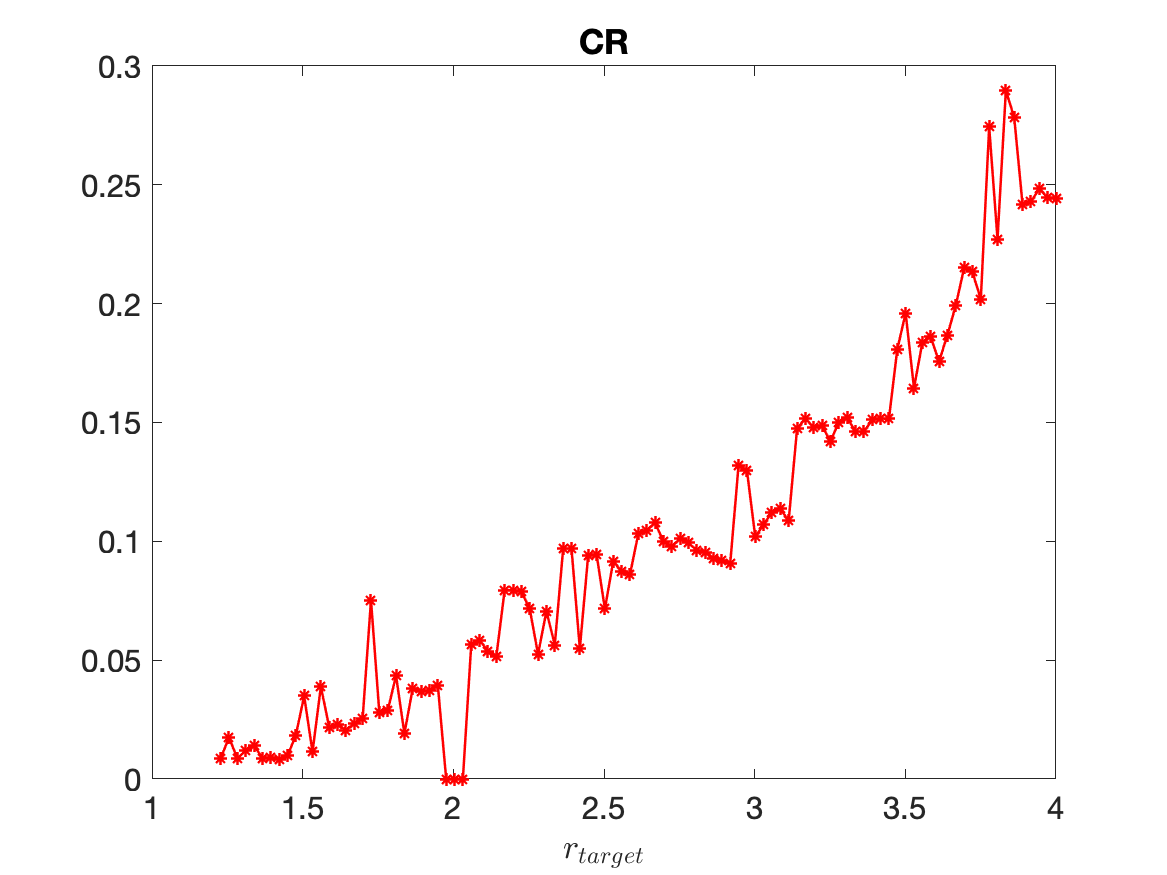}{\figWidth}};
   \draw(7.5,0.0) node[anchor=south west,xshift=-15pt,yshift=-8pt] {\trimfig{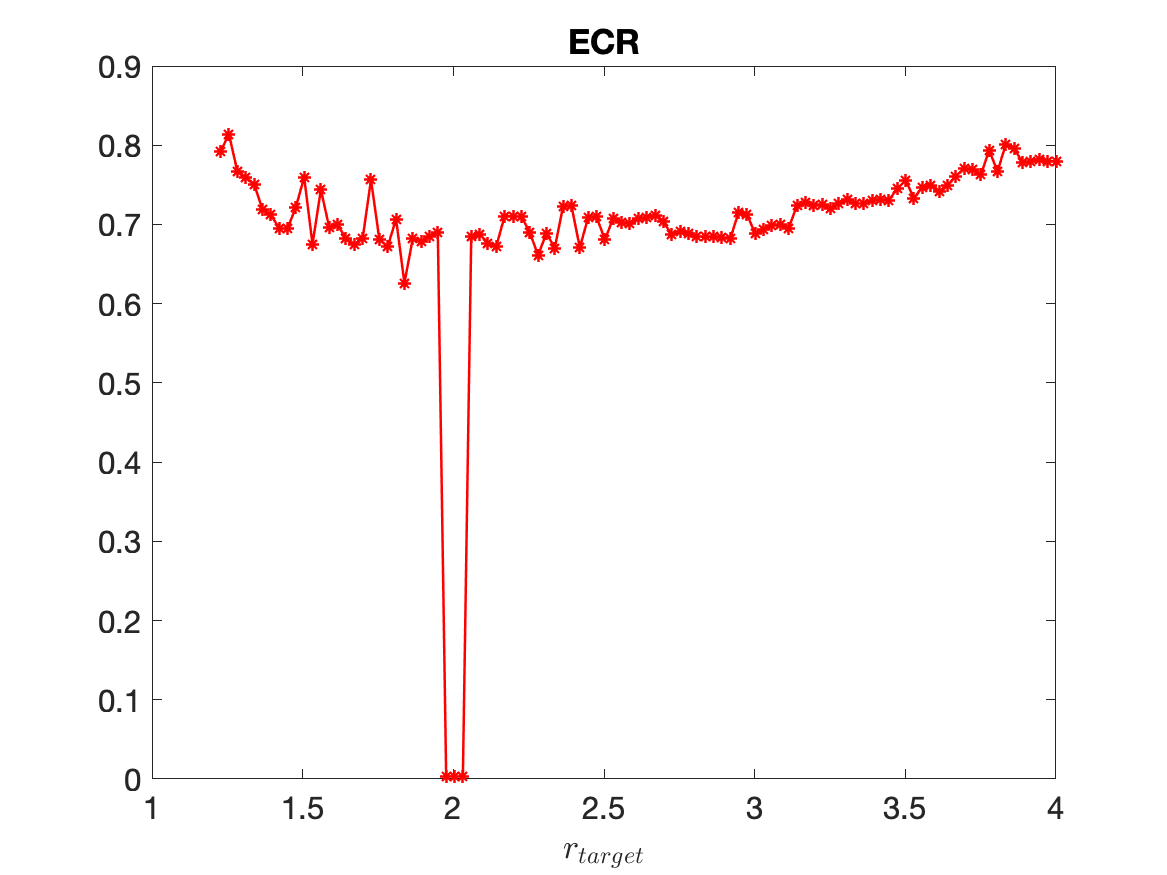}{\figWidth}};

%    \orderTwoVcycleThreeLevel{xshift=1.7cm,yshift=4.75cm,scale=1}{$\nu\!=\!2$}{$l(\tilde{r})$};
  \VcycleThreeLevel{xshift=1.7cm,yshift=4.75cm,scale=1}{$\nu\!=\!2$}{${\lmax}(\rTarget)$};  

   \ellipseLabel{11.5}{4.2}{Good ECR}{violet}{2.2}{.35};

   % \circleVortex{10.6}{4.25}{*fix me* Good ECR over a range of $\tilde{r}$}; 

   \circleLabel{10.65}{.7}{direct solver at $r=2$}{DarkGreen}{.3};

    % title:
   \draw(7.5,5.9)  node[draw=\colourTwoD,very thick,fill=\colourTwoD!20,anchor=south,yshift=-4pt,inner sep=2.5pt]
             {\small \CC[factor-$r$], \DIM[1], \LEV[multi]}; 

   \begin{scope}[xshift=11.5cm,yshift=1.25cm]
     \draw(0,0.0) node[anchor=south west,xshift=-15pt,yshift=-8pt] {\trimfiga{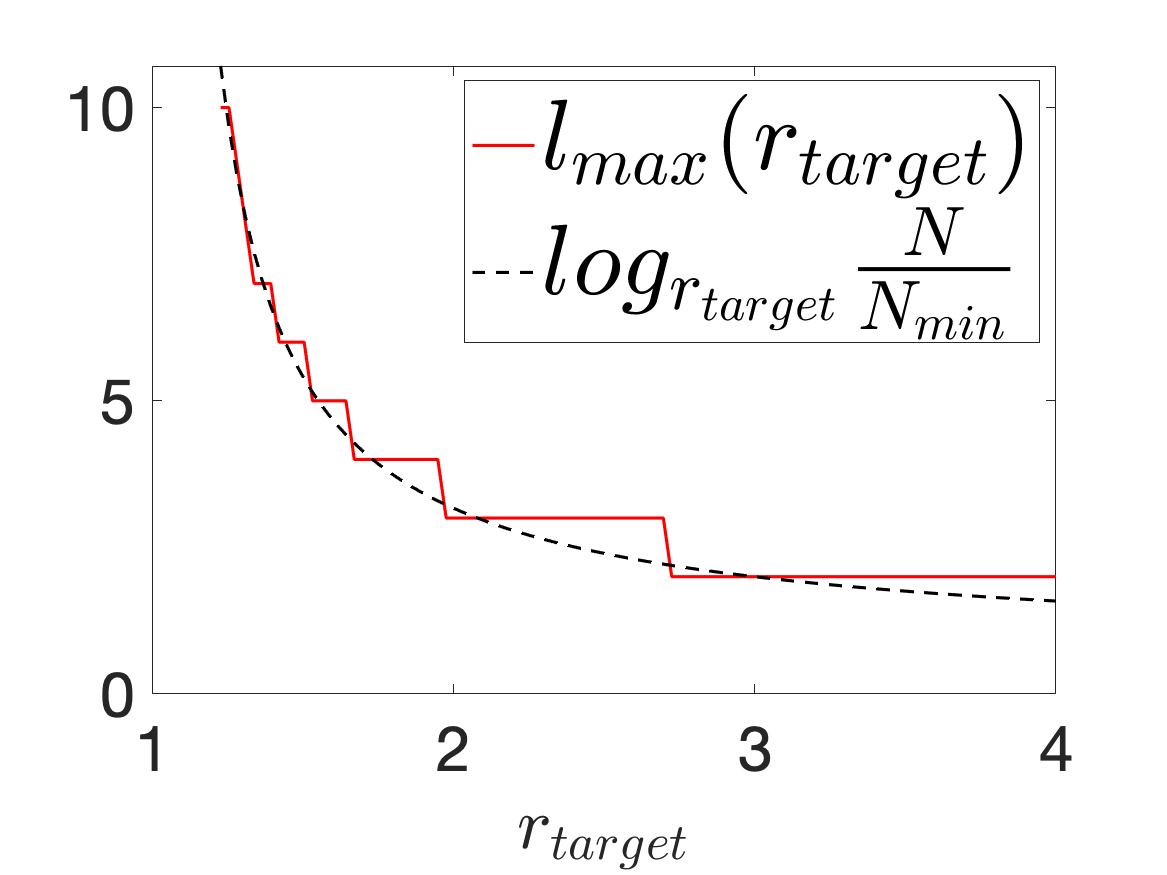}{\figWidtha}};
     %\draw(2,1) node[anchor=south] {\footnotesize $\lmax(\tilde{r})$};
    \end{scope}

% grid:
  %\draw[step=1cm,gray] (0,0) grid (16,6);
\end{tikzpicture}
\end{center}
\caption{CR and ECR versus $\rTarget$.
  \DIM[1], $\CYC{V}[1, 1]$
  cycle with red-black smoothing ($\omega = 1$; and non-Galerkin coarse-level operators),
  $r_l \approx \rTarget$ ($N = 64,\ \Nmin = 8$).
  A relatively good ECR is obtained over a wide range of $\rTarget$.
}
  \label{fig:CR_r_l_1d}
\end{figure}
}

% ------------------------------------------------------------------------------------------------------------------------
\subsection{General factor-$r$ coarsening in two dimensions}  \label{sec:generalCoarseningTwoDimensions}

The extension of general factor-$r$ coarsening to more space dimensions is straightforward. 
The number of grid points on coarse levels is chosen in a similar way to one dimension.
Tensor product linear interpolation is used to define the coarse-to-fine interpolation operator
and the fine-to-coarse restriction operator is then formed from the adjoint of the interpolation operator. 

{
\newcommand{\figWidth}{6cm}% height 
\newcommand{\trimfig}[2]{\trimh{#1}{#2}{.0}{.0}{.0}{.0}}
\begin{figure}[h!]
\begin{center}
\begin{tikzpicture}[scale=1]
  \useasboundingbox (0,0.5) rectangle (15,6.25);  % set the bounding box (so we have less surrounding white space)

   \draw(0.0,0.0) node[anchor=south west,xshift=-15pt,yshift=-8pt] {\trimfig{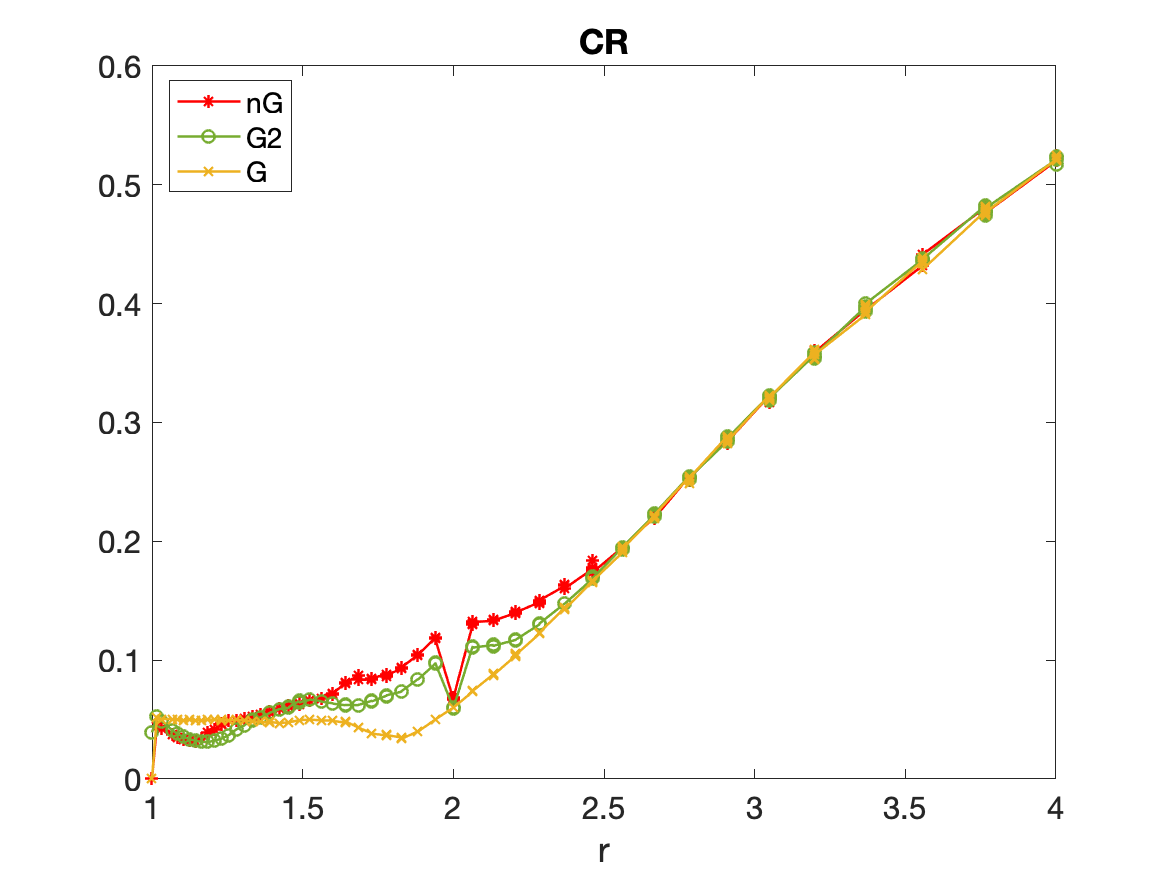}{\figWidth}};
   \draw(7.5,0.0) node[anchor=south west,xshift=-15pt,yshift=-8pt] {\trimfig{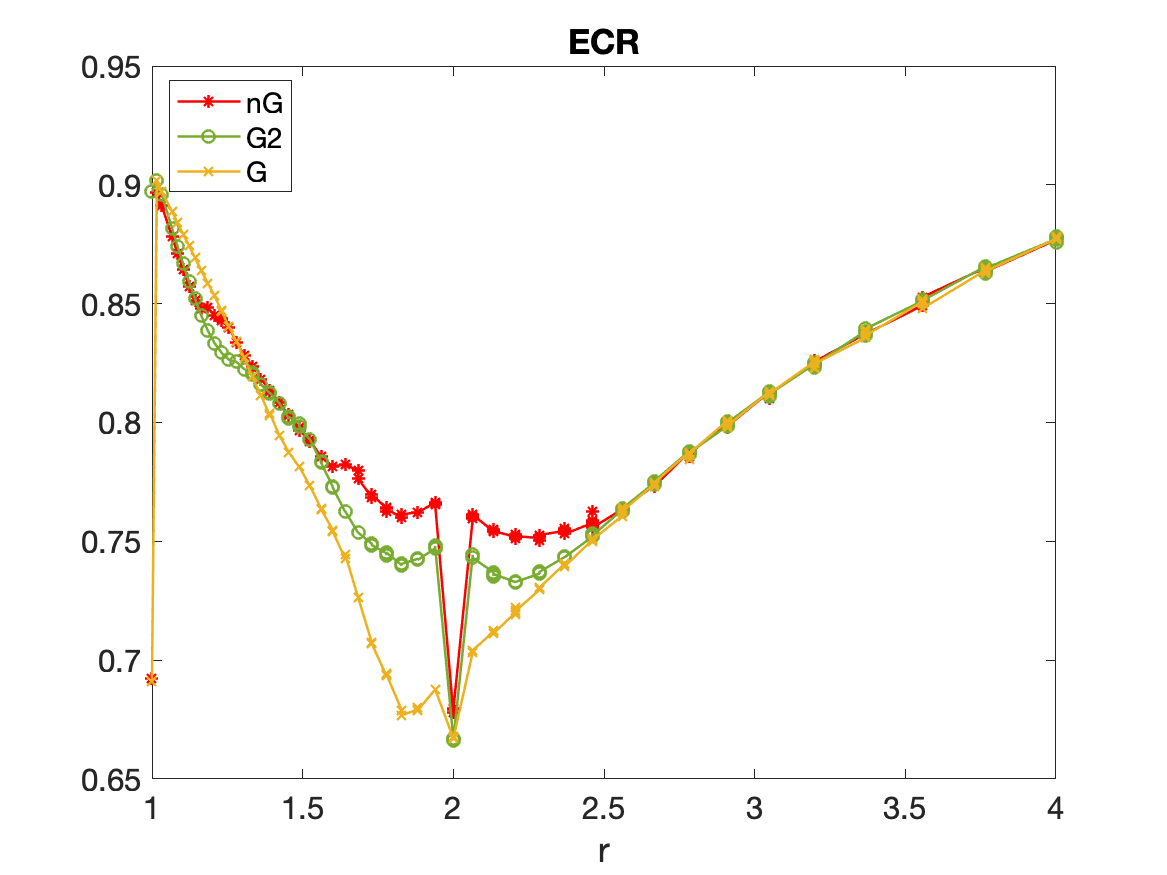}{\figWidth}};

   %\orderTwoVcycleTwoLevel{xshift=1.7cm,yshift=4.75cm,scale=1}{$\nu\!=\!2$}; this is std-C
   \VcycleTwoLevel{xshift=10.5cm,yshift=4.75cm,scale=1}{$\nu\!=\!2$};
   %\circleVortex{10.35}{1}{Opt for $r<2$}; 2-level opt does not mean anything!

    % title:
   \draw(7.5,5.8)  node[draw=\colourTwoD,very thick,fill=\colourTwoD!20,anchor=south,yshift=-4pt,inner sep=2.5pt]
             {\small \CC[factor-$r$], \DIM, \LEV}; 

% grid:
%\draw[step=1cm,gray] (0,0) grid (16,6);
\end{tikzpicture}
\end{center}
\caption{CR and ECR versus $r$.
	\DIM,
   $(h,\ H)$ \LEV\ $V[1, 1]$ cycle with red-black smoothing
  ($\omega = 1$), $\frac{H}{h}=r$ ($N = 64$).
(nG=$\LNAH$, G=$\LGGH$, G2=G($r=2$).)
  }
  \label{fig:CR_r_2d}
\end{figure}
}
Figure~\ref{fig:CR_r_2d} shows results for factor-$r$ coarsening in two dimensions for two-level $V[1,1]$ cycles.
Besides the non-Galerkin and Galerkin coarse-level operators, we also show results for a ``fixed Galerkin'' operator which uses the Galerkin operator from standard coarsening with a factor $r=2$ (``G2'').
As in one dimension, the convergence rates increase nearly monotonically with increasing coarsening factor $r$.
The effective convergence rate, on the other hand, shows a fairly broad minimum near $r=2$. 
The results for the Galerkin, fixed-Galerkin and non-Galerkin coarse-level operators are
similar with the Galerkin one giving generally better results.

Results from multi-level $V[1,1]$ cycles with factor-$r$ coarsening are shown in Figure~\ref{fig:CR_r_l_2d}.
The number of levels is chosen as a function of $\rTarget$ following formula~\eqref{eq:numberOfLevelsVersusR}.
Results are shown for the non-Galerkin coarse-level operators, that is, the
same stencil is used on every level; graphs for the Galerkin coarse-level operators are expected to be similar in form.
The multi-level convergence rates, which roughly mimic the two-level results, increase almost monotonically as $r$ increases.
The effective convergence rate shows a fairly broad minimum around $\rTarget=2$ with
a dip at exactly $\rTarget=2$, largely due to the work-units being discontinuously small at this value
(see the discussion in Section~\ref{Sec: rC_J}).
As in one dimension, we can again conclude that there is a range of coarsening factors near $r=2$
where the multigrid scheme behaves well, and there is no obvious optimal value of the coarsening factor based on the ECR, except for the special behavior at $r=2$.
% \old{Following our discussion in Section \ref{Sec: rbC_3d} with regard to red-black coarsening, there is a particular interest in examine the case $r \approx \sqrt{2}$. However, from Figure \ref{fig:CR_r_l_2d} it does not appear to be optimal there.}
We are going to further analyze the question of optimal coarsening factor in the next Section~\ref{Sec: rC_J}.
 
{
\newcommand{\figWidth}{6cm}% height 
\newcommand{\trimfig}[2]{\trimh{#1}{#2}{.0}{.0}{.0}{.0}}
\newcommand{\figWidtha}{2.5cm}% height 
\newcommand{\trimfiga}[2]{\trimhb{#1}{#2}{.0}{.0}{.0}{.0}}
\begin{figure}[h!]
\begin{center}
\begin{tikzpicture}[scale=1]
  \useasboundingbox (0,0.5) rectangle (16,6.5);  % set the bounding box (so we have less surrounding white space)

   \draw(0.0,0.0) node[anchor=south west,xshift=-15pt,yshift=-8pt] {\trimfig{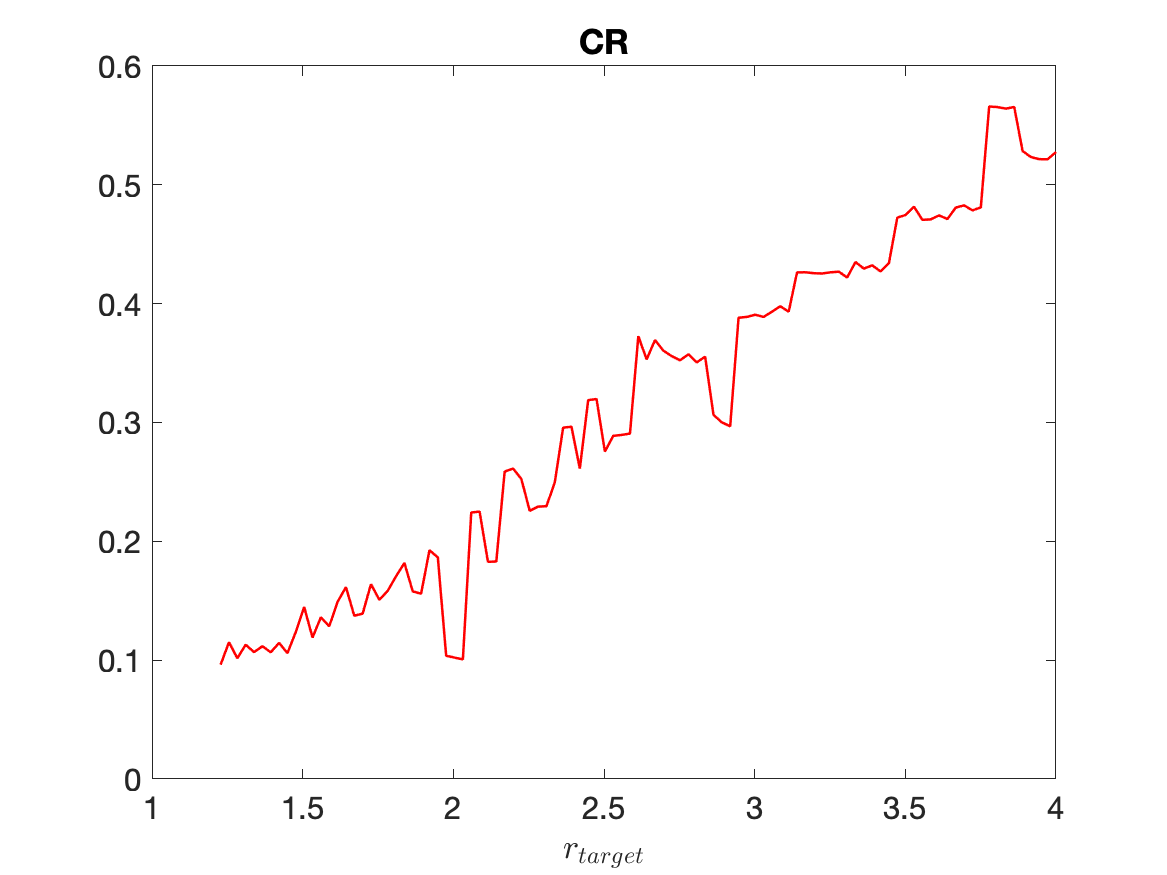}{\figWidth}};
   \draw(8.0,0.0) node[anchor=south west,xshift=-15pt,yshift=-8pt] {\trimfig{ECR2d_Rcoarsening_l_RB_gamma1w10}{\figWidth}};

  % \orderTwoVcycleThreeLevel{xshift=1.7cm,yshift=4.75cm,scale=1}{$\nu\!=\!2$}{$l(\tilde{r})$}; this is std-C! Maybe just the cycle label??
  \VcycleThreeLevel{xshift=1.7cm,yshift=4.75cm,scale=1}{$\nu\!=\!2$}{${\lmax}(\rTarget)$};

   % \circleVortex{10.6}{3}{*fix me* Good ECR over a range of $\tilde{r}$}; 
   \ellipseLabel{11.1}{3.25}{Good ECR}{violet}{1.2}{.6};

    % title:
   \draw(8,6)  node[draw=\colourTwoD,very thick,fill=\colourTwoD!20,anchor=south,yshift=-4pt,inner sep=2.5pt]
             {\small  \CC[factor-$r$], \DIM, \LEV[multi]}; 
 
   \begin{scope}[xshift=4.95cm,yshift=.75cm]\draw(0,0.0) node[anchor=south west,xshift=-15pt,yshift=-8pt] {\trimfiga{l_r}{\figWidtha}};
     %\draw(1.5,1) node[anchor=south] {\footnotesize ${\lmax}(\tilde{r})$};
    \end{scope}
% grid:
%\draw[step=1cm,gray] (0,0) grid (16,6);
\end{tikzpicture}
\end{center}
 \caption{CR and ECR versus $\rTarget$.
 	%\textbf{\red **FIX ME**} 
   %\textbf{\red Question: What happens to the figure if we use an optimal $\omega$ ?  }
   \DIM, $\CYC{V}[1, 1]$ cycle with red-black smoothing
   ($\omega = 1$; and non-Galerkin coarse-level operators), $r_l \approx \rTarget$ ($N = 64,\ \Nmin = 8$).
   A relatively good ECR can be obtained for a wide range of $\rTarget \in(1.5,2.5)$.
 }
 \label{fig:CR_r_l_2d}
\end{figure}
}

%------------------------------------------------------------------------------------------------------------------------
\newcommand{\ThetaLowr}{\Theta^{\rm low}_r}
\newcommand{\ThetaHighr}{\Theta^{\rm high}_r}
\newcommand{\JacobiSymbol}{\hat{S}_h}% Jacobi smoother symbol *new* 
\subsection{Local Fourier smoothing analysis of factor-$r$ coarsening with a Jacobi smoother} \label{Sec: rC_J}

In this section we use local Fourier analysis to better understand the computational results
on coarsening by a general factor $r$ obtained previously in Sections~\ref{sec:generalCoarseningOneDimension} and~\ref{sec:generalCoarseningTwoDimensions}.
To simplify the analysis we consider a multigrid cycle that uses a Jacobi smoother with a relaxation parameter $\omega$
combined with an ideal coarse-level correction.
We choose to analyze the $\omega$-Jacobi smoother instead of the red-black smoother for
a number of reasons.
The Jacobi smoother is easier to analyze since the Fourier modes are its eigenfunctions.
It is straightforward to choose an optimal relaxation parameter $\omega^*(r)$ for the Jacobi smoothing rate,
and thus we can have a fairer comparison of the convergence rates at different $r$.
With Jacobi the smoothing rate generally gives a better estimate of the overall multigrid convergence rate compared to the corresponding estimates with a red-black smoother, and thus Jacobi smoothing analysis will be easier to compare with computational results.
% \kP{and its smoothing rate and multigrid convergence rate are more close to each other than those of the red-black smoother.}
At the same time we expect the use of the Jacobi smoother will still cover most of the key features of factor-$r$ coarsening.

We study the model problem of Poisson's equation in $d$ dimensions on the infinite grid $\grid G_h$ \eqref{G_h_d}.
Let $L_h = -\Delta_h$ be the standard  second-order accurate approximation to the negative Laplacian in $d$ space dimensions.
The Jacobi smoothing operator with a relaxation parameter $\omega$ is given by
\begin{equation}
     S_h = I_h - \frac{\omega h^2}{2d} L_h. 
\end{equation}
The Fourier symbol of $S_h$, corresponding to $\phi_h (\cdot,\ \thv),\ \thv = [\th_1,\ \cdots,\ \th_d]^T \in \Theta \equiv [-\pi,\pi)^d$, is 
\begin{equation}
   \JacobiSymbol(\thv) = 1 - 2\omega \xi_{\thv},
\end{equation}
where 
\begin{equation}
    \xi_{\thv} \eqdef \frac{1}{d} \sum_{k=1}^d \sin^2 \frac{\th_k}{2}.
\end{equation}
\commentB{In order to determine the smoothing factor, as an estimate of the multigrid convergence rate, 
consider an $(h,\ rh)$ 2-level cycle 
with coarsening by a factor $r$.
According to the coarse grid, the sets of low and high frequencies are given simply by 
}
Consider a coarse grid 
\ba \label{G_H_r}
\grid G_{rh} = \Big\{  \xv =\jv\, rh:~ \jv \in \Integer^d \Big\}
\ea
on the first coarse level.
We will estimate the multigrid convergence rate with factor-$r$ coarsening from the
smoothing factor $\mu$ of $\omega$-Jacobi. 
For an $(h,\ rh)$ two-level cycle, % because for a 3-level cycle, for example, we have a quite different Theta-low.
the sets of low and high frequencies are given by
\ba
     \ThetaLowr \eqdef \big[-\frac{\pi}{r},\ \frac{\pi}{r} \big)^d, \quad
     \ThetaHighr = \Theta \setminus \ThetaLowr.
\ea
The fraction of the volume of frequency space occupied by $\ThetaLowr$ is thus $\frac{1}{r^d}$, which tends to one as $r\rightarrow 1$ and tends to zero as $r\rightarrow \infty$.
The smoothing factor for $\omega$-Jacobi is determined from the worst-case convergence rate over the high frequencies: 
\ba
   \mu(r,\ \omega) \eqdef & \sup_{\thv \in \ThetaHighr} |\JacobiSymbol(\thv)|
         = \sup_{\xi \in [\zeta,\ 1]} |1 - 2\omega \xi|
         = \max \big\{ \,  |1-2\omega\zeta|,\ |1-2\omega| \,  \big\},
%- %          \quad \geq\ \mu(\omega^*),
\ea
where 
\begin{equation}
    \zeta(r) \eqdef \frac{1}{d}\sin^2{\frac{\pi}{2r}}.
\end{equation}
An optimal relaxation parameter $\omega^*$ is chosen to minimize the smoothing factor $\mu$ over $\omega$:
\ba \label{mu_w_J_r}
&  \omega^*(r) = \frac{1}{1+\zeta(r)},
\ea
which leads to the optimal smoothing factor 
\ba
 &  \mu^*(r)  \eqdef \mu(r,\ \omega^*(r)) = \frac{1-\zeta(r)}{1+\zeta(r)}.  \label{eq:muStarfactorR}
\ea
Figure~\ref{fig:sr_J_r} (left) shows the form of $\mu^*(r)$ as a function of $r$ in different dimensions $d$. $\mu^*(r)$ is a monotone increasing function for $r>1$.
For a multigrid algorithm with $\nu$ smoothing steps per cycle, 
a convergence rate might be expected to be approximately $(\mu^*)^\nu$.
Figure~\ref{fig:sr_J_r} (right) shows this expected rate for the typical case of $\nu=2$ smoothing steps per cycle.
%These curves have a similar form to the previous computational results in Figures~\ref{fig:CR_r_l_1d}, \ref{fig:CR_r_2d} and~\ref{fig:CR_r_l_2d}.

{
\newcommand{\figWidth}{6cm}% height 
\newcommand{\trimfig}[2]{\trimh{#1}{#2}{.0}{.0}{.0}{.0}}
\begin{figure}[h!]
\begin{center}
\begin{tikzpicture}[scale=1]
  \useasboundingbox (0,0.5) rectangle (15,6.);  % set the bounding box (so we have less surrounding white space)

   \draw(0.0,0.0) node[anchor=south west,xshift=-15pt,yshift=-8pt] {\trimfig{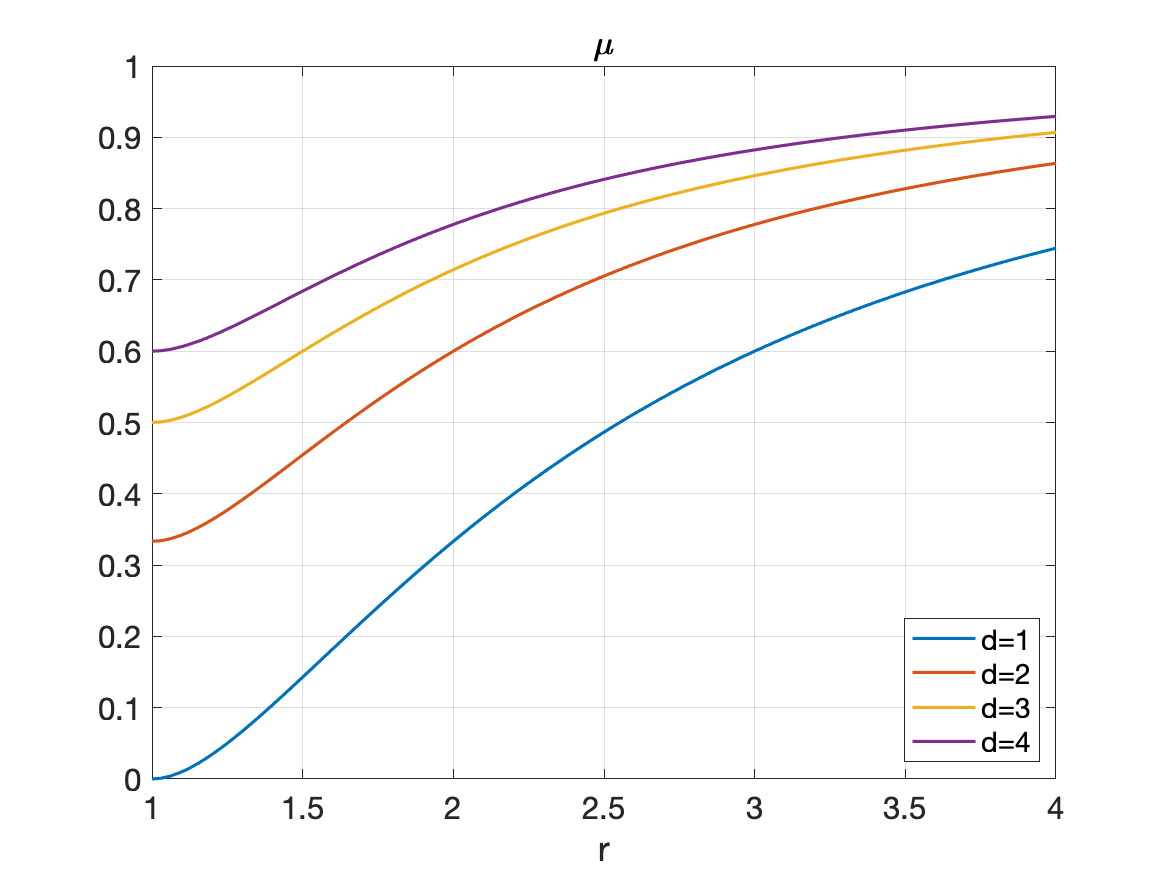}{\figWidth}};
   \draw(7.5,0.0) node[anchor=south west,xshift=-15pt,yshift=-8pt] {\trimfig{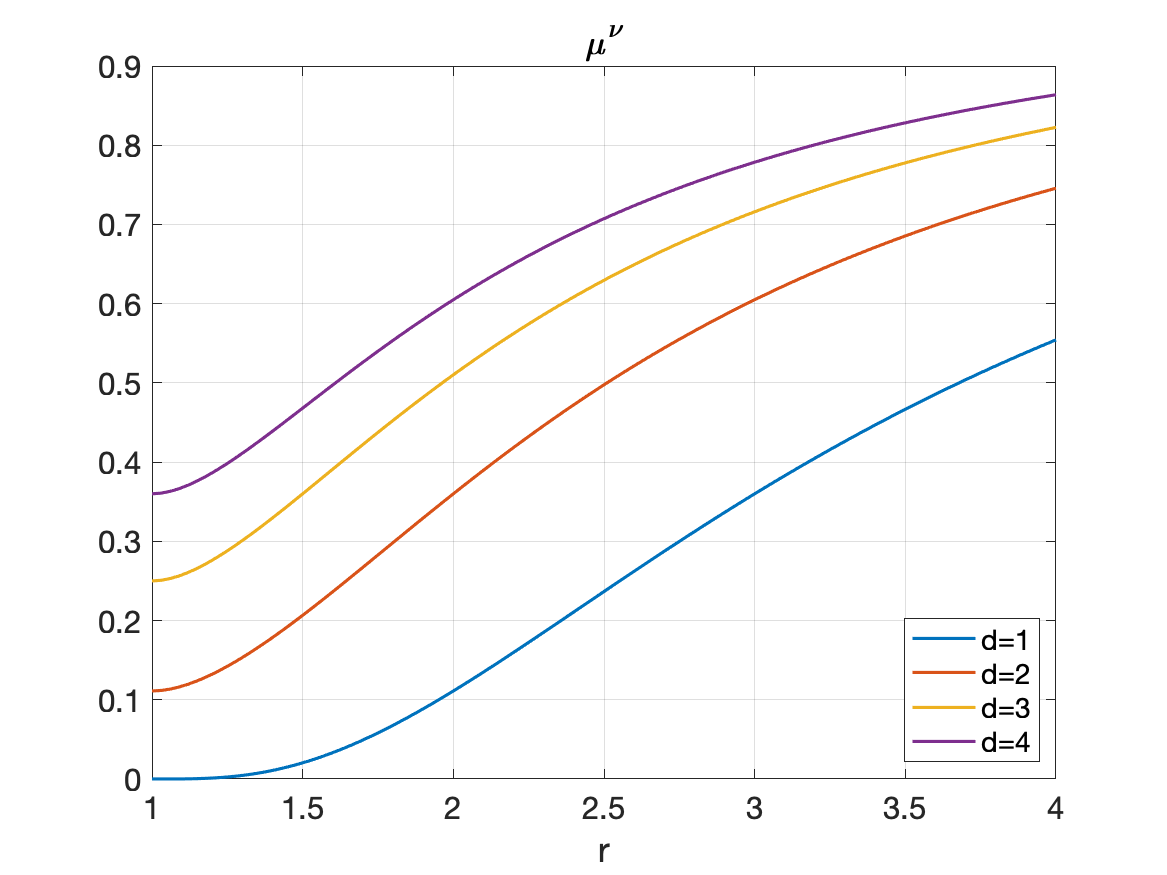}{\figWidth}};
% grid:
% \draw[step=1cm,gray] (0,0) grid (16,6);
\end{tikzpicture}
\end{center}
\caption{Smoothing rate $\mu$ (left) and estimate of CR $\mu^{\nu}$ (right, $\nu=2$)
  versus $r$ for V and W cycles with $\omega^*$-Jacobi smoothing and factor-$r$ coarsening.}
  \label{fig:sr_J_r}
\end{figure}
}

We focus on the two-dimensional case ($d=2$) and
to compare the estimates so far to computations, we
compute the actual convergence rates of a multigrid algorithm with 
%$N=64$ grid points on the finest level and 
Jacobi smoothing, with $\omega^*$ as given in (\ref{mu_w_J_r}). 
Computed CRs using $V[1,1]$ and $W[1,1]$ multi-level cycles are compared
to the theoretical curves from~\eqref{eq:muStarfactorR} in Figure~\ref{fig:cr_J_2d_r}.
As might be expected, the W cycle results compare quite well with the theory except for some deviation as $r$ tends to one.
The V cycle results compare well with the theory for $r>2$ but there is a more pronounced deviation for $r<2$ where the computed convergence rate levels off to a value of $CR \approx 0.35$. 
%This deviation from theory is at this point not fully understood but is likely attributed to the transfer operators which have to transfer higher and higher frequencies as $r$ approaches one.
\kl{
%{\blue I am not sure about this:} Me neither.. That's why it says it's not fully understand..
This deviation from the theoretical estimate comes from the distinction between the actual and ideal coarse-level corrections which is more dominant as $r$ approaches one, but at this point is not fully understood. It is likely that, for a V cycle the deviation mostly comes from the error in the coarse-level solves, while for a W cycle (or $\gamma \geq 2$) the deviation mostly comes from the transfer errors since higher and higher frequencies have to be transferred as $r$ approaches one.}
{
\newcommand{\figWidth}{6cm}% height 
\newcommand{\trimfig}[2]{\trimh{#1}{#2}{.0}{.0}{.0}{.0}}
\begin{figure}[h!]
\begin{center}
\begin{tikzpicture}[scale=1]
  \useasboundingbox (0,0.5) rectangle (15,6.);  % set the bounding box (so we have less surrounding white space)

   \draw(0.0,0.0) node[anchor=south west,xshift=-15pt,yshift=-8pt] {\trimfig{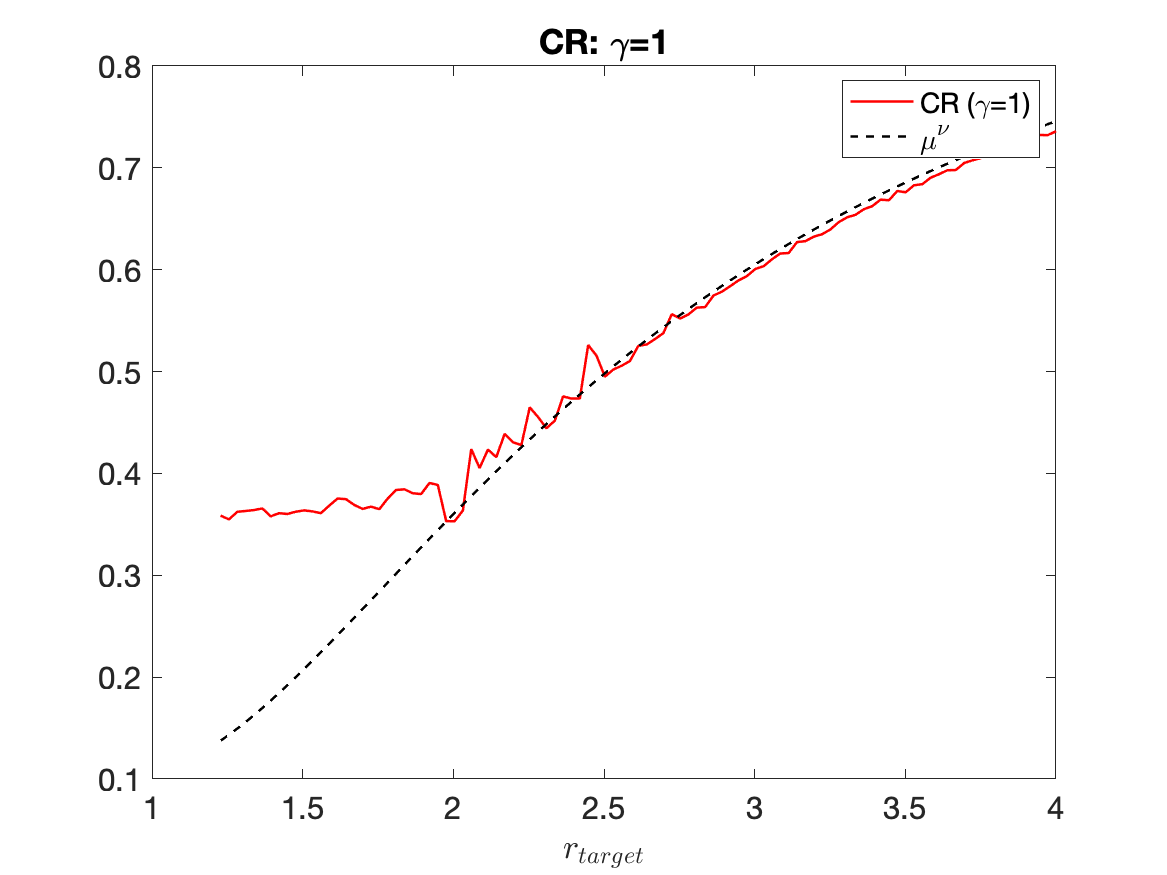}{\figWidth}};
   \draw(7.5,0.0) node[anchor=south west,xshift=-15pt,yshift=-8pt] {\trimfig{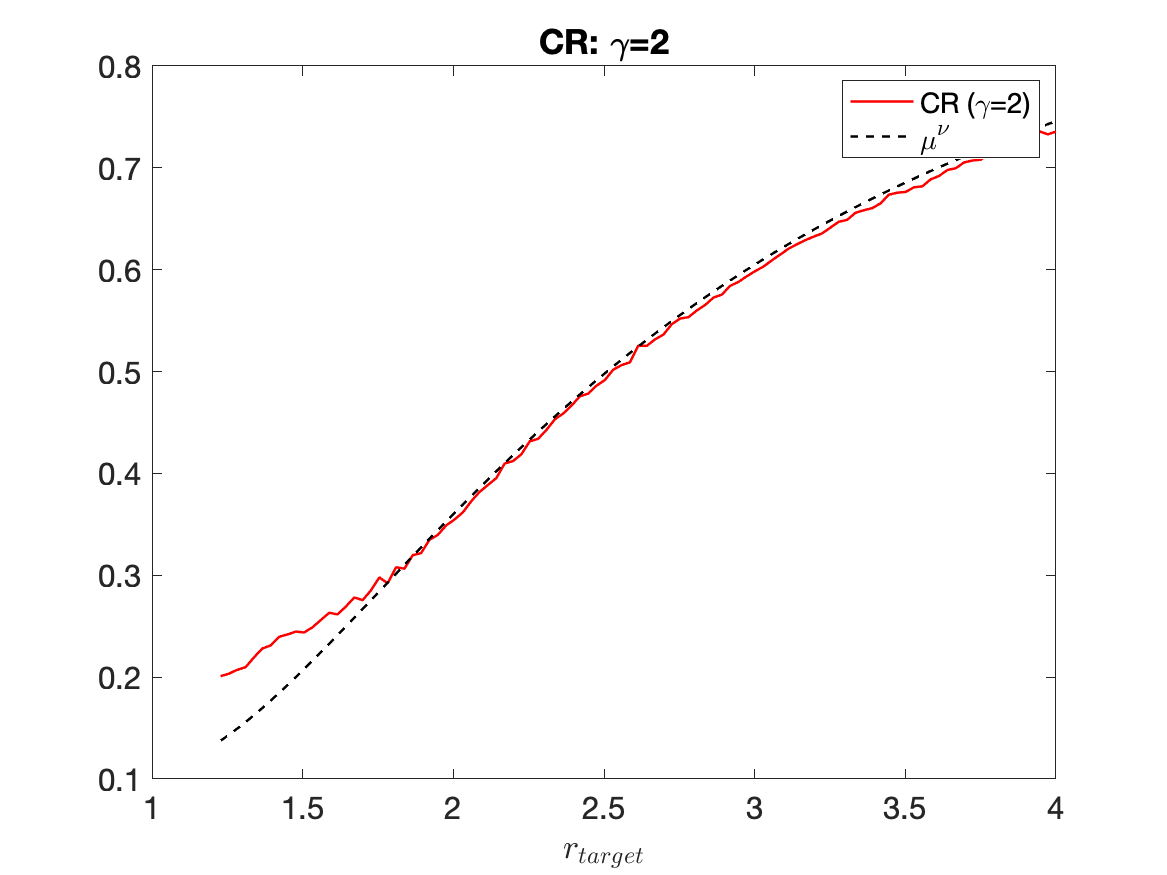}{\figWidth}};
% grid:
% \draw[step=1cm,gray] (0,0) grid (16,6);
\end{tikzpicture}
\end{center}
  \caption{CR versus $\rTarget$, in comparison with the estimate $\mu^{\nu}$.
  	\DIM, V ($\gamma=1$, left) and W ($\gamma=2$, right) multi-level cycles with $\omega^*$-Jacobi smoothing ($\nu=2$) and factor-$r$ coarsening (and non-Galerkin coarse-level operators; $N = 64,\ \Nmin = 8$).
  }
  \label{fig:cr_J_2d_r}
\end{figure}
}

 From the derivation of this simple example we can partly answer the question raised in Section \ref{Sec: rbC_3d}, with regard to the idea of a factor-$2^{\frac{1}{d}}$ coarsening in $d$ dimensions. As we can see, in two dimensions $\sqrt{2}$ coarsening does not give convergence rates comparable to red-black coarsening. The reason is that the premise of $2^{\frac{1}{d}}$ coarsening is to have equal volume of $\Theta^{low}$ and $\Theta^{high}$, which results in average-case optimization; for fastest convergence, on the other hand, we need \textsl{worst-case} optimization.
 From this comparison we can conclude that the fast convergence with red-black coarsening mainly comes from its structural alignment with the smoother.

\medskip
\newcommand{\ESR}{\text{ESR}}
\newcommand{\WU}{\text{WU}}

To get an idea of which coarsening factor $\rTarget$ gives the fastest algorithm we compare
  the effective convergence rates as a function of $\rTarget$. To do this
  we require an estimate of the computational work when coarsening by a factor of $r$. 
\commentB{
  Now let's go back to the question of optimal coarsening ratio $\rTarget$ to optimize the the effective convergence rate, for which 
%To compare effective convergence rates between the theory and computational results 
we require an estimate of the computational work when coarsening by a factor of $r$. } 
In particular, for the interpolation and restriction operators respectively,
we base the work-units on using linear interpolation and its adjoint. 
Note that at exactly $r=2$ (or factors of 2) the computational cost is somewhat reduced due to the alignment of the coarse grid points with half the fine grid points.
The work-units of a general multi-level $\gamma$-cycle with factor-$r$ coarsening is approximately given by 
(see \ref{sec: WU})
\ba
 \WU[r;\ \gamma,\ \nu] \approx 
    (\nu+3) \sum_{\ll = 0}^{\lmax-1} (\frac{\gamma}{r^d})^{\ll-1}
 %= (\nu+3) \frac{1-(\frac{\gamma}{r^d})^{\lmax}}{1-\frac{\gamma}{r^d}}
 \approx
 \begin{cases}
    (\nu+3) \frac{r^d}{r^d -\gamma}                                         & \text{for $\gamma/r^d < 1$}, \\[10pt]
    (\nu+3)\ \lmax                                      & \text{for $\gamma/r^d = 1$}, \\[10pt]
    (\nu+3) \frac{ ( \frac{\gamma}{r^d} )^{\lmax} }{ \frac{\gamma}{r^d} -1}  & \text{for $\gamma/r^d > 1$},
 \end{cases}
\ea
where $\lmax$ is chosen as \eqref{eq:numberOfLevelsVersusR}. 
Note the change in behavior in the WU depending on whether $\gamma < r^d$ or $\gamma \geq r^d$. 
This estimate for the computation work, and approximating the convergence rate by
$
  \rho \approx (\mu^*)^{\nu},
$
lead to the definition of the 
\textit{effective smoothing rate},
\begin{equation}
   \ESR \eqdef (\mu^*)^{\frac{\nu}{WU}}, \label{eq:ESR}
\end{equation}
as an estimated effective convergence rate for a multigrid cycle with $\nu$ smoothing steps per cycle.
This estimate for the ECR is graphed in Figure~\ref{fig:esr_J_r} for V ($\gamma=1$) and W ($\gamma=2$) cycles with $\nu=2$.
From equation~\eqref{eq:ESR} we can estimate the optimal coarsening factor, given the number of space dimensions $d$;
these appear at the minima in the \ESR~in Figure~\ref{fig:esr_J_r}.
For a V cycle, the optimal coarsening factor is always less than two,
in particular for $d=2$ the minimizer $r \approx 1.55$.
%Interestingly for $d=2$ it is close to $\sqrt{2}$, the coarsening ratio for red-black coarsening in two dimensions.
For a W cycle the optimal coarsening ratio, $r\approx 1.97$ is close to 2 in two dimensions
while it has a value $r\approx 1.79$ in three dimensions.

{
\newcommand{\figWidth}{6cm}% height 
\newcommand{\trimfig}[2]{\trimh{#1}{#2}{.0}{.0}{.0}{.0}}
\begin{figure}[h!]
\begin{center}
\begin{tikzpicture}[scale=1]
  \useasboundingbox (0,0.5) rectangle (15,6.);  % set the bounding box (so we have less surrounding white space)

   \draw(0.0,0.0) node[anchor=south west,xshift=-15pt,yshift=-8pt] {\trimfig{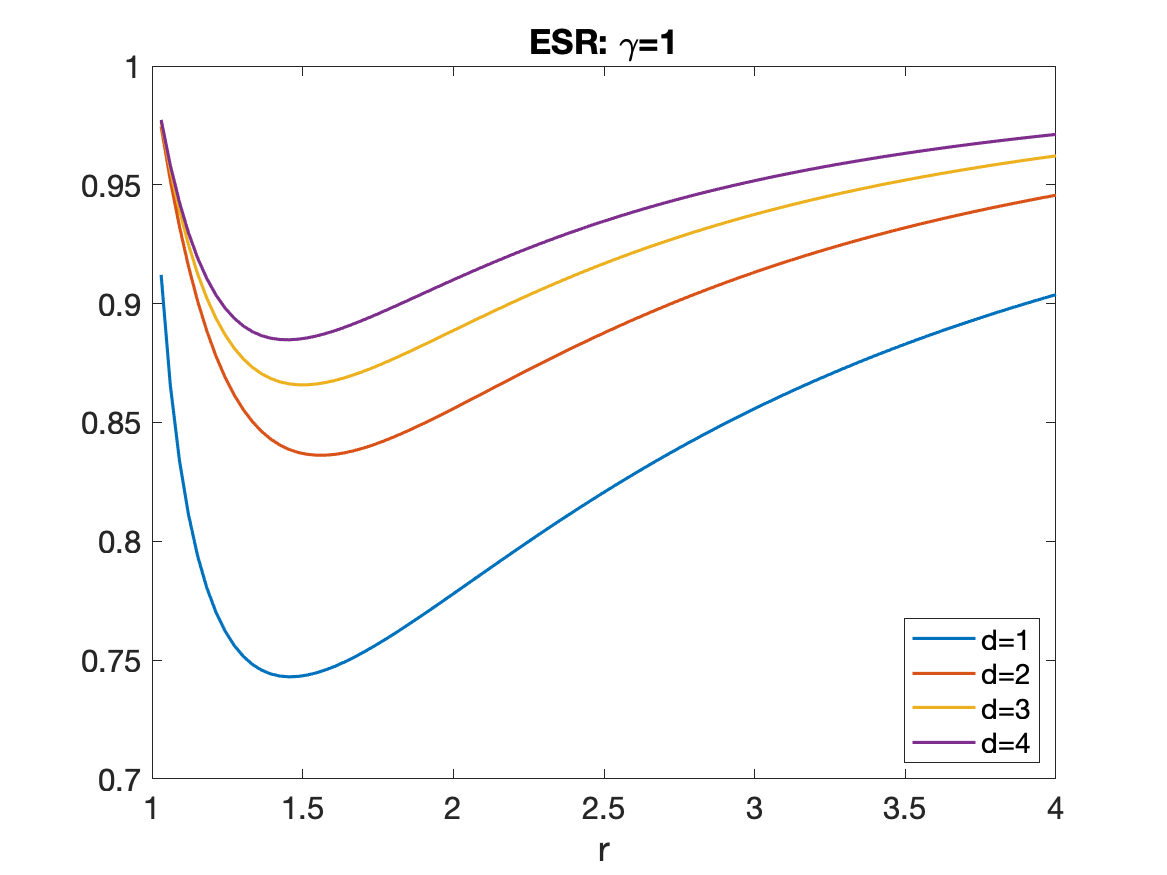}{\figWidth}};
   \draw(7.5,0.0) node[anchor=south west,xshift=-15pt,yshift=-8pt] {\trimfig{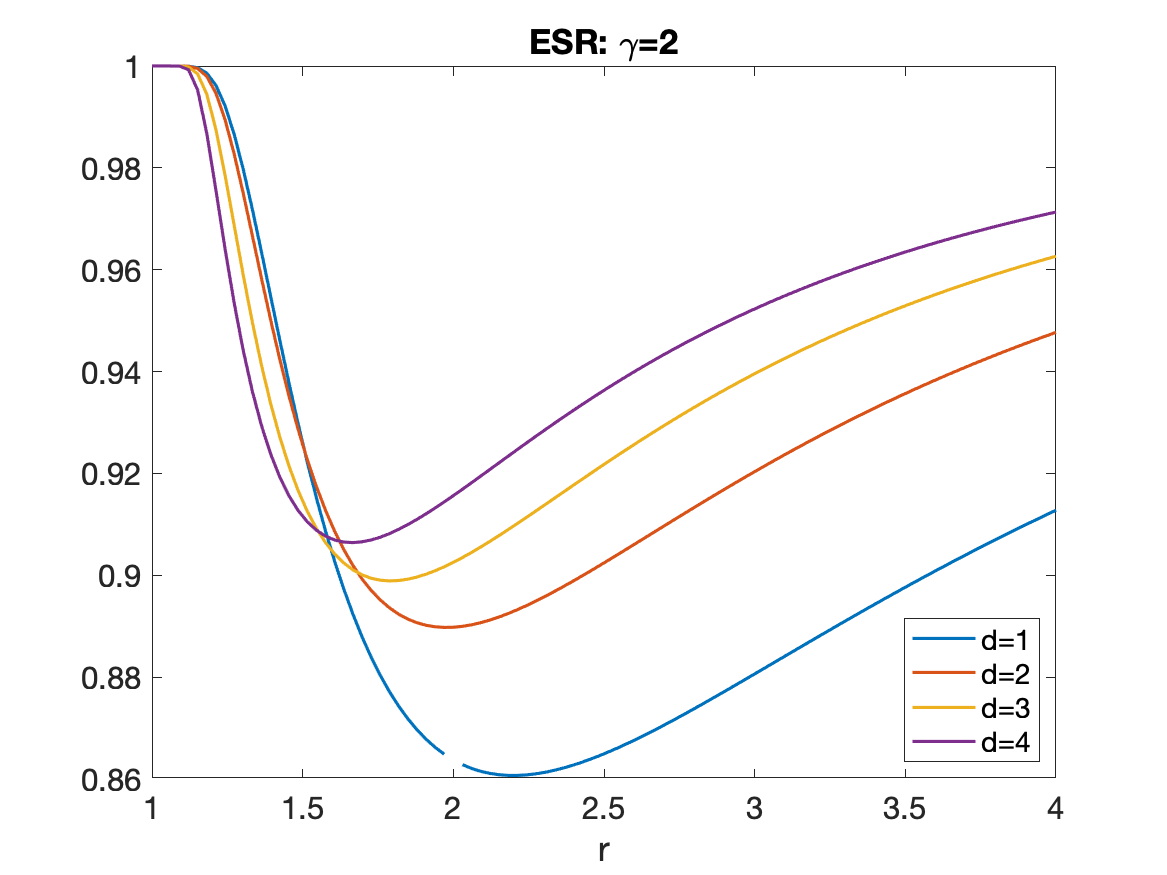}{\figWidth}};
% grid:
% \draw[step=1cm,gray] (0,0) grid (16,6);
\end{tikzpicture}
\end{center}
  \caption{ESR (estimate of the ECR) versus $r$. V ($\gamma=1$, left) and W ($\gamma=2$, right) cycles with $\omega^*$-Jacobi smoothing ($\nu=2$) and factor-$r$ coarsening.}
  \label{fig:esr_J_r}
\end{figure}
}

{
\newcommand{\figWidth}{6cm}% height 
\newcommand{\trimfig}[2]{\trimh{#1}{#2}{.0}{.0}{.0}{.0}}
\begin{figure}[h!]
\begin{center}
\begin{tikzpicture}[scale=1]
  \useasboundingbox (0,0.5) rectangle (15,6.);  % set the bounding box (so we have less surrounding white space)

   \draw(0.0,0.0) node[anchor=south west,xshift=-15pt,yshift=-8pt] {\trimfig{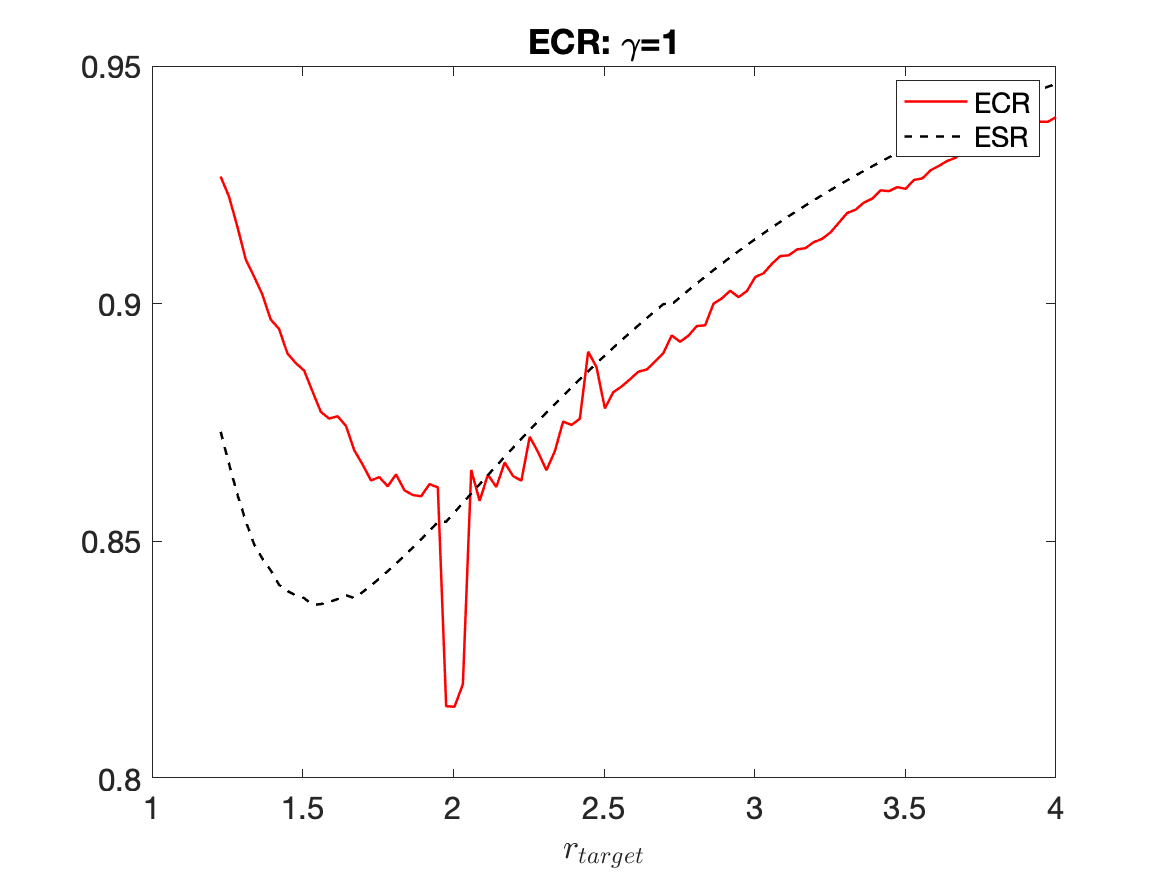}{\figWidth}};
   \draw(7.5,0.0) node[anchor=south west,xshift=-15pt,yshift=-8pt] {\trimfig{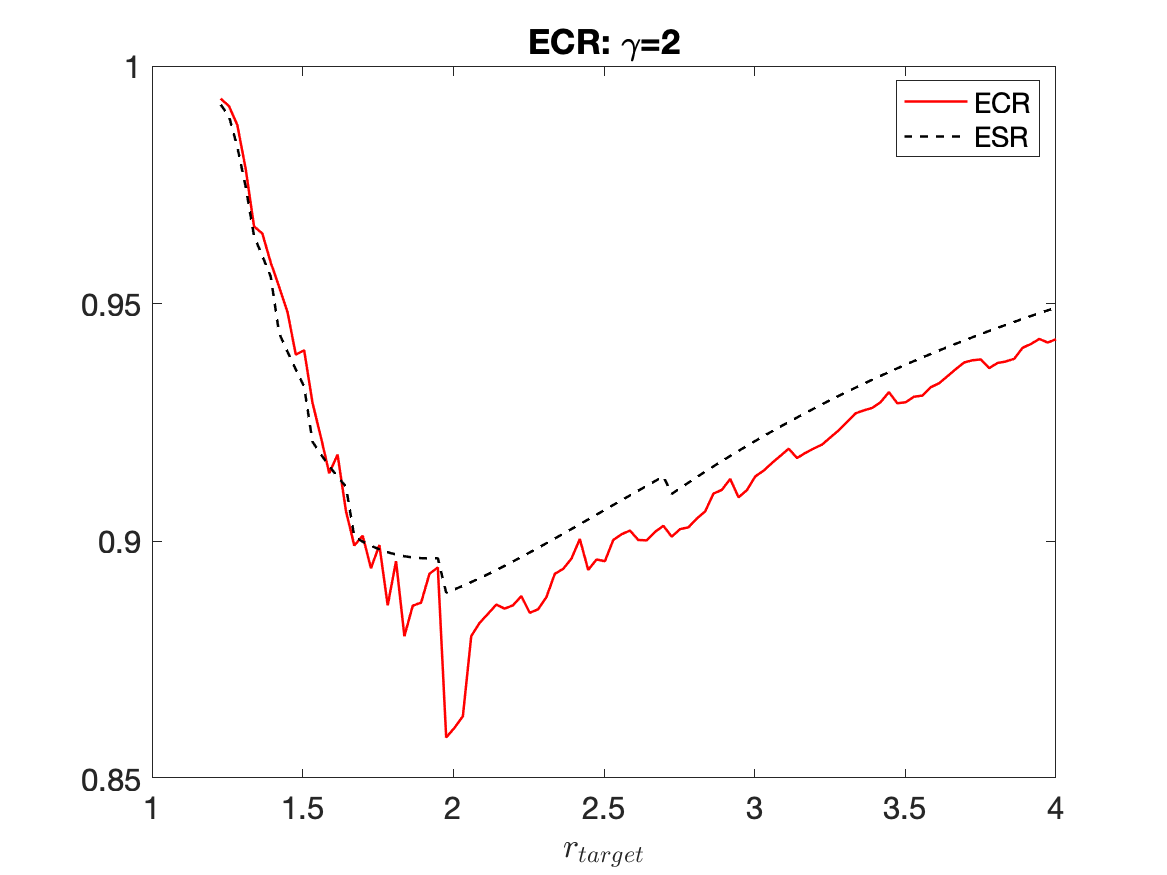}{\figWidth}};
% grid:
% \draw[step=1cm,gray] (0,0) grid (16,6);
\end{tikzpicture}
\end{center}
  \caption{ECR versus $\rTarget$, in comparison with the estimate ESR.
  	\DIM, V ($\gamma=1$, left) and W ($\gamma=2$, right) multi-level cycles with $\omega^*$-Jacobi smoothing ($\nu=2$) and factor-$r$ coarsening (and non-Galerkin coarse-level operators; $N = 64,\ \Nmin = 8$).
  }
  \label{fig:ecr_J_2d_r}
\end{figure}
}

In Figure~\ref{fig:ecr_J_2d_r} the estimated effective convergence rates \ESR, are compared with ECRs from actual computations for multigrid in two dimensions using $\omega^*$-Jacobi smoothing 
(the corresponding CRs appear in Figure~\ref{fig:cr_J_2d_r}).
We can see that the ESR gives a fairly good estimate for the ECR for a W cycle.
The comparison is not as good for the V cycle when $\rTarget<2$,  but this is to be expected
due to the deviation in the convergence rates between the theory and computation
as shown in Figure~\ref{fig:cr_J_2d_r}. 
Also, note that the singularly good ECR at $\rTarget = 2$ comes from the fact that the WU is singularly low when the coarse grid is embedded in the fine grid.

Finally, let us compare the convergence rates in Figure~\ref{fig:cr_J_2d_r}
and the effective convergenve rates in Figure~\ref{fig:ecr_J_2d_r} 
with Jacobi smoothing (with $\omega^*$), and the results in Figure~\ref{fig:CR_r_l_2d} with
red-black smoothing ($\omega=1$) in two dimensions. 
We can see that the red-black smoother, as usual, gives faster
convergence than Jacobi smoother, but the CRs and ECRs
seem to vary with $\rTarget$ in a very similar manner. 
Thus local Fourier smoothing analysis based on the $\omega$-Jacobi smoother
appears to reasonably describe the basic features of multigrid with coarsening by a general factor $r$.

\bigskip
\section{Conclusions}\label{Sec: conclusion}

In this paper, we have considered several nonstandard coarsening strategies for geometric multigrid, both for the grids and the operators on coarse levels. 

In the first part, we discussed the use of lower-order accurate coarse-level operators in multigrid algorithms for higher-order discretizations to the model problem. We showed though LFA results that, in particular by adjusting the value of a smoothing parameter $\omega$, we can usually achieve better convergence rates with second-order coarse-level operators. The resulting multigrid cycles are thus much more efficient.

In the second part, for 
the standard second-order accurate discretization to 
the two-dimensional Poisson's equation, we proposed red-black coarsening as an optimal strategy, to be paired with red-black smoothing. We presented detailed local Fourier analysis for the resulting multigrid cycles, in particular as a generalization from the one-dimensional case.
Numerical results show consistency with the theory and good convergence rates for general multi-level cycles.

In the third part, by simply discretizing the domain with a larger grid-spacing,
an algorithm for geometric multigrid cycles that coarsens by a general factor has been presented and analyzed. 
%In particular, we found that approximately coarsening by a factor $\sqrt{2}$ in two dimensions in this way does not yield such good convergence rate as red-black coarsening. 
This algorithm works in any dimension and brings some much needed flexibility in multi-level grids for complicated geometry,
and we've shown that good convergence rates are retained for a wide range of coarsening factors.

%Future Work

%\bigskip
\clearpage
%-----------------------------------------------------------------------------------------------------------
%theory-----------------------------------------------------------------------------------------------------
\appendix

\section{Further discussions and theoretical considerations}\label{Sec: theory}
In \ref{Sec: rb_1d} and \ref{Sec: rb_2d}, we record some of the details of the local Fourier analysis for one and two dimensional multigrid with standard coarsening.
\ref{Sec: rbC_LFA} concerns red-black coarsening, corresponding to Section~\ref{M2: rbC}. 
A three-level local Fourier analysis is given in \ref{Sec: rbC_3_LFA}, and some more general Laplace operators are considered in \ref{Sec: rbC_anisotropic} and \ref{Sec: rbC_4o}.
In~\ref{Sec: cyclic_2d}, the algorithm in Section~\ref{Sec: rbC} is reformulated as red-black reduction; and in~\ref{Sec: connections} we discuss the connections of the proposed multigrid cycles by red-black coarsening with reduction methods and some of the ideas from algebraic multigrid.
With regard to Section~\ref{M1: orderC}, we give the local Fourier analysis for multigrid for the fourth-order accurate discrete Laplacian in \ref{Sec: orderC_LFA}, with discussion on the order of accuracy of the coarse-level correction.
In \ref{sec: WU} we give a WU estimate for a general multi-level cycle.

\subsection{One-dimensional local Fourier analysis}\label{Sec: rb_1d}
This section supplements Section~\ref{sec:redBlack1d}, concerning the one-dimensional discrete Poisson's equation with \eqref{L_h_1}, and the two-level cycle with red-black smoothing and standard coarsening.

The matrix symbol \eqref{S_RB} of the red-black smoother on $\mathbb{E}_{\th}$ \eqref{E_1d} is derived as follows.
Consider the red-black smoothing operator's action on the Fourier modes, $\th \in \Theta$:
\begin{equation}\label{S_R phi}
\begin{aligned}
S^R_h\ \phi_h (\cdot,\ \th) &=\begin{cases}
\symF S^J_h({\th}) \ \phi_h (\cdot,\ \th), & \text{on}\quad  \grid G_h^R\\
\phi_h (\cdot,\ \th), & \text{on}\quad   \grid G_h^B
\end{cases}\\
&= 
\frac{1}{2}\left[(1 + \symF S^J_h({\th})) \ \phi_h (\cdot,\ \th)+(1 - \symF S^J_h({\th})) \ \phi_h (\cdot,\ \bar{\th})  \right],
\end{aligned}
\end{equation}
where the last equation follows from (\ref{eqn_minus});
and similarly
\begin{equation}\label{S_B phi}
S^B_h\ \phi_h (\cdot,\ \th) = \frac{1}{2}\left[(1 + \symF S^J_h({\th}))  \ \phi_h (\cdot,\ \th) -(1 -\symF S^J_h({\th}))\ \phi_h (\cdot,\ \bar{\th})  \right],
\end{equation}
where
\begin{equation}
\symF S^J_h({\th}) = 1-\frac{\omega h^2}{2} \symF L_h({\th}) = 1-2\omega \xi_{\th}.
\end{equation}
We can see that the red-black smoothing operator $S_h$ is invariant on each subspace $\mathbb{E}_{\th}$, with the matrix representation
$
\symE S(\th) = \symE S^B (\th)\ \symE S^R(\th),
$
where
\begin{equation}
\begin{aligned}
\symE S^R(\th) &= 
\frac{1}{2}\begin{bmatrix}
1+\symF S^J_h({\th}) & 1-\symF S^J_h(\bar{\th}) \\
1-\symF S^J_h({\th}) & 1+\symF S^J_h(\bar{\th})
 \end{bmatrix}
 = \mathbf{I} - \omega \begin{bmatrix}
 1 & -1\\ -1 & 1
 \end{bmatrix} \begin{bmatrix}
 \xi_{\th} \\ & 1-  \xi_{\th}
 \end{bmatrix},\\ 
\symE S^B (\th) &= 
\frac{1}{2}\begin{bmatrix}
1+\symF S^J_h({\th}) & -(1-\symF S^J_h(\bar{\th})) \\
-(1-\symF S^J_h({\th})) & 1+\symF S^J_h(\bar{\th})
 \end{bmatrix}
  = \mathbf{I} - \omega \begin{bmatrix}
 1 & 1\\ 1 & 1
 \end{bmatrix} \begin{bmatrix}
 \xi_{\th} \\ & 1-  \xi_{\th}
 \end{bmatrix}.
 \end{aligned}
\end{equation}

Then the smoothing factor, as defined in (\ref{mu}), is given by
\begin{equation}\label{mu_1d}
\mu = \rho^{\frac{1}{\nu}} \big( Q_h^H S_h^{\nu}\big)
= \sup_{\th \in \Theta^{\text{low}}} { \rho^{\frac{1}{\nu}} \big(\symE Q\ \symE S^{\nu}(\th) \big)}
\end{equation}
for $\nu$ smoothing steps per cycle, 
where $\symE Q = \begin{bmatrix} 0\\ & 1 \end{bmatrix}$ is the matrix symbol of the ideal coarse-level correction $Q_h^H$ on $\mathbb{E}_{\th}$.

\subsubsection{Cyclic reduction}\label{Sec: cyclic_1d}

We now relate the two-level cycle with standard coarsening $H=2h$ to one step of \textsl{cyclic reduction}.
Consider the $N \times N$ matrix corresponding to the discrete (negative) Poisson's equation
$\mathbf{A}_h\ \mathbf{x} = \mathbf{b}$
with periodic boundaries (instead of an infinite domain, for simplicity of notation), on a grid of size $N$:
\begin{equation}
\mathbf{A}_h = \begin{bmatrix}
2 & -1  & & & &   -1\\
-1 & 2 & -1 \\
 & -1 & 2 & -1 \\
 & & \ddots  & \ddots & \ddots \\
 & & & -1 & 2 & -1 \\
-1  & &  & &  -1 & 2
\end{bmatrix}.
\end{equation}
Cyclic reduction rearranges the grid points into odd-even (or red-black) ordering as 
\begin{equation}\label{system_rb}
\tilde{\mathbf{A}}_h\ \tilde{\mathbf{x}} = \tilde{\mathbf{b}},
\end{equation}
with
\begin{equation}
\tilde{\mathbf{x}} = \begin{bmatrix}
\mathbf{x}_R \\ \mathbf{x}_B
\end{bmatrix},\quad
\tilde{\mathbf{b}} = \begin{bmatrix}
\mathbf{b}_R \\ \mathbf{b}_B
\end{bmatrix},
\end{equation}
and
\begin{equation}
\tilde{\mathbf{A}}_h = \begin{bmatrix}
2\mathbf{I} & -\mathbf{L}^T \\
- \mathbf{L} & 2\mathbf{I}
\end{bmatrix},
\end{equation}
where $\mathbf{I}$ and 
\begin{equation}\label{LL}
\mathbf{L} = \begin{bmatrix}
1  & 1 \\
& 1 & 1 \\
& & \ddots & \ddots \\
& & & 1 & 1\\
1 & & & & 1
\end{bmatrix}
\end{equation}
are of size $\frac{N}{2} \times \frac{N}{2}$. 
Applying Gaussian elimination we can factor the permuted matrix as
\begin{equation}
\tilde{\mathbf{A}}_h  = \hat{\mathbf{L}} \hat{\mathbf{A}} \hat{\mathbf{L}}^T
\end{equation}
with
\begin{equation}
\hat{\mathbf{L}} = \begin{bmatrix}
\mathbf{I} \\ -\frac{1}{2}\mathbf{L} & \mathbf{I}
\end{bmatrix},\quad
\hat{\mathbf{A}} = \begin{bmatrix}
2\mathbf{I} \\
& \frac{1}{2} \mathbf{A}_H
\end{bmatrix},
\end{equation}
where $\mathbf{A}_H$ is of the same form as $\mathbf{A}_h$ but of size $\frac{N}{2} \times \frac{N}{2}$.

The system (\ref{system_rb}) is equivalent to
\begin{equation}
\hat{\mathbf{A}} \hat{\mathbf{L}}^T\ \tilde{\mathbf{x}} = \hat{\mathbf{L}}^{-1}\tilde{\mathbf{b}},
\end{equation}
where
\begin{equation}
\hat{\mathbf{L}}^{-1} = \begin{bmatrix}
\mathbf{I} \\ \frac{1}{2}\mathbf{L} & \mathbf{I}
\end{bmatrix}.
\end{equation}
Or more specifically,
\begin{equation}\label{eqn_cyclic}
\begin{cases}
\mathbf{x}_R = \frac{1}{2} (\mathbf{L}^T \mathbf{x}_B + \mathbf{b}_R);\\
\frac{1}{4} \mathbf{A}_H\ \mathbf{x}_B = \frac{1}{4} (\mathbf{L} \mathbf{b}_R + 2 \mathbf{b}_B),
\end{cases}
\end{equation}
where the $\frac{1}{4}$ of the second equation comes from the scaling $\frac{1}{H^2} = \frac{1}{4} \frac{1}{h^2}$. 

\begin{theorem}
	The system of equations (\ref{eqn_cyclic}) coming from cyclic reduction is equivalent to an $(h,\ 2h)$ two-level multigrid cycle with one post red-smoothing step with $\omega = 1$. In particular, the second equation is equivalent to the coarse-level correction, while the first equation corresponds exactly to one post red-smoothing step. Moreover, the two-level algorithm converges in one cycle.
\end{theorem}

The explanation of the above theorem is as follows. Given any current approximation to the solution in red-black ordering
\begin{equation}
\tilde{\mathbf{x}}^{(k)} = \begin{bmatrix}
\mathbf{x}_R^{(k)}  \\ \mathbf{x}_B^{(k)} 
\end{bmatrix},
\end{equation}
the coarse-level correction, written explicitly, consists of the following steps: 
the residual 
\begin{equation}
\tilde{\mathbf{d}}^{(k)} = \begin{bmatrix}
\mathbf{d}_R^{(k)}  \\ \mathbf{d}_B^{(k)} 
\end{bmatrix}
= \tilde{\mathbf{b}} - \tilde{\mathbf{A}}_h\ \tilde{\mathbf{x}}^{(k)}
\end{equation}
is restricted to the coarse grid as
\begin{equation}
\tilde{\mathbf{d}}_H^{(k)} = \frac{1}{4} \begin{bmatrix}
\mathbf{L} & 2\mathbf{I} 
\end{bmatrix}\ \tilde{\mathbf{d}}^{(k)}
=  \frac{1}{4} \begin{bmatrix}
\mathbf{L} & 2\mathbf{I} 
\end{bmatrix} \tilde{\mathbf{b}} -  \frac{1}{4} \mathbf{A}_H \mathbf{x}_B^{(k)} ;
\end{equation}
then the error equation 
\begin{equation}\label{eqn_B}
\frac{1}{4} \mathbf{A}_H \hat{\mathbf{v}}_H^{(k)} = \tilde{\mathbf{d}}_H^{(k)}
\end{equation}
is solved on the coarse grid; and finally, the error estimate on the coarse grid $\hat{\mathbf{v}}_H^{(k)}$ is interpolated back to the fine grid as
\begin{equation}
\tilde{\mathbf{v}}^{(k)} = \begin{bmatrix}
\hat{\mathbf{v}}_R^{(k)}  \\ \hat {\mathbf{v}}_B^{(k)} 
\end{bmatrix} 
= \frac{1}{2} \begin{bmatrix}
\mathbf{L} \\ 2\mathbf{I} 
\end{bmatrix}\ \hat{\mathbf{v}}_H^{(k)},
\end{equation}
in particular, $\hat {\mathbf{v}}_B^{(k)} = \hat{\mathbf{v}}_H^{(k)}$. Thus (\ref{eqn_B}) is equivalent to
\begin{equation}
\frac{1}{4} \mathbf{A}_H \hat{\mathbf{v}}_B^{(k)} = 
\frac{1}{4} \begin{bmatrix}
\mathbf{L} & 2\mathbf{I} 
\end{bmatrix} \tilde{\mathbf{b}} -  \frac{1}{4} \mathbf{A}_H \mathbf{x}_B^{(k)},
\end{equation}
which by linearity is equivalent to the second equation of (\ref{eqn_cyclic}). This discussion not only gives another explanation of the fact that the $(h,\ 2h)$ two-level cycle with $\omega=1$ converges after 1 post red-smoothing, but more importantly, it presents us a way to construct the (Galerkin) coarse-level operator as well as transfer operators by cyclic (red-black) reduction. We shall see that this generalizes readily to red-black reduction for higher dimensions.

\subsection{Two-dimensional local Fourier analysis} \label{Sec: rb_2d}
%\subsubsection{Red-black smoothing analysis in two dimensions}
This section corresponds to Section~\ref{Sec: rb_2d_std}, concerning the two-dimensional discrete Poisson's equation with \eqref{L_h_2}, and the two-level cycle with red-black smoothing and standard coarsening.

Analogous to (\ref{S_R phi}) and (\ref{S_B phi}) in one dimension, the red-black smoothing operator satisfies
\begin{equation}\label{S_hat_rb_2d}
\begin{cases}
S^R_h\ \phi_h (\cdot,\ \thv) =
\frac{1}{2}\left[(1 + \symF S^J_h({\thv})) \ \phi_h (\cdot,\ \thv)+(1 - \symF S^J_h({\thv})) \ \phi_h (\cdot,\ \bar{\thv})  \right],\\
S^B_h\ \phi_h (\cdot,\ \thv) = \frac{1}{2}\left[(1 + \symF S^J_h({\thv}))  \ \phi_h (\cdot,\ \thv) -(1 -\symF S^J_h({\thv}))\ \phi_h (\cdot,\ \bar{\thv})  \right],
\end{cases}
\end{equation}
where
\begin{equation}
\symF S^J_h({\thv}) = 1-\frac{\omega h^2}{4} \symF L_h({\thv}) = 1-2\omega \xi_{\thv}.
\end{equation}
This follows form the equation 
\begin{equation}\label{eqn_minus_2d}
\phi_h(\jv h,\ \bar{\thv}) = (-1)^{j_1\pm j_2} \phi_h(\jv h,\ \thv),\quad \jv=[j_1,\ j_2]^T \in \Integer^2,
\end{equation}
similar to (\ref{eqn_minus}).
Thus we have the invariance of the red-black smoother on the two-dimensional subspaces $\mathbb{E}_{\thv}$ \eqref{E_2d},
%\begin{equation}
%\mathbb{E}_{\thv} = \text{span} \{\phi_h (\cdot,\ \thv),\quad \phi_h (\cdot,\ \bar{\thv})\},\quad \thv \in \Theta,
%\end{equation}
on which the smoother has matrix representation
$
\symE S(\thv) = \symE S^B (\thv)\ \symE S^R(\thv)
$
with
\begin{equation}
\begin{aligned}
\symE S^R(\thv) &= 
\frac{1}{2}\begin{bmatrix}
1+\symF S^J_h({\thv}) & 1-\symF S^J_h(\bar{\thv}) \\
1-\symF S^J_h({\thv}) & 1+\symF S^J_h(\bar{\thv})
 \end{bmatrix}
 = \mathbf{I} - \omega \begin{bmatrix}
 1 & -1\\ -1 & 1
 \end{bmatrix} \begin{bmatrix}
 \xi_{\thv} \\ & 1-  \xi_{\thv}
 \end{bmatrix},\\ 
\symE S^B (\thv) &= 
\frac{1}{2}\begin{bmatrix}
1+\symF S^J_h({\thv}) & -(1-\symF S^J_h(\bar{\thv})) \\
-(1-\symF S^J_h({\thv})) & 1+\symF S^J_h(\bar{\thv})
 \end{bmatrix}
  = \mathbf{I} - \omega \begin{bmatrix}
 1 & 1\\ 1 & 1
 \end{bmatrix} \begin{bmatrix}
 \xi_{\thv} \\ & 1-  \xi_{\thv}
 \end{bmatrix}.
 \end{aligned}
\end{equation}
This gives \eqref{S_hat_RB_2d}.

Consider the coarse grid
\ba\label{G_H_2d_std}
    \grid G_H = \Big\{  \xv =\jv\, H:~ \jv=[j_1,\ j_2]^T \in \Integer^2 \Big\}.
\ea
by standard coarsening with $H = 2h$. 
%Note that $\grid G_H$ is embedded in $\grid G_h$ and consists of $\frac{1}{4}$ of the grid points. 
The $(h,\ 2h)$ two-level cycle has the four-dimensional eigen-subspaces $\mathbb{F}_{\boldsymbol{\theta}}$ as given in (\ref{F}).
In particular, following from (\ref{eqn_minus_2d}),
\ba
    \phi_h (\xv,\ \th_1,\ \th_2) &= \phi_h (\xv,\ \thbar_1,\ \thbar_2) = \phi_h (\xv,\ \thbar_1,\ \th_2) = \phi_h (\xv,\ \th_1,\ \thbar_2) \\
    &= \phi_H (\xv,\ 2\th_1,\ 2\th_2),\quad \forall \xv \in \grid G_H,\quad [\theta_1,\  \theta_2]^T \in \ThetaLowStd. 
\ea

We now finish the local Fourier analysis of the two-level cycle on $\mathbb{F}_{\thv}$ by introducing the coarse-level correction $K_h^H$. 
Consider the full-weighting restriction and linear interpolation operators (\ref{transfer_2d}).
%with $I_h^H = (I^h_H)^* = \frac{1}{4} (I^h_H)^T$. 
For each $\thv \in \ThetaLowStd$, 
%on $\grid G_H$,
\begin{equation}\label{transfer_2dstd_eqn}
\begin{cases}
I_h^H \PhiVectorF = \phi_H(\cdot,\ 2\thetav)\ \symE I_h^H(\boldsymbol{\theta}),\\
I^h_H\ \phi_H(\cdot,\ 2\thetav) =  \PhiVectorF \symE I^h_H (\boldsymbol{\theta}),
\end{cases}
\end{equation}
where
\begin{equation}\label{transfer_hat_2d}
\symE I_h^{H}(\thv) = \begin{bmatrix}
\symF I_h^H(\th_1,\ \th_2) &
\symF I_h^H({\bar{\th}_1},\ {\bar{\th}_2}) &
\symF I_h^H({\bar{\th}_1},\ \th_2) &
\symF I_h^H(\th_1,\ {\bar{\th}_2})
\end{bmatrix},\quad \symE I^h_{H}(\thv) = \symE I_h^{H}(\thv)^T,
\end{equation}
with 
\begin{equation}
\symF I_h^H(\th_1,\ \th_2) =( 1- \sin^2 \frac{\th_1}{2}) ( 1- \sin^2 \frac{\th_2}{2}) 
,\quad [\th_1,\ \th_2]^T \in \Theta.
\end{equation}

On $\grid G_H$, consider both the non-Galerkin coarse-level operator (\ref{L_H_2d})
with the same stencil as the fine-level operator $L_h$, with eigenvalue
\begin{equation}
\symF L_H(2\thv) = \frac{4}{H^2} (\sin^2 \th_1 + \sin^2 \th_2);
\end{equation}
and the Galerkin coarse-level operator
\begin{equation}\label{L_2h_G}
L_H =  I_h^H L_h I^h_H = 
\frac{1}{4H^2} \begin{bmatrix}
-1 & -2 & -1\\
-2 & 12 & -2\\
-1 & -2 & -1
\end{bmatrix}_H,
\end{equation}
with eigenvalue
\begin{equation}
\symF L_H(2\thv) = \symE I_h^H(\thv)\ \symE L_{\mathbb{F}}(\thv)\ \symE I^h_H (\thv),
\end{equation}
where
\begin{equation}
\symE L_{\mathbb{F}}(\thv) = \begin{bmatrix}
\symF L_h(\th_1,\ \th_2)\\ & \symF L_h(\bar \th_1,\ \bar \th_2) \\
& & \symF L_h(\bar \th_1,\ \th_2) \\
& & & \symF L_h(\th_1,\ \bar \th_2) 
\end{bmatrix} 
\equiv  \begin{bmatrix}
\symE L(\th_1,\ \th_2) \\ &  \symE L(\bar \th_1,\ \th_2) 
\end{bmatrix}
\end{equation}
is the matrix symbol of $L_h$ on $\mathbb{F}_{\thv}$. We have 
$L_H\ \phi_H (\cdot,\ 2\thv) = \symF L_H(2\thv)\ \phi_H (\cdot,\ 2\thv).$

Thus on each $\mathbb{F}_{\thv}$ ($\thv \neq \mathbf{0}$), the coarse-level correction
$
K_h^H 
$
has the matrix representation 
\begin{equation}\label{K_hat_2}
\symE K (\thv) = \mathbf{I} - \symF L_H^{-1}(2\thv)\
\symE I^h_H (\thv)\ \symE I_h^H(\thv)\ \symE L_{\mathbb{F}}(\thv).
\end{equation}
With the Galerkin $L_H$, $K_h^H$ is a projector.

Finally, we have the following theorem.
\begin{theorem}[red-black smoothing and standard coarsening in two dimensions - LFA]\label{Th_2d_std}
The $(h,\ 2h)$ two-level iteration operator $M_h^H$ in two dimensions, with the red-black smoother, standard coarsening and transfer operators \eqref{transfer_2d}, has its matrix representation on each $\mathbb{F}_{\thv}$ (\ref{F}) ($\thv \neq \mathbf{0}$) 
\begin{equation}
\symE M_{\mathbb{F}}(\thv) = \symE S_{\mathbb{F}}^{\nu_2}(\thv) \ \symE K(\thv)\ \symE S_{\mathbb{F}}^{\nu_1}(\thv),
\end{equation}
where $\symE S_{\mathbb{F}}(\thv)$ is given by (\ref{S_hat_F}) and $\symE K(\thv)$ is given by (\ref{K_hat_2}).
\end{theorem}

\medskip
\subsection{Red-black coarsening}\label{Sec: rbC_LFA}
In this section we discuss some further theoretical considerations of red-black coarsening presented in Section \ref{M2: rbC}.
We look at the three-level local Fourier analysis in \ref{Sec: rbC_3_LFA}, which is of theoretical interest.
We consider an anisotropic Poisson's equation in \ref{Sec: rbC_anisotropic}, and a fourth-order accurate discretization in \ref{Sec: rbC_4o}.
We discuss the connections of red-black coarsening with reduction methods as well as algebraic multigrid in \ref{Sec: cyclic_2d} and \ref{Sec: connections}.
\subsubsection{Three-level red-black coarsening}\label{Sec: rbC_3_LFA}
\begin{figure}[h!]
	\newcommand{\figWidth}{.495\linewidth}
	\begin{center}
%		\plotFig{Theta_2d_rb2}{\figWidth}\\
		\begin{tikzpicture}[scale=.5]
\useasboundingbox (-2,-.5) rectangle (8,8);
\draw
(0,0) -- (8,0) -- (8,8) -- (0,8) -- (0,0);
\draw (4,0) -- (4,8);
\draw (0,4) -- (8,4);
\draw[densely dotted, thick, orange]  (0,4) -- (4,0) -- (8,4) -- (4,8) -- (0,4);
\draw[thick, red]  (2,2) -- (6,2) -- (6,6) -- (2,6) -- (2,2);
\draw[fill] (0,0) circle (1pt) node[align=left,   below] {$-\pi$}
							   node[align=left,   left] {$-\pi$};
\draw[fill] (8,0) circle (1pt) node[align=right, below] {$\pi$};
\draw[fill] (0,8) circle (1pt) node[align=right, left] {$\pi$};
\draw[fill] (4,0) circle (1pt) node[align=center, below] {$0$} 
							   node[align=center, below=10pt] {$\theta_1$};
\draw[fill] (0,4) circle (1pt) node[align=center, left] {$0$} 
							   node[align=center, left=10pt] {$\theta_2$} ;
\draw[fill] (4,4) circle (1pt) node[align=center, red, above] {$\Theta^{(2)}$};
\draw[red,fill=red] (2,4) circle (1pt) node[align=left, red, left] {$-\frac{\pi}{2}$};
\draw[red,fill=red] (6,4) circle (1pt) node[align=left, red, right] {$\frac{\pi}{2}$};
\draw[red,fill=red] (4,2) circle (1pt) node[align=left, red, below] {$-\frac{\pi}{2}$};
\draw[red,fill=red] (4,6) circle (1pt) node[align=left, red, above] {$\frac{\pi}{2}$};
\end{tikzpicture} 
	\end{center}
	\caption{$\Theta$, 2D, 3-level, red-black coarsening}
	\label{fig: Theta_2d_rb2}
\end{figure}
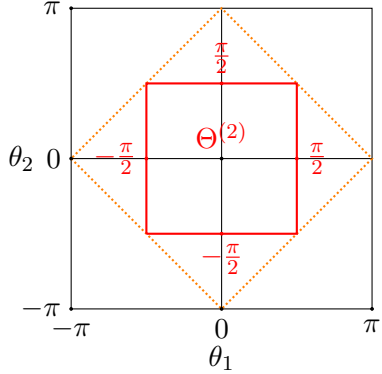

It is instructive to focus our interest on the $(h,\ H,\ 2h)$ three-level cycle with red-black coarsening for a moment, which is to be compared with the $(h,\ 2h)$ two-level cycle with standard coarsening in \ref{Sec: rb_2d}. Note that the coarsest grid $\grid G_{2h}$ is the same, and $\Theta_{rb}^{(2)} = [-\frac{\pi}{2},\ \frac{\pi}{2})^2$ is the same as the
$\ThetaLowStd$ for standard coarsening (see Figure \ref{fig: Theta_2d_rb2}). Moreover, we have the
eigen-subspaces of the three-level cycle as $\mathbb{E}^{(2)}_{\thv} = \mathbb{F}_{\thv}$, since
$(\sqrt{2}\mathbf{U})^{-1} \overline{\sqrt{2}\mathbf{U} \thv} = [\bar{\th}_1,\ \th_2]^T.$

The three-level $\gamma$-cycle iteration operator is $M_h = S_h^{\nu_2}\ K_h\ S_h^{\nu_1}$, with the
coarse-level correction
\begin{equation}
K_h \equiv \tilde{K}_h^H = I_h - I_H^h \tilde{L}_H^{-1} I_h^H L_h,
\end{equation}
where the coarse-level solve on $\grid G_H$ is approximated by $\gamma$ of the $(H,\ 2h)$ cycles (with initial guess 0) as
\begin{equation}\label{L_H_approx}
\tilde{L}_H^{-1} = \big(I_H - (M_H^{2h})^{\gamma}\big)\ L_H^{-1}.
\end{equation}
Its matrix representation on $\mathbb{F}_{\thv}$ is thus given by
\begin{equation}\label{M_3level}
\symE M_{\mathbb{F}}(\thv) = \symE S_{\mathbb{F}}^{\nu_2}(\thv) \ \tilde{\symE K}(\thv)\ \symE S_{\mathbb{F}}^{\nu_1}(\thv),\quad
\tilde{\symE K}(\thv) = \mathbf{I} - \symE I^h_{H,\ \mathbb{F}} (\thv)\  \tilde{\symE L}_H^{-1}(\thv)\ \symE I_{h,\ \mathbb{F}}^H(\thv)\ \symE L_{\mathbb{F}}(\thv)
\end{equation}
where $\tilde{\symE L}_H^{-1}(\thv)$ is the $2 \times 2$ matrix representation of the approximate coarse-level solve $\tilde{L}_H^{-1}$ on $\mathbb{E}_{\sqrt{2}\mathbf{U}\thv}$, and
\begin{equation}
\symE I_{h,\ \mathbb{F}}^{H}(\thv) = \begin{bmatrix}
\symF I_h^H(\th_1,\ \th_2) & \symF I_h^H(\bar \th_1,\ \bar \th_2) \\
& & \symF I_h^H(\bar \th_1,\ \th_2) & \symF I_h^H(\th_1,\ \bar \th_2) 
\end{bmatrix},\quad 
\symE I^h_{H,\ \mathbb{F}}(\thv) = \symE I_{h,\ \mathbb{F}}^{H}(\thv)^T,
\end{equation}
with
\begin{equation}
\symF I_h^H(\th_1,\ \th_2) =1- \frac{1}{2}(\sin^2 \frac{\th_1}{2} + \sin^2 \frac{\th_2}{2}) 
,\quad [\th_1,\ \th_2]^T \in \Theta.
\end{equation}
If we compare (\ref{M_3level}) to (\ref{M_2}) and (\ref{K_hat_2}) corresponding to the $(h,\ 2h)$ cycle with standard coarsening, the main difference is that the modes in $ \mathbb{E}_{\th_1,\ \th_2}$ and $\mathbb{E}_{\bar{\th}_1,\ \th_2}$ are dealt with differently by $\tilde{L}_H^{-1}\ (\approx {L}_H^{-1}) $, as opposed to by a single coarse-level solve $L_{2h}^{-1}$.

{
\newcommand{\figWidth}{6cm}% height 
\newcommand{\trimfig}[2]{\trimh{#1}{#2}{.0}{.0}{.0}{.0}}
\begin{figure}[h!]
\begin{center}
\begin{tikzpicture}[scale=1]
  \useasboundingbox (0,0.5) rectangle (15,13);  % set the bounding box (so we have less surrounding white space)
 % V ----------------------
   \draw(0.0,6.5) node[anchor=south west,xshift=-15pt,yshift=-8pt] {\trimfig{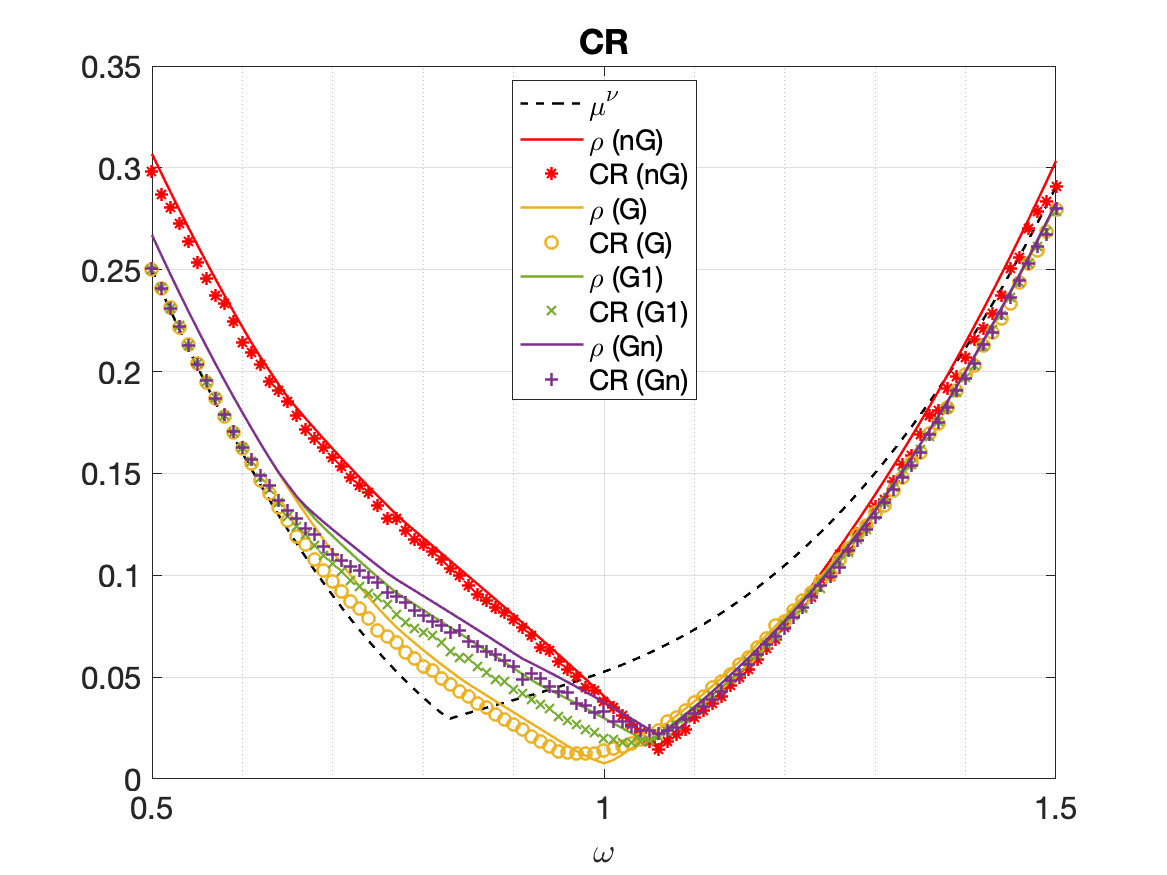}{\figWidth}};
   \draw(7.5,6.5) node[anchor=south west,xshift=-15pt,yshift=-8pt] {\trimfig{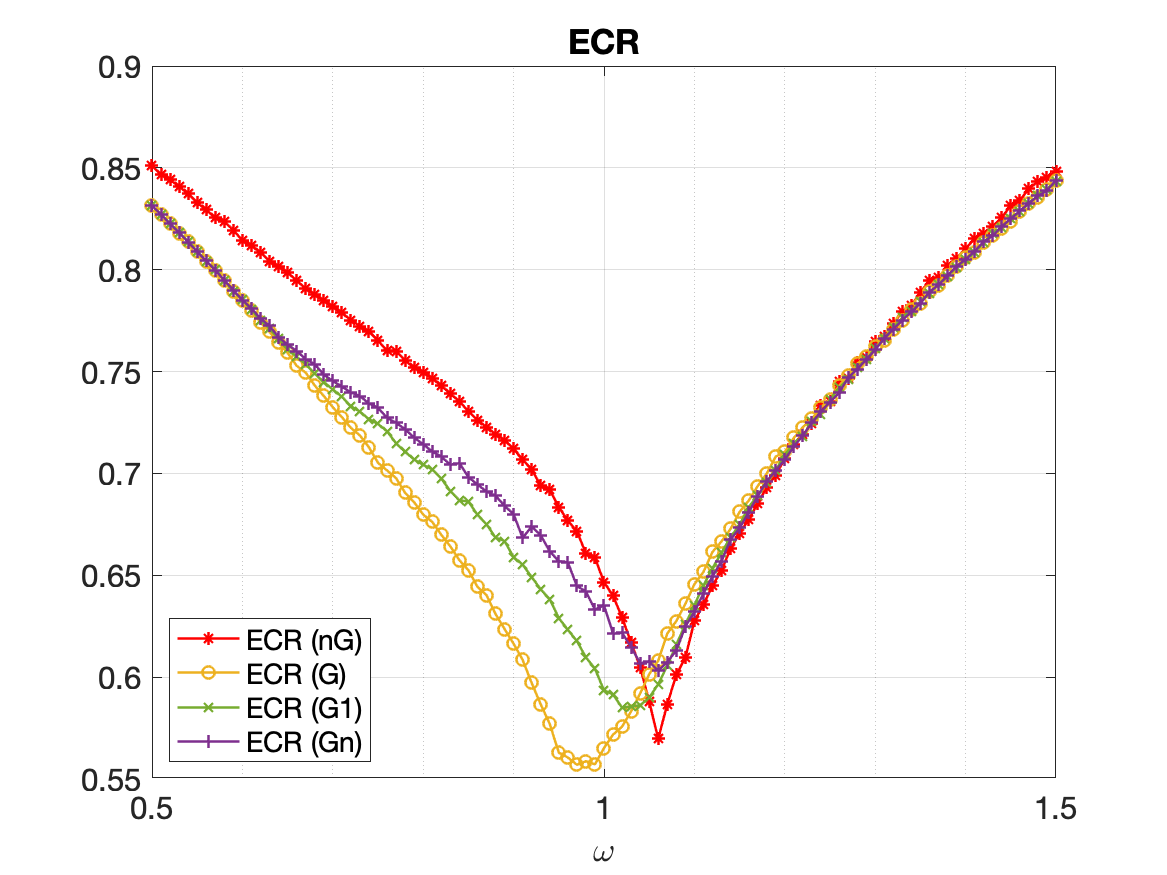}{\figWidth}};

   %  label:
   % \circleVortex{11.5}{7.8}{optimal $\omega \ne 1$};
   \circleLabel{11.5}{7.3}{optimal $\omega \ne 1$}{DarkGreen}{.25};

   % title:
   \draw(7.5,12.4)  node[draw=\colourTwoD,very thick,fill=\colourTwoD!20,anchor=south,yshift=-4pt] {\small \DIM, \ORD, \LEV[3], \CYC{V} cycle, \Cc[RB]}; 

   % Combined cartoons
   \orderTwoVcycleThreeLevelRotated{xshift=7cm,yshift=9.5cm,scale=1}{$\nu\!=\!2$}{};
   %\orderTwoVcycleThreeLevelRotated{xshift=5.5cm,yshift=11.8cm,scale=1}{$\nu\!=\!2$}{};
   %\orderTwoVcycleThreeLevelRotated{xshift=11.25cm,yshift=11.5cm,scale=1}{$\nu\!=\!2$}{};

 % W ----------------------
   \draw(0.0,0.0) node[anchor=south west,xshift=-15pt,yshift=-8pt] {\trimfig{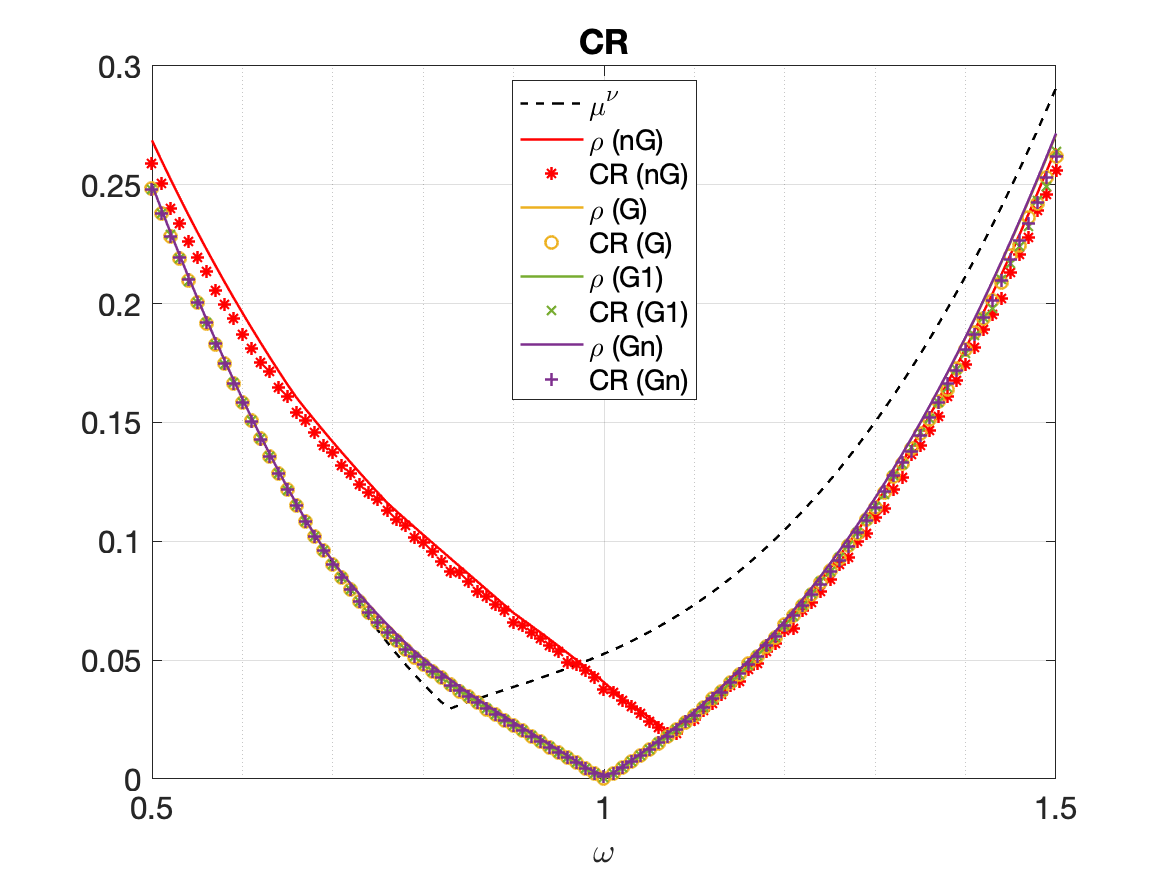}{\figWidth}};
   \draw(7.5,0.0) node[anchor=south west,xshift=-15pt,yshift=-8pt] {\trimfig{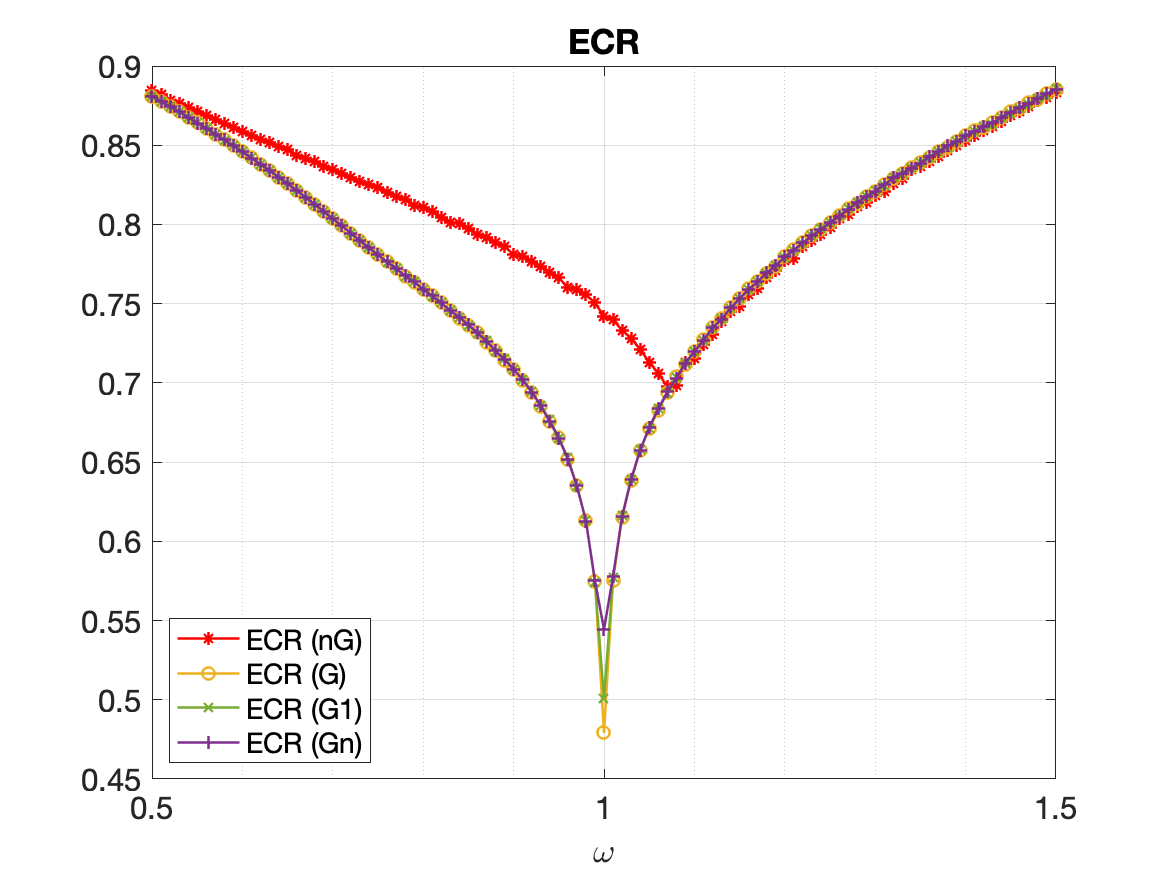}{\figWidth}};
   
   %  label:
   %\circleVortex{11.65}{.9}{optimal $\omega\!=\! 1$};
   \circleLabel{11.65}{.9}{optimal $\omega\!=\! 1$}{DarkGreen}{.25};

   % title:
   \draw(7.5,5.9)  node[draw=\colourTwoD,very thick,fill=\colourTwoD!20,anchor=south,yshift=-4pt] {\small \DIM, \ORD, \LEV[3], \CYC{W} cycle, \Cc[RB]}; 
   
   % Combined cartoons
    \orderTwoWcycleThreeLevelRotated{xshift=7cm,yshift=2.5cm,scale=1}{$\nu\!=\!2$};
   %\orderTwoWcycleThreeLevelRotated{xshift=5.5cm,yshift=4.8cm,scale=1}{$\nu\!=\!2$};
   %\orderTwoWcycleThreeLevelRotated{xshift=12.75cm,yshift=2.25cm,scale=1}{$\nu\!=\!2$};
   
% grid:
% \draw[step=1cm,gray] (0,0) grid (15,13.5);
\end{tikzpicture}
\end{center}
\caption{CR and ECR versus $\omega$. \DIM, \ORD, 
  \LEV[3] cycles with red-black smoothing and
  \CC[red-black],
  \commentB{
  \wdh{{\red FIX ME}
    non-Galerkin (nG$=\ngcg{2}$) and Galerkin (G$=\gcg{2}{2}{2}$, G1$=\gcg{?}{?}{?}$, Gn$=\gcg{?}{?}{?}$) 
  }}
\kl{non-Galerkin (nG=$\LNAH$) and Galerkin (G=$\LGGH$, G1$=\LGAAH$, Gn) }
coarse-level operators ($N = 32$).
Top: $\CYC{V}[1, 1]$ cycle, bottom: $\CYC{W}[1, 1]$ cycle.
}
  \label{fig:CR_3level_2d_rb}
\end{figure}
}

The results for the three-level $V[1, 1]$ and $W[1, 1]$
RBC cycles are shown in Figure~\ref{fig:CR_3level_2d_rb}; these should be compared to the two-level RBC results in Figure~\ref{fig:CR_2level_2d_rb} and the two-level SC results in Figure~\ref{fig:redBlackSmootherStandardCoarsening1d2d} (bottom).
The best convergence rate (CR) for the $V[1, 1]$ RBC cycle is $CR\approx .01$,
while the best convergence rate (CR) for the $W[1, 1]$ RBC cycle is $CR\approx .001$. 
%{\red (Get accurate numbers)}.
These rates are substantially better than the corresponding rate with standard coarsening (Figure~\ref{fig:redBlackSmootherStandardCoarsening1d2d}, bottom) which is  $CR\approx .03$. 
It is noted that the shape of the CR graphs for the RBC W cycle with Galerkin coarse-level operators (Figure~\ref{fig:CR_3level_2d_rb}, lower left),
is very similar to the corresponding two-level RBC result in Figure~\ref{fig:CR_2level_2d_rb}. This is in agreement with our
observation that the convergence of the three-level RBC cycle should approach that of the two-level RBC cycle as the
three-level coarse-level correction is made more accurate (i.e. increasing $\gamma$).

A comparison of the effective convergence rates (ECR) shows that the three-level $W[1, 1]$
has an optimal $ECR\approx .47$ (achieved at $\omega=1$) while for the $V[1, 1]$
cycle $ECR \approx .56$ (achieved at $\omega \approx 0.97$), with Galerkin coarse-level operators.
For the V cycle, the true Galerkin coarse-level operator gives the best ECR results.
For the W cycle the true Galerkin is still the best, although the G1 and Gn coarse-level operators are very good as well.

\subsubsection{Anisotropy}\label{Sec: rbC_anisotropic}
If a bit of generalization is added to the two-dimensional Poisson's equation by adding anisotropy to the Laplacian, it is straightforward to adapt the proposed multigrid algorithm with a re-scaling of the two dimensions. Consider the continuous operator
\begin{equation}
L = -(c_1 \partial_x^2 + c_2 \partial_y^2)
\end{equation} 
where  $c_1$ and $c_2$ are constants with $c_1 + c_2 = 2$. ($L = -\Delta$ if $c_1 = c_2 = 1$.) We assume the asymmetry is moderate in that the sizes of the coefficients $c_1$ and $c_2$ are still comparable, for in the case of extreme anisotropy it would be more reasonable to consider semi-coarsening rather than red-black coarsening. (This situation of moderate anisotropy is very common in a composite grid with change-of-variables suited for general geometry.) In this case, the problem can be transformed back to the isotropic Poisson's equation, by either a change of variables, or setting the fine grid with grid-spacing $(h_1,\ h_2)$ with 
%$\frac{h_1^2}{h_2^2} = \frac{c_1}{c_2}.$
\begin{equation}\label{h_c}
\frac{h_1^2}{h_2^2} = \frac{c_1}{c_2}.
\end{equation}
Given $(h_1,\ h_2)$, red-black coarsening as before yields a coarse grid with $H_1 = H_2 = \sqrt{h_1^2+h_2^2}$. Note that although this coarse grid is not physically orthogonal, if we consider it as discretization on the domain in a set of new variables $(r,\ s)$ from a continuous linear mapping $\mathbf{V}$, the operator $L$ is still diagonal:
\begin{equation}
-L = \partial_{\mathbf{x}}^T \mathbf{A} \partial_{\mathbf{x}} = 
\partial_{\mathbf{r}}^T \mathbf{V A V}^T \partial_{\mathbf{r}},
\end{equation}
where $\mathbf{A} = \begin{bmatrix}
c_1\\ & c_2
\end{bmatrix}$, and $\mathbf{V A V}^T$ is diagonal if (\ref{h_c}) is satisfied. 
Thus the coarse-level operator, Galerkin or non-Galerkin, can be constructed the same way as the isotropic case.

%------------------------------------------------------------------------------------------
{
\newcommand{\figWidth}{6cm}% height 
\newcommand{\trimfig}[2]{\trimh{#1}{#2}{.0}{.0}{.0}{.0}}
\begin{figure}[h!]
\begin{center}
\begin{tikzpicture}[scale=1]
  \useasboundingbox (0,0.5) rectangle (15,6.);  % set the bounding box (so we have less surrounding white space)

   \draw(4.0,0.0) node[anchor=south west,xshift=-15pt,yshift=-8pt] {\trimfig{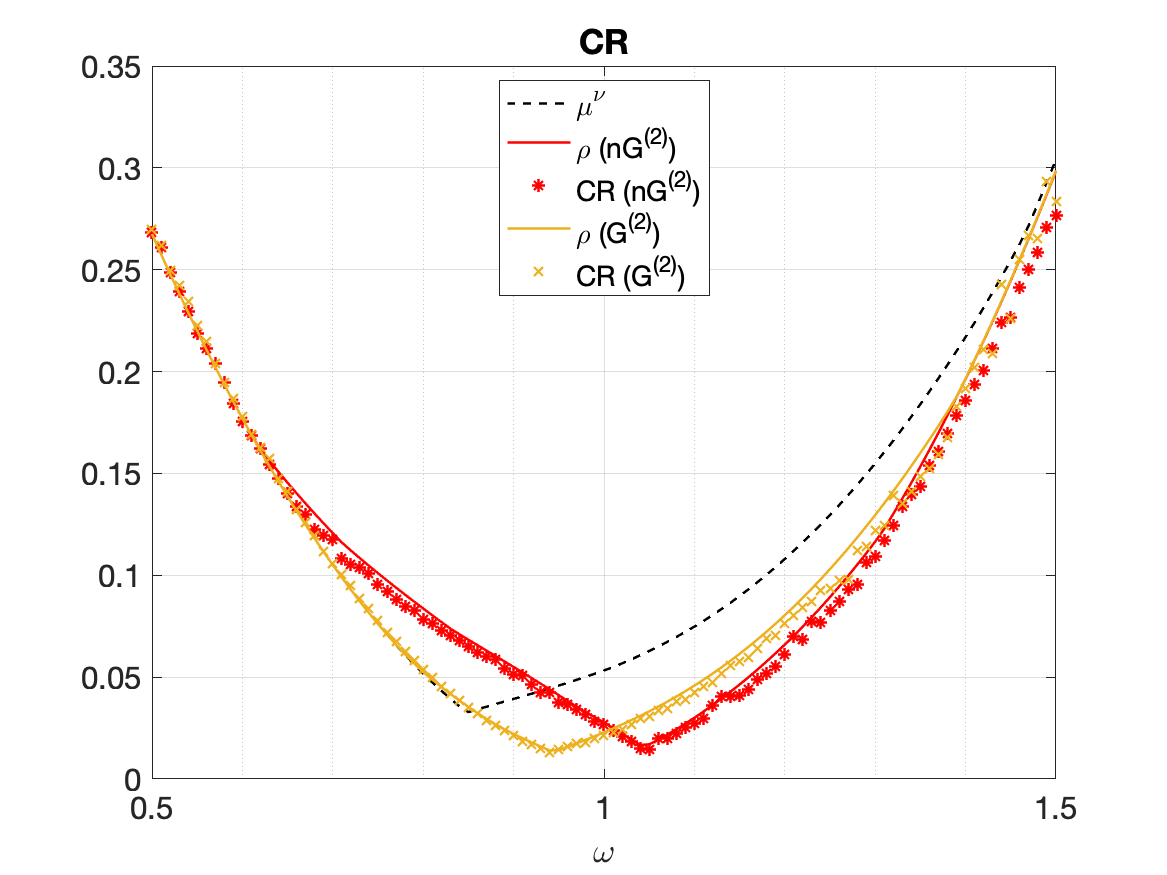}{\figWidth}};
   %\circleVortex{7.8}{0.8}{$ECR=0.6$};
   \circleLabelLeft{7.8}{0.7}{$ECR=0.6$}{purple}{.3};
   %\draw(7.5,0.0) node[anchor=south west,xshift=-15pt,yshift=-8pt] {\trimfig{ECR_2drt_4o_l1_nu2}{\figWidth}};
   
   \VcycleTwoLevel{xshift=6.5cm,yshift=4.7cm,scale=1}{$\nu\!=\!2$}{};
   \draw(8, 6)  node[draw=\colourTwoD,very thick,fill=\colourTwoD!20,anchor=south,yshift=-4pt,inner sep=2.5pt] {\small \DIM, \ORD[4], \LEV, \Cc[RB]}; 
  
% grid:
% \draw[step=1cm,gray] (0,0) grid (16,6);
\end{tikzpicture}
\end{center}
  \caption{CR versus $\omega$. \DIM, \ORD[4], \LEV\ $V[1, 1]$ cycle with red-black Gauss-Seidel smoothing and \CC[red-black], 2nd-order accurate non-Galerkin (nG$^{(2)} = \LNH{2}$) and Galerkin (G$^{(2)} = \LGqH{2}$) coarse-level operators ($N = 32$).}
  \label{fig:CR_2level_2d_4o_rb}
\end{figure}
}
{
\newcommand{\figWidth}{6cm}% height 
\newcommand{\trimfig}[2]{\trimh{#1}{#2}{.0}{.0}{.0}{.0}}
\begin{figure}[h!]
\begin{center}
\begin{tikzpicture}[scale=1]
  \useasboundingbox (0,0.5) rectangle (15,6);  % set the bounding box (so we have less surrounding white space)

   \draw(0.0,0.0) node[anchor=south west,xshift=-15pt,yshift=-8pt] {\trimfig{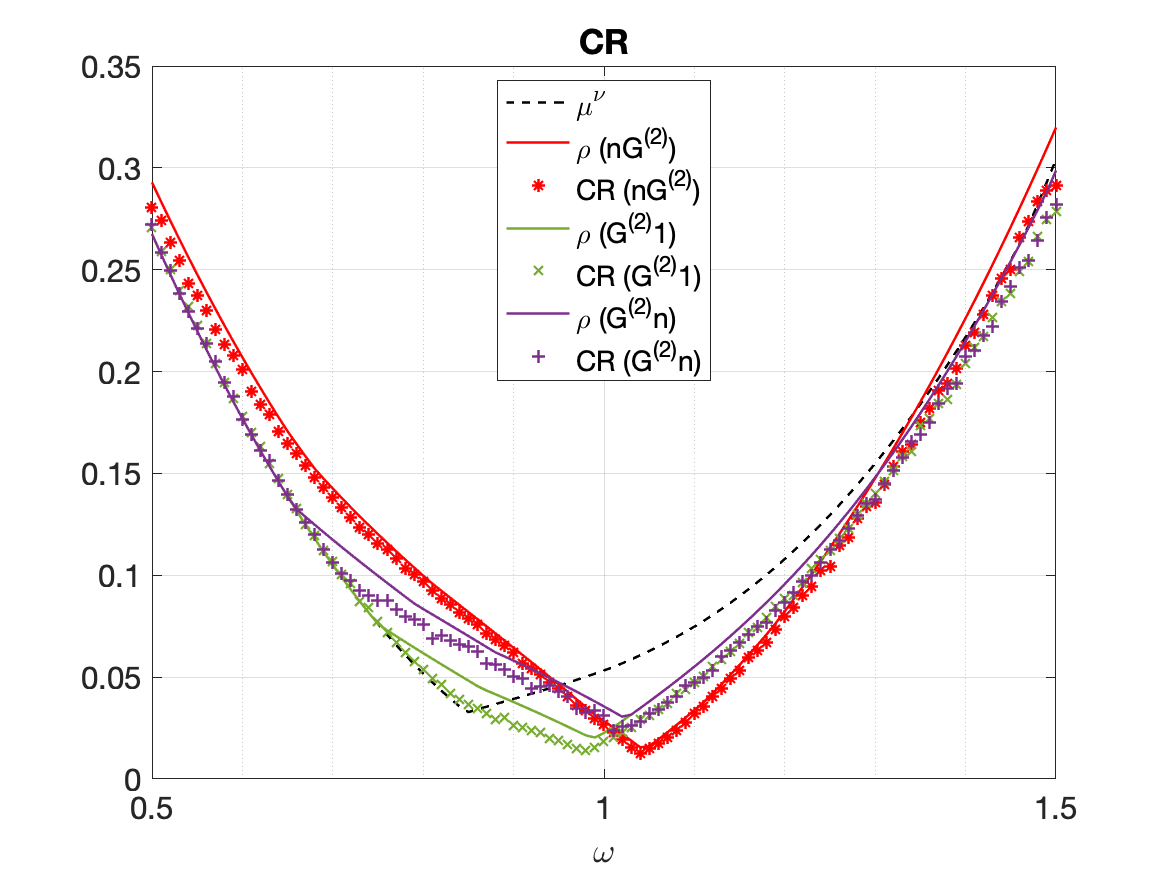}{\figWidth}};
  % \draw(7.5,6.0) node[anchor=south west,xshift=-15pt,yshift=-8pt] {\trimfig{ECR_2drt_4o_l2gamma10_nu2}{\figWidth}};
   %\circleVortex{4.5}{0.7}{$ECR=0.56$};
   \circleLabel{4.5}{0.7}{$ECR=0.56$}{purple}{.3};
   \VcycleThreeLevel{xshift=2.5cm,yshift=4.7cm,scale=1}{$\nu\!=\!2$}{};
   
   \draw((7.5,0.0) node[anchor=south west,xshift=-15pt,yshift=-8pt] {\trimfig{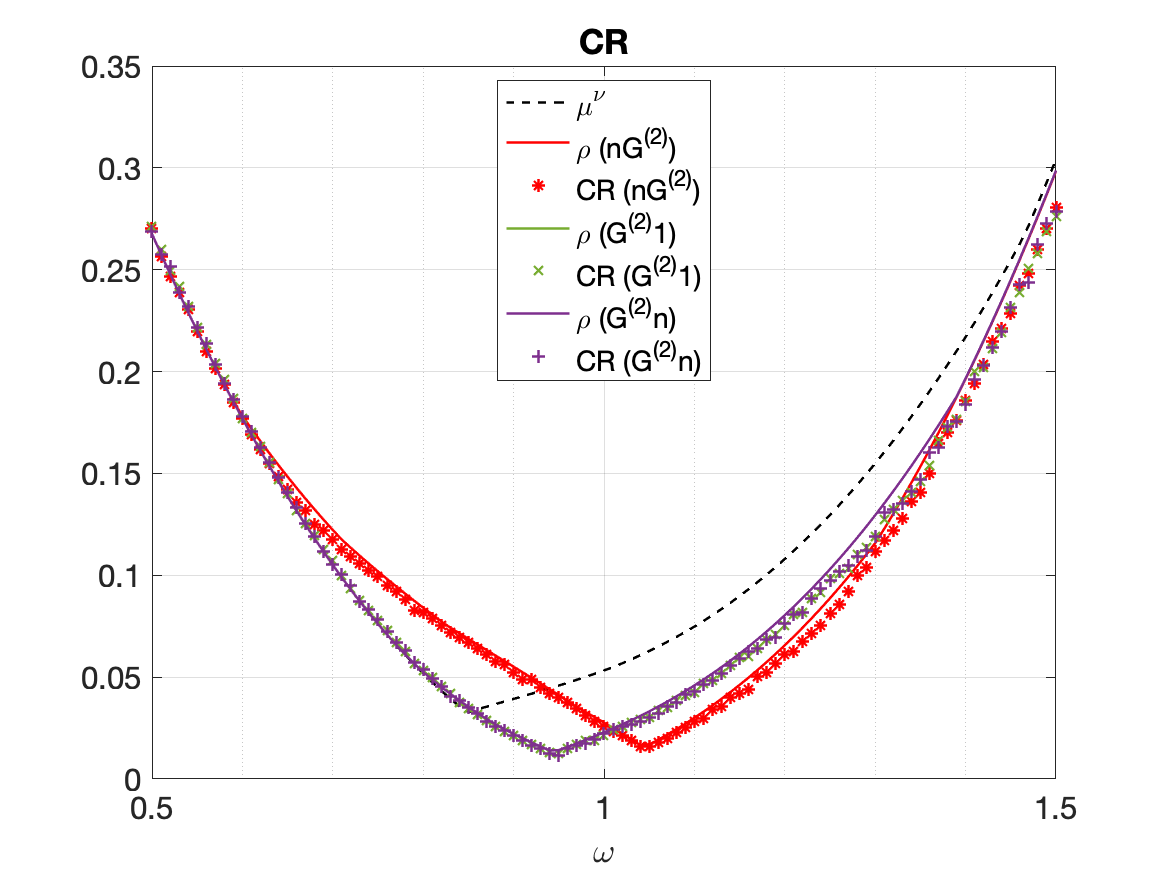}{\figWidth}};
  % \draw(7.5,0.0) node[anchor=south west,xshift=-15pt,yshift=-8pt] {\trimfig{ECR_2drt_4o_l2gamma20_nu2}{\figWidth}};
   %\circleVortex{11.3}{0.7}{$ECR=0.67$};
   \circleLabelLeft{11.3}{0.7}{$ECR=0.67$}{purple}{.3};
   \WcycleThreeLevel{xshift=10cm,yshift=4.7cm,scale=1}{$\nu\!=\!2$}{}; 
   
   \draw(7.5,5.8)  node[draw=\colourTwoD,very thick,fill=\colourTwoD!20,anchor=south,yshift=-4pt,inner sep=2.5pt] {\small \DIM, \ORD[4], \LEV[3], \Cc[RB]};  
   
% grid:
% \draw[step=1cm,gray] (0,0) grid (16,6);
\end{tikzpicture}
\end{center}
  \caption{CR versus $\omega$. \DIM, \ORD[4], \LEV[3]\ cycles with red-black Gauss-Seidel smoothing and \CC[red-black], 
  	2nd-order accurate non-Galerkin (nG$^{(2)} = \LNH{2}$) and 
  	Galerkin (G$^{(2)}$1=$\LGAH{2}$, G$^{(2)}$n) 
  	coarse-level operators ($N = 32$).
  Left: $\CYC{V}[1, 1]$ cycle, right: $\CYC{W}[1, 1]$ cycle.}
  \label{fig:CR_3level_2d_4o_rb}
\end{figure}
}
\subsubsection{Fourth-order accuracy}\label{Sec: rbC_4o}
Consider the fourth-order accurate finite difference discretization of the minus Laplacian (\ref{L_h_2d_4o}).
Following the idea in Section \ref{M1: orderC},
on the coarse grid $\grid G_H$, instead of the fourth-order accurate operator, we consider the second-order accurate $L_H^{(2)}$, both Galerkin and non-Galerkin, constructed from the second-order fine-level operator $L_h^{(2)}$ as in (\ref{L_h_2}). This is more computationally efficient, and as we can see in Figure \ref{fig:CR_2level_2d_4o_rb} as well as Figure \ref{fig:CR_3level_2d_4o_rb}, the convergence rates are not hindered by this choice, and the ECRs are comparable to those for the second-order accurate difference equation.

\subsubsection{Red-black reduction in two dimensions}\label{Sec: cyclic_2d}

Consider the matrix of the system for the two-dimensional discrete Poisson's equation with periodic boundaries, on a square grid with $N$ rows and $N$ columns, of size $N^2 \times N^2$ with row-wise ordering:
\begin{equation}
\mathbf{A}_h = \begin{bmatrix}
\mathbf{T} & -\mathbf{I}  & & & & -\mathbf{I}\\
-\mathbf{I} & \mathbf{T} & -\mathbf{I} \\
& -\mathbf{I} & \mathbf{T} & -\mathbf{I} \\
& & \ddots  & \ddots & \ddots \\
& & & -\mathbf{I} & \mathbf{T} & -\mathbf{I} \\
-\mathbf{I}  & &  & &  -\mathbf{I} & \mathbf{T}
\end{bmatrix},
\end{equation}
in which each $N \times N$ diagonal block is
\begin{equation}
\mathbf{T} = \begin{bmatrix}
4 & -1  & & & &   -1\\
-1 & 4 & -1 \\
& -1 & 4 & -1 \\
& & \ddots  & \ddots & \ddots \\
& & & -1 & 4 & -1 \\
-1  & &  & &  -1 & 4
\end{bmatrix}.
\end{equation}
Analogous to the one-dimensional case, cyclic reduction rearranges the block-cyclic $\mathbf{A}_h$ into odd and even block-rows, and eliminates the odd ones. This would corresponds to semi-coarsening by a factor of 2. Instead, we consider red-black ordering (and again consider a grid point red if its indexes add to an odd number) and get a permuted matrix as
\begin{equation}
\tilde{\mathbf{A}}_h = \begin{bmatrix}
4\mathbf{I} & -\tilde{\mathbf{L}}^T \\
- \tilde{\mathbf{L}} & 4\mathbf{I}
\end{bmatrix}.
\end{equation}
Here $\mathbf{I}$ and $\tilde{\mathbf{L}}$ are of size $\frac{N^2}{2} \times \frac{N^2}{2}$, and 
\begin{equation}
\tilde{\mathbf{L}} = \begin{bmatrix}
\mathbf{L}  & \mathbf{I} &  & & \mathbf{I} \\
\mathbf{I} & \mathbf{L}^T & \mathbf{I} \\
& \ddots & \ddots & \ddots \\
& & \mathbf{I} & \mathbf{L} & \mathbf{I}\\
\mathbf{I} & &  &  \mathbf{I} & \mathbf{L}^T
\end{bmatrix},
\end{equation}
in which each block is of size $\frac{N}{2} \times \frac{N}{2}$, and $\mathbf{L}$ is the same as (\ref{LL}). By one red-black reduction step, we eliminate the red points by factoring the matrix as $\tilde{\mathbf{A}}_h = \hat{\mathbf{L}} \hat{\mathbf{A}} \hat{\mathbf{L}}^T$ with
\begin{equation}
\hat{\mathbf{L}} = \begin{bmatrix}
\mathbf{I} \\ -\frac{1}{4}\tilde{\mathbf{L}} & \mathbf{I}
\end{bmatrix},\quad
\hat{\mathbf{A}} = \begin{bmatrix}
4\mathbf{I} \\
& \mathbf{A}_H
\end{bmatrix},
\end{equation}
in which each block is of size $\frac{N^2}{2} \times \frac{N^2}{2}$. Thus by
\begin{equation}
 \hat{\mathbf{A}} = \hat{\mathbf{L}}^{-1} \tilde{\mathbf{A}}_h \hat{\mathbf{L}}^{-T},\quad
 \hat{\mathbf{L}}^{-1} = \begin{bmatrix}
 \mathbf{I} \\ \frac{1}{4}\tilde{\mathbf{L}} & \mathbf{I}
 \end{bmatrix},
\end{equation}
we can get the (scaled) coarse-level matrix $\mathbf{A}_H$ on $\grid G_H = \grid G_h^B$ 
%(\ref{G_r2}) 
corresponding to the Galerkin coarse-level operator (\ref{L_G1}); as well as the restriction
$
\frac{1}{8} \begin{bmatrix}
\tilde{\mathbf{L}} & 4\mathbf{I} 
\end{bmatrix}
$
and interpolation
$
\frac{1}{4} \begin{bmatrix}
\tilde{\mathbf{L}} \\ 4\mathbf{I} 
\end{bmatrix}
$
operators corresponding to (\ref{transfer_r2}). Then following the same argument as in \ref{Sec: cyclic_1d}, we observe once again the equivalence:
\begin{theorem}[red-black reduction in two dimensions]
One step of red-black reduction in two dimensions is equivalent to a two-level multigrid cycle with red-black coarsening and one post red-smoothing step with $\omega = 1$. Moreover, the two-level algorithm converges in one cycle.
\end{theorem}

As we shall see in \ref{Sec: reduction}, red-black reduction is simply an equivalent formulation of what is called total reduction \cite{reduction}.
\subsubsection{Connections and discussions}\label{Sec: connections}

In this section we relate the proposed multigrid cycles with red-black coarsening to reduction methods and algebraic multigrid algorithms. In particular, from the starting point of a two-level algorithm for the simple 5-point stencil of the two-dimensional Laplacian, we discuss the connections and differences of perspectives when the methods are applied to more general settings, such as multiple levels and other difference operators.

\subsubsubsection{Total Reduction}\label{Sec: reduction}
Reduction methods for constant-coefficient difference operators are formulated and summarized in \cite{reduction}, in which the elimination process is replaced by applications of `conjugate' difference stars. Here we give a brief summary of the ideas, and in particular, identify red-black reduction as another formulation of total reduction in two dimensions.

For a given difference star (difference operator in stencil notation) in two dimensions
$S = [S_{ij}]$,
define its `conjugate' as
$\bar{S} = [(-1)^{i+j} S_{ij}]$.
In particular, consider the model problem with the difference star 
\begin{equation}
S_h = \begin{bmatrix}
& -1 \\
-1 & 4 & -1\\
& -1
\end{bmatrix}_h
\end{equation}
on $\grid G_h$ given by (\ref{G_h_d}) with $d=2$.
Then
\begin{equation}
\bar{S}_h = \begin{bmatrix}
& 1 \\
1 & 4 & 1\\
& 1
\end{bmatrix}_h.
\end{equation}
Note that the `conjugate' stars can be written as $S = S^+ + S^-$, and $\bar{S} = S^+ - S^-$, where
\begin{equation}
S^+ = \begin{bmatrix}
&  \\
&  4 &\\
& 
\end{bmatrix},\quad
S^- = \begin{bmatrix}
& -1 \\
-1 & &  -1\\
& -1
\end{bmatrix}.
\end{equation}
Applying to the difference star its `conjugate', the resulting difference star 
\begin{equation}
\bar{S}_h S_h = S_h\bar{S}_h = \begin{bmatrix}
& & -1 \\
& -2 & & -2 \\
-1 & & 12 & & -1\\
& -2 & & -2 \\
& & -1
\end{bmatrix}_h \equiv
\begin{bmatrix}
-1 & -2 & -1\\
-2 & 12 & -2\\
-1 & -2 & -1
\end{bmatrix}_H \equiv S_H
\end{equation}
is invariant on the coarse grid $\grid G_H = \grid G_h^B$ (\ref{G_r2}), $H = \sqrt{2}h$. Then an equation 
%$S_h x_h = b_h$ 
\begin{equation}
S_h\ x_h = b_h
\end{equation}
is reduced to
\begin{equation}\label{S_x_B}
S_H\ x_B = \bar{S}_h\ b_h,
\end{equation}
and given $x_B$,
\begin{equation}\label{S_x_R}
S^+_h\ x_R = S^-_h\ x_R + b_R
\end{equation}
with $x_h = [x_R,\ x_B]^T$, $b_h = [b_R,\ b_B]^T$ in red-black ordering. Note that $S^-_h\ x_R$ only involves $x_B$. This is application of one step of total reduction.  

It is straightforward to see that the above is completely equivalent to red-black reduction as formulated in \ref{Sec: cyclic_2d} by equaling
\begin{equation}
S_h= h^2 L_h,\quad S_H = 4H^2 L_H = 8h^2 L_H,\quad \bar{S}_h = 8 I_h^H,
\end{equation}
where $L_h$ is given by (\ref{L_h_2}), $L_H$ is the Galerkin coarse-level operator given by (\ref{L_G1}), and $I_h^H$ is the interpolation operator given by (\ref{transfer_r2}). In addition, (\ref{S_x_R}) corresponds to post red-smoothing. Thus once again, one step of total, or red-black reduction is equivalent to the two-level multigrid cycle with red-black coarsening and the Galerkin coarse-level operator, applied to the 5-point stencil (\ref{L_h_2}).

Incidentally, partial reduction, in contrast to total reduction, applies the `conjugate' operator only in $x$ or $y$ direction as
$\bar{S}^x = [(-1)^i S_{ij}]$ or $\bar{S}^y = [(-1)^j S_{ij}]$, which is equivalent to block-cyclic reduction, and is related to multigrid cycles with semi-coarsening.

The same total (red-black) reduction step, or equivalently $\bar{S}_H$, can be applied on the reduced system (\ref{S_x_B}); and so on, recursively until the size of the problem is proper for a direct solver. We can relate $m$ steps of total (red-black) reduction with an $(m+1)$-level multigrid V-cycle with red-black coarsening, in that the sets of variables of the reduced systems correspond to the grids on the coarse levels by red-black coarsening. 

However, there are some major differences between multigrid and reduction methods. Multigrid is among the error-correction algorithms, and the equations to be solved on coarse levels are approximate error equations, instead of reduced sets of equations or applications of `conjugate' difference operators. One important feature of the reduction methods is that all the operators such as $\bar{S}$, $S^-$ come from the difference operator $S$; or in the matrix language in \ref{Sec: cyclic_2d}, operators such as $\mathbf{A}_H$ and $\tilde{\mathbf{L}}$ come directly from the matrix $\mathbf{A}$. That is, in the equivalent two-level cycle, the transfer operators as well as the coarse-level operator are all constructed based on the fine-level operator. When our multigrid cycle with red-black coarsening are extended to multiple levels (as in Section \ref{Sec: rbC_multilevel}), since the coarse-level error equations need only to be solved approximately, we have the flexibility to keep using the simple transfer operators $I_h^H$ and $I_H^h$ without having to consider $\bar{S}_H$. Thus the resulting Galerkin coarse-level operators can be much simpler (with short-range stencils; coarse-level matrices still enjoy enough sparsity). Moreover, we have the flexibility to choose more naive Galerkin-like coarse-level operators or even non-Galerkin ones. The story is similar for different fine-level operators such as an anisotropic Laplacian (\ref{Sec: rbC_anisotropic}) or discretizations of higher-order accuracy (\ref{Sec: rbC_4o}).
In short, the different components in a geometric multigrid algorithm can be somewhat independent from each other, making use of the information from the continuous problem and the geometric grids.

It is also to be noted that total reduction and multigrid with red-black coarsening have the same difficulty of extending to three dimensions.

\subsubsubsection{Algebraic Multigrid}\label{Sec: AMG}
As with reduction methods, algebraic multigrid algorithms work directly with the discrete system of equations, even in a setting with geometrically based problems. Unlike reduction methods, which rely on the uniformity of the difference operator, algebraic multigrid algorithms locally adapt the coarsening strategy based on the connectivity of the nonzeros in the matrix. The `grid' $\grid G_h$ in algebraic multigrid corresponds to the variables or nodes of the directed graph associated with the matrix $A_h$. The algorithm sets up a C/F(coarse/fine)-splitting of the grid $\grid G_h = \grid G_h^C \cup \grid G_h^F$, and the set of C-variables is chosen as the coarse grid $\grid G_H = \grid G_h^C$. The variational property is usually assumed in algebraic multigrid, that is, the transfer operators are adjoint to each other ($I_H^h = (I_h^H)^*$), and the the coarse-level operator is constructed based on Galerkin recursion ($A_H = I_h^H A_h I_H^h$). Thus after coarsening, the only operator needed to be constructed is the interpolation $I_H^h$. With simple pointwise smoothers (often used in C/F-ordering), algebraic multigrid algorithms construct coarsening and interpolation following two basic principles: first, for fast convergence, the coarse-level correction should work well on algebraically smooth errors, that is, those components of error which the smoother is inefficient in reducing; second, for computational efficiency, the proportion of the C-variables in a given grid, as well as the complexity of the coarse-level operator, should be limited. 

Algebraic multigrid theory gives uniform (two-level) convergence estimates for certain classes of matrices. Of course, the more restricted and well-behaved the class of matrices is, the better the convergence estimates. Algebraic multigrid for sparse, symmetric positive (semi-)definite matrices, in particular those matrices coming from discretizations of elliptic partial differential equations, is very well studied, and the application of algebraic multigrid algorithms on these matrices are very efficient.
%In particular, those matrices coming from discretizations of elliptic partial differential equations enjoy many nice properties that make the application of algebraic multigrid algorithms very efficient.  

In the simplest case, consider the second-order discretization of the model problem (\ref{MP}) on the uniform Cartesian grid (\ref{G_h_d}) in two dimensions. For the first coarse level, standard algebraic multigrid coarsening strategy based on strong F-to-C direct connectivity leads to exactly the coarse grid (\ref{G_r2}) given by red-black coarsening. This comes readily from the shape and symmetry of the 5-point stencil (\ref{L_h_2}). In this special case the interpolation can be constructed so that the two-level multigrid cycle corresponds to a direct solver, if for instance, the range of the coarse-level correction operator aligns with the nullspace of the smoother, as illustrated in Section \ref{Sec: rbC}. As pointed out in \cite{MGtext},
as the process is preceded to coarser levels, the algebraic multilevel V-cycle is actually equivalent to total reduction. Such a multilevel direct solver is not very practical, however, due to growth of the complexity of the coarse-level operators by Galerkin recursion. One of the heuristic principles in constructing the interpolation for general matrices is to approximate the limit case of direct solvers by clever interplay of the interpolation and the smoother \cite{MGtext}. 

Although with red-black smoothing, the two-level cycle for the model problem with the 5-point stencil (\ref{L_h_2}) constructed by algebraic multigrid is identical with the one constructed by red-black coarsening (with the Galerkin coarse-level operator) in Section \ref{Sec: rbC}, since the motivation is different, after the second level, the coarse-level operators (including the coarse-level grids) and interpolation constructed by algebraic multigrid will deviate from the multilevel cycles proposed in Section \ref{Sec: rbC_multilevel},
and similarly for problems with other fine-level operators.

The fact that the coarsening and interpolation in algebraic multigrid algorithms are dependent on the fine-level operator (matrix), on the one hand helps ensure fast convergence, but on the other makes it harder to control the complexity of the coarse-level operators.
The interplay of the smoother and the coarse-level correction (in particular, the interpolation), and the balancing of convergence rates with computational costs, lie in the heart of any multigrid algorithm. In geometric multigrid, more focus is put on designing a better smoother, while in algebraic multigrid the stress is on the coarse-level correction (that is, coarsening and interpolation) \cite{MGtext}. For geometric multigrid, a low convergence rate is harder to achieve, while in algebraic multigrid the control of computational cost needs more attention, at least for well-behaved classes of problems. 

When we have a continuous problem, in particular a Poisson-like problem in mind, however, basing the algorithm solely on the discretized matrices throws away some potentially useful information, especially when many of the nice properties of the matrices to be made use of come directly from the continuous operators.
Following the spirit of algebraic multigrid, on the other hand, in geometric multigrid we can tailor the coarsening strategy more carefully according to the fine-level operator and the smoother, while still enjoy the simplicity of coarsening and transfers on coarser levels, as well as the flexibility of coarse-level operators.

%variable coefficients
%higher dimensions

\medskip
\subsection{Fourth-order accuracy}\label{Sec: orderC_LFA}
In this section we give more details of local Fourier analysis for the fourth-order discretization of the two-dimensional model problem as discussed in Section~\ref{M1: orderC}.
We also give a heuristic argument why operators of low-order accuracy can be used on the coarse levels. 
For simplicity here we consider standard coarsening. This idea has also been applied in \ref{Sec: rbC_4o} with red-black coarsening.

The fourth-order accurate negative Laplacian in two dimensions $L_h$ as given by (\ref{L_h_2d_4o}) has the Fourier symbol
%eigenvalues
\begin{equation}\label{L_hat_2d_4o}
\symF L_h(\thetav) = \frac{8}{h^2} \left(\bar{\xi}_{\thv} +  \frac{1}{3}\bar{\xi^2_{\thv}}\right), 
\end{equation}
where
\begin{equation}
\bar{\xi}_{\thv} \equiv \frac{1}{2} \left(\sin^2 \frac{\th_1}{2} + \sin^2 \frac{\th_2}{2}\right),\qquad
\bar{\xi^2_{\thv}} \equiv  \frac{1}{2} \left(\sin^4 \frac{\th_1}{2} + \sin^4 \frac{\th_2}{2}\right).
\end{equation}
The local Fourier analysis for this fourth-order accurate problem mimics the discussion in Section~\ref{Sec: rb_2d_std} and \ref{Sec: rb_2d}. 

For the $\omega$-red-black Gauss-Seidel smoother, (\ref{S_hat_rb_2d}) still holds, with $\symF S^J_h$ replaced by
\begin{equation}
\symF S^{GS}_h(\thv) = -\frac{60(1-\frac{1}{\omega})-32( \cos \th_1 + \cos \th_2) + e^{i2\th_1} + e^{i2\th_2}}{60\frac{1}{\omega} + e^{-i2\th_1} + e^{-i2\th_2}}
\end{equation}
as the symbol of (\ref{S_GS}).

%As in Section~\ref{Sec: rb_2d_std}, we consider an $(h,\ H)$ two-level cycle with the coarse grid $\grid G_H$ given by (\ref{G_H_2d_std}), and the same transfer operators (\ref{transfer_2d}).
%On the coarse grid $\grid G_H$, we compare the use of second-order operators (\ref{L_H_2d}) and (\ref{L_H_2G_2d}) to the performance with the fourth-order coarse-level operators  (\ref{L_H_4n_2d}) or (\ref{L_H_4G}).

The cubic interpolation operator and its adjoint restriction operator, considered in Section~\ref{Sec: orderC_t}, are given by
\begin{align}%\label{transfer_2d4o}
I_h^H &= \frac{1}{1024}
\begin{bmatrix}
1 & 0 & -9 & -16  & - 9 & 0 & 1\\
0 & 0 & 0 & 0  & 0 & 0 & 0\\
-9 & 0 & 81 & 144 & 81 & 0 & -9\\
-16 & 0 & 144 & 256 & 144 & 0 & -16\\
-9 & 0 & 81 & 144 & 81 & 0 & -9\\
0 & 0 & 0 & 0  & 0 & 0 & 0\\
1 & 0 & -9 & -16  & - 9 & 0 & 1\\
\end{bmatrix}_h^H, \label{interp_2d4o} \\
I^h_H &= \frac{1}{256} \left]
\begin{matrix}
1 & 0 & -9 & -16  & - 9 & 0 & 1\\
0 & 0 & 0 & 0  & 0 & 0 & 0\\
-9 & 0 & 81 & 144 & 81 & 0 & -9\\
-16 & 0 & 144 & 256 & 144 & 0 & -16\\
-9 & 0 & 81 & 144 & 81 & 0 & -9\\
0 & 0 & 0 & 0  & 0 & 0 & 0\\
1 & 0 & -9 & -16  & - 9 & 0 & 1\\
\end{matrix}\right[^h_H, \label{restr_2d4o}
\end{align} 
and (\ref{transfer_2dstd_eqn}) and (\ref{transfer_hat_2d}) hold with
\begin{equation}
\symF I_h^H(\th_1,\ \th_2) = [ 1- \sin^4 \frac{\th_1}{2} (3-2\sin^2 \frac{\th_1}{2} ) ] [ 1- \sin^4 \frac{\th_2}{2} (3-2\sin^2 \frac{\th_2}{2} ) ] 
,\quad [\th_1,\ \th_2]^T \in \Theta.
\end{equation}

\subsubsection{Heuristic argument: lower-order accurate coarse-level operators}\label{Sec: orderC_heuristic}
In this section, we seek to argue heuristically from local Fourier analysis why using lower-order accurate operators $L_H$ on coarse levels, in particular replacing (\ref{L_H_4n_2d}) or (\ref{L_H_4G}) by (\ref{L_H_2d}) or (\ref{L_H_2G_2d}) in the coarse-level correction operator $K_h^H$, does not generally hinder the convergence rate $\rho$ of the multigrid cycle. 
On the other hand, since in practice the coarse-level problems are solved only approximately (for instance by using multigrid recursively (as in \eqref{L_H_approx} for example), a related consideration is to ask how accurate the coarse-level solve $\tilde{L}_H^{-1}$ needs to be to maintain comparable rate of convergence as two-level cycles with exact solve on the coarse level. 
The discussion in this section follows and extends from the work of Hemker\cite{Hemker90}.

In general, consider a PDE of order $m = 2$ and its difference approximation ($L_h$) of order $p = 4$. Consider an $(h,\ H)$ two-level cycle to solve the difference equation with smoothing factor $\mu$ 
%(as defined in (\ref{mu_2_2})) 
and $\nu$ smoothing sweeps per cycle. Then its asymptotic convergence rate $\rho \equiv \rho(M_h^H) = \rho\big( K_h^H S_h^{\nu}\big)$ 
%(as defined in (\ref{rho_2d_std})) 
is comparable to $\mu^{\nu} = \rho\big( Q_h^H S_h^{\nu}\big)$, which corresponds to a cycle with an ideal coarse-level correction $Q_h^H$. In reality the coarse-level correction $K_h^H = I_h - I^h_H L^{-1}_H I_h^H L_h$ is used instead of $Q_h^H$. We discuss the approximations that can be made on $L^{-1}_H$ so that the convergence rate $\rho$ is still comparable to the reference rate $\mu^{\nu}$. 
To do this, we are going to consider the matrix representation $\symE K(\thv)\ \symE S^{\nu}(\thv)$\footnote{In this section we omit the subscript $\mathbb{F}$ of the $4\times 4$ matrices.} of the operator $K_h^H S_h^{\nu}$ on the eigenspaces $\mathbb{F}_{\thv}\ (|\thv| \neq 0)$, with 
%$\symE K (\thv)$ given by (\ref{K_hat_2})).
\begin{equation}\label{K_hat_hat}
\symE K (\thv) = \mathbf{I} - \symF L_H^{-1}(2\thv)\
\symE I^h_H (\thv)\ \symE I_h^H(\thv)\ \symE L(\thv).
\end{equation}

Rewrite the eigenspace (\ref{F}) as
\begin{equation}
\mathbb{F}_{\thv} = \text{span} \{\phi_h (\cdot,\ \thv),\quad  \bar{\mathbb{F}}_{\thv}\},\quad
 \thv \in \Theta^{low},
\end{equation}
where we denote the high-frequency subspace corresponding to $\thv \in \Theta^{low}$ as 
\begin{equation}\label{F_bar}
\bar{\mathbb{F}}_{\thv} \equiv \text{span} \{\phi_h (\cdot,\ \thbar_1,\ \thbar_2) ,\quad \phi_h (\cdot,\ \thbar_1,\ \th_2),\quad  \phi_h (\cdot,\ \th_1,\ \thbar_2) \}.
\end{equation}
Thus we rewrite the $4\times 4$ matrix representations on $\mathbb{F}_{\boldsymbol{\theta}}$ into a $2\times 2$ block structure. For instance, the matrix of $L_h$ on $\mathbb{F}_{\boldsymbol{\theta}}$ can be written as
\begin{equation}
\symE L(\thv) = \begin{bmatrix}
\symF L_h(\thv) \\ &  \symE L_{\bar{\mathbb{F}}}(\thv)
\end{bmatrix},\quad
\symE L_{\bar{\mathbb{F}}}(\thv) =  \begin{bmatrix}
\symF L_h(\bar \th_1,\ \bar \th_2) \\
 & \symF L_h(\bar \th_1,\ \th_2) \\
 & & \symF L_h(\th_1,\ \bar \th_2) 
\end{bmatrix},
\end{equation}
where $\symF L_h$ is given by (\ref{L_hat_2d_4o}). Note that this notation 
makes it straightforward to extend the following analysis to higher dimensions.
%can make the following analysis be easily extended to higher dimensions.

We now consider the perturbations of $\symE K (\thv)$ that result from the approximations discussed before. It is of particular interest, as we shall illustrate later, to have an estimate as $|\thv| \to 0$.

%\subsubsection{Order of transfers}
Theoretically, define the `ideal' restriction $q_h^H$ and interpolation $q^h_H$ so that
\begin{equation}
\begin{cases}
q_h^{H} \PhiVectorF = \phi_{H}(\cdot,\ 2\thv)\ \symE q_h^{H}(\thv),\\
q^h_{H}\ \phi_{H}(\cdot,\ 2\thv) = \PhiVectorF \symE q^h_{H} (\thv)
\end{cases}
\end{equation}
on $\grid G_H$,
with
\begin{equation}
\symE q_h^{H}(\thv) = \begin{bmatrix}
1 & \mathbf{0}
\end{bmatrix},\quad \symE q^h_{H}(\thv) = \symE q_h^{H}(\thv)^T.
\end{equation}
We have 
\begin{equation}\label{Q_bar}
\symE q^h_H (\thv)\ \symE q_h^H(\thv) = \bar{\symE Q} \eqdef \mathbf{I} - \mathbf{Q} = \begin{bmatrix}
1 \\ & \mathbf{O}\end{bmatrix}.
\end{equation}
(Here $\mathbf{O}$ denotes the $3\times 3$ zero matrix.)
The transfer operators (\ref{transfer_2d}) between the fine and coarse grids have order of accuracy $t=2 \leq p$,
while the operators (\ref{restr_2d4o}) and (\ref{interp_2d4o}) have order of accuracy $t=4 \leq p$,
that is,
\begin{equation}
\symE I_h^H(\thv) = \symE q_h^H(\thv) + O(|\thv|^t),
\end{equation}
as can be readily seen from (\ref{transfer_hat_2d}).
Thus 
\begin{equation}
\symE I^h_H (\thv)\ \symE I_h^H(\thv) = \begin{bmatrix}
1 + O(|\thv|^t)  & O(|\thv|^t) \\ O(|\thv|^t)  & O(|\thv|^{2t}) 
\end{bmatrix} = \bar{\symE Q} + O(|\thv|^t).
\end{equation}

Now we consider the approximations on $L_H^{-1}$.
%\subsubsection{Lower-order $L_H^{-1}$}\label{Sec: L_H_lo}
First, suppose the coarse-level problem is of order $q = 2 \leq p$, we have the matrix representation of $L^{-1}_H  L_h$ on $\mathbb{F}_{\thv}\ (|\thv| \neq 0)$ as
\ba \label{L_H_lo}
\symF L_H^{-1}(2\thv)\ \symE L(\thv) = 
\begin{bmatrix}
	1 + O(|\thv|^q) \\ & O(|\thv|^{-m})
\end{bmatrix}.
\ea
%
%\subsubsection{Approximate $\tilde{L}_H^{-1}$} \label{Sec: L_H_approx}
Secondly, suppose on the coarse level the error equation $L_H v_H = d_H$ is solved approximately (via iteration) by $ \tilde{v}_H = \tilde{L}_H^{-1} d_H$ so that the residual on this level
\begin{equation}\label{L_H_inv_approx}
\tilde{d}_H = (I_H - L_H\tilde{L}_H^{-1})\ d_H
\end{equation}
is reduced by a factor $\sigma$. Then we make a relative error $\e_{\thv}$ from $\symF L_H^{-1}(2\thv)$ to $(1 + \e_{\thv})\ \symF L_H^{-1} (2\thv)$ with $| \e_{\thv}| \leq \sigma$.
Note that here in $\e_{\thv}$ we characterize those approximations that are not particularly dependent on $\thv$, unlike the lower-order approximation in (\ref{L_H_lo}).
%
%\subsubsection{Convergence}
Combining the above two sources of perturbations (\ref{L_H_lo}) and (\ref{L_H_inv_approx}), we can estimate (\ref{K_hat_hat}) as
\ba 
\symE K (\thv) &= \mathbf{I} - (1 + \e_{\thv})
\begin{bmatrix}
	1 + O(|\thv|^t)  & O(|\thv|^t) \\ O(|\thv|^t)  & O(|\thv|^{2t}) 
\end{bmatrix} 
\begin{bmatrix}
	1 + O(|\thv|^q) \\ & O(|\thv|^{-m})
\end{bmatrix}\\
&= \symE Q -  \Big( \e_{\thv} \bar{\symE Q} + \begin{bmatrix}
	O(|\thv|^t)+ O(|\thv|^q) & O(|\thv|^{t-m}) \\ O(|\thv|^t) & O(|\thv|^{2t-m})
\end{bmatrix} \Big)\\
&= \symE Q -  \symE E(\thv), \label{K_hat_est}
\ea
where
\ba \label{E_hat_est}
\symE E(\thv) \eqdef
\begin{bmatrix}
 \e_{\thv}+ O(|\thv|^t)+ O(|\thv|^q) & O(|\thv|^{t-m}) \\ O(|\thv|^t) & O(|\thv|^{2t-m})
\end{bmatrix},
\ea
$|\e_{\thv}| \leq \sigma$.
Note that under the assumptions made for multigrid convergence it is necessary that 
%$2t-m \geq 0$, 
$t \geq \frac{m}{2}$~\cite{Hemker90}.
\commentB{
\textcolor{red}{Thus for convergence it is necessary that $\sigma = o(1)$ (wrong) so that the coarse-level correction $K_h^H$ is a $o(1)$ perturbation of the ideal $Q_h^H$:}
\begin{equation}
\symE K(\thv) = \symE Q - \symE E(\thv),\quad \|\symE E (\thv)\| = o(1) wrong,\ |\thv| \to 0.
\end{equation}
}

For the asymptotic convergence rate of the multigrid cycle $\rho = \rho(M_h^H) = \rho(K_h^H S_h^{\nu})$, consider the operator $K_h^H S_h^{\nu}$. On each $\mathbb{F}_{\thv}$ its matrix representation
\begin{equation}
\symE K(\thv)\ \symE S^{\nu}(\thv) = \symE Q\ \symE S^{\nu}(\thv) - \symE E(\thv)\ \symE S^{\nu}(\thv). 
\end{equation}
The first part gives the $\rho(Q_h^H S_h^{\nu}) = \mu^{\nu}$ when $\sup$ over $\thv \in \Theta^{low}$, giving the reference convergence rate. 
For simplicity, we can roughly characterize the symbol of the smoother as 
\begin{equation}\label{S_hat_rough}
\symE S(\thv) \approx \begin{bmatrix}
1\\  & \mu
\end{bmatrix},
\end{equation} 
for $\thv \in \Theta^{low}$. Then we have the perturbation 
\begin{equation}
\symE E(\thv)\ \symE S^{\nu}(\thv) \approx \begin{bmatrix}
\sigma + O(|\thv|^t)+ O(|\thv|^q) & O(|\thv|^{t-m})\mu^{\nu} \\ O(|\thv|^t) & O(|\thv|^{2t-m})\mu^{\nu}
\end{bmatrix}
\end{equation}
So we have the rough estimate
\begin{equation}
\rho(\symE K(\thv)\ \symE S^{\nu}(\thv)) \approx \max\{\sigma + O(|\thv|^t)+ O(|\thv|^q),\quad \mu^{\nu}\},\quad |\thv| \to 0,
\end{equation}
the maximum of the reduction ratios of low-frequency and high-frequency modes, and one would need $\sigma$ to be approximately the size as $\mu^{\nu}$, as discussed in Section~\ref{M1: orderC}. 
(Note that the above estimates only apply for $|\thv|$ small, since higher-order terms in $|\thv|$ are omitted.)
Thus we've established that $\rho(\symE K(\thv)\ \symE S^{\nu}(\thv)) \approx \mu^{\nu}$ for $|\thv|$ small.

\smallskip
Actually, the perturbation $\symE E(\thv)\ \symE S^{\nu}(\thv)$ of the iteration operator is small for $\thv$ across the full range $\Theta^{low}$, and thus we have comparable convergence $\rho \approx \mu^{\nu}$. Heuristically, for $|\thv| \to 0$ the perturbation is small because $\symE E(\thv)$ is small shown by \eqref{E_hat_est}, while for $\thv$ away from $\mathbf{0}$ the smoother actually does better than \eqref{S_hat_rough}, and has the the property that $\symE S(\thv)$ is `small', as illustrated in Figure~\ref{fig:rho_S_4o}\footnote{Note that in the illustrations we switch the smoother from RB-GS to RB-J, just so that we have symmetry in $\th_1$ and $\th_2$ and the pictures look more pleasant. This change is irrelevant to the arguments we make.} (left). As we can see, since the smoother does so good a job when $\thv$ is away from $\mathbf{0}$, it is when $|\thv|$ is small where the coarse-level correction needs to work the hardest. The smoothing factor, corresponding to the ideal coarse-level correction, is illustrated in Figure~\ref{fig:rho_S_4o} (right), where we note that $\sup_{\thv \in \Theta^{low}} \rho (\symE Q \symE S) =  \mu$. In accordance with \eqref{K_hat_est}, the situation is very similar with the actual coarse-level correction $K$ since where it matters, that is when $|\thv|$ is small, $K$ behaves very similar as $Q$. This is illustrated in 
Figure~\ref{fig:rho_M_4o}\footnote{Note that we illustrate a matrix symbol by its spectral radius. However, we stress again that the spectral radii do not multiply.}, where $\rho(\symE M(\thv)) = \rho(\symE K(\thv)\ \symE S(\thv))$ (right) is to be compared with $\rho(\symE Q\ \symE S(\thv))$ in 
Figure~\ref{fig:rho_S_4o} (right). 

{
	\newcommand{\figWidth}{5cm}% height 
	\newcommand{\trimfig}[2]{\trimh{#1}{#2}{.0}{.0}{.0}{.0}}
	\begin{figure}[h!]
		\begin{center}
			\begin{tikzpicture}[scale=1]
			\useasboundingbox (0,0.5) rectangle (15,5.);  % set the bounding box (so we have less surrounding white space)
			
			\draw(-0.5,0.0) node[anchor=south west,xshift=-15pt,yshift=-8pt] {\trimfig{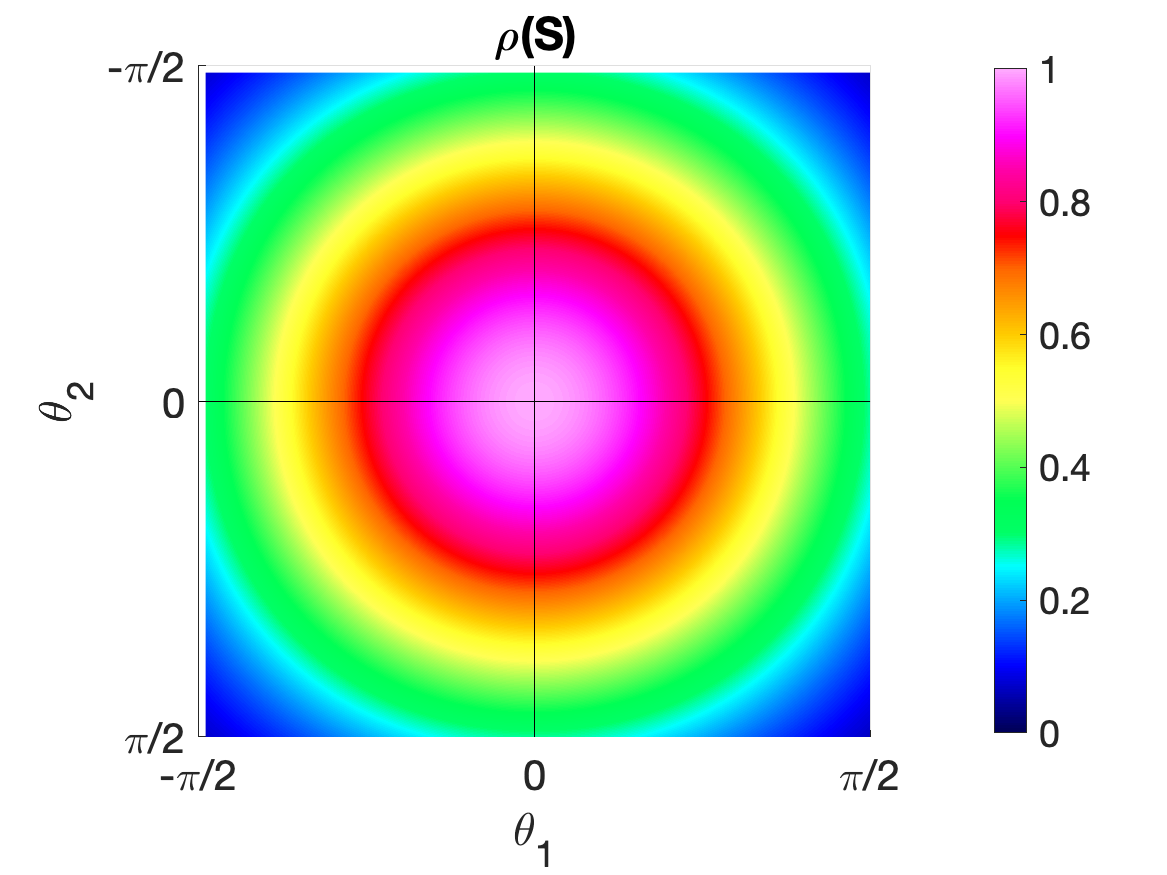}{\figWidth}};
			\draw(7.5,0.0) node[anchor=south west,xshift=-15pt,yshift=-8pt] {\trimfig{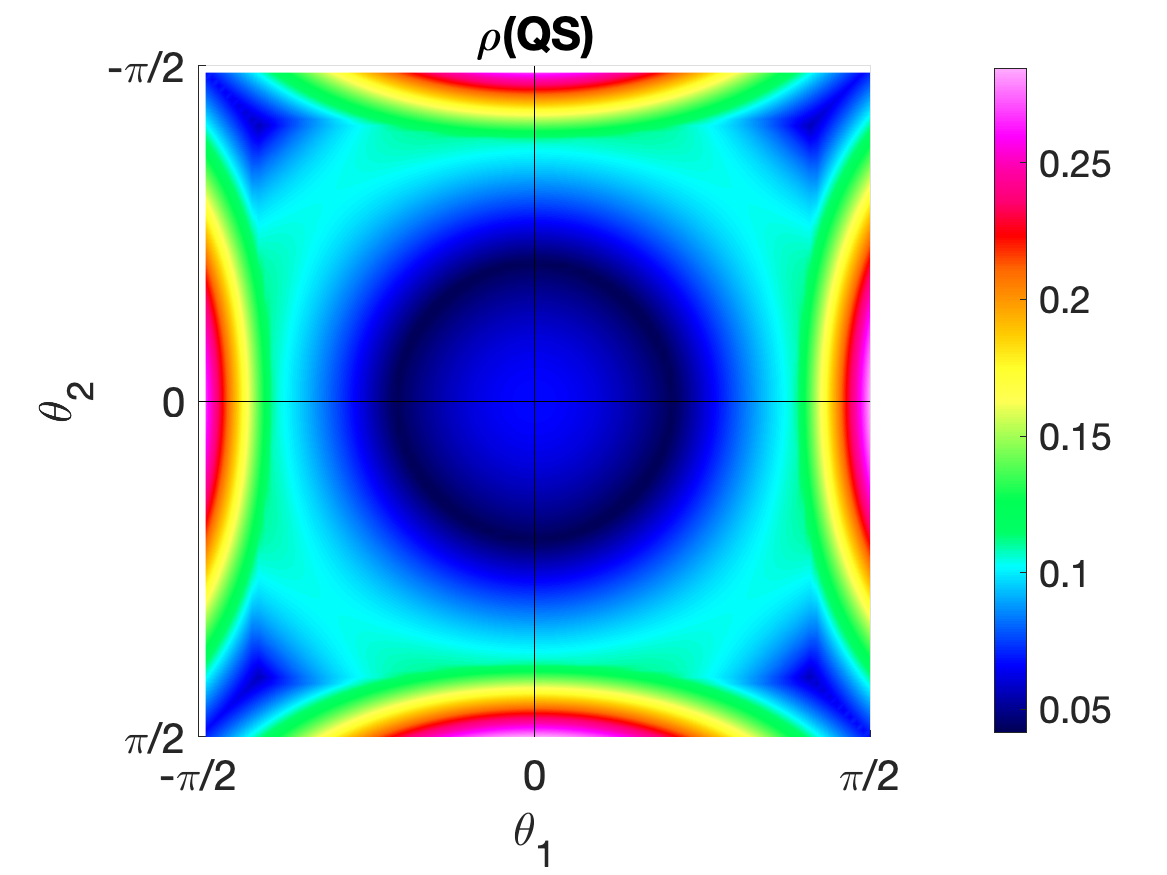}{\figWidth}};
			% grid:
			% \draw[step=1cm,gray] (0,0) grid (16,5);
		
			\end{tikzpicture}
		\end{center}
		\captionof{figure}{Illustration of the red-black-Jacobi
			 smoother (with $\omega=1$, $\nu=1$) on $\Theta^{low}$, 2D, 4th-order accurate, standard coarsening. Left: $\rho(\symE S(\thv))$, right: $\rho(\symE Q\ \symE S(\thv))$. }
		\label{fig:rho_S_4o}
	\end{figure}
}
{	\newcommand{\figWidth}{5.2cm}% height 
	\newcommand{\trimfig}[2]{\trimh{#1}{#2}{.05}{.225}{.0}{.0}}
	% -- for colour bar ----
	\newcommand{\cbWidth}{.2cm}% colour bar width
	\newcommand{\cbHeight}{4cm}% colour bar height
	\newcommand{\xcb}{.5cm}% colour bar lower left corner
	\newcommand{\ycb}{-.2cm}% colour bar lower left corner
	\setlength{\ycbTop}{\ycb+\cbHeight}% colour bar top label position
	\setlength{\ycbMid}{\ycb+\cbHeight*\real{.5}}% colour bar top label position
	\newcommand{\trimfigcb}[3]{\includegraphics[width=#2, height=#3, clip, trim=17cm 2.35cm 1.65cm 2.35cm]{#1}}
	% horizontal colour bar: \colourBar[shiftCommands]{cMin}{cMax}
	\newcommand{\colourBarHorizontal}[3]{
		\begin{scope}[#1]
			\draw (\xcb,\ycb) node[anchor=south west,xshift=0.25cm,yshift=.5cm,rotate=-90] {\trimfigcb{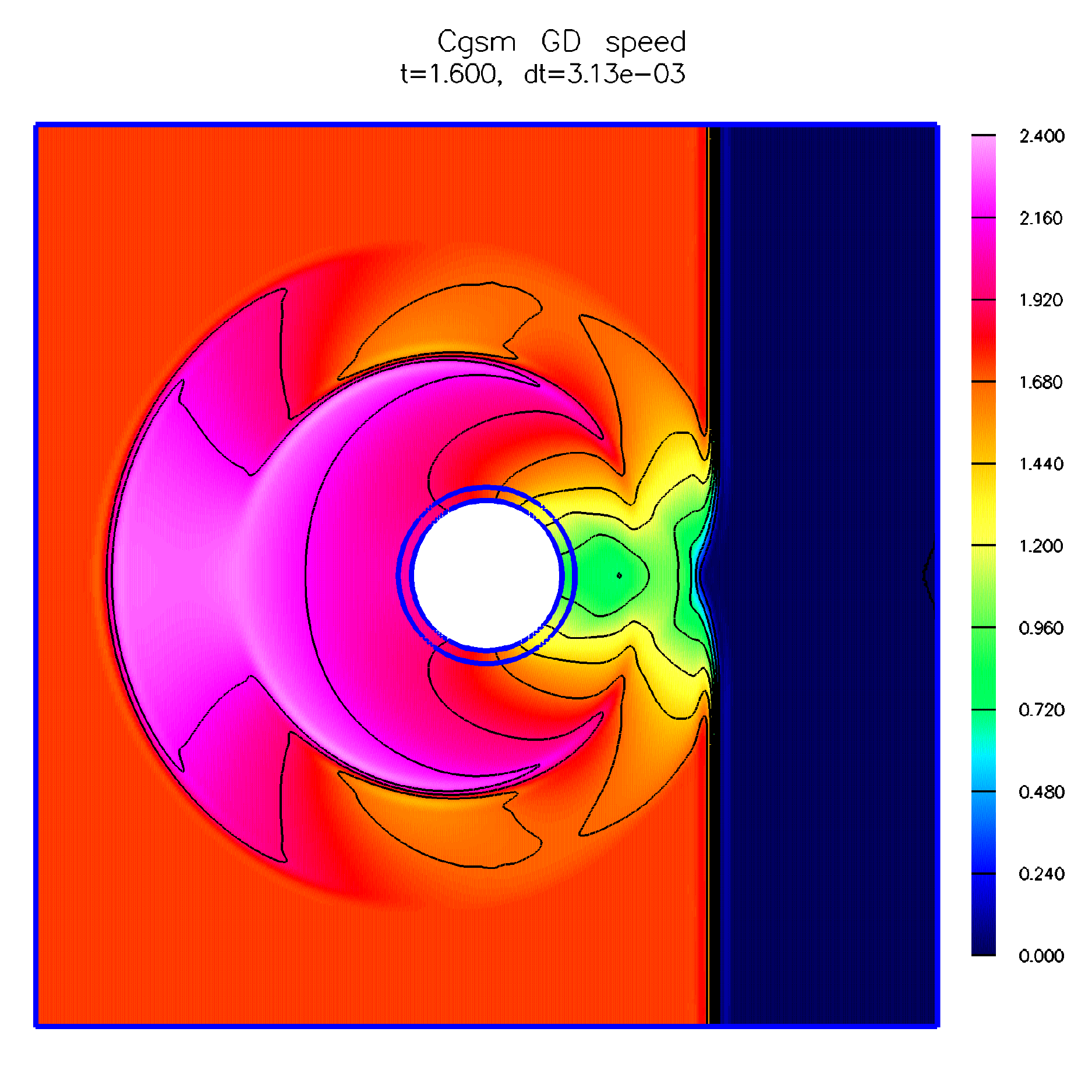}{\cbWidth}{\cbHeight}};
			\draw (0.35,0.25) node[anchor=north,xshift=+3pt,yshift=+2pt] {\scriptsize $#2$};
			\draw (5.25,0.25) node[anchor=north,xshift=+0pt,yshift=+2pt] {\scriptsize $#3$};
		\end{scope}
	}% end colour bar 
	\begin{figure}[t]
		\begin{center}
			\begin{tikzpicture}[scale=1]
			\useasboundingbox (0,0.5) rectangle (15,22);  % set the bounding box (so we have less surrounding white space)
			
			\draw(0,17) node[anchor=south west,xshift=-15pt,yshift=-8pt] {\trimfig{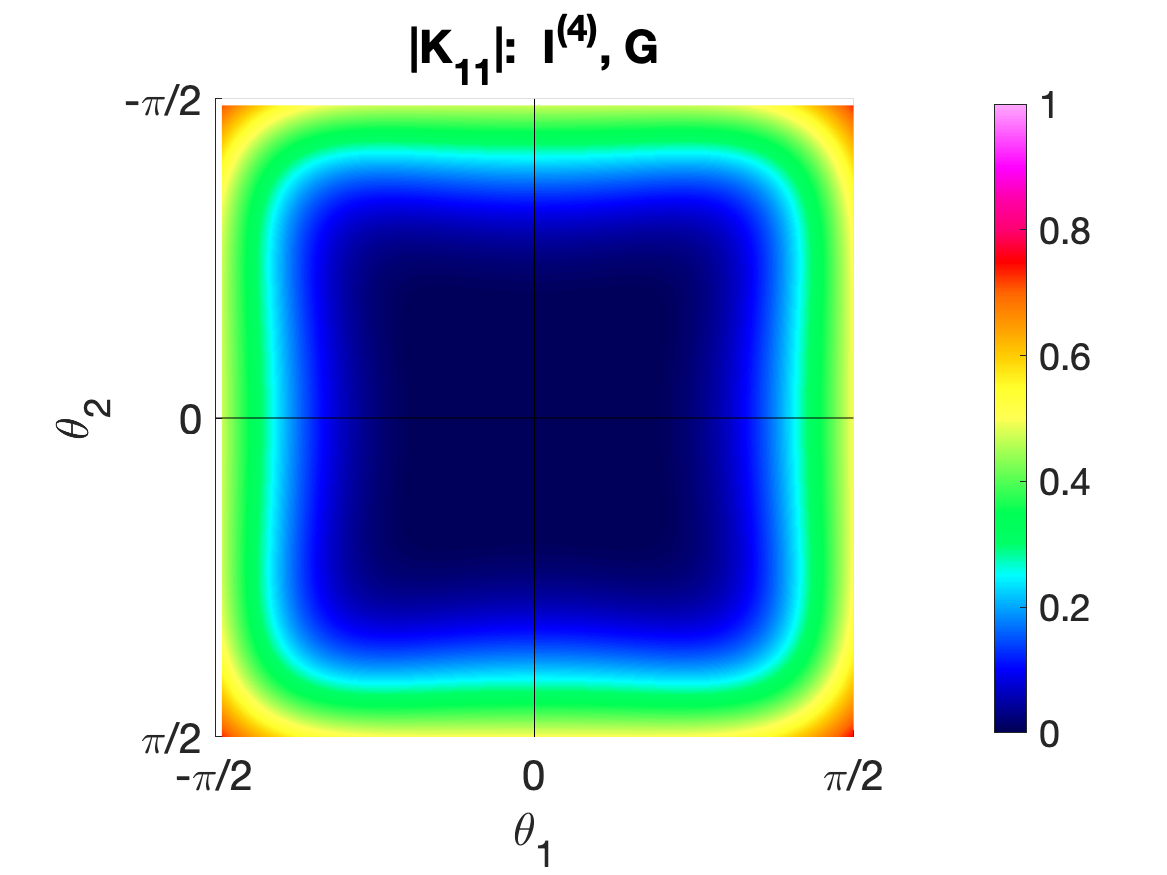}{\figWidth}};
			\draw(7.5,17) node[anchor=south west,xshift=-15pt,yshift=-8pt] {\trimfig{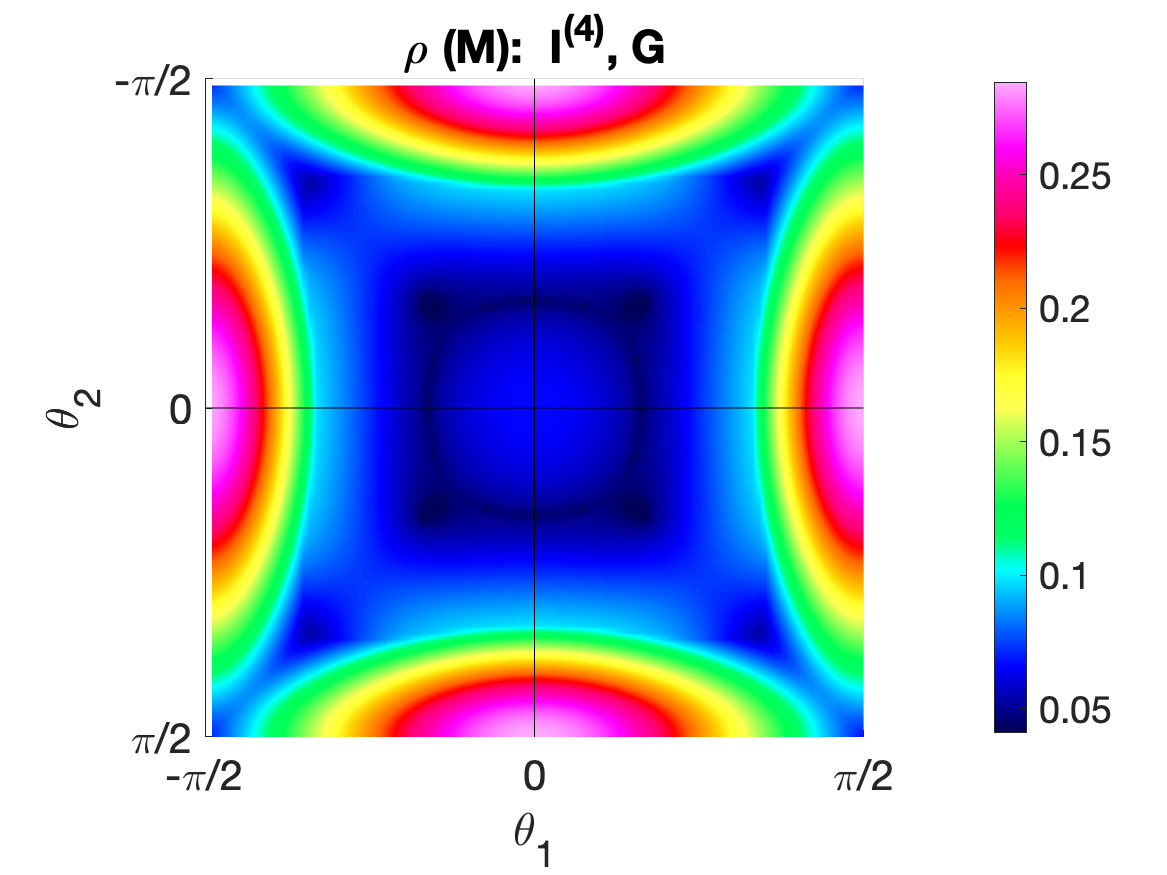}{\figWidth}};
			
			\draw(0,11.5) node[anchor=south west,xshift=-15pt,yshift=-8pt] {\trimfig{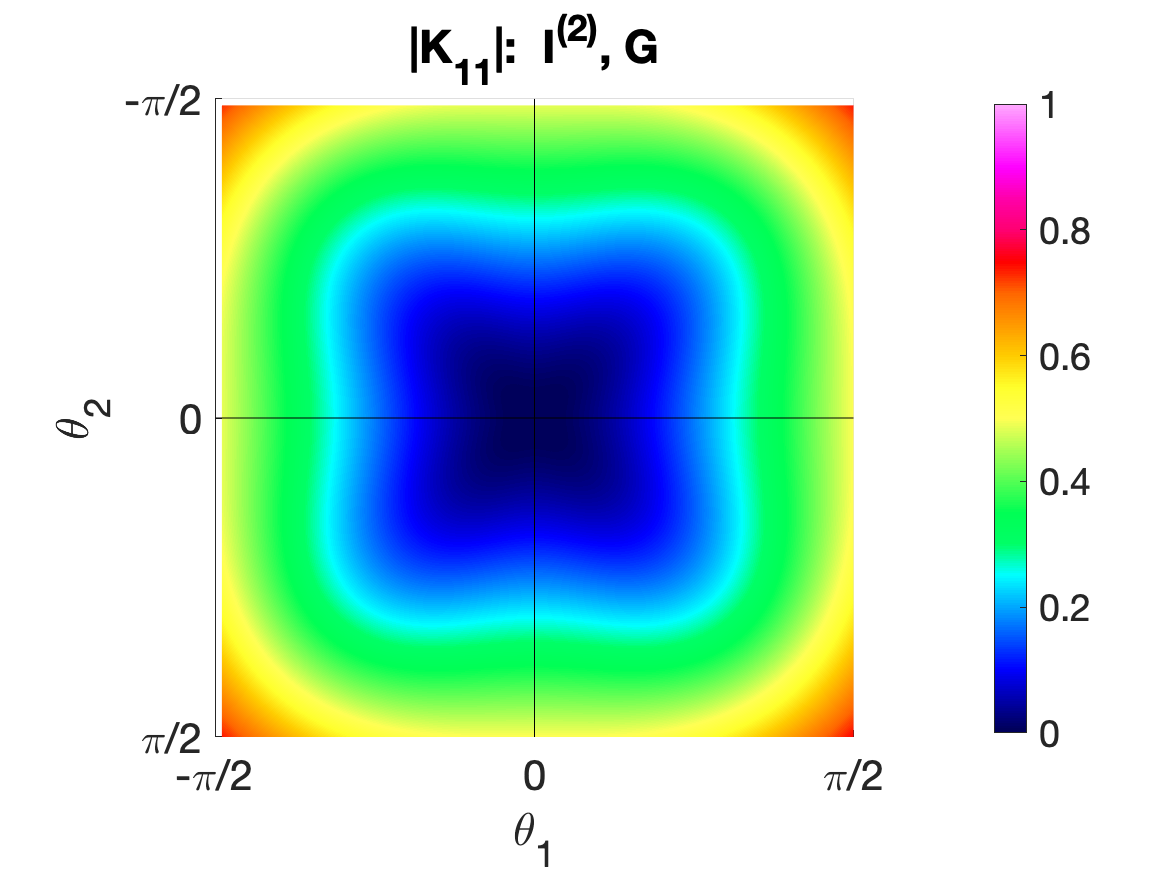}{\figWidth}};
			\draw(7.5,11.5) node[anchor=south west,xshift=-15pt,yshift=-8pt] {\trimfig{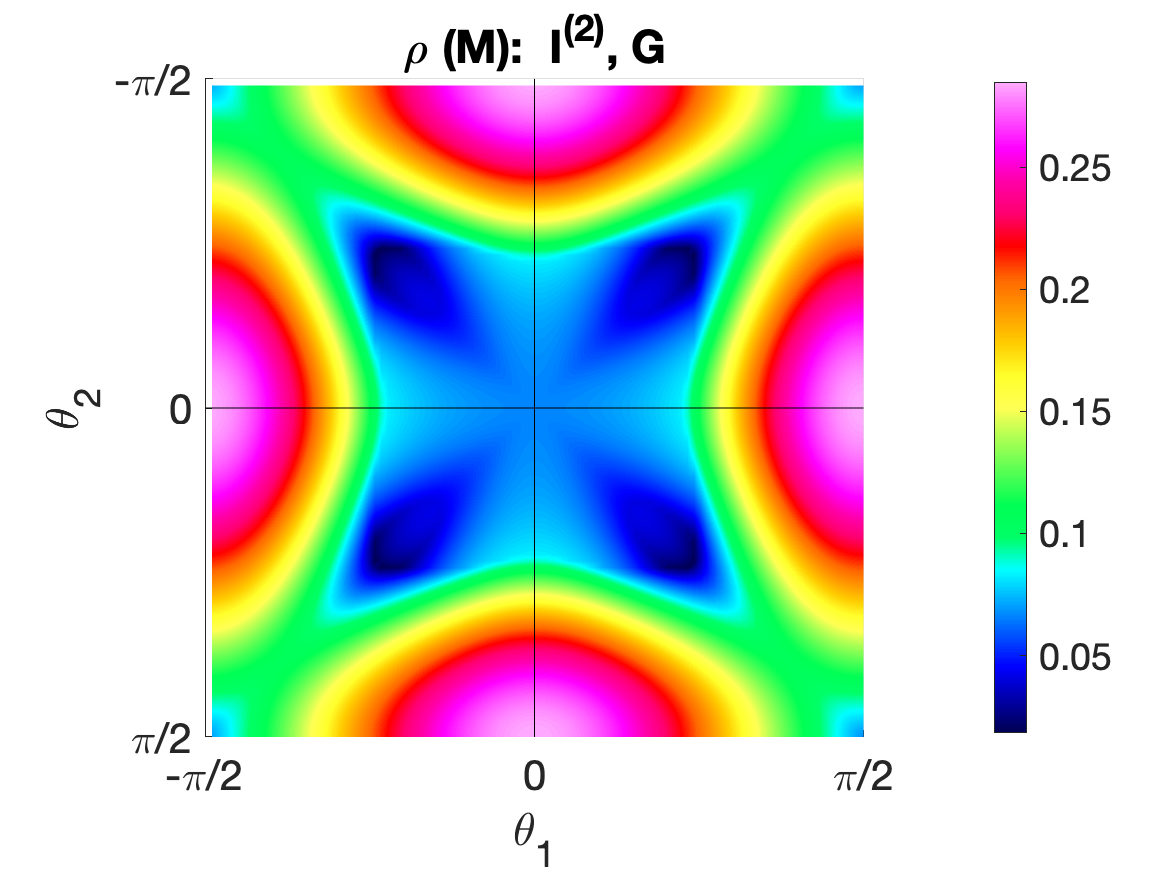}{\figWidth}};
			
			\draw(0,6) node[anchor=south west,xshift=-15pt,yshift=-8pt] {\trimfig{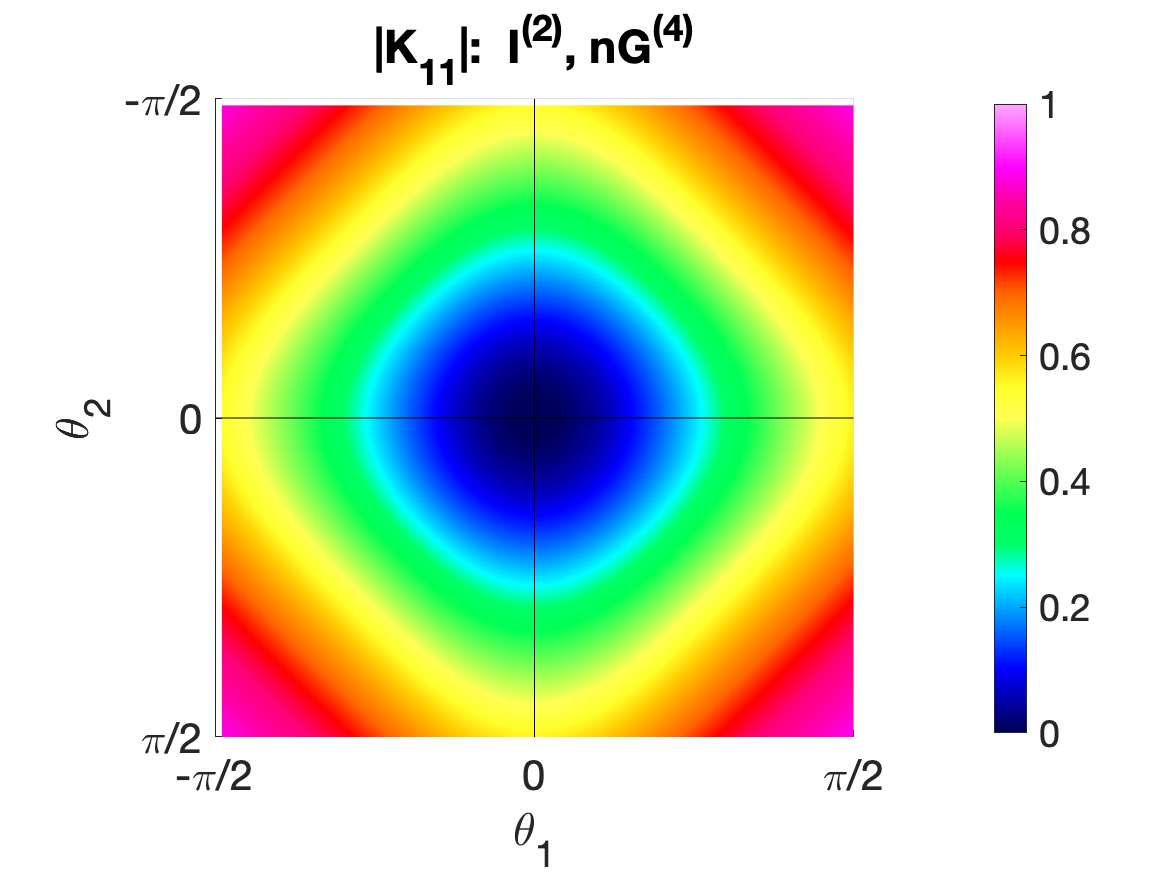}{\figWidth}};
			\draw(7.5,6) node[anchor=south west,xshift=-15pt,yshift=-8pt] {\trimfig{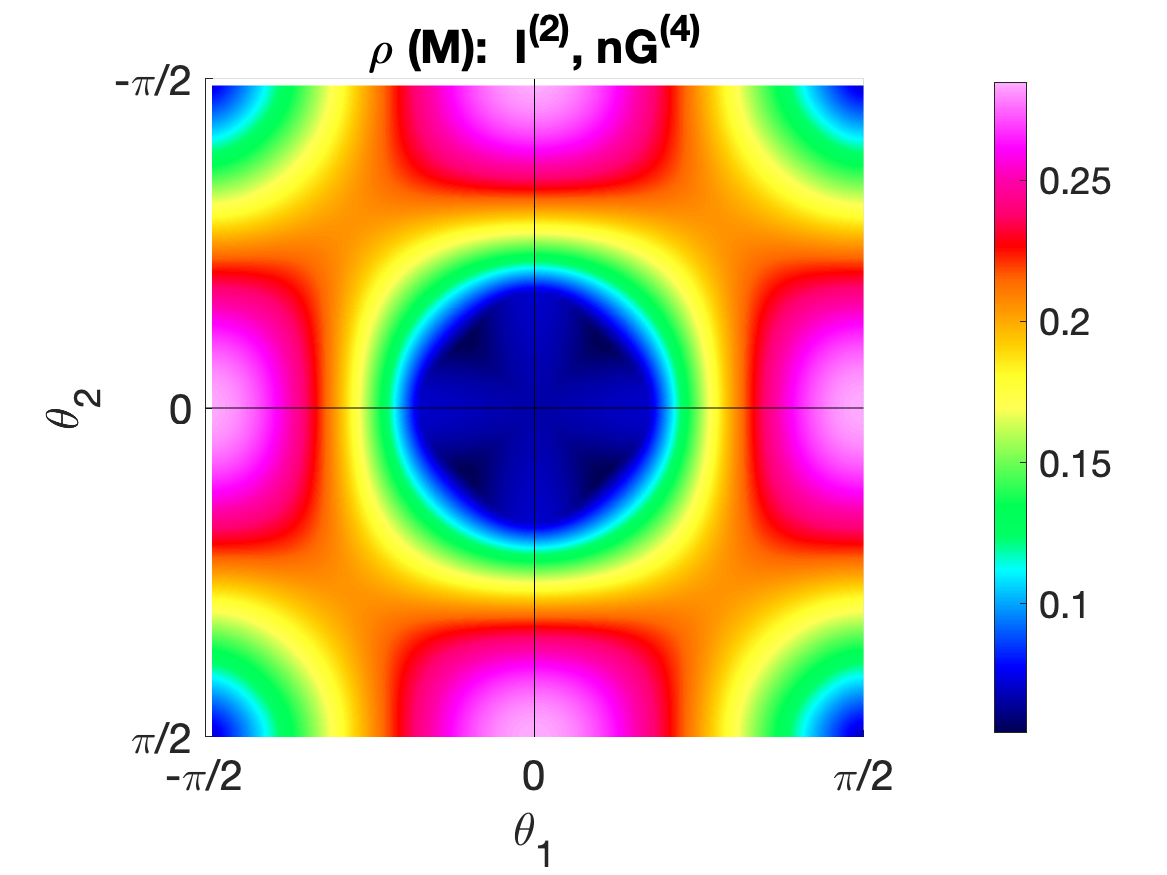}{\figWidth}};
			
			\draw(0,.5) node[anchor=south west,xshift=-15pt,yshift=-8pt] {\trimfig{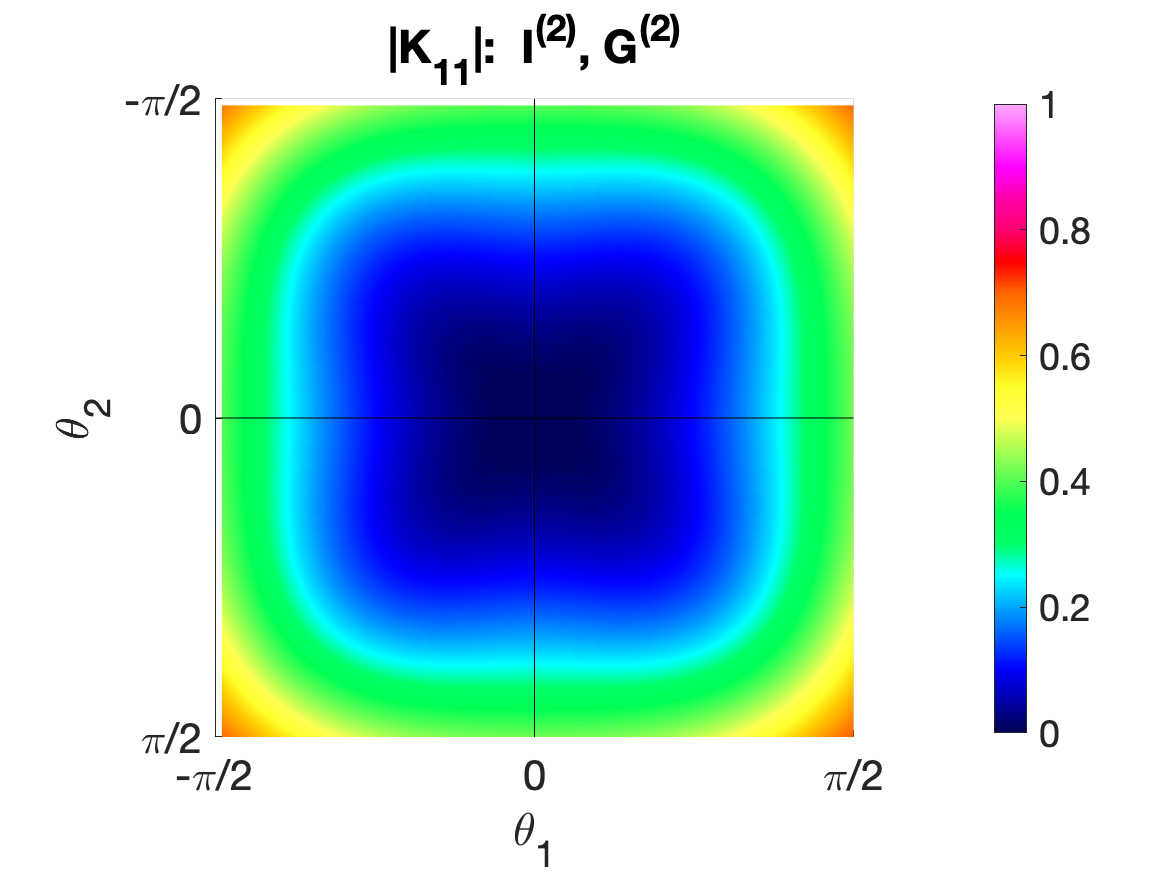}{\figWidth}};
			\draw(7.5,.5) node[anchor=south west,xshift=-15pt,yshift=-8pt] {\trimfig{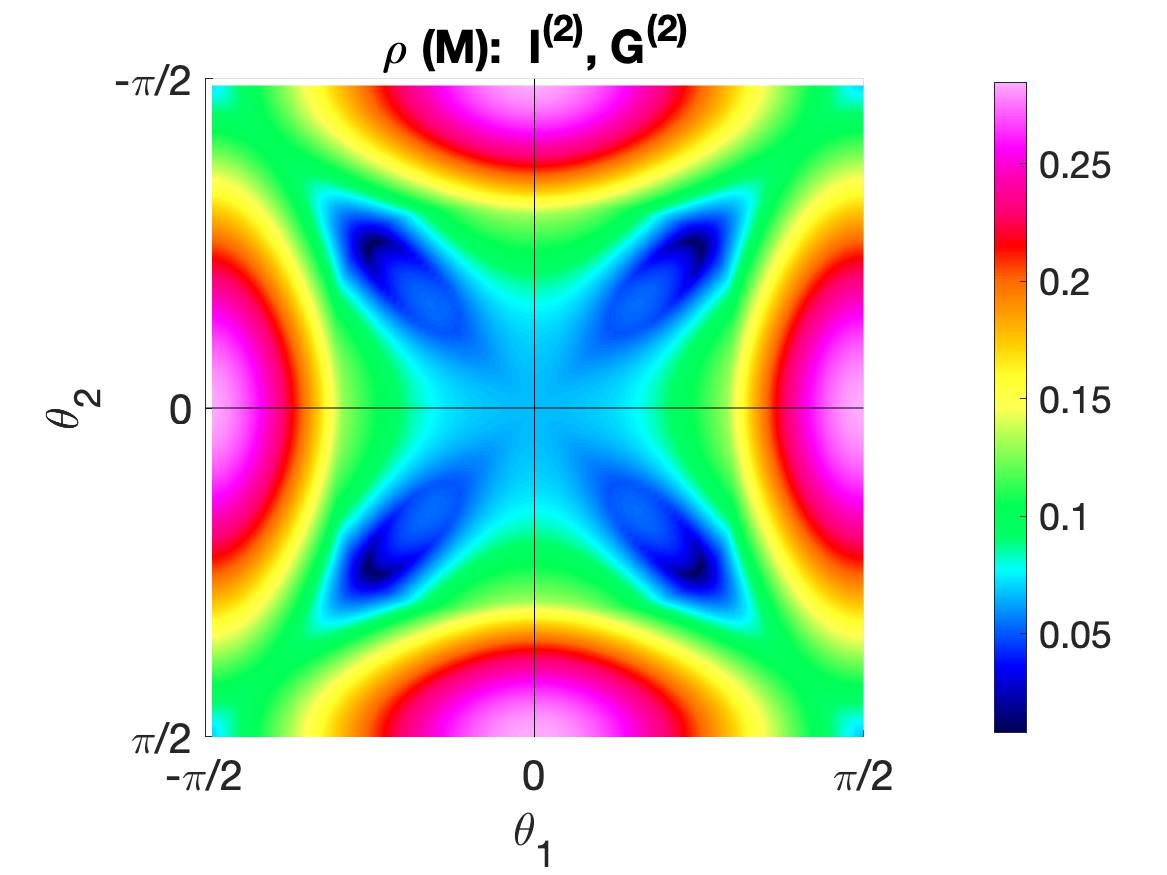}{\figWidth}};
			
			\colourBarHorizontal{xshift=.0cm,yshift=0cm}{0}{1};
			\colourBarHorizontal{xshift=7.5cm,yshift=0cm}{0}{0.3};
			
			% grid:
			 %\draw[step=1cm,gray] (0,0) grid (16, 22);
			\end{tikzpicture}
		\end{center}
		\captionof{figure}{Illustration of the coarse-level correction (with various transfer operators (I) and coarse-level operators) and the iteration operator (with red-black-Jacobi smoothing, $\omega=1$, $\nu=1$) on $\Theta^{low}$, 2D, 4th-order accurate, standard coarsening. Left: $|\symE K_{11}(\thv)|$, right: $\rho(\symE M(\thv)) = \rho(\symE K(\thv)\ \symE S(\thv))$. }
		\label{fig:rho_M_4o}
	\end{figure}
}

Figure~\ref{fig:rho_M_4o} (left) is meant to illustrate the fact that $K$ is close to $Q$ when $|\thv|$ is small, as in \eqref{K_hat_est}. Note that the ideal coarse-level correction $Q$ is 1 for high-frequency modes and 0 for low-frequency modes. As for the actual coarse-level correction $K$, 
$\rho(\symE K(\thv)) = 1$ always; and we focus on the behavior of $K$ on the low-frequency modes.
Analogous to the smoothing factor corresponding to $\rho(\symE Q\ \symE S^{\nu}(\thv))$ with ideal coarse-level correction, we illustrate the coarse-level correction $K$ by
\begin{equation}
\rho(\symE K(\thv)\ \bar{\symE Q}) = |\symE K_{11}(\thv)|,
\end{equation}
where $\bar{\symE Q}$, as defined in \eqref{Q_bar}, represents an `idealized' smoother that eliminates all of the high-frequency component of the error.
(For the ideal coarse-level correction $Q$, in comparison, $\symE Q_{11} = 0$.)
We can see that the various choices of coarse-level corrections $K$ are all very small (and similar) near $\thv = \mathbf{0}$; where they are different is when $\thv$ is near the corners, for which the smoother is very good. The Galerkin coarse-level correction operators only behave badly near the tiny corners, and are more similar to the ideal $Q$.
Thus we can see the iteration operator $M_h$ has its worst modes at neither extremes.

\subsection{Work-units estimates}\label{sec: WU}
In this section we give a rough estimate of the work-units (WU) of the multigrid cycles for the model problem \eqref{MP_h} in $d$ dimensions with second-order accuracy.
Consider a fine-level problem of size $N^d$ with $L_h$, and a general $(\lmax +1)$-level cycle.
Assume evenly distributed grid points on each level, 
and the number of grid points of level $\ll$ to be $N_{\ll}^d$, $\ll=0,\ 1,\ \cdots,\ \lmax$, $N_0 \equiv N$. 
Suppose a coarsening factor so that
\begin{equation}\label{N_d_C}
\frac{N_{\ll}^d}{N_{\ll+1}^d} = r^d,\quad \ll = 0,\ 1,\ \cdots,\ \lmax.
\end{equation}
For example, for standard coarsening $r=2$; for red-black coarsening (Section~\ref{M2: rbC}) $r^d=2$; for general factor coarsening (Section~\ref{M3: rC}) it is the same $r$ here.
Suppose, per usual, the number of levels is determined by the size of the problem $N^d$ as
\begin{equation}
\lmax(r) \approx \log_r \frac{N}{\Nmin},
\end{equation} 
%\eqref{eq:numberOfLevelsVersusR},
and $N_{\lmax} \approx \Nmin$.

Note that for the second-order accurate discretization in $d$ dimension, one work unit is defined as the work (FLOPs) for one Jacobi smoothing step of $L_h$ with the $(2d+1)$-point stencil, which we estimate as
approximately $(2d+1) N^d$.

On a level $\ll$ with grid size $N_{\ll}^d$, we estimate the FLOPs for a single smoothing sweep to be about $(2d+1) N_{\ll}^d$,
that is, about $\frac{N_{\ll}^d}{N^d} = \frac{1}{r^{d\ll}}$ work units
(even though the Galerkin operators may have fuller or wider stencils), 
and the work for computing the residual to be about the same. 
We estimate the transfers between levels ${\ll}$ and ${{\ll}+1}$ to be also about 
$\frac{1}{r^{d\ll}}$ work units for both standard coarsening and red-black coarsening.
On the other hand, one subtlety arises that the interpolation involves more work units when the coarse grid is not embedded in the fine grid.
If $\grid G_{h_{\ll+1}} \subset \grid G_{h_{\ll}}$, one only needs to interpolate grid points ${\grid G_{h_{\ll}}}\setminus {\grid G_{h_{\ll+1}}}$ (and some of these points has even smaller interpolation stencils according to the location of the fine-grid point with respect the coarse grid).
The restriction operator, chosen as the adjoint of the interpolation operator, requires roughly the same amount of work units as the interpolation.
%(because the transpose conserves the sparsity).
So in general for coarsening factor as in \eqref{N_d_C},
we estimate the WU of the transfers between levels $\ll$ and $\ll+1$ to be about twice as the WU in the case of standard coarsening (or red-black coarsening). 
On the coarsest grid with size $N_{\lmax}^d \approx \Nmin^d$, the direct (or iterative) solver for the banded system is estimated (somewhat generously)
to cost about 
%\kP{$2 N_{min}^{2d-1}$ } or 
\kP{$(2d+1) \Nmin^{d+1}$ }
FLOPs. 

Summing the WU on each level we have for a $\gamma$-cycle with $\nu$ smoothing sweeps per cycle, 
\begin{align}
\text{WU}[r^d;\ \gamma,\ \nu] &= 
%\frac{2}{2d+1} \Nmin^{d-1} \frac{\gamma^{{\lmax}-1}}{r^{d\ \lmax}} 
 \Nmin \frac{\gamma^{{\lmax}-1}}{r^{d\ \lmax}} 
+ (\nu+1+T) \sum_{{\ll}=0}^{{\lmax}-1} (\frac{\gamma}{r^d} )^{\ll}\\
&\approx (\nu+1+T)
\begin{cases}
\frac{r^d}{r^d-\gamma}\quad &\gamma < r^d,\\
\lmax \quad & \gamma = r^d,\\
\frac{r^d}{\gamma-r^d} \big(\frac{\gamma}{r^d} \big)^{\lmax} \quad &\gamma < r^d,
\end{cases} 
\end{align}
where $T=2$, unless the coarse grid is embedded in the fine grid (standard or red-black coarsening) in which case $T=1$. 
Note that the number of levels is assumed to be large enough so that the cost of coarsest-level solve is comparably small;
and also it is more reasonable to stay in the case $\gamma < r^d$ rather than $\gamma \geq r^d$.

%\bigskip
\clearpage
\bibliography{henshaw,henshawPapers,mg}
\bibliographystyle{elsart-num}
\end{document}